\documentclass[12pt]{amsbook}

\usepackage{amsmath, amssymb, amsfonts, amsthm, amscd}
\usepackage{manfnt}
\usepackage{mathtools}
\usepackage{fullpage}
\usepackage[all]{xy}
   \SelectTips{cm}{10}
\usepackage{txfonts}
\usepackage{amsmidx} 
\usepackage{url}

\newtheorem{thm}{Theorem}[section]
\newtheorem{lemma}[thm]{Lemma}
\newtheorem{prop}[thm]{Proposition}
\newtheorem{cor}[thm]{Corollary}

\newtheorem{conj}[thm]{Conjecture}
\newtheorem{hypothesis}[thm]{Hypothesis}
\newtheorem{prop-conj}[thm]{Proposition-Conjecture}

\theoremstyle{definition}
\newtheorem{defn}[thm]{Definition}

\theoremstyle{remark}
\newtheorem{rmk}[thm]{Remark}

\theoremstyle{remark}
\newtheorem{notation}[thm]{Notation}

\theoremstyle{remark}
\newtheorem{question}[thm]{Question}

\theoremstyle{remark}
\newtheorem{eg}[thm]{Example}

\newcommand{\Q}{\mathbb{Q}}
\newcommand{\Qb}{\overline{\mathbb{Q}}}
\newcommand{\Z}{\mathbb{Z}}
\newcommand{\CC}{\mathbb{C}}
\newcommand{\RR}{\mathbb{R}}
\newcommand{\bS}{\mathbb{S}}
\newcommand{\Ql}{\mathbb{Q}_{\ell}}
\newcommand{\Qlb}{\overline{\mathbb{Q}}_\ell}
\newcommand{\Elb}{\overline{E_{\lambda}}}
\newcommand{\af}{\mathbf{A}_F}
\newcommand{\afp}{\mathbf{A}_{F'}}

\DeclareMathOperator{\Hom}{Hom}
\DeclareMathOperator{\Ext}{Ext}
\DeclareMathOperator{\End}{End}
\DeclareMathOperator{\Aut}{Aut}
\DeclareMathOperator{\Ind}{Ind}
\DeclareMathOperator{\BC}{BC}
\DeclareMathOperator{\Sym}{Sym}
\DeclareMathOperator{\Ad}{Ad}
\DeclareMathOperator{\rec}{rec}
\DeclareMathOperator{\tr}{tr}
\DeclareMathOperator{\im}{im}
\DeclareMathOperator{\coker}{coker}
\DeclareMathOperator{\Res}{Res}

\DeclareMathOperator{\gr}{gr}

\DeclareMathOperator{\pr}{pr}
\DeclareMathOperator{\Rep}{Rep}
\DeclareMathOperator{\Spec}{Spec}
\DeclareMathOperator{\id}{id}

\newcommand{\gal}[1]{\Gamma_{#1}} % an absolute Galois group
\newcommand{\Gal}{\mathrm{Gal}} % a relative Galois group
\newcommand{\Bdr}{\mathrm{B}_{\mathrm{dR}}}

\newcommand{\Bht}{\mathrm{B}_{\mathrm{HT}}}
\newcommand{\Ddr}{\mathrm{D}_{\mathrm{dR}}}
\newcommand{\Dcris}{\mathrm{D}_{\mathrm{cris}}}
\newcommand{\Dht}{\mathrm{D}_{\mathrm{HT}}}

\newcommand{\into}{\hookrightarrow}
\newcommand{\onto}{\twoheadrightarrow}

\newcommand{\mc}{\mathcal}
\newcommand{\mf}{\mathfrak}
\newcommand{\mr}{\mathrm}
\newcommand{\mbf}{\mathbf}
\newcommand{\mbb}{\mathbb}

\newcommand{\rl}{\rho_{\lambda}}

\newcommand{\tG}{\widetilde{G}}
\newcommand{\tT}{\widetilde{T}}
\newcommand{\tZ}{\widetilde{Z}}

\newcommand{\tpi}{\tilde{\pi}}
\newcommand{\tomega}{\tilde{\omega}}
\newcommand{\tmu}{\tilde{\mu}}
\newcommand{\tnu}{\tilde{\nu}}

\newcommand{\bs}{\backslash}

\makeindex{t} %index of terms
\makeindex{s} %index of symbols
\begin{document}
\frontmatter

\title{Variations on a theorem of Tate}
\author{Stefan Patrikis}
\address{Princeton University, Department of Mathematics\\ Fine Hall \\ Washington Road\\ Princeton, NJ 08544}
\curraddr{Department of Mathematics \\  MIT \\ Cambridge, MA 02139}
\email{patrikis@math.mit.edu}
\urladdr{\url{http://math.mit.edu/~patrikis/}}
\thanks{Partially supported by an NDSEG fellowship and NSF grant DMS-1303928.}



\date{}

\begin{abstract}
Let $F$ be a number field. These notes explore Galois-theoretic, automorphic, and motivic analogues and refinements of Tate's basic result that continuous projective representations $\mathrm{Gal}(\overline{F}/F) \to \mathrm{PGL}_n(\mathbb{C})$ lift to $\mathrm{GL}_n(\mathbb{C})$. We take special interest in the interaction of this result with algebraicity (for automorphic representations) and geometricity (in the sense of Fontaine-Mazur). On the motivic side, we study refinements and generalizations of the classical Kuga-Satake construction. Some auxiliary results touch on: possible infinity-types of algebraic automorphic representations; comparison of the automorphic and Galois ``Tannakian formalisms"; monodromy (independence-of-$\ell$) questions for abstract Galois representations.
\end{abstract}

\maketitle
%\dedicatory{}

\tableofcontents

\mainmatter

\chapter{Introduction}
\section{Introduction}
Let $F$ be a number field, with $\gal{F}= \Gal(\overline{F}/F)$ its Galois group relative to a choice of algebraic closure. Tate's theorem that $H^2(\gal{F}, \Q/\Z)$ vanishes (see \cite[Theorem 4]{serre:DSsurvey}) encodes one of the basic features of the representation theory of $\gal{F}$: since $\coker(\mu_{\infty}(\CC) \to \CC^\times)$ is uniquely divisible, $H^2(\gal{F}, \CC^\times)$ vanishes as well, and therefore all obstructions to lifting projective representations $\gal{F} \to \mr{PGL}_n(\CC)$ vanish. More generally, as explained in \cite{conrad:dualGW}, for any surjection $H' \onto H$ of linear algebraic groups over $\Qlb$ with central torus kernel, any homomorphism $\gal{F} \to H(\Qlb)$ lifts to $H'(\Qlb)$. The simplicity of Tate's theorem is striking: replacing $\Q/\Z$ by some finite group $\Z/m$ of coefficients, the vanishing result breaks down.

Two other groups, more or less fanciful, extend $\gal{F}$ and conjecturally encode the key structural features of, respectively, automorphic representations and (pure) motives. The more remote automorphic Langlands group, $\mc{L}_F$, at present exists primarily as heuristic, whereas the motivic Galois group $\mc{G}_F$ would take precise form granted the standard conjectures, and in the meantime can be approximated by certain unconditional substitutes\footnote{Namely, Deligne's theory of absolute Hodge cycles, or Andr\'{e}'s theory of motivated cycles; see \S \ref{motivatedcycles}.}. In either case, however, we can ask for an analogue of Tate's theorem and attempt to prove unconditional results. The automorphic analogue is simpler to formulate, amounting to studying the fibers of the functorial transfer ${}^L \tG \to {}^L G$ where $\tG$ and $G$ are connected reductive $F$-groups, $\tG$ formed from $G$ by extending the (possibly non-connected) center $Z_G$ to a central subgroup $\tZ$, i.e. $\tG= (G \times \tZ)/Z_G$, so that $\tZ/Z_G$ is a torus. \index{s}{$G$}\index{s}{$\tG$}For simplicity, we always take the $F$-groups $G$ (and $\tZ$) to be split, but the results should also hold in the quasi-split case.\footnote{To be specific, the precise description of the fibers should be a feature of the quasi-split case; see Conjecture \ref{locmultone} and the subsequent discussion. The fact that the fibers are merely non-empty should be completely general.}  We show (Proposition \ref{llgeneralization}) that cuspidal representations of $G(\af)$ do indeed lift to $\tG(\af)$. The analogue is false for $\mc{G}_F$, and \textit{a priori} there is really no reason to believe that `lifting problems' for the three groups $\gal{F}, \mc{L}_F$, and $\mc{G}_F$ should qualitatively admit the same answer: a general automorphic representation has seemingly no connection with either $\ell$-adic representations or motives; and a general $\ell$-adic representation has no apparent connection with classical automorphic representations or motives. Two quite distinct kinds of transcendence-- one complex, one $\ell$-adic-- prevent these overlapping theories from being reduced to one another. A key problem, therefore, is to identify the overlap and ask how the lifting problem behaves when restricted to this (at least conjectural) common ground; it is only in posing this refined form of the lifting problem that the relevant structures emerge.   

So we take a detour to discuss notions of algebraicity for automorphic representations. Weil laid the foundation for this discussion in his paper \cite{weil:characters}, by showing that for Hecke characters (automorphic representations of $\mr{GL}_1(\af)$), an integrality condition on the archimedean component suffices to imply algebraicity of the coefficients of the corresponding $L$-series (i.e. algebraicity of its Satake parameters). Waldschmidt (\cite{waldschmidt:typeA}) later proved necessity of this condition. Weil, Serre, and others also showed that a subset, the `type $A_0$' Hecke characters (Definition \ref{AA0def}), moreover give rise to compatible systems of Galois characters $\gal{F} \to \Qlb^\times$ (and, via the theory of CM abelian varieties, motives underlying these compatible systems). The general feeling since has been that the most obvious analogue of Weil's type $A_0$ condition should govern the existence of associated $\ell$-adic representations. \index{t}{type $A_0$ Hecke character}To be precise, let $\pi$ be an automorphic representation of our $F$-group $G$, and let $T$ be a maximal torus of $G$. We will use the terminology `$L$-\textit{algebraic}' of \cite{buzzard-gee:alg} as the general analogue of type $A_0$ characters. That is, fixing at each $v \vert \infty$ an isomorphism $\iota_v \colon \overline{F}_v \xrightarrow{\sim} \CC$, we can write (in Langlands' normalization of \cite{langlands:archllc}) the restriction to $W_{\overline{F}_v}$ of its $L$-parameter as
\[
\rec_{v}(\pi_v) \colon z \mapsto \iota_v(z)^{\mu_{\iota_v}} \bar{\iota}_v(z)^{\nu_{\iota_v}} \in T^\vee(\CC).
\] 
with $\mu_{\iota_v}, \nu_{\iota_v} \in X^\bullet(T)_{\CC}$ and $\mu_{\iota_v}-\nu_{\iota_v} \in X^\bullet(T)$ (here and throughout, $\bar{\iota}_v$ denotes the complex conjugate of $\iota_v$).\index{s}{$\mu_{\iota_v}, \nu_{\iota_v}$} Unless there is risk of confusion, we will omit reference to the embedding $\iota_v$, writing $\mu_v= \mu_{\iota_v}$, etc.\index{s}{$\mu_v, \nu_v$} 
\begin{defn}\label{Lalgdefintro}
$\pi$ is $L$-algebraic if for all $v \vert \infty$, $\mu_v$ and $\nu_v$ lie in $X^\bullet(T)$. \index{t}{L-algebraic automorphic representation}
\end{defn}
The na\"{i}ve reason for focusing on this condition is that, for $G= \mr{GL}_n$, it lets us see the Hodge numbers of the (hoped-for) corresponding motive. Clozel (\cite{clozel:alg}), noticing that cohomological representations need not satisfy this condition, but have parameters $\mu_v, \nu_v \in \rho + X^\bullet(T)$, where $\rho$ is the half-sum of the positive roots,\footnote{For some choice of Borel containing $T$. Such a choice was implicit in defining a dual group, with its Borel $B^\vee$ containing a maximal torus $T^\vee$, and the above archimedean $L$-parameters.} studied this alternative integrality condition for $G= \mr{GL}_n$. Buzzard and Gee (\cite{buzzard-gee:alg}) have recently elaborated on this latter condition (`\textit{C-algebraic}') and its relation with $L$-algebraicity. \index{t}{C-algebraic automorphic representation}For $\mr{GL}_1$, the two notions $C$ and $L$-algebraic coincide, and in general they are both useful and distinct generalizations of Weil's type $A_0$ condition.

\index{t}{type $A$ Hecke character}But Weil's `integrality condition,' which he calls `type $A$' (Definition \ref{AA0def}), is more general than the type $A_0$ condition, and its analogues in higher rank seem to have been largely neglected.\footnote{One substantial use of type $A$ but not $A_0$ Hecke characters since the early work of Weil and Shimura occurs in the paper \cite{blasius-rogawski:motiveshmfs}, in which Blasius-Rogawski associate motives to certain tensor products of Hilbert modular representations.} It is of independent interest to resuscitate this condition, but we moreover find it essential for understanding the interaction of descent problems and algebraicity, including our original lifting question. By this we mean the following problem: given connected reductive (say quasi-split) $F$-groups $H$ and $G$ with a morphism of $L$-groups ${}^L H \to {}^L G$, and given a cuspidal $L$-algebraic representation $\pi$ of $G(\af)$ which is known to be in the image of the associated functorial transfer, is $\pi$ the transfer of an $L$-algebraic representation? The answer is certainly `no,' even for $\mr{GL}_1$, but it fails to be `yes' in a tightly constrained way. I believe the most useful general notion is the following:\index{t}{$W$-algebraic automorphic representation}
\begin{defn} We say that $\pi$ is $W$-\textit{algebraic} if for all $v \vert \infty$, $\mu_v$ and $\nu_v$ lie in $\frac{1}{2}X^\bullet(T)$. 
\end{defn}
This is more restrictive than Weil's type $A$ condition (which allows twists $|\cdot|^r$ for $r \in \Q$ as well), and it is easy to concoct examples of automorphic representations with algebraic Satake parameters that are not even twists of $W$-algebraic representations. It is one of our guiding principles, suggested by the (archimedean) Ramanujan conjecture, that all such examples are degenerate, and up to unwinding these degeneracies, all representations with algebraic Satake parameters should be built up from $W$-algebraic pieces. Moreover, for $L$-homomorphisms ${}^L H \to {}^L G$ with kernel equal to a central torus, the distinction between $W$- and $L$-algebraicity underlies the obstructions to preserving $L$-algebraicity in the lift. Some initial evidence comes from the tensor product $\mr{GL}_n \times \mr{GL}_n' \to \mr{GL}_{nn'}$ (Proposition \ref{tensordescent}). Specializing to the context of Hilbert modular forms, we then discuss $W$-algebraicity and the $\mr{GL}_2 \times \mr{GL}_2 \xrightarrow{\boxtimes} \mr{GL}_4$ transfer (its fibers and algebraicity properties) in detail, moreover linking it to the Galois side.\footnote{Various interesting examples arise, such as the construction of $\ell$-adic Galois representations which are pure but have no geometric twists.} Focusing finally on a transfer ${}^L \tG \to {}^L G$ as in Tate's lifting theorem, we prove the following algebraic refinement of the already-mentioned lifting result:
\begin{prop}\label{cmliftintro}
Let $F$ be CM field, and for simplicity let $G$ be a split semi-simple $F$-group, and let $\tZ$ be a torus. Let $\pi$ be a cuspidal representation of $G(\af)$. Assume $\pi_\infty$ is tempered.
\begin{enumerate}
\item If $\pi$ is $L$-algebraic, then there exists an $L$-algebraic lift to a cuspidal automorphic representation $\tpi$ of $\tG(\af)$.
\item If $\pi$ is $W$-algebraic, then there exists a $W$-algebraic lift $\tpi$.
\end{enumerate}
\end{prop}
In \S \ref{cmdescentsection} we develop a conjectural framework that allows this result to be extended to all totally imaginary $F$. See page \pageref{descentintro} of this introduction. In contrast, over totally real fields, we have the following:
\begin{prop}\label{totrealliftintro}
Now suppose either $F$ is totally real, with $\pi$ as before, or $F$ is arbitrary, but that the central character $\omega_\pi$ admits a finite-order extension to a Hecke character of $\tZ(\af)$. Continue to assume $\pi_\infty$ is tempered.
\begin{enumerate}
\item If $\pi$ is $L$-algebraic, then it admits a $W$-algebraic lift $\tpi$.
\item Assume $F$ is totally real, and for the `only if' direction of the following statement assume Conjecture \ref{locmultoneintro} below. Then the images of $\mu_v$ and $\nu_v$ under $X^\bullet(T) \to X^\bullet(Z_G)$ lie in $X^\bullet(Z_G)[2]$, and $\pi$ admits an $L$-algebraic lift if and only if these images are independent of $v \vert \infty$. 
\end{enumerate}
\end{prop}
To give a few examples, when $F$ is totally real there is no obstruction to finding $L$-algebraic lifts for $G$ a simple split group of type $A_{2n}, E_6, E_8, F_4$ or $G_2$. The natural question to ask, then, is whether for other groups there actually are counter-examples. Corollary \ref{counterexamples} applies limit multiplicity formulas (as in \cite{clozel:limmult}) to produce many discrete-series examples generalizing the basic case of mixed-parity Hilbert modular forms.\footnote{For the relationship of these results with Arthur's conjectural construction, which requires modification, of a morphism $\mc{L}_F \to \mc{G}_F$, see Remark \ref{arthurgripe2}.} The Conjecture \ref{locmultoneintro} mentioned here is a problem in local representation theory that plays a key role in the comparison of $L$-packets (both local and global) on $G$ and $\tG$; we require it to deduce a description of the fibers of the transfer ${}^L \tG \to {}^L G$.  In its most general form, Adler-Prasad conjectured (\cite{adler-prasad}) the following for $G$ and $\tG$ quasi-split:\index{t}{Adler-Prasad conjecture}
\begin{conj}[\cite{adler-prasad} Conjecture $2.6$]\label{locmultoneintro}
Let $v$ be a place of $F$, and let $G$ and $\tG$ be quasi-split groups, related as above. Then for any irreducible admissible (smooth for finite $v$, Harish-Chandra module for infinite $v$) representation $\tpi_v$ of $\tG(F_v)$, the restriction $\tpi_v|_{G(F_v)}$ decomposes with multiplicity one.
\end{conj}
The motivation for this conjecture is the uniqueness of Whittaker models for quasi-split groups. It is now known (by Adler-Prasad and others) for pairs $(\tG, G)= (\mr{GL}_n, \mr{SL}_n)$ or $(GU(V), U(V))$ where $V$ is a symplectic or orthogonal space over $F_v$, as well as for generic $\tpi_v$.\footnote{For archimedean $v$, it is obvious in many cases, and in this context it should be accessible in general.} We also note that the assumptions on temperedness at infinity are not so serious as might appear: for a cuspidal automorphic representation of $\mr{GL}_n(\af)$, $W$-algebraic implies (unconditionally) tempered at infinity. Certainly some results for non-tempered $\pi$ are also possible, especially now that Arthur (\cite{arthur:classical}) has proven his conjectures for classical groups. 

The Galois question parallel to these refined automorphic lifting results has been raised by Brian Conrad in \cite{conrad:dualGW}. In that paper, Conrad addresses lifting problems of the form
\[
\xymatrix{
& H'(\Qlb) \ar[d] \\
\gal{F} \ar@{-->}[ru]^{\tilde{\rho}} \ar[r]^{\rho} & H(\Qlb) \\
},
\]
where $H' \onto H$ is a surjection of linear algebraic groups with central kernel. He discusses existence (a local-global principle), ramification control, and $\ell$-adic Hodge theory properties, and the results are comprehensive, except for one question (see Remark $1.6$ and Example $6.8$ of \cite{conrad:dualGW}):
\begin{question}\label{conradquestion}\index{t}{Galois lifting problem}
Suppose that the kernel of $H' \onto H$ is a torus. If $\rho$ is geometric, when does there exist a geometric lift $\tilde{\rho}$?
\end{question}
This, in the case $H= G^\vee$, $H'= \tG^\vee$, is the natural Galois analogue of Propositions \ref{cmliftintro} and \ref{totrealliftintro}; indeed, it provided much of the motivation to understand those automorphic questions. By `geometric,' we mean almost everywhere unramified and de Rham at places above $\ell$.\footnote{As part of a general definition, we should also include reductive algebraic monodromy group-- the Zariski closure of the image-- but that is not essential for this question.} 
Conrad's Example $6.8$ shows the answer is at least `not always;' he produces a character $\hat{\psi} \colon \gal{L} \to \Qlb^\times$ over certain CM fields $L$ (with $F$ the totally real subfield) such that $\Ind_{L}^F (\hat{\psi})$ reduces to a geometric projective representation that has no geometric lift. We have an essentially general, but unfortunately conditional, solution to this problem. Having assumed that certain ($\mr{GL}$-valued) Galois representations are either automorphic or motivic, or merely satisfy certain Hodge-Tate weight symmetries (Hypothesis \ref{HTsym}) that are consequences of these assumptions, we prove (Theorem \ref{anymonodromylift} and remarks following):
\begin{thm}\label{anymonodromyliftintro}
Let $F$ be totally imaginary, and let $\rho \colon \gal{F} \to H(\Qlb)$ be a geometric representation satisfying Hypothesis \ref{HTsym}. Then $\rho$ admits a geometric lift $\tilde{\rho} \colon \gal{F} \to H'(\Qlb)$.
\end{thm}
There is an analogue in the totally real case (Corollary \ref{fullmonodromytotreal} and Remark \ref{liftsupp}), which is an exact parallel of Proposition \ref{totrealliftintro}. It is somewhat surprising that the Galois theory is not more complicated than the automorphic theory: roughly speaking, the automorphic input of temperedness is missing, and the obstruction to finding a geometric lift, say for the simple yet decisive case $\mr{GL}_n \to \mr{PGL}_n$, seems to be a question of $\frac{1}{n}$-integrality rather than $\frac{1}{2}$-integrality; this first appearance is, however, a red herring, and $W$-algebraicity remains the condition of basic importance on the Galois side as well. 

We also include a local version of this result; it uses the same argument but requires no additional assumption on `Hodge-Tate symmetry.' For $K/\Q_{\ell}$ finite, we can ask the same sorts of lifting questions with $\gal{K}$ in place of $\gal{F}$. A theorem of Wintenberger (\cite{wintenberger:relevement}), in the case of $H' \onto H$ a central isogeny, generalized by Conrad to $H' \onto H$ with kernel of multiplicative type, asserts that for any $\ell$-adic Hodge theory property $\mathbf{P}$ (i.e. Hodge-Tate, crystalline, semi-stable, or de Rham\index{t}{$\ell$-adic Hodge theory property $\mathbf{P}$}) a $\rho$ of type $\mathbf{P}$ admits a type $\mathbf{P}$ lift if and only if $\rho$ restricted to the inertia group $I_K$ admits a Hodge-Tate lift. This need not hold in the isogeny case, but we can complete this story by showing it holds unconditionally for central torus quotients: see Corollary \ref{HTlift}. This question too was suggested by Conrad, and a more difficult variant will appear in \S $3.6$ of \cite{chai-conrad-oort:cm}.

Having developed the automorphic and Galois sides in parallel, and having found the same answer to the two algebraic variants of Tate's lifting problem, we can use known results on the existence of automorphic Galois representations to construct, unconditionally, infinite families of Galois representations over totally real fields that fail to lift: \S \ref{spineg} treats $\mr{GSpin}_{2n+1} \to \mr{SO}_{2n+1}$ in the regular case (where on the automorphic side we try to lift $L$-algebraic, discrete series at infinity, representations $\pi$ of $\mr{Sp}_{2n}(\af)$ to $\mr{GSp}_{2n}(\af)$). In other cases, we can use the conjectural Fontaine-Mazur-Langlands correspondence for $\mr{GL}_n$-- we specify what we mean by this in Conjecture \ref{FML} below-- to develop parallel results for other groups. Here is an example from \S \ref{descentcomparison}:
\begin{prop}\label{BGegintro}
Let $F$ be any number field. Assume the Fontaine-Mazur-Langlands conjecture (\ref{FML}) and Conjecture \ref{cmdescent}. Then any geometric $\rho \colon \gal{F} \to \mr{PGL}_n(\Qlb)$ is automorphic, in the weak sense that there exists an $L$-algebraic automorphic representation $\pi$ of $\mr{SL}_n(\af)$ such that $\rho$ corresponds almost everywhere locally to $\pi$.\footnote{This side-steps the following subtlety, which is at the heart of the failure of multiplicity one for $\mr{SL}_n(\af)$: two projective representations of a finite group can be element-by-element conjugate but not globally conjugate. See Example \ref{locglobconj} and the discussion in \S \ref{automorphynotions}.}  
\end{prop}
For $F$ imaginary, we can invoke Theorem \ref{fullmonodromylift}. For $F$ totally real, a rather delicate descent argument is required. Similarly, we can deduce (see Proposition \ref{BGeg}) from Fontaine-Mazur-Langlands some cases of the converse, providing evidence for a conjecture of Buzzard and Gee (Conjecture $3.2.1$ of \cite{buzzard-gee:alg}). 

More generally, we are led to ask whether descent questions on the automorphic and Galois sides obey the same formalism: if an automorphic representation $\pi$ of $G(\af)$ is $L$-algebraic with associated Galois representation $\rho_\pi \colon \gal{F} \to {}^L G(\Qlb)$ (as conjectured in \cite{buzzard-gee:alg}, for instance), then we can ask two parallel questions for a given $L$-homomorphism ${}^L H \to {}^L G$:
\begin{itemize}
\item Is $\pi$ the functorial transfer of an automorphic representation of $H$?
\item Does $\rho_{\pi}$ lift through ${}^L H(\Qlb)$?
\end{itemize}
Closely related to these are the parallel questions:
\begin{itemize}
\item If $\pi$ descends to $H$, then does it have an $L$-algebraic descent?
\item If $\rho_{\pi}$ lifts through ${}^L H(\Qlb) \to {}^L G(\Qlb)$, then does it have a geometric factorization? 
\end{itemize}
We have already discussed some cases (arising from Tate's lifting question) of this general problem. Here is another example:
\begin{prop}
Continue to assume Fontaine-Mazur-Langlands and Conjecture \ref{cmdescent}. Let $\Pi$ be a cuspidal $L$-algebraic representation of $\mr{GL}_n(\af)$, and suppose that $\rho_\Pi \cong \rho_1 \otimes \rho_2$, where $\rho_i \colon \gal{F} \to \mr{GL}_{n_i}(\Qlb)$. Then there exist cuspidal automorphic representations $\pi_i$ of $\mr{GL}_{n_i}(\af)$ such that $\Pi= \pi_1 \boxtimes \pi_2$.
\end{prop}
The interest of this result is that examples generated from mixed-parity Hilbert modular forms show that the $\pi_i$ need not be $L$-algebraic; we observe in passing that in that context ($F$ totally real, $n_1=n_2=2$), certain cases of this result are, allowing a finite base-change, unconditional (via potential automorphy theorems). The structural point here is that the two Tannakian formalisms seem to be compatible in this tensor product example. In contrast:
\begin{prop}
Let $F$ be totally real. Let $\rho \colon \gal{F} \to \mr{PGL}_2(\Qlb)$ be geometric but have no geometric lift, and suppose $\Ad(\rho) \colon \gal{F} \to \mr{SO}_3(\Qlb)$ satisfies the (potential automorphy) hypotheses of Corollary $4.5.2$ of \cite{blggt:potaut}. Then there exists a lift $\tilde{\rho} \colon \gal{F} \to \mr{GL}_2(\Qlb)$ and, after a totally real base change $F'/F$, a mixed-parity Hilbert modular representation $\pi$ on $\mr{GL}_2/F'$, such that we have matching of local parameters
\[
\Sym^2(\tilde{\rho})|_{\gal{F'}} \sim_w \Sym^2(\pi) \otimes \chi
\]
for some finite order, necessarily non-trivial, Hecke character $\chi$ of $F'$. In contrast, restricting to $L'=F' K$ where $K/\Q$ is a quadratic imaginary extension in which $\ell$ is inert, we can find a lift $\tilde{\rho}$ such that $\Sym^2(\tilde{\rho})$ corresponds to $\BC_{L'/F'}(\Sym^2(\pi))$.

Conversely, starting with a mixed-parity $\pi$, we can produce such a $\chi$ and $\tilde{\rho}$.
\end{prop}
This proposition gives a precise example of the Fontaine-Mazur-Langlands correspondence breaking down past the $L$-algebraic/geometric boundary. It also quantifies a difference between descent problems (in this case, characterizing the image of $\Sym^2$) on the automorphic and Galois sides.\footnote{More simply, such examples exist for $n \colon \mr{GL}_1 \to \mr{GL}_1$.}
In all, these examples motivate a comparison of the images of $r \colon {}^L H \to {}^L G$ on the automorphic and Galois sides, when $r$ is an $L$-morphism with central kernel. The most optimistic expectation (for $H$ and $G$ quasi-split) is that if $\ker(r)$ is a central torus, then the two descent problems for ($L$-algebraic) $\Pi$ and (geometric) $\rho_\Pi$ are equivalent; whereas if $\ker(r)$ is disconnected, there is an obstruction to the comparison, that nevertheless can be killed after a finite base-change.

Despite these comparisons of the Galois and automorphic formalisms, the theory we develop here is deeply inadequate. Note that Proposition \ref{BGegintro} (as in the Buzzard-Gee conjecture) only makes claims about weak (almost everywhere local) equivalence of automorphic and Galois representations for groups other than $\mr{GL}_n$.\index{t}{weak equivalence of automorphic representations} The following two simple questions point toward just how much we are missing:
\begin{question}\label{modularlifting}
\begin{enumerate}
\item Suppose that $\rho \colon \gal{F} \to G^\vee(\Qlb)$ is weakly-equivalent to an automorphic representation $\pi$ of $G(\af)$. Do there exist weakly-equivalent lifts $\tilde{\rho} \colon \gal{F} \to \tG^\vee(\Qlb)$ and $\tpi$ on $\tG(\af)$?\footnote{We do not return to this question in these notes, but certain cases are manageable, and will be explained in a subsequent paper.}
\item When does a weakly compatible system of representations $\rho_{\ell} \colon \gal{F} \to G^\vee(\Qlb)$ lift to a weakly compatible system $\tilde{\rho}_{\ell} \colon \gal{F} \to \tG^\vee(\Qlb)$? (For the precise definition of `weakly compatible system,' see Definition \ref{wcsystem} below.)
\end{enumerate}
\end{question} 
The difficulties in answering the first question are reflected in the second, which in general has a negative answer; see Example \ref{locglobconj} for the prototypical counterexample, which is familiar from the study of automorphic multiplicities and endoscopy. The only way I know to address the second problem involves both assuming and proving much more: granting that the $\rho_{\ell}$ are $\ell$-adic avatars of a (pro-algebraic) representation of the motivic Galois group $\mc{G}_{F, E}$ of motives over $F$ with coefficients in a number field $E$ (implicit here are specified $E$-forms of $G^\vee$ and $\tG^\vee$), and showing that this representation lifts to $\tG^\vee$, possibly after enlarging $E$. We have therefore come full-circle to the motivic lifting-problem raised at the beginning of this introduction. This is far more difficult than the corresponding automorphic and Galois questions, but we can unconditionally treat certain examples.

In order to pose the question precisely, we need a category of `motives' in which the Galois formalism of $\mc{G}_F$ is unconditional. There are two common constructions, one based on Deligne's theory (\cite{DMOS}) of absolute Hodge cycles, the other based on Andr\'{e}'s theory (\cite{andre:motivated}), much inspired by Deligne's work, of motivated cycles. We work with motivated cycles, since the inclusions
\[
\text{algebraic cycles} \subseteq \text{motivated cycles} \subseteq \text{absolute Hodge cycles} \subseteq \text{Hodge cycles}
\]
more or less imply that results about motivated cycles (eg, `Hodge cycles are motivated') follow \textit{a fortiori} for absolute Hodge cycles. Note that the Hodge conjecture asserts that each $\subseteq$ is an equality; Deligne's `espoir' (\cite[0.10]{deligne:valeurs}) asserts that the last $\subseteq$ is an equality. We let $\mc{M}_{F, E}$ denote the category of motives for motivated cycles over $F$ with coefficients in $E$ (see \S \ref{motivatedcycles} and \ref{motiveswithcoefficients}). This is a semi-simple $E$-linear Tannakian category, and, choosing an embedding $F \into \CC$, we can associate (via the $E$-valued Betti fiber functor) a (pro-reductive) motivic Galois group $\mc{G}_{F, E}$. We will prove a lifting result of the form
\[
\xymatrix{
& \mr{GSpin}(V_E) \ar[d] \\
\mc{G}_{F, E} \ar@{-->}[ur]^{\tilde{\rho}} \ar[r]^{\rho} & \mr{SO}(V_E)
}
\]
for certain $\rho$ arising from degree $2$ primitive cohomology $V_{\Q}$ of a hyperk\"{a}hler variety (see Definition \ref{HK}) over $F$ (or, rather, for $\rho$ having analogous formal properties); the extension of scalars $V_E= V_{\Q} \otimes_{\Q} E$ is essential for the lifting result to hold. 

The starting-point for this refined motivic lifting result is work of Andr\'{e} (\cite{andre:hyperkaehler}) on the `motivated' theory of hyperk\"{a}hlers; indeed, his results imply a version of this lifting statement with $F$ replaced by a large but quantifiable (with considerable work: see Theorem $8.4.3$ of \cite{andre:hyperkaehler}) finite extension $F'/F$. Our contribution is the arithmetic descent from $F'$ to $F$, and for this we use a variant of a technique introduced by Ribet (\cite{ribet:Qcurves}) to study so-called `$\Q$-curves,' elliptic curves over $\Qb$ that are isogenous to all of their $\gal{\Q}$-conjugates. Passing from a softer geometric statement (compare: `a variety over $F$ has a point over some finite extension') to a precise arithmetic version (compare: `when does a variety over $F$ have a point over $F$?') typically requires some deep input; in this case, Faltings's isogeny theorem (\cite{faltings:endlichkeit}) does the hard work. The method requires a case-by-case analysis (depending on the motivic group of the transcendental lattice), and I have decided not to pursue all the cases here, but here is a partial result (see Theorem \ref{motiviclift}, Theorem \ref{nongenericoddmotiviclift}, and Proposition \ref{gettingtired} and following for more cases and more precise versions):
\begin{thm}\label{HKliftintro}
Let $(X, \eta)$ be a polarized variety over $F$ satisfying Andr\'{e}'s conditions $A_k, B_k, B_k^+$ (see \S \ref{HKsetup}). For instance, with $k=1$, $X$ can be a hyperk\"{a}hler variety with $b_2(X)>3$. Suppose that the transcendental lattice $T_{\Q} \subset V_{\Q}= Prim^{2k}(X_{\CC}, \Q)(k)$ satisfies either
\begin{itemize}
\item $\End_{\Q-Hodge}(T_{\Q})= \Q$; or
\item $T_{\Q}$ is odd-dimensional.
\end{itemize}
For simplicity, moreover assume that $\det V_{\ell}= \det T_{\ell}=1$ as $\gal{F}$-representations.\footnote{This can be achieved in each case after a quadratic, independent-of-$\ell$ extension on $F$. It is only so we can work with $\mr{SO}(V_{\ell})$ rather than $\mr{O}(V_{\ell})$, but since our abstract Galois-lifting results apply to non-connected groups (see Theorem \ref{anymonodromyliftintro}, for instance), this hypothesis is not essential.} Then after some finite scalar extension to $V_E= V_{\Q} \otimes_{\Q} E$ there is a lifting of the motivic Galois representation $\rho^V \colon \mc{G}_{F} \to \mr{SO}(V_{\Q})$:
\[
\xymatrix{
& \mr{GSpin}(V_E) \ar[d] \\
\mc{G}_{F, E} \ar[ur]^{\tilde{\rho}} \ar[r]^{\rho^V} & \mr{SO}(V_E).
}
\]
Moreover, this lifting gives rise to a weakly compatible system of lifts $\tilde{\rho}_{\lambda} \colon \gal{F} \to \mr{GSpin}(V_{\lambda})$ on $\lambda$-adic realizations, for all finite places $\lambda$ of $E$.
\end{thm}
For more cases in which $T_{\Q}$ is even-dimensional, see \S \ref{Teven}; the omitted cases should yield to similar methods. In particular, we obtain in some cases a positive answer to the second part of Question \ref{modularlifting}; this compatibility of $\lambda$-adic realizations is \textit{not} automatic for Andr\'{e}'s motives (nor for absolute Hodge motives), so we need to exploit an explicit description of the lift. Moreover, the excluded case $b_2(X)=3$ of Theorem \ref{HKliftintro} is related to a more general result for potentially abelian motives: we prove such a lifting result across an arbitrary surjection $H' \to H$ with central torus kernel: see Proposition \ref{taniyamalift} and Lemma \ref{b2=3}. 

Our next example seems to be a novel result even in complex algebraic geometry, where many authors have studied variants of the Kuga-Satake construction for Hodge structures of `$K3$-type' (see, for example, \cite{morrison:KSabsurface}, \cite{voisin:kuga-satake}, \cite{galluzzi:KSabfourfold}). Generalizing the well-known case (\cite{morrison:KSabsurface}) of $H^2$ of an abelian surface-- which has a Hodge structure of $K3$ type, so falls within the ken of classical Kuga-Satake theory\footnote{The paper \cite{galluzzi:KSabfourfold} treats the case of a particular family of abelian four-folds, where a constraint on the Mumford-Tate group allows one to extract a Hodge structure of $K3$-type; as we show, though, the lifting phenomenon is completely general.}-- we prove an analogous $\mr{GSpin} \to \mr{SO}$ motivic lifting result for $H^2$ of any abelian variety (see Corollary \ref{AVKS}). Moreover, we can describe the corresponding motive (in the spin representation) as a Grothendieck motive, without assuming the Standard Conjectures.

Note that Serre has asked (\cite[8.3]{serre:motivicgalois}) whether homomorphisms $\mc{G}_{\overline{F}} \to \mr{PGL}_2$ lift to $\mr{GL}_2$. This is closely related to the question, also raised by Serre, of whether the derived group of $\mc{G}_{\overline{F}}$ is simply-connected, but Serre notes that questions of this sort do not have obvious conjectural answers, even if we assume we are in \textit{le paradis motivique}. Our abstract Galois lifting results, as well as the handful of motivic examples, suggest and provide evidence for the following sharpening and generalization of Serre's question, that this motivic lifting phenomenon is utterly general:
\begin{conj}\label{motivicliftconjecture}
Let $F$ and $E$ be number fields, and let $H' \onto H$ be a surjection, with central torus kernel, of linear algebraic groups over $E$. Suppose we are given a homomorphism $\rho \colon \mc{G}_{F, E} \to H$. Then if $F$ is imaginary, there is a finite extension $E'/E$ and a homomorphism $\tilde{\rho} \colon \mc{G}_{F, E'} \to H_{E'}$ lifting $\rho \otimes_E E'$. If $F$ is totally real, then such a lift exists if and only if the Hodge number parity obstruction of Corollary \ref{fullmonodromytotreal} vanishes.
\end{conj}
This is consistent with the Fontaine-Mazur conjecture, Serre's expectation for $\mc{G}_{\overline{F}}$, and Theorem \ref{anymonodromyliftintro}, and indeed the three of them taken together imply this conjecture (mimic the proof of Proposition \ref{taniyamalift}).  As we have seen, the cases of potentially abelian motives, hyperk\"{a}hler motives, and abelian varieties provide some evidence both for this conjecture and for Serre's original question. In current work-in-progress, we study examples of this motivic lifting phenomenon for (unlike the examples in this paper) motives not lying in the Tannakian category generated by abelian varieties and Artin motives. There is clearly vast terrain waiting to be discovered here, some of which is not entirely hostile to exploration. 

A few other general themes recur throughout these notes, in automorphic, Galois-theoretic, and motivic variants. The first is the systematic exploitation of `coefficients,' and the general principle that our arithmetic objects will naturally have coefficients in CM fields.\index{t}{CM coefficients} In \S \ref{cmdescentsection}, we use this principle to formulate (and prove, for regular representations; see Proposition \ref{cmdescent}) a conjectural generalization of Weil's result that a type $A$ Hecke character of a number field $F$ descends, up to finite order twist, to the maximal CM subfield $F_{cm}$; this generalization asserts that the \textit{infinity-type} of an $L$-algebraic cuspidal automorphic representation of $\mr{GL}_n(\af)$, over a totally imaginary field $F$, necessarily descends to $F_{cm}$. \label{descentintro} In \S \ref{hodgesymmetry}, we use this principle, which in the motivic context is more or less the Hodge index theorem, to establish for Andr\'{e}'s motives the Hodge-Tate weight symmetries needed for the Galois lifting Theorem \ref{anymonodromyliftintro}. Finally, Corollary \ref{regularoverCMfields} gives an example of how having CM coefficients can be exploited even in the study of abstract compatible systems of $\ell$-adic representations.\footnote{Since the writing of this paper, the same principles have been applied in \cite{patrikis-taylor:irr} to establish new results on the potential automorphy of regular motives, and the irreducibility of automorphic Galois representations.}

This last result is based on the following abstract independence-of-$\ell$ result:\footnote{See Proposition \ref{tracezero} for a more general statement.}
\begin{prop}\label{introtracezero}
Let $\rho_{\lambda} \colon \gal{F} \to \mr{GL}_n(\Elb)$ be a compatible system of irreducible, Lie-multiplicity-free representations of $\gal{F}$ with coefficients in a number field $E$, so that $\rho_{\lambda} \cong \Ind_{L^{\lambda}}^F(r_{\lambda})$ for some Lie-irreducible representation $r_{\lambda}$ of $\gal{L^{\lambda}}$, where the number field $L^{\lambda}$ \textit{a priori} depends on $\lambda$. Then the set of places of $F$ having a split factor in $L^{\lambda}$ is independent of $\lambda$. If we further assume that the $L^{\lambda}/F$ are Galois, then $L^{\lambda}$ is independent of $\lambda$.
\end{prop}
By Lie-multiplicity-free, we mean multiplicity-free after all finite restrictions. \index{t}{Lie-multiplictiy free, or LMF}One application is a converse to a theorem of Rajan (Theorem $4$ of \cite{rajan:sm1}), also generalizing a result of Serre (Corollaire $2$ to Proposition $15$ of \cite{serre:cebotarev}). See \S \ref{monodromy} for further discussion, as well as the aformentioned application (Corollary \ref{regularoverCMfields}), which is a weak Galois-theoretic shadow of the automorphic Proposition \ref{cmdescent}. In \S \ref{monodromy}, we also record a variant for number fields of a result of Katz (for affine curves over a finite field; see \cite{katz:duke}), which in particular clarifies the place in the general theory occupied by Lie-multiplicity-free (or more specifically, Lie-irreducible) representations:
\begin{prop}
Let $\rho \colon \gal{F} \to \mr{GL}_n(\Qlb)$ be an irreducible representation. Then either $\rho$ is induced, or there exists $d \vert n$, a Lie irreducible representation $\tau$ of dimension $n/d$, and an Artin representation $\omega$ of dimension $d$ such that $\rho \cong \tau \otimes \omega$. This decomposition is unique up to twisting by a finite-order character. Consequently, any (irreducible) $\rho$ can be written in the form
\[
\rho \cong \Ind_{L}^F (\tau \otimes \omega)
\]  
for some finite $L/F$ and irreducible representations $\tau$ and $\omega$ of $\gal{L}$, with $\tau$ Lie-irreducible and $\omega$ Artin.
\end{prop}
This is a very handy result, which despite its basic nature does not seem to be widely-known. Some of our descent arguments in \S \ref{hyperlift} rely on it. 

Finally, time and again our arguments are buttressed by some explicit knowledge of the constraints on infinity-types of automorphic representations, and on Hodge-Tate weights of Galois representations. More generally, we find that asking finer structural questions about the interaction of functoriality and algebraicity naturally leads to existence (and non-existence) problems that are often dissociated from functoriality.\footnote{As we will see in \S \ref{review}, this theme is present even in the proof of Tate's original vanishing theorem $H^2(\gal{F}, \Q/\Z)=0$.} We certainly have many more questions than answers (scattered through \S's \ref{GL1}-\ref{GaloisHMFs}), although fortunately our main results depend largely on detailed study of $\mr{GL}_1$ (\S \ref{GL1}), which provides both the technical ingredients for later arguments and the motivation for higher-rank results (namely, after a close reading of Weil's \cite{weil:characters}, the definition of $W$-algebraicity and the key descent principle of Proposition \ref{cmdescent}). In this latter respect, I do not believe that $\mr{GL}_1$ has yielded all of its fruit-- see the discussion surrounding Corollary \ref{HCstructurethm}-- but the extreme difficulty of establishing the (non-)existence of automorphic representations with given infinity-types, even in qualitative form, necessarily tempers further conjecture.

For a reader interested only in certain aspects of these notes, I hope that the table of contents is a clear reference to the points of interest. 
\section{What is assumed of the reader: background references}
Because this monograph studies all three vertices-- and, to a far lesser extent, the edges-- of the mysterious Galois-automorphic-motivic triangle, the reader will need some limited familiarity with all three subjects. In this section, we indicate what background will be assumed, give some useful references, and also explain that while the background required may be broad, it is not terribly deep; indeed, the arguments of this paper are elementary, once some basic definitions are assimilated. Let us deal individually with each of these three topics.
\subsection{Automorphic representations}
Preliminary to the study of automorphic representations, one must have some familiarity with the theory of reductive algebraic groups. For simplicity, we always restrict to split groups, so it is sufficient to know the theory over algebraically closed fields. Many arguments in this paper rest on the manipulation of root data; Springer's survey \cite{springer:corvallis} and his book \cite{springer:lag} (especially chapters 7-10) are very clear, and provide more than enough background. We also assume familiarity with the algebraic representation theory of (connected) reductive groups (i.e., elementary highest-weight theory).

We use a little of the representation theory of reductive groups over local fields. Familiarity with unramified representations of split $p$-adic groups and the Satake correspondence (\cite{gross:satake} is a beautiful guide; \cite[III]{cartier:p-adicrep} is a thorough treatment of unramified representations; \cite{casselman:notes} is the best general introduction to the theory of admissible smooth representations), and with some aspects of the formalism of the archimedean local Langlands correspondence will suffice. As for the global theory of automorphic representations themselves, and the theory of the L-group, the standard and more than adequate reference is \cite{borel:L}; \cite[\S 9-11]{borel:L} would be particularly useful background reading, sketching the archimedean theory and, crucially for our purposes, explaining the desired, in some cases proven, connection between central characters and L-parameters.
\subsection{$\ell$-adic Galois representations}\label{ladicrefs}
Apart from some very elementary notions, we only require some familiarity with the formal, non-geometric, aspects of $\ell$-adic Hodge theory, which is to say the study of $\ell$-adic representations of $\gal{K}$ for a finite extension $K/\Q_{\ell}$ (when studied in a purely local context, $\ell$ is traditionally replaced by $p$, as in '$p$-adic Hodge theory'). An excellent overview, with references to detailed proofs, is \cite[\S I-II]{berger:introp-adicrep}. A `text-book-style' reference (very useful but not quite in final form) for everything we do is \cite{brinon-conrad:cmi}.\footnote{It is also highly recommended to visit Laurent Berger's website to check on the status of the course-notes from his course at IHP in the Galois Trimestre of 2010.} Finally, we should remark to a reader new to $p$-adic Hodge theory that the subject has been much in flux in the last 10 years, and that one should take heed of recent conceptual advances (eg \cite{beilinson:derivedderham}, \cite{bhargav:p-adicderdR}, \cite{fargues-fontaine:curve}, \cite{scholze:p-adichodgetheory}) before reading \textit{too} deeply into the `classical' theory. 

\subsection{Motives}
Here we recommend some familiarity with the notion of motive, as conceived by Grothendieck, and with Grothendieck's Standard Conjectures on algebraic cycles.\index{t}{Standard Conjectures} Kleiman's articles \cite{kleiman:algcycles} and \cite{kleiman:standardconjectures} provide a clear, careful, and concise introduction. Certainly familiarity with the various cohomological realizations (Betti, de Rham, $\ell$-adic), and the relations between them, of a smooth projective variety is necessary; the first section, `Review of Cohomology,' of Deligne's article \cite[\S I.1]{DMOS} will provide quick and easy orientation for someone new to the subject. We will also require (some of) the theory of Tannakian categories; we will give a brief introduction that should be enough for a reader to follow all of our arguments, but a thorough treatment is \cite{deligne-milne}. Finally, although not necessary for the present work, the excellent book \cite{andre:introduction} provides a broad survey of the theory of motives, from its inception to more recent developments.

\subsection{Connecting the dots}
These notes rely on the systematic transfer of intuition back and forth between these three areas; it may be helpful, then, to collect in precise form the towering conjectures that dominate this conceptual landscape. The following three-part conjecture, and the consequence of combining parts \ref{FMT} and \ref{FMLsub}, summarize the principal problems in the field.\index{t}{geometric Galois representation}
\begin{conj}[Fontaine-Mazur-Langlands-Tate]\label{FML}
Let $\rho \colon \gal{F} \to \mr{GL}_n(\Qlb)$ be an irreducible geometric Galois representation. Recall from the introduction that this means $\rho$ is almost everywhere unramified, and is de Rham at all places above $\ell$. Then:
\begin{enumerate}
\item (Fontaine-Mazur)\index{t}{Fontaine-Mazur conjecture} There exists a smooth projective variety $X/F$ and an integer $r$ such that $\rho$ is isomorphic to some sub-quotient of $H^j(X_{\overline{F}}, \Qlb)(r)$. 
\item\label{FMT} (Fontaine-Mazur-Tate)\index{t}{Fontaine-Mazur-Tate conjecture} Such a $\rho$ is \textit{motivic} in the sense that it is cut out by $\Qlb$-linear combinations of (homological) algebraic cycles. More precisely, for any embedding $\iota \colon \Qb \into \Qlb$, there is an idempotent project $e$ in the algebra $C^0_{hom}(X, X)_{\Qb}$ (see \cite[\S 1.3.8]{kleiman:algcycles}, where this is denoted $\mc{A}(X)$) of self-correspondences of $X$ with $\Qb$-coefficients, such that
\[
\rho \cong e (H^*(X)(r)_{\Qb}) \otimes_{\Qb, \iota} \Qlb
\]
as $\gal{F}$-representations.
\item\label{FMLsub} (Fontaine-Mazur-Langlands)\index{t}{Fontaine-Mazur-Langlands conjecture} There exists a cuspidal automorphic representation $\Pi$ of $\mr{GL}_n(\af)$ such that for almost all $v$ unramified for $\rho$, the eigenvalues of $\rho(fr_v)$ agree with the Satake parameters of $\Pi_v$, viewed in $\Qlb$ via the composition $\iota_{\ell} \circ \iota_{\infty}^{-1}(\rec_v(\Pi_v))$.\footnote{A stronger version moreover asks that for all finite places $v$,
\[
\mr{WD}(\rho|_{\gal{F_v}})^{fr-ss} \cong \iota_{\ell} \circ \iota_{\infty}^{-1}(\rec_v(\Pi_v))
\]} Moreover for $\iota \colon F \into \CC$, the archimedean Langlands parameter 
\[
\rec_v(\Pi_v)|_{\overline{F_v}^\times} \colon W_{\overline{F_v}} \to \mr{GL}_n(\CC),
\]
landing in a maximal torus $T_n$, is of the form $z \mapsto z^{\mu_{\iota}} \bar{z}^{\nu_{\iota}}$ for $\mu_{\iota}, \nu_{\iota} \in X_{\bullet}(T_n)$, and for any $\tau \colon F \subset F_v \into \Qlb$, the Hodge-Tate co-character $\mu_{\tau}$ of $\rho|_{\gal{F_v}}$ is (conjugate to) $\mu_{\iota^*_{\infty, \ell}(\tau)}$.\footnote{$\iota^*_{\infty, \ell}(\tau)$ is the pullback of $\tau$ to an archimedean embedding via $\iota_{\ell}$, $\iota_{\infty}$; see \S \ref{notation}.}

Conversely, given a cuspidal automorphic representation $\Pi$ of $\mr{GL}_n(\af)$ with integral archimedean Langlands parameters, there exists an irreducible geometric Galois representation $\rho_{\Pi} \colon \gal{F} \to \mr{GL}_n(\Qlb)$ (depending on $\iota_{\infty}, \iota_{\ell}$) satisfying the above compatibilities.
\end{enumerate}
\end{conj}
\begin{rmk}
\begin{enumerate}
\item Note that the Standard Conjectures are implicit in part \ref{FMT} of Conjecture \ref{FML}. We will rarely work directly with motives for homological equivalence in these notes, substituting instead Andr\'{e}'s category of motivated motives. Keeping one's faith in the Standard Conjectures, but wanting to assume less at the outset, one could substitute the algebra of motivated correspondences $C^0_{mot}(X, X)$ (Definition \ref{motivatedcycle} and following) for $C^0_{hom}(X, X)$ in Conjecture \ref{FML}. 
\item This is not a historically faithful presentation of these conjectures; for our purposes in these notes, however, this formulation is convenient.
\item Part \ref{FMLsub} implies a similar correspondence between semi-simple (not necessarily irreducible) geometric Galois representations and (suitably algebraic) isobaric automorphic representations of $\mr{GL}_n(\af)$. Alternatively, one can begin from this more general conjecture and deduce that cuspidal automorphic representations must correspond to irreducible Galois representations using standard properties of Rankin-Selberg L-functions.
\end{enumerate}
\end{rmk}
The Fontaine-Mazur-Tate conjecture challenges us, given a geometric $\ell$-adic representation $\rho_{\ell}$, to produce a motive with $\rho_{\ell}$ as $\ell$-adic realization, and in particular to produce a family of $\ell'$-adic realizations (for all $\ell'$) that are `compatible' with $\rho_{\ell}$. It is often convenient to abstract this notion of compatible system of Galois representations, as an intermediary between the isolated $\rho_{\ell}$ and the robust motive. 
\begin{defn}\label{wcsystem}\index{t}{weakly compatible system of $\lambda$-adic, or $\ell$-adic, representations}
Let $F$ and $E$ be number fields, and let $N$ be a positive integer. A rank $N$ \textit{weakly compatible system of} $\lambda$-adic (or, informally, $\ell$-adic) \textit{representations of} $\gal{F}$ \textit{with coefficients in} $E$ is a collection 
\[
\mc{R}= \left( \{\rho_{\lambda}\}_{\lambda}, S, \{Q_v(X)\}_{v \not \in S} \right),
\]
consisting of:
\begin{enumerate}
\item for each finite place $\lambda$ of $E$, a continuous semi-simple geometric representation 
\[
\rho_{\lambda} \colon \gal{F} \to \mr{GL}_N(\overline{E_{\lambda}});
\]
\item a finite set of places $S$ of $F$, containing the infinite places, such that for all $v \not \in S$ and for all $\lambda$ of different residue characteristic from $v$, $\rho_{\lambda}|_{\gal{F_v}}$ is unramified;
\item for all $v \not \in S$, a polynomial $Q_v(X) \in E[X]$ such that for all $\lambda$ of different residue characteristic from $v$, $Q_v(X)$ is the characteristic polynomial of $\rho_{\lambda}(fr_v)$.
\end{enumerate}
We will sometimes use a similar notion where $\mr{GL}_N$ is replaced by a, for simplicity, split connected reductive group $H$ over the number field of coefficients $E$. In this case, repeat the above definition verbatim, except replace the characteristic polynomial of $\rho_{\lambda}(fr_v)$ with the analogous point of the space of Weyl-invariant functions on a maximal torus in $H$, i.e. the image of $\rho_{\lambda}(fr_v) \in H(\overline{E_{\lambda}})$ under the Chevalley map induced by the map on coordinate rings
\[
E[H] \supset E[H]^H \xrightarrow[\mr{res}]{\sim}E[T]^W,
\]
where $T$ is an $E$-split maximal torus, and $W$ is the Weyl group of $(H, T)$.\index{t}{weakly compatible system of $\lambda$-adic representations valued in a reductive group}
\end{defn}
We will sometimes impose a \textit{purity} hypothesis on $\mc{R}$:
\begin{defn}\label{purewcs}
Fix an integer $w$. For any finite place $v$ of $F$, let $q_v$ denote the order of the residue field of $F$ at $v$. We say that a weakly compatible system $\mc{R}$ is \textit{pure} of weight $w \in \Z$ if for a density one set of places $v$ of $F$, each root $\alpha$ of $Q_v(X)$ in $\overline{E}$ satisfies $|\iota(\alpha)|^2= q_v^w$ for all embeddings $\iota \colon \overline{E} \into \CC$.\index{t}{pure weakly compatible system}
\end{defn}
Another standard variant of the definition of weakly compatible system would also specify that the $\ell$-adic Hodge numbers are suitably `independent of $\ell$,' but we will not require this; see for instance \cite[\S 5.1]{blggt:potaut}. Finally, an elementary argument (see, eg, \cite[Lemma 2.1.5]{clozel-harris-taylor}) shows that any continuous $\overline{E_{\lambda}}$-representation of $\gal{F}$ (a compact group) takes values in some $\mr{GL}_N(E')$ for some finite extension $E'$ of $E_{\lambda}$.

\section{Notation}\label{notation}
For a number field $F$, we always choose an algebraic closure $\overline{F}/F$ and let $\gal{F}$ denote $\Gal(\overline{F}/F)$.\index{s}{$\gal{F}$} We write $C_F= \af^\times/F^\times$ for the idele class group of $F$.\index{s}{$C_F$} 

If $L/F$ is a finite extension inside $\overline{F}$ and $W$ a representation of $\gal{L}$ ($= \Gal(\overline{F}/L)$), we abbreviate $\Ind_{\gal{L}}^{\gal{F}}(W)$ to $\Ind_L^F(W)$. For $g \in \gal{F}$, we write $(gW)$ for the conjugate representation of $g \gal{L} g^{-1}$.

We fix separable closures $\Qb$ and $\Qlb$ of $\Q$ and $\Q_\ell$ for all $\ell$. We fix throughout embeddings $\iota_\ell \colon \Qb \into \Qlb$ and $\iota_{\infty} \colon \Qb \into \CC$.\index{s}{$\iota_{\ell}$}\index{s}{$\iota_{\infty}$} An archimedean embedding $\iota \colon F \into \CC$ thus induces an $\ell$-adic embedding $\tau_{\ell, \infty}^*(\iota)= \iota_\ell \circ \iota_{\infty}^{-1} \circ \iota \colon F \into \Qlb$; and similarly an $\ell$-adic embedding $\tau \colon F \into \Qlb$ induces $\iota_{\infty, \ell}^*(\tau) \colon F \into \CC$.\index{s}{$\tau_{\ell, \infty}^*(\iota)$, $\tau^*(\iota)$}\index{s}{$\iota_{\infty, \ell}^*(\tau)$, $\iota^*(\tau)$} If there is no risk of confusion, we just write $\tau^*(\iota)$ and $\iota^*(\tau)$. These embeddings will be invoked (often implicitly) whenever we associate automorphic forms and Galois representations.

For a connected reductive $F$-group $G$, we construct dual and $L$-groups $G^\vee$ and ${}^L G$ over $\Qb$\index{s}{$G^\vee$, ${}^L G$} (having chosen a maximal torus, Borel, and splitting, although we will only ever make the maximal torus explicit), and then use $\iota_\ell$ and $\iota_\infty$ to regard the dual group over $\Qlb$ or $\CC$ as needed. 

If $v$ is a place of $F$, we denote by $\rec_v$\index{s}{$\rec_v$}\index{t}{local reciprocity map} the local reciprocity map from irreducible admissible smooth representations ($v$ finite) or irreducible admissible Harish-Chandra modules ($v$ infinite) to (frobenius semi-simple) representations of the Weil(-Deligne) group of $F_v$. We use this in the unramified and archimedean cases, and for $\mr{GL}_n$. For the group $\mr{GL}_1$, this is local class field theory, normalized so that uniformizers correspond to geometric frobenii, which we denote $fr_v$.

It will be convenient to have a short-hand for the assertion that an automorphic representation and Galois representation `correspond' under the conjectural global Langlands correspondence. We will elaborate on a number of related notions in \S \ref{automorphynotions}, but for now we give two that will come up frequently. Recall that we have fixed embeddings $\iota_{\ell}$ and $\iota_{\infty}$ of $\Qb$ into $\Qlb$ and $\CC$, respectively. If $\pi$ is an automorphic representation of $G(\af)$, and $\rho \colon \gal{F} \to {}^L G(\Qlb)$ is a continuous, almost everywhere unramified, representation, then we write $\rho \sim_w \pi$ \index{s}{$\sim_w$}\index{t}{correspondence between automorphic representations and Galois representations} if for almost every unramified place $v$ of $F$ (at which we may assume $\rho$ is unramified), $\rec_v(\pi_v) \colon W_{F_v} \to {}^L G(\CC)$ can be realized (up to $G^\vee$-conjugation) over $\Qb \xrightarrow{\iota_{\infty}} \CC$, and that then
\[
\iota_{\ell} \circ \iota_{\infty}^{-1} \left( \rec_v(\pi_v) \right) \sim \left( \rho|_{\gal{F_v}} \right)^{\mr{fr-ss}},
\]
where $\sim$ here denotes $G^\vee(\Qlb)$-conjugacy, and `fr-ss' means we replace $\rho(fr_v)$ with its semi-simple part. More concretely, we restrict to places $v$ at which both $\pi$ and $\rho$ are unramified, so that $\rec_v(\pi_v)$ and $\rho|_{\gal{F_v}}$ are both just given by their evaluations at $fr_v$; after using $\iota_{\ell} \circ \iota_{\infty}^{-1}$ to regard them as defined over the same field, we ask that these two elements be $G^\vee(\Qlb)$-conjugate.

As in the Fontaine-Mazur-Langlands Conjecture \ref{FML}, we also often want to compare the archimedean component of $\pi$ with the Hodge-Tate weights of $\rho$; when $\pi$ is L-algebraic (see Definition \ref{Lalgdefintro}), and $\rho$ is de Rham, we may write, for all $\iota \colon F \into \CC$, the archimedean Langlands parameter $\rec_v(\pi_v)|_{\overline{F}_v^\times}$ in the form 
\[
z \mapsto \iota(z)^{\mu_{\iota}}\bar{\iota}(z)^{\nu_{\iota}}
\]
for $\mu_{\iota}, \nu_{\iota} \in X_{\bullet}(T^\vee)$. Then write $\rho \sim_{w, \infty} \pi$ \index{s}{$\sim_{w, \infty}$} if $\rho \sim_w \pi$, and moreover for all such $\iota \colon F \into \CC$, and for all $\tau \colon F \subset F_v \into \Qlb$, the $\tau$-labeled Hodge-Tate co-character (see \S \ref{labeledweights}) of $\rho|_{\gal{F_v}}$ is (conjugate to) $\mu_{\iota^*_{\infty, \ell}(\tau)}$. 

By a CM field we mean as usual a quadratic totally imaginary extension of a totally real field; these, and their real subfields, are the number fields on which complex conjugation is well-defined, independent of the choice of complex embedding.\index{t}{CM field} For any number field, we write $F_{cm}$\index{s}{$F_{cm}$} for the maximal subfield of $F$ on which complex conjugation is well-defined; thus it is the maximal CM or totally real subfield, depending on whether $F$ contains a CM subfield. We also write $\Q^{cm}$\index{s}{$\Q^{cm}$} for the union of all CM extensions of $\Q$ inside $\Qb$.

For any topological group $A$, we denote by $A^D$\index{s}{$(\cdot)^D$, with argument a topological group} the group $\Hom_{cts}(A, \mr{S}^1)$ ($\mr{S}^1$ is the unit circle); when $A$ is abelian and locally compact, we topologize this as the Pontryagin dual.

For any ground field $k$ (always characteristic zero for us) with a fixed separable closure $k^{s}$, and a (separable) algebraic extension $K/k$ inside $k^{s}$, we write $\widetilde{K}$\index{s}{$\widetilde{(\cdot)}$, with argument a field} for the normal (Galois) closure of $K$ over $k$ inside $k^{s}$.

We always denote base-changes of schemes by sub-scripts: thus, $X_{k'}$ is the base-change $X \times_{\Spec k} \Spec k'$ for a $k$-scheme $X$ and a base-extension $\Spec k' \to \Spec k$. Notation such as $M \otimes_{k} k'$ will be reserved (see \S \ref{motiveswithcoefficients}) to describe extension of scalars from an object $M$ of a $k$-linear abelian category to one of a $k'$-linear abelian category; when there is no risk of confusion we will sometimes denote this also by $M_{k'}$.

Finally, we will signal either a significant change in running hypotheses, or an essential yet easily-overlooked point, with the `dangerous bend' symbol \dbend \index{s}{\dbend} 

\section{Acknowledgments}
Some of this work was carried out with the support of an NDSEG fellowship, and final revisions were made with the support of NSF grant DMS-1303928. I thank the Institut Henri Poincar\'{e}, which I visited for the $2010$ Trimestre Galois, and the Oxford University mathematics department, which I visited in $2011-2012$, for their hospitality. I of course am grateful to the Princeton math department, both for its official support, and for its singular camaraderie.  

I thank the anonymous referee(s) for many helpful comments improving the readability of this paper. It also gives me great pleasure to thank the following people. My debt to Brian Conrad, who wrote the paper (\cite{conrad:dualGW}) that catalyzed much of the present work, will be clear. He also provided very helpful comments on an earlier draft. I am grateful to Peter Sarnak for his encouragement and several enjoyable conversations about this material, and to Christopher Skinner for reading an earlier draft. Most of all, I thank Andrew Wiles, my advisor, for his ongoing support and patient encouragement over the last few years.

%\chapter{Foundations \& examples}
\chapter{Foundations \& examples}
This chapter discusses the motivating examples and technical ingredients that underly the general arguments of Chapter \ref{2}. After (\S \ref{review}) recalling foundational lifting results of Tate, Wintenberger, and Conrad, and (\S \ref{l-adichodge}) setting up the tools we will need from $\ell$-adic Hodge theory, we undertake a detailed discussion (\S \ref{GL1}) of Hecke characters and $\ell$-adic Galois characters, especially with reference to their possible infinity-types and Hodge-Tate weights. In \S \ref{formal} we first explain an abstract principle that is important for `doing Hodge theory with coefficients;' in \S \ref{cmdescentsection} we apply this to generalize Weil's result on descent of type $A$ Hecke characters. Drawing further inspiration from \cite{weil:characters}, we then discuss (\S \ref{walgsection}) what seems to be the most useful generalization of the type $A$ condition to higher-rank groups, which we call $W$-algebraic representations.  The subsequent sections (\S \ref{HMFs} and \S \ref{GaloisHMFs}) discuss the simplest non-abelian example, that of $W$-algebraic Hilbert modular forms, their associated Galois representations, and the first interesting cases of Conrad's geometric lifting question (Question \ref{conradquestion}); although certain results in these sections are superseded by the general theorems of Chapter \ref{2}, we can also prove much more refined statements in the Hilbert modular case, as well as some complementary results about the $\mr{GL}_2 \times \mr{GL}_2 \xrightarrow{\boxtimes} \mr{GL}_4$ functorial lift.\footnote{Some of the results of \S \ref{GaloisHMFs}, among other things, have been independently obtained by Tong Liu and Jiu-Kang Yu in their preprint \cite{liu-yu:GO4}.} Another case where more precise results are accessible with the current technology is for certain automorphic representations on symplectic groups, and their associated (orthogonal and spin) Galois representations (\S \ref{spineg}).

\section{Review of lifting results}\label{review}
In this section, we review the basic Galois lifting results due to Tate, Wintenberger, and Conrad. We then highlight some of the basic problems that remain unaddressed by these foundational results. The two main elements of Tate's proof will come up repeatedly throughout these notes, so they are worth making explicit: he requires information coming both from the `automorphic-Galois correspondence' and from the bare automorphic theory, which in our work amounts to some question of what infinity-types automorphic representation can have.\footnote{A very simple aspect of this extremely difficult problem, that is.} In Tate's proof, the former is global class field theory-- in the form of the local-global structure of the Brauer group-- and the latter takes the form of precise results about the structure of the connected component of the idele class group $C_F$.  We will give a slightly different proof that emphasizes the continuity with some of our other arguments, especially Lemma \ref{Heckelift}.
\begin{thm}[Tate; see {\cite[Theorem 4]{serre:DSsurvey}}] \label{tate'stheorem}
Let $F$ be a number field\footnote{As Tate observes, the same holds for global function fields, but we are not concerned with this case.}. Then $H^2(\gal{F}, \Q/\Z)=0$.
\end{thm}\index{t}{Tate's theorem}
\proof
It suffices to prove $H^2(\gal{F}, \Q_p/\Z_p)=0$ for all primes $p$, and then an easy inflation-restriction argument shows we may assume that $F$ contains $\mu_p$. Since $H^2(\gal{F}, \Q_p/\Z_p)$ is $p$-power torsion, to show it is zero we may instead show that multiplication by $p$ is injective, or equivalently that the boundary map
\[
H^1(\gal{F}, \Q_p/\Z_p) \xrightarrow{\delta} H^2(\gal{F}, \Z/p) \cong \mr{Br}(F)[p]
\]
is surjective (the identification with the Brauer group is possible since $\mu_p \subset F$). First we note that the corresponding claim holds for the completions $F_v$: since $\mr{Br}(F_v)[p] \cong \Z/p$ (for $v$ finite; for $v$ archimedean the argument is even simpler), we only need the corresponding (at $v$) boundary map $\delta_v \colon H^1(\gal{F_v}, \Q_p/\Z_p) \to H^2(\gal{F_v}, \Z/p)$ to be non-zero, and thus only that multiplication by $p$ on $H^1(\gal{F_v}, \Q_p/\Z_p)$ should not be surjective. But a character $F_v^\times \to \Q_p/\Z_p$ is a $p^{th}$-power if and only if it is trivial on $\mu_p(F_v) \subset F_v^\times$, and there are obviously characters not satisfying this condition. Note also that for $\phi_v \in H^1(\gal{F_v}, \Q_p/\Z_p)$, the image $\delta_v(\phi_v)$ depends only on the restriction $\phi_v|_{\mu_p(F_v)}$.

For global $F$, let $\alpha$ be a $p$-torsion element of $\mr{Br}(F) \subset \oplus_v \mr{Br}(F_v)$, with local components $\alpha_v$. By the local theory, for all $v$ we have characters $\phi_v \colon F_v^\times \to \Q_p/\Z_p$ such that $\delta_v(\phi_v)= \alpha_v$, and the collection of $\alpha_v$ is equivalent to the collection of restrictions $\phi_v|_{\mu_p(F_v)}$, i.e. to the corresponding character $\phi \colon \mu_p(\af) \to \Q_p/\Z_p$. The fact that the $\alpha_v$ arise from a global Brauer class $\alpha$ implies that $\phi$ factors through $\mu_p(F) \bs \mu_p(\af)$, and we need only produce a finite-order extension to a Hecke character $\tilde{\phi} \colon C_F \to \Q/\Z$. The archimedean classes $\alpha_v$ are trivial for $v$ complex, so Lemma \ref{Heckelift} implies the existence of such an extension (note that it suffices to produce, as in the Lemma, an extension to a character $\tilde{\phi} \colon C_F \to \Q/\Z$, since we can then project down to $\Q_p/\Z_p$ and still have an extension of $\phi$).
\endproof
\begin{eg}
As mentioned in the introduction, Tate's result certainly breaks down with finite coefficients. For example, consider the projectivization of the $\ell$-adic Tate module of an elliptic curve over a totally real field $F$. The obstruction to lifting to $\mr{SL}_2(\Qlb)$ lives in $H^2(\gal{F}, \Z/2)$, and it vanishes if and only if the $\ell$-adic cyclotomic character admits a square-root, which cannot happen for $F$ totally real.
\end{eg}
%\item Let $S$ be a finite set of places of $F$, and let $\Gamma_{F, S}$ be Galois group of the maximal extension of $F$ unramified outside $S$. Then $H^2(\Gamma_{F, S}, \Q/\Z)$ can be non-zero.

The papers \cite{conrad:dualGW} and \cite{wintenberger:relevement} both study the problem of lifting Galois representation $\gal{F} \to H(\Qlb)$ through a surjection of linear algebraic groups $H' \to H$ with central kernel. The problem naturally breaks into two cases: that of finite kernel (an isogeny), and that of connected (torus) kernel. It is helpful to contrast these cases using an algebraic toy model:
\begin{eg}\label{toy}
Let $H' \onto H$ be a surjection of \textit{tori} over an algebraically closed field of characteristic zero, with kernel $D$, and suppose we are given an algebraic homomorphism $\rho \colon \mathbb{G}_m \to H$. Via the (anti-) equivalence of categories between diagonalizable algebraic groups and abelian groups (taking character groups), we see that the obstruction to lifting $\rho$ as an algebraic morphism lives in $\Ext^1(X^\bullet(D), \Z)$; in particular, when $D$ is a torus, all such $\rho$ lift. For instance, we can fill in the dotted arrow in the diagram
\[
\xymatrix{
& \mathbb{G}_m^2 \ar[d]^{(w, z) \mapsto w^r z^s}\\
\mathbb{G}_m \ar@{-->}[ur] \ar[r]_{z \mapsto z^n} & \mathbb{G}_m
}
\]
precisely when $gcd(r,s) \vert n$. The kernel of $\mathbb{G}_m^2 \to \mathbb{G}_m$ is connected precisely when $gcd(r, s)=1$, so there is no obstruction in that case, and in general there is an explicit congruence obstruction.
\end{eg}

Since it is the basis of all that follows, I recall the application of Theorem \ref{tate'stheorem} to lifting through central torus quotients (see Proposition $5.3$ of \cite{conrad:dualGW}).
\begin{prop}[$5.3$ of \cite{conrad:dualGW}]\label{tatelift}
Let $H' \onto H$ be a surjection of linear algebraic groups over $\Qlb$ with kernel a central torus $S^\vee$.\footnote{I use this awful notation to remain consistent with the `dual picture' that we will later use to think about these questions.} Then any continuous representation $\rho \colon \gal{F} \to H(\Qlb)$ lifts to $H'(\Qlb)$, i.e. we can fill in the diagram
\[
\xymatrix{
& H'(\Qlb) \ar[d] \\
\gal{F} \ar@{-->}[ru]^{\tilde{\rho}} \ar[r]_{\rho} & H(\Qlb).
}
\]
\end{prop}
\proof
By induction, we may assume $S^\vee= \mathbb{G}_m$. There is an isogeny complement $H'_1$ in $H'$ to $S^\vee$: that is, we still have a surjection $H'_1 \onto H$, but now $H'_1 \cap S^\vee$ is finite,\footnote{The simplest example is $\mr{SL}_n \subset \mr{GL}_n$.} and thus equal to $\mu_{n_0} \subset \mathbb{G}_m$ for some $n_0$. For any integer $n$ divisible by $n_0$, we can enlarge $H'_1$ to $H'_n:= H'_1 \cdot S^\vee[n]$, which now surjects onto $H$ with kernel $S^\vee[n] \cong \mu_n$. These isogenies yield obstruction classes $c_n \in H^2(\gal{F}, \Z/n)$ (as $\gal{F}$-module, the $\mu_n \subset S^\vee[n]$ is of course trivial) that are compatible under the natural maps $\Z/n \to \Z/n'$ for $n \vert n'$. Tate's theorem (\cite{serre:DSsurvey}) tells us that 
\[
\varinjlim_{n} H^2(\gal{F}, \Z/n)= H^2(\gal{F}, \Q/\Z)=0,
\]
so for sufficiently large $n$, there exists a lift $\tilde{\rho} \colon \gal{F} \to H'_n(\Qlb)$.
\endproof
\begin{rmk}\label{aeunr}
\begin{itemize}
\item Note that crucial to lifting here is the ability to enlarge the coefficients: if $H' \onto H$ is a morphism of groups over $\Ql$ and $\rho$ lands in $H(\Ql)$, then we only obtain a lift to $H'(\Ql(\mu_n))$ for sufficiently large $n$.
\item For example (and similarly in general; compare Lemma \ref{twist}), in the case of $\mr{GL}_n \to \mr{PGL}_n$, this proof produces lifts with finite-order determinant. If $\rho$ is Hodge-Tate, it is possible that no such lift is also Hodge-Tate. These non-Hodge-Tate Galois lifts have no parallel in the `toy model' of Example \ref{toy}, but we will use them as a stepping-stone toward finding geometric lifts, just as in the automorphic context we can choose an initial, possibly non-algebraic, lift, and then modify it to something algebraic (see Proposition \ref{cmextension}). By contrast, in the `motivic' version of this problem (see \S \ref{hyperlift}), the motivic Galois group does not have the extra flexibility to admit such `non-algebraic' lifts.
\item If $\rho$ is almost everywhere unramified, then it is easy to see that $\tilde{\rho}$ is almost everywhere unramified; see Lemma $5.2$ of \cite{conrad:dualGW}.
\end{itemize}
\end{rmk}
The problem of lifting through the isogeny $H'_1 \to H$ (or any isogeny $H' \to H$) is taken up in \cite{wintenberger:relevement}. Before describing this work, we state results of Wintenberger and Conrad related to lifting $p$-adic Hodge theory properties through central quotients.\footnote{Wintenberger treated isogenies; the general case is due to Conrad.} 
\begin{thm}[Corollary $6.7$ of \cite{conrad:dualGW}]\label{conradliftingP}
Let $H' \onto H$ be a surjection of linear algebraic groups over $\Qlb$ with central kernel of multiplicative type. For $K/\Ql$ finite, let $\rho \colon \gal{K} \to H(\Qlb)$ be a representation satisfying a basic $\ell$-adic Hodge theory property\footnote{Hodge-Tate, de Rham, semistable, crystalline} $\mbf{P}$. Provided that $\rho$ admits a lift $\tilde{\rho} \colon \gal{K} \to H'(\Qlb)$, it has a lift satisfying $\mbf{P}$ if and only if $\rho|_{I_K}$ has a Hodge-Tate lift $I_K \to H'(\Qlb)$.
\end{thm}
\begin{rmk}
When the kernel of $H' \to H$ is a central torus, Proposition \ref{tatelift} shows that $\rho$ always has some lift. We will see later (Corollary \ref{HTlift}) that it even always has a Hodge-Tate lift, so there is no obstruction to finding a lift $\tilde{\rho} \colon \gal{K} \to H'(\Qlb)$ satisfying $\mbf{P}$. Even simple cases of this theorem yield very interesting results: for instance, applied to $H'= \mr{GL}_{n_1} \times \mr{GL}_{n_2} \to H \subset \mr{GL}_{n_1 n_2}$, for $H$ the image of the exterior tensor product, it shows that if a tensor product of $\gal{K}$-representations satisfies $\mathbf{P}$, then up to twist they themselves satisfy $\mathbf{P}$. We will use this example in Proposition \ref{tensorformalism}.
\end{rmk}

We now sketch the proof of Wintenberger's main global result. This provides an occasion to simplify and generalize the arguments using subsequent progress in $p$-adic Hodge theory; it also serves to make clear what these methods do and do not prove, and to set up a contrast with the quite different techniques that we will use.
\begin{thm}[$2.1.4$ of \cite{wintenberger:relevement}]\label{wintenberger}
Let $H' \xrightarrow{\pi} H$ be a central isogeny of linear algebraic groups over $\Q$. Let $S$ be a finite set of non-archimedean places of $F$. Then there exist two extensions of number fields $F'' \supset F' \supset F$, and a finite set of finite places $S'$ of $F'$ such that for any prime number $\ell$ and any representation $\rho_{\ell} \colon \gal{F} \to H(\Ql)$ satisfying
\begin{itemize}
\item $\rho_{\ell}$ has good reduction\footnote{Meaning unramified for places not dividing $\ell$, and crystalline for places above $\ell$.} outside $S$;
\item For all $v \vert \ell$, the one-parameter subgroups $\mu_v \colon \mathbb{G}_m \to H_{\CC_{F_v}}$ giving the Hodge-Tate structure of $\rho_{\ell}|_{\gal{F_v}}$ lift to $H'_{\CC_{F_v}}$;\footnote{When some local lift $I_{F_v} \to H'(\Ql)$ exists, this condition is equivalent to the existence of a local Hodge-Tate lift. When $\rho_{\ell}|_{I_{F_v}}$ is moreover crystalline, or semistable, it is equivalent to the existence of a crystalline, or semistable, lift.}
\end{itemize}
then the restriction $\rho_{\ell}|_{\gal{F'}}$ admits a geometric lift $\rho'_{\ell}$, unramified outside $S'$, to $H'(\Ql)$; if $\rho_{\ell}$ is crystalline (resp. semistable) at places above $\ell$, then the lift may be chosen crystalline (resp. semistable). Moreover, the restriction $\rho'_{\ell}|_{\gal{F''}}$ is unique, i.e. independent of the initial choice of lift $\rho'_{\ell}$.
\end{thm}
\proof
We follow Wintenberger's arguments in detail, simplifying where the technology allows. First we explain some generalities: since the map on $\Ql$-points $\pi_\ell \colon H'(\Ql) \to H(\Ql)$ need not be surjective, we have an initial obstruction $\mc{O}^1(\rho_{\ell})$ in $\Hom(\gal{F}, H^1(\gal{\Ql}, \ker \pi))$ given by the composition
\[
\gal{F} \xrightarrow{\rho_{\ell}} H(\Ql)/\im(\pi_{\ell}) \into H^1(\gal{\Ql}, \ker \pi(\Qlb)).
\]
We have the choice of killing this obstruction by enlarging $\Ql$ to some finite extension, or by restricting $F$; following Wintenberger, we will do the former. Having dealt with this obstruction, we will face the more serious lifting obstruction $\mc{O}^2(\rho_{\ell})$ in $H^2(\gal{F}, \ker \pi (\Ql))$.

With these preliminaries aside, we construct the field $F'$, first dealing with all but the (finite) set of $\ell$ lying below a place in $S$. Let $a(\pi)$ be the annihilator of $\ker \pi$. The field $F'$ (and set of places $S'$) will be defined, independent of $\ell$, in the following three steps:
\begin{itemize}
\item Let $F_1$ be the maximal abelian extension of exponent $a(\pi)$, unramified outside $S$. For all $\ell$ and $\rho_{\ell}$, the class $\mc{O}^1(\rho_{\ell})$ dies after restriction to $\gal{F_1}$: outside of $S_{\ell}= S \cup \{v \vert \ell\}$, it is of course unramified, and if $v \vert \ell$ but $v \not \in S$ (i.e. $\rho_{\ell}|_{\gal{F_v}}$ is crystalline), then Th\'{e}or\`{e}me $1.1.3$ of \cite{wintenberger:relevement}) shows the existence of a crystalline lift, valued in $H'(\Ql)$, of $\rho_{\ell}|_{I_{F_v}}$.\footnote{Since weakly-admissible is now known to be equivalent to admissible, Wintenberger's argument applies without restriction on the Hodge-Tate weights, or the ramification of the local field $F_v/\Ql$.} 
\item Let $F_2/F_1$ be a totally imaginary extension containing $\Q(\zeta_{a(\pi)})$ and such that for all places $v'$ above a place $v$ in $S$, $a(\pi)$ divides the local degree $[F_{2, v'}: F_{1, v}(\zeta_{a(\pi)})]$. If we write $S_2$ for the set of places of $F_2$ above $S$, $\ell$, or $a(\pi_{\ell})$ (the annihilator of $\ker \pi(\Ql)$), then this implies local triviality of the obstruction classes, as long as we assume $S$ contains no places above $\ell$. That is,
\[
\mc{O}^2(\rho_{\ell}) \in \ker \left( H^2(\Gamma_{F_2, S_2}, \ker \pi(\Ql)) \to \bigoplus_{v \in S_2 \cup S_{\infty}} H^2(\gal{F_v}, \ker \pi(\Ql)) \right).
\]
At $v'$ above $S$, the fact that $F_2$ contains $\zeta_{a(\pi)}$ makes the local obstruction class a Brauer obstruction, which is killed over $F_2$ by the assumption $a(\pi) \vert [F_{2, v'}: F_{1, v}(\zeta_{a(\pi)})]$. At places not in $S$ but above $a(\pi_{\ell})$ there is no obstruction, since unramified reprsentations always lift (granted that $\mc{O}^1(\rho_{\ell})$ is trivial). Finally, at places above $\ell$, which we have assumed for now do not lie in $S$, Wintenberger's local result implies the existence of a (crystalline, even) lift.\footnote{On the full decomposition group $\gal{F_v}$, this is the Proposition of $1.2$ of \cite{wintenberger:relevement}; to extend from $I_{F_v}$ to $\gal{F_v}$, one must use the fact that a crystalline lift on inertia is unique.}
\item Let $F_3$ be the Hilbert class field of $F_2$. By global duality, the above kernel is Pontryagin dual to
\[
\ker \left( H^1(\Gamma_{F_2, S_2}, \ker \pi(\Ql)^*) \to \bigoplus_{v \in S_2 \cup S_{\infty}} H^1(\Gamma_{F_2, S_2}, \ker \pi (\Ql)^*) \right)
\]
The Galois module $\ker \pi(\Ql)^*$ is trivial since $F_2$ contains $\zeta_{a(\pi)}$, so this is just the space of homomorphisms $\gal{F_2} \to \ker \pi (\Ql)$ that are unramified everywhere and trivial at the places $S_2$. Restriction of the class $\mc{O}^2(\rho_{\ell})$ to $F_3$ corresponds to (via global duality for $F_3$ now) the *transfer* on the Pontryagin dual of the $H^1$'s, so by the Hauptidealsatz (trivialization of ideal classes upon restriction to the Hilbert class field), $\mc{O}^2(\rho_{\ell})$ dies upon restriction to $F_3$.
\end{itemize}
We conclude that, independent of $\ell$ not below places of $S$, and of $\rho_{\ell}$ satisfying the hypotheses of the theorem, there exist lifts over the number field $F_3$. To deal with $\ell$ for which the local representations are not crystalline, we apply Proposition $0.3$ of \cite{wintenberger:relevement} to find a finite extension $F(\ell)/F$, unramified outside $S_{\ell}$, and such that all $\rho_{\ell}|_{\gal{F(\ell)}}$ will lift.\footnote{This is a `soft' result. It is easy to see that for fixed $\ell$, all $\rho_{\ell}$ unramified outside $S$ lift after restriction to some $F(\ell)$. The point of Wintenberger's theorem is that, with the $\ell$-adic Hodge theory hypotheses, $F(\ell)$ can be taken independent of $\ell$.} Enlarging $F_3$ by the composite of the $F(\ell)$ (for $\ell$ below places of $S$), we obtain a number field $F_4$ and a set of places $S_4$ equal to all places above $S$ or $a(\pi)$, such that for \textit{all} $\ell$, any $\rho_{\ell}$ as in the theorem has a lift with good reduction outside $S_4$ after restriction to $\gal{F_4}$. 

The final step of the argument makes one final extension to show that these lifts can be taken crystalline (or semistable) whenever $\rho_{\ell}$ is, and even unramified outside all primes above $S_{\ell}$ (including those above $a(\pi)$, although this is not necessary for the theorem's conclusion).\footnote{We omit the details. See Lemme $2.3.5$ and following of \cite{wintenberger:relevement}.} This will be the number field $F'$, and the set of places $S'$ can then be taken to be those above $S$. Finally, all such lifts $\rho'_{\ell}$ (with good reduction outside $S'$) are the same over the maximal abelian extension of exponent $a(\pi)$, unramified outside $S'$, of $F'$; this is our $F''$.
\endproof 

\dbend
\begin{rmk}
We can apply the theorem to a system of representations $\rho_{\ell}$ all having good reduction outside $S$, and all having liftable Hodge-Tate cocharacters (which in the `motivic' case we may assume all arise from a common underlying geometric Hodge structure), to obtain a common number field $F'$, with finite set of places $S'$, and a further extension $F''$ such that the whole system lifts to $F'$ (with good reduction outside $S'$), and such that all the lifts are unique over $F''$.  The crucial limitation of this result is that, even if the $\rho_{\ell}$ are assumed to be the system of $\ell$-adic realizations of some motive (for absolute Hodge cycles, as in \cite{wintenberger:relevement}, or for motivated cycles), the proof of the theorem produce lifts $\rho'_{\ell}$ that are neither (proven to be) motivic nor even weakly compatible (almost everywhere locally conjugate). So, even with the strong uniqueness properties ensured by Wintenberger's theorem, it is still not at all clear how to lift a compatible (or even motivic) system of $\ell$-adic representations to another compatible (or motivic) system. 
\end{rmk}
Here is the simplest cautionary example of how weakly compatible systems need not lift to weakly compatible systems; it arises from an example due to Serre (see \cite{larsen:locglobconj}) involving projective representations that are everywhere locally, but not globally, conjugate:\index{t}{compatible system of $\ell$-adic representations valued in linear algebraic group}
\begin{eg}\label{locglobconj}
For any positive integer $n$, let $H_n$ denote the Heisenberg group of upper-triangular unipotent elements in $\mr{GL}_3(\Z/n\Z)$. This is a non-split extension
\[
1 \to \Z/n \to H_n \to \Z/n \times \Z/n \to 1,
\]
and writing $A$ and $B$ for lifts of the generators of the quotient copies of $\Z/n$, and $Z$ for a generator of the center, we have the commutation relation $[A, B]= Z$. Now let $\zeta$ be a primitive $n^{th}$ root of unity; for any $\alpha \in \Z/n$, we can then define a representation $\rho_{\alpha} \colon H_n \to \mr{GL}_n(\CC)$ by, in the standard basis $e_i$, $i=1, \ldots, n$, and abusively letting $e_0$ denote $e_n$ as well, 
\begin{align*}
A(e_i)&= e_{i-1} \\
B(e_i)&= \zeta^{(i-1)\alpha}e_i \\
Z(e_i)&= \zeta^{\alpha} e_i
\end{align*}
for $i= 1, \ldots, n$. One can check that for all $\alpha \in (\Z/n\Z)^{\times}$, the associated projective representations $\overline{\rho_{\alpha}} \colon H_n \to \mr{PGL}_n(\CC)$ are everywhere locally conjugate. Nevertheless, for distinct $\alpha$, they are not globally conjugate: if they were, then the corresponding $\mr{GL}_n(\CC)$-representations would be twist-equivalent. The implications for our problem are the following: take $\overline{\rho_{\alpha}}$ and $\overline{\rho_{\beta}}$ for $\alpha \neq \beta$, and form a weakly-compatible system $\rho_{\ell} \colon \gal{F} \onto H_n \to \mr{PGL}_n(\Qlb)$ ($H_n$ certainly arises as Galois group of a finite extension of number fields.) in which some $\rho_{\ell}$ are formed from $\overline{\rho_{\alpha}}$ and others from $\overline{\rho_{\beta}}$. There is then no system of weakly-compatible lifts $\tilde{\rho}_{\ell} \colon \gal{F} \to \mr{GL}_n(\Qlb)$, since any character of $\gal{F}$ is trivial on $Z$ (viewed as an element of $\Gal(F^{ab}F'/F)$, where $F'/F$ is the given $H_n$-extension; note that $\Gal(F^{ab}F'/F^{ab})$ is generated by $Z$). Also note that we can even produce such examples where the lifts to $\mr{GL}_n(\CC)$ have the same determinant:
\[
\det(\rho_{\alpha}) \colon
\begin{cases}
A \mapsto (-1)^{n-1}; \\
B \mapsto \text{$1$ if $n$ is odd, or if $\alpha$ and $n$ are even; $-1$ otherwise};\\
Z \mapsto 1.
\end{cases}
\]
Thus, even weakly compatible systems $\gal{F} \to \mr{PGL}_n(\Qlb) \times \Qlb^\times$ need not have compatible lifts through the isogeny $\mr{GL}_n \to \mr{PGL}_n \times \mathbb{G}_m$.

To connect this example more explicitly with endoscopy, note that $\mr{Cent}_{\mr{PGL}_n(\CC)}(\overline{\rho_{\alpha}})$ is the (cyclic order $n$) subgroup generated by $\overline{\rho_{\alpha}}(A)$, whereas $\rho_{\alpha}$ itself is irreducible. These centralizers partially govern the multiplicities of discrete automorphic representations in Arthur's conjectures. Explicitly, $\rho_{\alpha}(ABA^{-1})= \zeta^{\alpha} B$.
\end{eg} 

\section{$\ell$-adic Hodge theory preliminaries}\label{l-adichodge}
\subsection{Basics}\label{l-adicbasics}
In this section, we recall some background and prove some simple lemmas in $\ell$-adic Hodge theory. Throughout, let $K/ \Ql$ be a finite extension. Choose an algebraic closure $\overline{K}/K$, and let $\CC_K$ denote the completion of $\overline{K}$.  Since $\gal{K}$ acts by isometries on $\overline{K}$, $\CC_K$ inherits a continuous $\gal{K}$-action. We will mainly use only the simplest of Fontaine's period rings, $\Bht$ and $\Bdr$, and even these only formally. Again, we refer the reader to the references in \S \ref{ladicrefs} for background, but here we quickly recall the little needed to follow future arguments. $\Bht$ is the graded $\CC_K$-algebra 
\[
\Bht= \bigoplus_{i \in \Z} \CC_K(i),
\]
\index{s}{$\Bht$}with the obvious continuous $\CC_K$-semilinear action of $\gal{K}$, $\CC_K(i)$ denoting a twist of $\CC_K$ by the $i^{th}$ power of the cyclotomic character. $\Bdr$\index{s}{$\Bdr$} is the fraction field of a complete DVR and $K$-algebra $\Bdr^+$ whose residue field is $\CC_K$. We do not recall the construction, but note that $\Bdr$ is a \textit{filtered} $K$-algebra (the filtration associated to the maximal ideal of $\Bdr^+$) with a continuous $K$-linear action of $\gal{K}$. Two of the fundamental results of the theory assert that $\gr^{\bullet} (\Bdr) \cong \Bht$ (as graded $\CC_K$-algebras with semi-linear $\gal{K}$-action) and $\Bdr^{\gal{K}}=\Bht^{\gal{K}}= K$. From this second point it follows that if $V$ is a finite-dimensional representation of $\gal{K}$ over $\Ql$, then 
\[
\Dht(V)= (\Bht \otimes_{\Ql} V)^{\gal{K}}
\]\index{s}{$\Dht$}
is a $K$-vector space of dimension at most $\dim_{\Ql}(V)$, and similarly
\[
\Ddr(V)= (\Bdr \otimes_{\Ql} V)^{\gal{K}}
\]\index{s}{$\Ddr$}
is a $K$-vector space of dimension at most $\dim_{\Ql}(V)$.
\begin{defn}
Let $V$ be a finite-dimensional representation of $\gal{K}$ over $\Ql$. $V$ is said to be \textit{Hodge-Tate} if $\dim_K \Dht(V)= \dim_{\Ql}(V)$, and \textit{de Rham} if $\dim_K \Ddr(V)= \dim_{\Ql}(V)$.\index{t}{Hodge-Tate Galois representation}\index{t}{de Rham Galois representation}
\end{defn}
If $V$ is de Rham, then it is Hodge-Tate, but not conversely.\footnote{The standard example is a non-split extension of $\gal{K}$- representations
\[
0 \to \Ql \to V \to \Ql(1) \to 0.
\]
Standard results on Galois cohomology of local fields show such a $V$ exists. It's not so hard to show $V$ is Hodge-Tate, but to show that it cannot be de Rham seems to require quite deep input; see for instance \cite[Example 6.3.5]{brinon-conrad:cmi}.
}

We will usually consider $\gal{K}$-representations $V$ with coefficients in some extension $E/\Ql$, typically either some unspecified finite extension, or $E= \Qlb$. Then $\Ddr(V)$ (mutatis mutandis for $\Dht$) is a filtered $K \otimes_{\Ql} E$-module; if $V$, viewed as a $\Ql$ representation by forgetting the $E$-linear structure, is de Rham, then $\Ddr(V)$ is a free $K \otimes_{\Ql} E$-module. To see this, we may assume $E$ is large enough to contain all embeddings of $K$ into $\Qlb$. Then for all $\Ql$-embeddings $\tau \colon K \into E$, we have orthogonal idempotents $e_{\tau} \in K \otimes_{\Ql} E$, giving rise to an isomorphism
\begin{align*}
&K \otimes_{\Ql} E \xrightarrow[\sim]{(e_{\tau}} \bigoplus_{\tau \colon K \into E} E\\
&x \otimes \alpha \mapsto \left( \tau(x)\alpha \right)_{\tau}.
\end{align*}
As $K \otimes_{\Ql} E$-module,
\[
\Ddr(V) \cong \bigoplus_{\tau} e_{\tau} \Ddr(V) \cong \bigoplus_{\tau} (\Bdr \otimes_{K, \tau} V)^{\gal{K}}
\]
(compare Lemma \ref{HTS}), and each space on the right-hand side has $E$-dimension at most $\dim_E(V)$. The sum of their $E$-dimensions is then at most $\dim_{\Ql}(K) \dim_E(V)$, so the $\Ql$-dimension of the right-hand side is at most $\dim_{\Ql}(K) \dim_{\Ql}(V)$, which is, since $V$ is de Rham, the $\Ql$-dimension of the left-hand side. Therefore equality holds everywhere, and it easily follows that $\Ddr(V)$ is free over $K \otimes_{\Ql} E \cong \oplus_{\tau} E$.

\begin{eg}
Note that even when $V$ is de Rham, other steps of the filtration need not be free over $K \otimes_{\Ql} E$. We sketch a basic example, omitting the details because they require introducing more about $\Bdr$ than we will subsequently need. Take any extension to a character $\psi$ of $\gal{K}$ of the character
\[
I_K \xrightarrow{\rec_K} \mc{O}_K^{\times} \xrightarrow{\tau} \Qlb^\times,
\]
where $\tau \colon K \into \Qlb$ is any fixed $\Ql$-embedding. If $K= \Ql$, so that there is only one such embedding, then $\psi|_{I_K}$ is just the restriction to inertia of the cyclotomic character, and $\Ddr(\psi)$ has filtration concentrated in degree -1. If $K \neq \Ql$, then for all embeddings $\tau' \neq \tau$, the filtered one-dimensional $\Qlb$-vector space $e_{\tau'}\Ddr(\psi)$ has its filtration concentrated in degree zero, whereas the $\tau$ piece is concentrated in degree -1. 
\end{eg}

It will be crucial for us to have a refined notion of Hodge-Tate weight that takes into account the different embeddings $\tau \colon K \into \Qlb$:
\begin{defn}\label{labeledHT} \index{t}{labeled Hodge-Tate weights}
Let $V$ be a de Rham representation of $\gal{K}$ on a $\Qlb$-vector space of dimension $d$, so that $\Ddr(V)$ is a free $K \otimes_{\Ql} \Qlb$-module of rank $d$. Then for each $\tau \colon K \into \Qlb$, we define the $d$-element multi-set of $\tau$\textit{-labeled Hodge-Tate weights}, $\mr{HT}_{\tau}(D)$,\index{s}{$\mr{HT}_{\tau}$} to be the collection of integers $h$ (with multiplicity) such that 
\[
\mathrm{gr}^h(e_\tau(\Ddr(V)) \neq 0.
\]
\end{defn}
We could also make this definition for $V$ having coefficients in a subfield $E\subset \Qlb$, possibly then having to enlarge $E$ in order to contain the embeddings of $K$ (and thus define the labeled Hodge-Tate weights). In \S \ref{formal}, we will discuss the elementary but significant implications of $E$ \textit{not} containing all embeddings of $K$.

\subsection{Labeled Hodge-Tate-Sen weights}\label{labeledweights}
Later we will work with $\gal{K}$-representations that are not even Hodge-Tate and will need a generalized notion of (labeled) Hodge-Tate weight. First, again let $V$ be a de Rham representation of $\gal{K}$ on a $\Qlb$-vector space. The identification of $\gr^\bullet(\Bdr)$ with $\mr{B}_{\mr{HT}}$ implies that the multiplicity of $q$ in $\mr{HT}_{\tau}(V)$ also equals $\dim_{\Qlb}(V \otimes_{\tau, K} \CC_K(-q))^{\gal{K}}$. We use this formulation to extend the definition of $\tau$-labeled weights to the Hodge-Tate case.\index{t}{labeled Hodge-Tate weights of a Hodge-Tate representation} Note that if $K'/K$ is any finite extension, and $\tau' \colon K' \into \Qlb$ is any embedding extending $\tau \colon K \into \Qlb$, then $\mr{HT}_{\tau'}(V|_{\gal{K'}})= \mr{HT}_{\tau}(V)$.

\index{t}{Hodge-Tate-Sen weights}Now we drop the assumption that $V$ is Hodge-Tate. For Sen's theory, which among other things introduces a generalized notion of Hodge-Tate-Sen weight, \cite[II.1.2]{berger:introp-adicrep} is a very brief, and very illuminating, overview; see the references there or \cite[Part IV]{brinon-conrad:cmi} for a detailed introduction. Sen's theory gives much more precise results than what we use here; the most important thing for our purposes is that it associates to any $\CC_K$-semi-linear representation $V$ of $\gal{K}$ a $\CC_K$-linear endomorphism $\Theta_{V}$ (the `Sen operator') of $V$ (see \cite[Corollary 15.1.6, Theorem 15.1.7]{brinon-conrad:cmi}).\index{t}{Sen operator} The eigenvalues of $\Theta_{V}$ are the `Hodge-Tate-Sen' weights; this terminology is justified by the fact that $V$ is Hodge-Tate\footnote{Strictly speaking, in \S \ref{l-adicbasics} we only defined Hodge-Tate representations on finite-dimensional $\Ql$-vector spaces. This is a special case of a more general notion for $\CC_K$-semi-linear $\gal{K}$-representations: see \cite[Definition 2.3.4]{brinon-conrad:cmi}.} if and only if $\Theta_{V}$ is semi-simple with integer eigenvalues (\cite[Exercise 15.5.4]{brinon-conrad:cmi}).

We will now be more precise about how to work with $\gal{K}$-representations with coefficients in $\Qlb$. We will also describe Sen operators for representations valued in a general linear algebraic group. For some finite extension $E/\Ql$ inside $\Qlb$, our (arbitrary) $\gal{K}$-representation $V$ descends to an $E$-linear representation $V_E$.  For any embedding $\iota \colon \Qlb \into \CC_K$, we have the Sen operator $\Theta_{\rho, \iota} \in \End_{\CC_K}(V \otimes_{\Qlb, \iota} \CC_K)$.  More precisely, let $K'= \iota(E) K$, so that $V_E \otimes_{E, \iota} \CC_K$ is a $\CC_K$-semilinear representation of $\gal{K'}$, to which we can in the usual way associate the Sen operator $\Theta_{\rho, \iota}$.\index{s}{$\Theta_{\rho, \iota}$}\index{t}{Sen operator} Note that the Sen operator of a $\CC_K$ semi-linear representation is insensitive to finite restriction, so this construction is independent of the choice of (sufficiently large) $K'$. Moreover, choosing a model $V_{E'}$ of $V$ with coefficients in some finite extension $E'/E$ inside $\Qlb$ also yields the same $\Theta_{\rho, \iota}$, since we don't change the $\CC_K$-semilinear representation of $\gal{K'}$ (this may be a bigger $K'$, having enlarged $E$). This allows the following Tannakian observation about Sen operators (\S 6 of \cite{conrad:dualGW}).
\begin{lemma}
Let $G$ be a linear algebraic group over $\Qlb$ with Lie algebra $\mf{g}$ and $\rho \colon \gal{K} \to G(\Qlb)$ a (continuous) representation. Then for each embedding $\iota \colon \Qlb \into \CC_K$, there exists a unique Sen operator $\Theta_{\rho, \iota} \in \mf{g} \otimes_{\Qlb, \iota} \CC_K$ such that for all linear representations $r \colon G \xrightarrow[/\Qlb]{} \mr{GL}_{V}$, $\mr{Lie}(r)(\Theta_{\rho, \iota})$ is equal to the previously-defined $\Theta_{r \circ \rho, \iota}$.
\end{lemma}
\proof
This follows from Tannaka duality for Lie algebras (a precise reference is \cite{harish-chandra:tannaka}). Namely, for each $\Qlb$-representation $G \to \mr{GL}_{V}$ we get an element $\Theta_{V, \iota} \in \mf{gl}_V \otimes_{\iota} \CC_K$, and these satisfy the functorial properties
\begin{itemize}
\item If $V_1 \xrightarrow{T} V_2$ is a $G$-morphism, then $\mr{Lie}(T) \circ \Theta_{V_1, \iota}= \Theta_{V_2, \iota} \circ \mr{Lie}(T)$;
\item $\Theta_{V_1 \oplus V_2, \iota}= \Theta_{V_1, \iota} \oplus \Theta_{V_2, \iota}$;
\item $\Theta_{V_1 \otimes V_2, \iota}= \Theta_{V_1, \iota} \otimes \mr{id}_{V_2} + \mr{id}_{V_1} \otimes \Theta_{V_2, \iota}$.
\end{itemize}
(If different fields $E$ and $K'$ as above are needed to define the $\Theta_{V_{i}, \iota}$, we can pass to a common extension and apply the remarks preceding the lemma.) Tannaka duality then implies that all $\Theta_{V, \iota}$ arise from a unique element $\Theta_{\rho, \iota} \in \mf{g}\otimes_{\iota} \CC_K$.
\endproof
\begin{rmk}
The necessary functorial properties of the Sen operators are not automatic. All the relevant statements are in \S 15 of \cite{brinon-conrad:cmi}.
\end{rmk}
As mentioned previously, $V \otimes_{\Qlb, \iota} \CC_K$ is Hodge-Tate if and only if $\Theta_{\rho, \iota}$ is semi-simple with integer eigenvalues. We will need the following comparison:
\begin{lemma}\label{HTS}
Suppose $\rho \colon \gal{K} \to \mathrm{Aut}_{\Qlb}(V)$ is Hodge-Tate. Any embedding $\iota \colon \Qlb \into \CC_K$ induces $\tau_{\iota} \colon K \into \Qlb$, and we then have
\[
\mr{HT}_{\tau_{\iota}}(\rho)= \{ \text{eigenvalues of $\Theta_{\rho, \iota}$}\}.
\]
\end{lemma}
\proof
Let $E$ and $K'$ be as above. The $\gal{K}$-representation $V_E$ is Hodge-Tate, and there is a natural isomorphism of graded $E \otimes_{\Ql} K'$-modules
\[
(V_E \otimes_{\Ql} \mr{B}_{HT})^{\gal{K}} \otimes_K K' \xrightarrow{\sim} (V_E \otimes_{\Ql} \mr{B}_{HT})^{\gal{K'}}.
\]
Restricting to the $q^{th}$ graded pieces, we get an $E \otimes_{\Ql} K'$-isomorphism
\[
(V_E \otimes_{\Ql} \CC_K(-q))^{\gal{K}} \otimes_K K' \xrightarrow{\sim} (V_E \otimes_{\Ql} \CC_K(-q))^{\gal{K'}},
\]
which in turn is an $E \otimes_{\Ql} K'$-isomorphism
\[
\bigoplus_{\tau \colon K \into E} (V_E \otimes_{\tau, K} \CC_K(-q))^{\gal{K}} \otimes_K K' \xrightarrow{\sim} \bigoplus_{\iota \colon E \into \CC_K} (V_E \otimes_{E, \iota} \CC_K(-q))^{\gal{K'}}.
\]
Projecting to the $\tau$-component, we obtain an $E \otimes_{\tau, K} K'$-isomorphism
\[
(V_E \otimes_{\tau, K} \CC_K(-q))^{\gal{K}} \otimes_K K' \xrightarrow{\sim} \bigoplus_{\iota \colon E \into \CC_K: \tau_{\iota}= \tau} (V_E \otimes_{E, \iota} \CC_K(-q))^{\gal{K'}}.
\]
Writing $m_{q, \tau}$ for the multiplicity of $q$ in $\mr{HT}_{\tau}(V)$, the left-hand side is a free $E \otimes_{\tau, K} K'$-module of rank $m_{q, \tau}$. Meanwhile, the $K'$-dimension of the $\iota$-factor of the right-hand side is by definition the multiplicity of $q$ as a Hodge-Tate weight of $V_E \otimes_{E, \iota} \CC_K$, hence is $q$'s multiplicity as an eigenvalue of $\Theta_{\rho, \iota}$. We deduce (applying the $\iota$-projection in $E \otimes_{\tau, K} K' \cong \prod_{\iota: \tau_{\iota}=\tau} K'$) that for all $\iota$ the eigenvalues of $\Theta_{\rho, \iota}$ match the multi-set $\mr{HT}_{\tau_{\iota}}(V)$.
\begin{rmk}
We will only apply this when $V$ is in fact de Rham, with the exception of Corollary \ref{HTlift}.
\end{rmk}
\subsection{Induction and $\ell$-adic Hodge theory}
We will later need a result on compatibility of Fontaine's functors with induction. The following lemma is certainly well-known, but I don't know of a proof in the literature, so I record some details. For a finite extension $K/\Ql$, we let (as is standard notation in the subject) $K_0 \subset K$ denote the maximal unramified sub-extension of $K$. 
\begin{lemma}\label{dRInd}
Let $L/K/\Q_\ell$ be finite, and let $W$ be a de Rham representation of $\gal{L}$.  Then $V= \mathrm{Ind}_L^K W$ is also de Rham, and $\Ddr(V)$ is the image under the forgetful functor $\mr{Fil}_L \to \mr{Fil}_K$ of $\Ddr(W)$.  Moreover, $\Ind_L^K W$ is crystalline if and only if $W$ is crystalline and $L/K$ is unramified.
\end{lemma}
\proof
This follows almost immediately from Frobenius reciprocity if one uses contravariant Fontaine functors: 
\[
\Ddr^*(V):= \Hom_{\gal{K}}(V, \Bdr) \cong \Hom_{\gal{L}}(W, \Bdr|_{\gal{L}})= \Ddr^*(W).
\]
Since $\Ddr^*(V) \cong \Ddr(V^*)$, and $\Ind_{L}^K(W^*) \cong \Ind_L^K(W)^*$, this shows that $\Ddr(V)$ is simply the filtered $K$-vector space underlying $\Ddr(W)$. Comparing dimensions, it is clear that $V$ is de Rham if and only if $W$ is. In the crystalline case, the same observation yields $\Dcris^*(V) \cong \Dcris^*(W)$, so
\[
\dim_{K_0} \Dcris(V)= \dim_{K_0} \Dcris(W)= \dim_{L_0} \Dcris(W) [L_0:K_0],
\]
and comparing dimensions we see that $V$ is crystalline if and only if $W$ is crystalline and $[L_0:K_0]=[L:K]$, i.e. $L/K$ is unramified.
\endproof

\section{$\mr{GL}_1$}\label{GL1}
A detailed analysis of the $\mr{GL}_1$ theory provides both motivation and technical tools for addressing the more general Galois and automorphic lifting problems. 

\subsection{The automorphic side}
We begin by reviewing some basic facts about Hecke characters $\psi \colon \af^\times/F^\times \to \CC^\times$. Any such $\psi$ is the twist by $|\cdot|_{\af}^r$, for some $r \in \RR$, of a unitary Hecke character. At each place $v$ of $F$, $\psi_v= \psi|_{F_v^\times}$ can be decomposed via projection to the maximal compact subgroup in $F_v^\times$; at $v \vert \infty$, this lets us write any unitary $\psi_v$ in terms of a given embedding $\iota_v \colon F_v \into \CC$ as
\[
\psi_v(x_v)= (\iota_v(x_v)/|\iota_v(x_v)|)^{m_v} |\iota_v(x_v)|_{\CC}^{i t_v},
\]
with $m_v= m_{\iota_v}$ an integer, and $t_v$ a real number.\index{s}{$m_v$}\index{s}{$t_v$}\index{t}{infinity-type of a Hecke character} To avoid cumbersome notation, we will sometimes omit reference to the embedding $\iota_v$. Hecke characters encode finite Galois-theoretic information as well as rather subtle archimedean information; the latter will be particularly relevant in what follows, and the basic observation (\cite{weil:characters}) is:
\begin{lemma}\label{infinitycriterion}
There is a (unitary, say) Hecke character $\psi$ of $F$ with archimedean parameters $\{m_v, t_v\}_{v \vert \infty}$, as above, if and only if for some positive integer $M$, all global units $\alpha \in \mc{O}_F^\times$ satisfy
\[
\left( \prod_{v \vert \infty} \left( \frac{\iota_v(\alpha)}{|\iota_v(\alpha)|} \right)^{m_v} |\iota_v(\alpha)|^{i t_v} \right)^M= 1.
\]
Equivalently, the inside product is trivial for all $\alpha$ in some finite-index subgroup of $\mc{O}_F^\times$.
\end{lemma}
Using the notation we have just established, we recall Weil's notions of type $A_0$ and type $A$ Hecke characters:\index{t}{type $A_0$ Hecke character}\index{t}{type $A$ Hecke character}
\begin{defn}\label{AA0def}
A Hecke character $\psi$ is said to be of type $A$ if $r \in \Q$ and $t_v=0$ for all $v \vert \infty$. $\psi$ is said to be type $A_0$ if $t_v=0$ and $\frac{m_v}{2}+r$ is an integer for all $v \vert \infty$. Equivalently, for all $v \vert \infty$, there exist integers $p_v$ and $q_v$ such that $\psi_v(z)= \iota_v(z)^{p_v} \bar{\iota}_v(z)^{q_v}$.  
\end{defn}
Then Weil (\cite{weil:characters}) shows:
\begin{lemma}[Weil]\label{weilalg}
Let $\chi$ be a type $A$ Hecke character locally trivialized by the modulus $\mf{m}$. $\chi$ then induces a character $\chi^{\mf{m}}$ of $I^{\mf{m}}$, the group of non-zero, prime to $\mf{m}$, fractional ideals, whose values are all algebraic numbers.
\end{lemma}
Later we will be interested in a higher-rank version of the type $A$ condition; in current discussions of algebraicity of automorphic representations, people tend to focus on analogues of type $A_0$. We recall a basic result of Weil-Artin (\cite{weil:characters}), which will later inspire much of our discussion in higher-rank as well.
\begin{lemma}[Weil-Artin]\label{weildescent}
Let $\psi$ be a type-$A$ Hecke character of a number field $F$. Write $F_{cm}$ for the maximal CM subfield of $F$. Then there exists a finite-order Hecke character $\chi_0$ of $F$ and a type-$A$ Hecke character $\psi'$ of $F_{cm}$ such that $\psi= \chi_0 \cdot \BC_{F/F_{cm}}(\psi')$. If $F_{cm}$ is totally real, then $\psi$ is a finite-order twist of some power $|\cdot|^r$, $r \in \Q$, of the absolute value.
\end{lemma}\index{t}{CM descent prototype}
\begin{rmk}
We emphasize that the true content of this result, which should fully generalize to higher-rank (see Proposition \ref{cmdescent}), is that the infinity-type of a type-$A$ Hecke character descends to the maximal CM subfield.
\end{rmk}
\proof
Weil leaves it as an exercise, so we give a proof. We then indicate a second proof, which we will apply more generally in Proposition \ref{cmdescent}. We may assume that $\psi$ is type $A_0$. First assume $F/\Q$ is Galois, and set $G= \Gal(F/\Q)$. We may regard $F$ as a subfield of $\CC$ (i.e., choose an embedding), so the above relation becomes
\[
\prod_{g \in G} g(\alpha)^{m_g}=1
\] 
for all $\alpha$ in some finite-index subgroup of $\mc{O}_F^\times$ and some integers $m_g$. Now take any embedding $\iota \colon F \into \CC$, using it to define a complex conjugation $c_{\iota}$ on $F$, and thus on $G$. In particular, when $F$ is totally imaginary we can rewrite the above as
\[
\prod_{G/c_\iota} g(\alpha)^{m_g} (c_{\iota} g)(\alpha)^{m_{c_\iota g}}=1.
\]
Applying $\log| \iota(\cdot)|$, we get (for all $\alpha$)
\[
\sum_{G/c_\iota} (m_g+ m_{c_\iota g}) \log|\iota(g(\alpha))|_{\CC}=0,
\]
which by the unit theorem is only possible when $m_g+m_{c_\iota g}$ equals a constant $w_\iota$ independent of $g$; the analogous argument in the totally real case shows that the $m_g$ themselves are independent of $g$, concluding that case of the proof. But $w_\iota$ is independent of $\iota$ as well,\footnote{For instance, whatever $\iota$, $\frac{1}{2}|G| w_\iota= \sum_{g \in G} m_g$.} hence $m_g+m_{c_\iota g}$ is independent both of $g$ and the choice of conjugation.  This implies that for any two $\iota, \iota'$, $m_g= m_{c_\iota c_{\iota'} g}$, and since $F_{cm}$ is precisely the fixed field (in $F$) of the subgroup generated by all products $c_{\iota}c_{\iota'}$, we are done.

Now let $F$ be arbitrary, and let $\psi$ be a type $A$ Hecke character of $F$. As usual letting $\tilde{F}$ denote the Galois closure, we know from the Galois case that $\psi|_{\tilde{F}}$ has infinity-type descending to $(\tilde{F})_{cm}$. Setting $H_1= \Gal(\tilde{F}/F)$ and $H_2= \Gal(\tilde{F}/(\tilde{F})_{cm})$, the values $m_{\iota_v}$ (for embeddings $\iota_v \colon \tilde{F} \into \CC$) are constant on $H_1$-orbits (by construction) and $H_2$-orbits (by the Galois case), hence on $H_1H_2$-orbits. But $F \cap (\tilde{F})_{cm}= F_{cm}$, so $H_1H_2= \Gal(\tilde{F}/F_{cm})$, and the infinity-type descends all the way down to $F_{cm}$, proving the lemma in general.

Later, we will generalize the following argument (see Proposition \ref{cmdescent} for more details) in the totally imaginary case: for all $\sigma \in \Aut(\CC)$, there is a Hecke character ${}^\sigma \psi$ whose finite part is ${}^\sigma \psi_f(x)= \sigma( \psi_f(x))$ and whose infinity-type is given by the integers $m_{\sigma^{-1}\iota_v}$ (the key but easy check is $F^\times$-invariance). We may assume $\psi$ is unitary, so ${}^c \psi= \psi^{-1}$, and therefore the fixed field $\Q(\psi_f) \subset \Qb$ of all $\sigma \in \Aut(\CC)$ such that ${}^\sigma \psi_f \cong \psi_f$ is contained in $\Q^{cm}$. The result follows from some Galois theory.

As pointed out to me by Brian Conrad, a third, more algebraic, approach is possible, where infinity-types of algebraic Hecke characters are related to (algebraic) characters of the connected Serre group; compare the discussion of potentially abelian motives in \S \ref{taniyamalift}. Of course, in all arguments, the unit theorem is the essential ingredient.
\endproof
We will have to analyze Hecke characters in the following setting: An automorphic representation $\pi$ of an $F$-group $G$ has a central character\index{t}{central character of an automorphic representation} $\omega_\pi \colon Z_G(F) \backslash Z_G(\mbf{A}) \to \CC^\times$, and we will want to understand the possible infinity-types of extensions
\[
\xymatrix{
 Z_G(F) \backslash Z_G(\mbf{A}) \ar[r]^-{\omega_\pi} \ar[d]  & \CC^\times \\
 \tZ(F) \backslash \tZ(\mbf{A}) \ar@{-->}[ru]_{\tilde{\omega}} & 
}
\]
where $\tZ$ is an $F$-torus containing $Z_G$.
The basic case to which this is reduced, for split groups at least, is where $Z_G= \mu_n$ and $\tZ= \mbf{G}_m$, where we have to extend a character from $\mu_n(F) \bs \mu_n(\mbf{A})$ to the full idele class group $F^\times \bs \mbf{A}^\times$.
\begin{lemma}\label{Heckelift}
Let $F$ be any number field. Given a continuous character $\omega \colon \mu_n(F) \bs \mu_n(\mbf{A}) \to \mr{S}^1$, fix embeddings $\iota_v \colon F_v \into \CC$ and write
\[
\omega_{\infty} \colon (x_v)_{v \vert \infty} \mapsto \iota_v(x_v)^{m_v}
\]
for some set of residue classes $m_v \in \Z/n\Z$. Then:
\begin{itemize} 
\item There exist type $A$ (unitary) Hecke character extensions $\tilde{\omega} \colon F^\times \bs \mbf{A}^\times \to \mr{S}^1$ if and only if the images $m_v \in \Z/n\Z$ depend only on the restriction $v|_{F_{cm}}$. There exists a finite-order extension $\tilde{\omega}$ if and only if for all complex $v \vert \infty$, the classes $m_v \in \Z/n\Z$ are all zero.
\item The extension $\tomega$ is unique up to the $n^{th}$ power of another type $A$ Hecke character, except in the Grunwald-Wang special case, where one first has two choices of extension to $C_F[n]$ (then the $n^{th}$ power ambiguity).
\item In particular, if $F$ is totally real or CM, then (unitary) type $A$ extensions of $\omega$ always exist. If $F$ is totally real, they are finite-order, and if $F$ is CM, a finite-order extension can be chosen if and only if $\omega|_{F_{\infty}^\times}$ is trivial.
\end{itemize}   
\end{lemma}
\proof
For the time being, let $F$ be arbitrary. The cokernel of $\mu_n(F) \bs \mu_n(\mbf{A}) \to C_F[n]$ is either trivial or order $2$ (the Grunwald-Wang special case), so we first choose an extension $\omega'$ to $C_F[n]$ (we will use this flexibility in Proposition \ref{llgeneralization}. Then Pontryagin duality gives an isomorphism
\[
{C_F^D}/{nC_F^D} \xrightarrow{\sim} (C_F[n])^D,
\]
so there exists some $\tomega \in C_F^D$ extending $\omega'$, and it is unique up to $n^{th}$ powers. Restricted to $F_\infty^\times$, we can write (implicitly invoking $\iota_v$)
\[
\tomega_\infty \colon (x_v)_{v \vert \infty} \to \prod_{v \vert \infty} (x_v/|x_v|)^{m_v} |x_v|_\CC^{i t_v}
\]
for integers $m_v$ and real numbers $t_v$. Then a type $A$ lift of $\omega'$ exists if and only if there is a Hecke character $\psi$ with infinity-type
\[
\psi_\infty \colon (x_v) \mapsto \prod_{v \vert \infty} (x_v/|x_v|)^{f_v} |x_v|_\CC^{i t_v/n}
\]
for some integers $f_v$ ($\tomega \psi^{-n}$ is then the desired type $A$ character). Following Weil, the existence of a Hecke character $\tomega$ corresponding to the infinity-data $\{m_v, t_v\}_v$ is equivalent to the existence of some integer $M$ such that for all $\alpha \in \mc{O}_F^\times$,
\[
\prod_{v \vert \infty}(\alpha/|\alpha|)^{m_v M} |\alpha|_v^{i t_vM}=1,
\]
and, similarly, such a $\psi$ can exist if and only if there exists $M' \in \Z$ such that
\[
\prod_{v \vert \infty}(\alpha/|\alpha|)^{f_v MM'n} |\alpha|_v^{i t_vMM'}=1
\]
for all $\alpha \in \mc{O}_F^\times$. Substituting, the left-hand side is $\prod_{v \vert \infty} (\alpha/|\alpha|)^{(f_v n-m_v)MM'}$, which can equal $1$ for all $\alpha \in \mc{O}_F^\times$ if and only if there exists a type $A$ Hecke character with infinity-type
\[
(x_v) \mapsto \prod_{v \vert \infty} (x_v/|x_v|)^{f_v n-m_v}.
\]
By Weil's classification of type $A$ Hecke characters, this is possible if and only if $f_v n-m_v$ depends only on $v|_{F_{cm}}$. This in turn is possible (for some choice of $f_v$) if and only if $m_v \in \Z/n\Z$ depends only on $v|_{F_{cm}}$. In particular, over CM and totally real fields $F$, there is no obstruction, so a type $A$ lift of $\omega$ always exists.

The claim about existence of finite-order extensions is a simple variant: if all $m_v$ ($v$ complex) are zero in $\Z/n\Z$, then there is a type $A$ lift $\tilde{\omega}$ with infinity-type
\[
\tomega_\infty \colon (x_v)_{v \vert \infty} \to \prod_{v \vert \infty} (x_v/|x_v|)^{m_v}
\]
for some integers $m_v$, all divisible by $n$ and only depending on $\iota_v|_{F_{cm}}$. But then there is also a type $A$ Hecke character $\psi$ of $F$ with infinity-type
\[
\psi_{\infty} \colon (x_v)_{v \vert \infty} \to \prod_{v \vert \infty} (x_v/|x_v|)^{m_v/n},
\]
so $\tilde{\omega} \psi^{-n}$ is a finite-order extension of $\omega$.

As for uniqueness, once a type $A$ lift of $\omega'$ is chosen, any other is (by the discussion at the start of the proof) a twist by the $n^{th}$ power of another Hecke character, which must itself clearly be of type $A$. This is then the ambiguity in extending $\omega$ itself, except in the Grunwald-Wang special case, in which there is the additional $\Z/2\Z$ ambiguity noted above.
\endproof
\begin{rmk}
\begin{itemize}
\item Of course any quasi-character of $\mu_n(\mbf{A})$ is unitary, and to understand all extensions it suffices (twisting by powers of $|\cdot|_{\mbf{A}}$) to understand unitary extensions-- hence the restriction to $\mr{S}^1$ instead of $\CC^\times$.
\item This raises a tantalizing question: certainly over a non-CM field $F$ we can produce characters $\omega$ that have no type $A$ extensions (although they certainly all have some extension to a Hecke character, by Pontryagin duality). What if $\omega$ actually arises as the central character of a suitably algebraic cuspidal automorphic representation (on $\mr{SL}_n(\mbf{A})$, for instance)? We return to this in \S \ref{cmdescentsection}.
\end{itemize}
\end{rmk}
A few consequences, either of the result or the method of proof, will follow; first, though, let us make a definition:
\begin{defn}
Let $F$ be any number field, and let $\psi \in C_F^D$ be a unitary Hecke character. We say that $\psi$ is \textit{of Maass type} if for all $v \vert \infty$, the restriction $\psi_v \colon F_v^\times \to \mr{S}^1$ has the form $x \mapsto |x|^{i t_v}$ for real numbers $t_v$.
\end{defn}
\begin{cor}\label{HCstructurethm} 
Let $\psi \colon C_F \to \CC^\times$ be a Hecke character of a number field $F$.
\begin{itemize}
\item If $F$ is totally real or CM, then $\psi$ can be decomposed as 
\[
\psi= \psi_{alg} \psi_{Maass} \psi_{fin} |\cdot|^w,
\]
where $w$ is the unique real number twisting $\psi$ to a unitary character, $\psi_{alg}$ is unitary type $A$, $\psi_{Maass}$ is of Maass type, and $\psi_{fin}$ is finite order. The last three characters are all unique up to finite-order characters.
\item If $\psi$ is of Maass type (after twisting to a unitary character), with $F$ arbitrary, then $\psi$ is `nearly divisible': for any $n \in \Z$ there exist Hecke characters $\chi$ of Maass type and $\chi_0$ of finite order such that $\chi^n \chi_0= \psi$. In particular, after a finite base-change $\psi$ is $n$-divisible.
\item Write $\mc{A}_F(1)$ for the space (topological group) of all Hecke characters of $F$. Suppose $F$ is CM, of degree $2s$ over $\Q$. Then there is an exact sequence
\[
1 \to \gal{F}^D \to \mc{A}_F(1) \to \RR \times \Z^s \times \Q^{s-1} \to 1.
\]
\end{itemize}
(As always, $\gal{F}^D$ denotes $\Hom_{cts}(\gal{F}, \mr{S}^1)$; it is the space of Dirichlet characters.) For the totally real subfield $F_+$ of $F$, we have a similar sequence
\[
1 \to \gal{F_+}^D \to \mc{A}_{F_+}(1) \to \RR \times \Q^{s-1} \to 1.
\]
\end{cor}
It is also possible (via the unit theorem) to compute all possible infinity-types of Hecke characters of any number field, but we will make no use of this calculation.\index{t}{infinity-types of Hecke characters, all possible} Roughly speaking, the new transcendentals that arise in the `mixed' case where algebraic and Maass parts cannot be separated are the arguments (angles) of fundamental units. It is very tempting to ask whether for CM (or totally real fields) algebraic and spherical (`Maass') infinity-types on higher-rank groups can twist together in a non-trivial way. If this question is too na\"{i}ve, is there any other higher-rank generalization of part $1$ of Corollary \ref{HCstructurethm}?

We include a couple other related useful results. For the first, note that in contrast with the situation for $\ell$-adic Galois characters (see Lemma \ref{Galoisroots} below), a general Hecke character cannot be written as an $n^{th}$-power up to a finite-order twist.\footnote{Consider, for example, a type $A$ Hecke character whose integral parameters at infinity are not divisible by $n$.}
\begin{lemma}\label{Heckeroots}
Let $\psi$ be a Hecke character of $F$, and suppose that $L/F$ is a finite Galois extension over which $\BC_{L/F}(\psi)= \psi \circ N_{L/F}$ is an $n^{th}$-power. Then up to a finite-order twist, $\psi$ itself is an $n^{th}$-power.
\end{lemma}
\proof
Let $\omega$ be a Hecke character of $L$ such that $\psi \circ N_{L/F}= \omega^n$. Then $\omega$ is $\Gal(L/F)$-invariant, and the next lemma shows that $\omega$ descends to $F$ up to a finite-order twist
\endproof
The next lemma completes the previous one, and also enables a slight refinement of a result of Rajan (see below):
\begin{lemma}\label{invariantHecke}
Let $L/F$ be a Galois extension of number fields, and let $\psi$ be a $\Gal(L/F)$-invariant Hecke character of $L$. Then there exists a Hecke character $\psi_F$ of $F$ and a finite-order Hecke character $\psi_0$ of $L$ such that $\psi= (\psi_F \circ N_{L/F}) \cdot \psi_0$.
\end{lemma}
\proof
We may assume $\psi$ is unitary, and we may choose, for all infinite places $w$ of $L$, embeddings $\iota_w \colon L_w \to \CC$ such that all $\iota_w$ for $w$ above a fixed place $v$ of $F$ restrict to the same embedding $\iota_v \colon F_v \to \CC$. $\Gal(L/F)$ acts transitively on the places $w \vert v$, so when we write
\[
\psi_w(x_w)= \left(\frac{\iota_w(x_w)}{|\iota_w(x_w)|} \right)^{m_w}\cdot |\iota_w(x_w)|^{i t_w},
\]
$\Gal(L/F)$-invariance implies that the $m_w$ and $t_w$ depend only on the place $v$ below $w$ (from now on in the proof, the embeddings $\iota_w$ will be implicit). We therefore denote these by $m_v$ and $t_v$. Lemma \ref{infinitycriterion} implies there is a Hecke character of $F$ with infinity-type given by the data $\{m_v, t_v\}_v$. Namely, there is an integer $M$ such that for all $\alpha \in \mc{O}_L^\times$,
\[
1= \left( \prod_{w \vert \infty} \left( \frac{\alpha}{|\alpha|} \right)^{m_w} |\alpha|_{\CC}^{i t_w} \right)^M.
\]
Restricting to $\alpha \in \mc{O}_F^\times$, this becomes
\[
1= \left( \prod_{v \vert \infty} \prod_{w \vert v}  \left( \frac{\alpha}{|\alpha|} \right)^{m_v} |\alpha|_{\CC}^{i t_v} \right)^M= \left( \prod_{v \vert \infty} \left( \frac{\alpha}{|\alpha|} \right)^{m_v} |\alpha|_{\CC}^{i t_v} \right)^{M \# \{w \vert v\}},
\]
($\#\{w \vert v\}$ is independent of $v$) which is simply the criterion for there to be a Hecke character of $F$ with infinity-type given by the collection $\{m_v, t_v\}_v$. Any such character has base-change differing from $\psi$ by a finite-order character, so we are done.
\endproof
\begin{rmk}
The aforementioned theorem of Rajan describes the image of solvable base-change\index{t}{image of a functorial transfer}{solvable base-change}: precisely, for a solvable extension $L/F$ of number fields and a $\Gal(L/F)$-invariant cuspidal automorphic representation $\pi$ of $\mr{GL}_n(\mathbf{A}_L)$, Theorem $1$ of \cite{rajan:solvable} asserts that there is a cuspidal representation $\pi_F$ of $\mr{GL}_n(\af)$ and a $\Gal(L/F)$-invariant Hecke character $\psi$ of $L$ such that
\[
\BC_{L/F}(\pi_0) \otimes \psi = \pi.
\]
Lemma \ref{invariantHecke} shows that this $\psi$ may be chosen finite-order; in particular, for some cyclic extension $L'/L$ (so $L'/F$ is still solvable), $\BC_{L'/L}(\pi)$ descends to $F$.
\end{rmk}
\subsection{Galois $\mr{GL}_1$}
We now discuss (continuous) Galois characters $\hat{\psi} \colon \gal{F} \to \Qlb^\times$. For the time being, $F$ is any number field. These are necessarily almost everywhere unramified, and we focus on $\ell$-adic Hodge theory aspects. The following is well-known, and is proven in \cite[Chapter III]{serre:ladic}:
\begin{thm}\label{FMLGL1}
Suppose that $\hat{\psi}$ is de Rham. Then for all places $v$ not dividing $\ell$, $\hat{\psi}|_{\gal{F_v}}$ assumes algebraic values in $\Qb \xrightarrow{\iota_{\ell}} \Qlb$, and there exists a type $A_0$ (i.e. $L=C$-algebraic) Hecke character $\psi \colon \mbf{A}_F/F^\times \to \CC^\times$ corresponding to $\hat{\psi}$ (via $\iota_\infty \colon \Qb \into \CC$). Moreover, $\hat{\psi}$ is motivic: the Fontaine-Mazur conjecture holds for $\mr{GL}_1/F$.
\end{thm}
\proof
We indicate how this follows from the arguments of \cite[III]{serre:ladic}; see Remark \ref{FMGL1} for the explanation of a simple case. Since $\hat{\psi}$ is Hodge-Tate at all $v \vert \ell$, it is (by a theorem of Tate: see \cite[III-A6]{serre:ladic}) locally algebraic in the sense of \cite[III-1.1]{serre:ladic}. Then the argument of \cite[III.2.3 Theorem 2]{serre:ladic} implies that, for a suitable modulus $\mathfrak{m}$, $\hat{\psi}$ is the $\ell$-adic Galois representation associated to an algebraic homomorphism $\mc{S}_{\mf{m}, \Qlb} \to \mathbb{G}_{m, \Qlb}$, where $\mc{S}_{\mf{m}}$ is the $\Q$-torus of \cite[II-2.2]{serre:ladic}. Up to isomorphism, this algebraic representation can be realized over some finite extension $E$ of $\Q$ inside $\Qlb$. Taking the `archimedean realization' as in \cite[II-2.7]{serre:ladic}, we obtain the desired type $A_0$ Hecke character. It then follows by Lemma \ref{weilalg} that $\hat{\psi}|_{\gal{F_v}}$ assumes algebraic values for $v \nmid \ell$. The Fontaine-Mazur conjecture for $\hat{\psi}$ follows from a theorem of Deligne (\cite[Proposition IV.D.1]{DMOS}), which shows the somewhat stronger statement that $\hat{\psi}$ is the $\ell$-adic realization of a motive for absolute Hodge cycles.
\endproof
\begin{rmk}\label{FMGL1}
Consider the simple case in which $F= \Q$. Then $\hat{\psi}|_{\gal{\Q_{\ell}}}$ has a single labeled Hodge-Tate weight; it therefore has the same (labeled) Hodge-Tate weights as an integer power $\omega_{\ell}^r$ of the $\ell$-adic cyclotomic character. It follows then from the above-cited theorem of Tate (\cite[III-A6]{serre:ladic}) and global class field theory that $\hat{\psi}\omega_{\ell}^{-r}$ is finite-order, hence is the $\ell$-adic realization (via a fixed isomorphism $\CC \xrightarrow{\sim} \Qlb$) of some finite-order Hecke character $\chi$. The Hecke character corresponding to $\hat{\psi}$ is then $\chi |\cdot|_{\mathbf{A}_{\Q}}^r$. It is also easy to show that $\hat{\psi}$ is motivic: $\omega_{\ell}^r$ is the $\ell$-adic realization of the Tate motive $\Q(r)$, and if $\hat{\psi}\omega_{\ell}^{-r}$ factors through a finite quotient $\Gal(F/\Q)$, then it is a suitable direct factor of $H^0(X_{\overline{\Q}}, \Qlb)$, for the zero-dimensional variety $X= \Res_{F/\Q}(\Spec F)$ whose $\Qb$-points are naturally indexed by embeddings $\tau \colon F \into \Qb$, with the $\gal{\Q}$-action on $H^0$ arising from permutation of the set of embeddings (for more details, see the discussion of Artin motives in \S \ref{homologicalmotives}). 

More generally, the idea behind establishing the Fontaine-Mazur conjecture in the abelian case is to find a sub-quotient of the cohomology of a CM abelian variety having `labeled Hodge numbers' matching those of $\hat{\psi}$; that this is possible essentially follows from the constraints on the $\infty$-type of the type $A_0$ Hecke character underlying $\hat{\psi}$ (as in Lemma \ref{weildescent}. If the CM abelian variety were actually defined over $F$, this would realize $\hat{\psi}$ inside some finite-order twist of its cohomology. The subtle part of the theorem is the need to have some control over the field of definition of the CM abelian variety. 
\end{rmk}
We make repeated use of the following well-known observation (see for instance Lemma $3.1$ of \cite{conrad:dualGW}). 
\begin{lemma}\label{Galoisroots}
Let $\chi \colon \gal{F} \to \Qlb^\times$ be an $\ell$-adic Galois character. For any non-zero integer $m$, there are characters $\chi_1, \chi_0 \colon \gal{F} \to \Qlb^\times$, with $\chi_0$ finite-order, such that $\chi= (\chi_1)^m \chi_0$.
\end{lemma}
As on the automorphic side, we want to consider Galois characters somewhat outside the algebraic range. As is customary, we write $\mr{HT}_{\tau}(\rho)$\index{s}{$\mr{HT}_{\tau}(\rho)$}\index{t}{labeled Hodge-Tate weights} to indicate the $\tau \colon F \into \Qlb$-labeled Hodge-Tate weights of an $\ell$-adic Galois representation $\rho$; see \S \ref{l-adichodge} for details, as well as for what we mean when we write non-integral $\tau$-labeled Hodge-Tate-Sen weights.
\begin{cor}\label{galchar}
Let $F$ be a totally imaginary field, and $n$ a non-zero integer. For all $\tau \colon F \into \Qlb$, fix an integer $k_{\tau}$. Then there exists a Galois character $\hat{\psi} \colon \gal{F} \to \Qlb^\times$ with $\mr{HT}_{\tau}(\hat{\psi})= \frac{k_\tau}{n}$ (respectively, with $\mr{HT}_{\tau}(\hat{\psi}) \equiv \frac{k_\tau}{n} \mod \Z$) if and only if 
\begin{enumerate}
\item $k_\tau$ depends only on $\tau_0:= \tau|_{F_{cm}}$ (respectively, only depends modulo $n$ on $\tau_0$);
\item and there exists an integer $w$ such that $k_{\tau_0}+ k_{\tau_0 \circ c}= w$ (respectively, $k_{\tau_0}+ k_{\tau_0 \circ c} \equiv w \mod n$) for all $\tau$, and $c$ the complex conjugation on $F_{cm}$.\footnote{If $F$ is Galois, we can rephrase this as $k_\tau+ k_{\tau \circ c}=w$ for all choices of complex conjugation $c$ on $F$.}
\end{enumerate}
\end{cor}
\proof
Suppose such a $\hat{\psi}$ exists, with weights $\frac{k_\tau}{n}$. $\hat{\psi}^n$ is geometric, and so (via $\iota_\ell$ and $\iota_{\infty}$) there exists a type $A_0$ Hecke character $\psi$ of $F$ corresponding to $\hat{\psi}^n$. For $v \vert \infty$ and $\iota_v \colon F_v \xrightarrow{\sim} \CC$, $\psi|_{F_v^\times}$ is given by
\[
\psi_v(x_v)= \iota_v(x_v)^{k_{\tau^*(\iota_v)}} \bar{\iota}_v(x_v)^{k_{\tau^*(\bar{\iota}_v)}},
\] 
where $\tau^*(\iota_v)= \tau_{\ell, \infty}^*(\iota_v)$ denotes the embedding $F \into \Qlb$ induced by $\iota_v$, $\iota_\ell$, and $\iota_\infty$. By purity for Hecke characters, there is an integer $w$ such that $k_{\tau^*(\iota_v)}+ k_{\tau^*(\bar{\iota}_v)}= w$; moreover, by Weil's descent result, $k_{\tau^*(\iota_v)}$ depends only on $\iota_v|_{F_{cm}}$. The constraint on the weights follows. 

Conversely, given a set of weights satisfying the purity constraint, we form a putative infinity-type $p_{\iota_v}= k_{\tau^*(\iota_v)}$ for a Hecke character of $F$; that this is in fact an achievable infinity-type follows from Weil. We form the associated geometric Galois character and then use the fact that, up to a finite-order twist, we can always extract $n^{th}$ roots of $\ell$-adic characters.

The $\mod n$ statement is a simple modification. For instance, to construct a character with given weights satisfying the purity constraint modulo $n$, proceed as follows. Choose a maximal set modulo $c$ of embeddings $\tau_0 \colon F_{cm} \into \CC$. Choose lifts to $\Q$ of the congruence classes $\frac{k_{\tau_0}}{n} \mod \Z$, and declare $k_\tau= k_{\tau_0}$ for all $\tau$ lying above these $\tau_0$. Then choose $w \in \Z$ lifting $k_{\tau_0}+ k_{\tau_0 \circ c}  \pmod n$, which is possible by hypothesis, and set $k_{\tau_0 \circ c}= w- k_{\tau_0}$.
\endproof
The same technique yields an easy example of the constraints on Galois characters over totally real fields:
\begin{lemma}\label{galtotreal}
Suppose $F$ is totally real, and $\hat{\psi} \colon \gal{F} \to \Qlb^\times$ is a character with all Hodge-Tate-Sen weights in $\Q$. Then all of these weights are equal. Conversely, for $x, d \in \Z$, there are global characters with all HTS weights equal to $\frac{x}{d}$.
\end{lemma}
\proof
For some integer $d$, all the HTS weights lie in $\frac{1}{d} \Z$, so $\hat{\psi}^d$ is geometric, and therefore corresponds to a type $A_0$ Hecke character. $F$ is totally real, so $\hat{\psi}^d$ must be, up to a finite-order twist, an integer power of the cyclotomic character.
\endproof
It is worth remarking that the most general answer to the question `what Galois characters $\hat{\psi} \colon \gal{F} \to \Qlb^\times$ exist' is essentially Leopoldt's conjecture.

Finally, we will need a lemma refining the construction of certain Galois characters over CM fields:
\begin{lemma}\label{galsqrt}
Let $F$ be a totally real field in which a prime $\ell \neq 2$ is unramified. Let $L=KF$ be its composite with an imaginary quadratic field $K$ in which $\ell$ is inert, and let $\psi$ be any unitary type $A$ Hecke character of $L$. Then:
\begin{enumerate}
\item There exists a Galois character $\hat{\psi} \colon \gal{L} \to \Qlb^\times$ such that $\hat{\psi}^2$ corresponds to $\psi^2$ (which is type $A_0$). 
\item Moreover, the frobenius eigenvalues $\hat{\psi}(fr_v)$ at all unramified $v$ lie in $\Q^{cm}$. 
\end{enumerate}
\end{lemma}
\proof
First, let $m_{\iota_v}$ as before denote the integers giving the infinity-type of $\psi$. 
Using the known algebraicity of $\psi$ (\cite{weil:characters}) and the fixed embeddings $\iota_\infty, \iota_\ell$, we can define the $\ell$-adic representation associated to $\psi^2$:
 \begin{align*}
(\widehat{\psi^2})_\ell \circ \mathrm{rec}_L^{-1} \colon &\mbf{A}_L^\times/(\overline{L^\times L_{\infty}^\times}) \to \Qlb \\
&(x_w) \mapsto \prod_{w \nmid \ell \infty} \psi_w(x_w)^2 \prod_{w \vert \ell} \left( \psi_w(x_w)^2 \prod_{\tau \colon L_w \into \Qlb} \tau(x_w)^{m_{\iota_{\infty, \ell}^*(\tau)}} \right).
\end{align*}
By the dual Grunwald-Wang theorem, to show that $(\widehat{\psi^2})_\ell$ is the square of some character $\gal{L} \to \Qlb^\times$, it suffices to check that for all places $w$ of $L$, the above character is locally on $L_w^\times$ a square. This in turn immediately reduces to seeing whether, for each $w \vert \ell$, the character $\chi_w \colon L_w^\times \to \Qlb^\times$ given by
\[
\chi_w(x_w)= \prod_{\tau \colon L_w \into \Qlb} \tau(x_w)^{m_{\iota_{\infty, \ell}^*(\tau)}} 
\]
is a square. Writing 
\[
L_w^\times= \langle \ell \rangle \times \mu_\infty(L_w^\times) \times (1+ \ell \mc{O}_w),
\]
it suffices to check on each component of this factorization. On the $\langle \ell \rangle$ factor, this is clear (choose a square root of $\chi_w(\ell)$). On the $1+ \ell \mc{O}_w$ factor, the $\ell$-adic logarithm, using our hypotheses that $\ell \neq 2$ is unramified, lets us define a single-valued square-root function, and thus extract a square root of $\chi_w|_{1+\ell \mc{O}_w}$. Now note that the product over $\tau \colon L_w \into \Qlb$ is a product over pairs of complex-conjugate embeddings extending a given $F_v \into \Qlb$, and $m_{\iota^*(\bar{\tau})}=m_{\overline{\iota^*(\tau)}}= -m_{\iota^*(\tau)}$ ($\psi$ is unitary).  $\mu_{\infty}(L_w^\times)$ is isomorphic to $\Z/(q-1)\Z$, where $q$ is the order of the residue field at $w$. Let $\zeta \mapsto 1$ give the isomorphism, for $\zeta$ a primitive $(q-1)^{st}$ root of unity. Complex conjugation must identify to multiplication by an element $r \in \Z/(q-1)$ satisfying $r^2 \equiv 1 \mod q-1$; in particular, $r$ is an odd residue class. Then $x \mapsto \tau(x)^{m_{\iota^*(\tau)}} \cdot (\tau \circ c) (x)^{-m_{\iota^*(\tau)}}$ takes $\zeta$ to an even power of $\tau(\zeta)$, hence $\chi_w|_{\mu_{\infty}(L_w^\times)}$ is a square.\footnote{Note also that if $\ell$ is split in $L/\Q$, $\chi_w$ cannot be a square when $m_{\iota_{\infty, \ell}^*(\tau)}$ is odd.} 

The second part of the lemma follows from the construction of $\hat{\psi}$ and the corresponding statement (\cite{weil:characters}) for the type $A$ Hecke character $\psi$. 
\endproof

\section{Coefficients: generalizing Weil's CM descent of type $A$ Hecke characters}\label{cmdescentsection}
In the coming sections we pursue higher-rank analogues of two aspects of Weil's paper \cite{weil:characters}. This section extends the important observation that type $A$ Hecke characters of a number field $F$ descend, up to a finite-order twist, to the maximal CM subfield $F_{cm}$ (see Lemma \ref{weildescent}). Of interest in its own right, this generalization also provides some of the intuition necessary for a general solution to Conrad's lifting question (Question \ref{conradquestion}). 

\subsection{Coefficients in Hodge theory}\label{formal}
The guiding principle that allows us to reinterpret, and correspondingly generalize, Weil's result is that careful attention to the `field of coefficients' of an arithmetic object can yield non-trivial information about its `field of definition,' or that of certain of its invariants.  We begin by recording in an abstract setting a lemma whose motivation is `doing Hodge theory with coefficients.' Let $k$ be a field of characteristic zero, and let $F$ and $E$ be extensions of $k$ with $F/k$ finite. Let $D$ be a filtered $F \otimes_k E$-module that is free of rank $d$. Let $E'$ be a Galois extension of $E$ large enough to split $F$ over $k$. For all ($k$-embeddings) $\tau \colon F \into E'$, we define as in Definition \ref{labeledHT} the $\tau$-labeled Hodge-Tate weights $\mr{HT}_{\tau}(D)$\index{s}{$\mr{HT}_{\tau}$}\index{t}{labeled Hodge-Tate weights} as follows: for such $\tau$ we have orthogonal idempotents $e_\tau \in F \otimes_k E'$ giving rise to projections
\begin{align*}
&F \otimes_k E' \xrightarrow[\sim]{(e_{\tau})} \prod_\tau E'\\
& x \otimes \alpha \mapsto \left( \tau(x) \alpha\right)_\tau.
\end{align*}
The projection $e_\tau(D \otimes_E E')$ is then a filtered $E'$-vector space of dimension $d$, and we define $\mr{HT}_{\tau}(D)$ to be the collection of integers $h$ (with multiplicity) such that 
\[
\mathrm{gr}^h(e_\tau(D \otimes_E E')) \neq 0.
\]
In \S \ref{l-adicbasics} we applied this formalism to the filtered $K \otimes_{\Ql} \Qlb$-module $\Ddr(V)$, when $V$ was a representation of $\gal{K}$ ($K/\Ql$ finite) on a $\Qlb$-vector space. In this section, however, we emphasize the general formalism: one should really keep in mind not $\ell$-adic Hodge theory, but rather the de Rham realization of a motive over $F$ with coefficients in $E$.
\begin{lemma}\label{hodgetheory}
$\mr{HT}_\tau(D)$ depends only on the $\Gal(E'/E)$-orbit of $\tau \colon F \into E'$. If $E/k$ is Galois, then $\mr{HT}_\tau(D)$ depends only on the restriction $\tau|_{\tau^{-1}(E)}$.
\end{lemma}
\proof
We decompose $F \otimes_k E$ into a product of fields $\prod E_i$, writing $q_i \colon F \otimes_k E \onto E_i$ for the quotient map. This yields filtered $E_i$-vector spaces $D_i$ for all $i$. Any $E$-algebra homomorphism $\tau \colon F \otimes_k E \to E'$ factors through $q_{i(\tau)}$ for a unique index $i(\tau)$, and then $\mr{HT}_{\tau}(D)$ is simply the multi-set of weights of $D_{i(\tau)}$. This implies the first claim, since the $\Gal(E'/E)$-orbit of $\tau$ is simply all embeddings $\tau'$ for which $i(\tau')= i(\tau)$.

Now we assume $E/k$ Galois and address the second claim. Having fixed a `reference' embedding $\tau_0 \colon F \into E'$, it makes sense to speak of $F \cap E:= F \cap \tau_0^{-1}(E)$ inside $E$. Obviously $E$ splits $F \cap E$ over $k$, so we can write
\[
F \otimes_k E \cong F \otimes_{F \cap E}(F \cap E \otimes_k E) \cong \prod_{\sigma} F \otimes_{F \cap E, \sigma} E,
\]
where the $\sigma$ range over all embeddings $F \cap E \into E$. These factors $F \otimes_{F \cap E, \sigma} E$ are themselves fields, since $E/k$ is Galois, and so this decomposition realizes explicitly the decomposition of $F \otimes_k E$ into a product of fields. In particular, by the first part of the Lemma, $\mr{HT}_{\tau}(D)$ depends only on $i(\tau)$, i.e. only on the $\sigma \colon F \cap E \into E$ to which $\tau$ restricts. 
\endproof
\begin{defn}\label{regular}
We call a filtered $F \otimes_k E$-module $D$ \textit{regular} if the multi-sets $\mr{HT}_{\tau}(D)$ are multiplicity-free.\index{t}{regular filtered module} 
\end{defn}
A very simple application that we will use later is:
\begin{cor}\label{notregular}
Let $D$ be a filtered $L \otimes_k E$-module (globally free) for some finite extension $L/F$, and suppose $L$ does not embed in $\tilde{E}$. Then the restriction of scalars (image under the forgetful functor) $\Res_{L/F}(D)$ is not regular.
\end{cor}

\subsection{CM descent}
We will see that Weil's result (Lemma \ref{weildescent}) is the conjunction of Lemma \ref{hodgetheory} with the fact that algebraic Hecke characters have CM `fields of coefficients.' This observation will lead us naturally to the desired higher-rank generalization.

The following notation will be in effect for the rest of the paper. Let $G$ be a connected reductive $F$-group. For each $v \vert \infty$ fix an isomorphism $\iota_v \colon \overline{F}_v \xrightarrow{\sim} \CC$. For $\pi$ an automorphic representation of $G(\af)$, we can write (in Langlands' normalization) the restriction to $W_{\overline{F}_v}$ of its $L$-parameter as
\[
\rec_{v}(\pi_v) \colon z \mapsto \iota_v(z)^{\mu_{\iota_v}} \bar{\iota}_v(z)^{\nu_{\iota_v}} \in T^\vee(\CC).
\] 
with $\mu_{\iota_v}, \nu_{\iota_v} \in X^\bullet(T)_{\CC}$ and $\mu_{\iota_v}-\nu_{\iota_v} \in X^\bullet(T)$. For $v$ imaginary, $\mu_{\bar{\iota}_v}= \nu_{\iota_v}$. Unless there is risk of confusion, we will omit reference to the embedding $\iota_v$, writing $\mu_v= \mu_{\iota_v}$, etc. Recall that, in the terminology of \cite{buzzard-gee:alg}, $\pi$ is $L$-algebraic if for all $v \vert \infty$, $\mu_v$ and $\nu_v$ lie in $X^\bullet(T)$; it is $C$-algebraic if $\mu_v$ and $\nu_v$ lie in $\rho + X^\bullet(T)$, where $\rho$ denotes the half-sum of the positive roots (with respect to our fixed Borel containing $T$ used to define a based root datum of $G$). 

Our starting point is a result (and, more generally, conjecture) of Clozel.  Take $G= \mr{GL}_n/F$ and, for simplicity, $F$ to be totally imaginary. We collect the data of $\pi$'s archimedean $L$-parameters as $M= \{\mu_\iota\}_{\iota}$, which we will loosely refer to as the `infinity-type' of $\pi$.\index{t}{infinity-type of an automorphic representation}  For $\sigma \in \Aut(\CC)$, define the action $^{\sigma} M= \{\mu_{\sigma^{-1} \iota}\}_{\iota}$. Recall that $\pi$ is said to be \textit{regular} if all roots of $\mr{GL}_n$ are non-vanishing on all of the co-characters $\mu_{\iota}$.\index{t}{regular automorphic representation for the group $\mr{GL}_n$}
\begin{thm}[Th\'{e}or\`{e}me $3.13$ of \cite{clozel:alg}]
Let $F$ be any number field, and suppose $\pi$ is a cuspidal, $C$- or $L$-algebraic automorphic representation of $\mr{GL}_n(\af)$ that is moreover regular. Then $\pi_f$ has a model over the fixed field\index{s}{$\Q(\pi_f)$} $\Q(\pi_f) \subset \Qb \subset \CC$ of all automorphisms $\sigma \in \Aut(\CC)$ such that ${}^{\sigma} \pi_f \cong \pi_f$, with $\Q(\pi_f)$ in fact a number field for $\pi$ C-algebraic. For each $\sigma \in \Aut(\CC)$ there is a cuspidal representation\index{s}{$^{\sigma} \pi$} $^{\sigma} \pi$ with finite part $^{\sigma} \pi_f= \pi_f \otimes_{\CC, \sigma} \CC$ and with infinity-type $^{\sigma} M$.
\end{thm}
More generally, Clozel conjectures this for any $C$-algebraic isobaric automorphic representation of $\mr{GL}_n(\af)$. Clozel's theorem and conjecture are stated for $C$-algebraic representations, but the $L$-algebraic analogue follows easily by twisting.
\begin{hypothesis}\label{algconj}
Throughout the rest of this section, we will assume that $\pi$ is an isobaric $C$- or $L$-algebraic automorphic representation of $\mr{GL}_n(\af)$ satisfying the conclusion of Clozel's theorem.
\end{hypothesis}
We first note an elementary but crucial refinement of Hypothesis \ref{algconj}:
\begin{cor}\label{automorphiccmcoefficients}
The field $\Q(\pi_f)$ is CM.
\end{cor}
\proof
It suffices to check for each factor in an isobaric decomposition, so we may assume $\pi$ is cuspidal. Consider a unitary twist $\pi |\cdot|^w$, with $w$ necessarily some half-integer $w$. Its field of definition is CM if and only if that of $\pi$ is. Therefore we may assume $\pi$ is unitary. The $L^2$ inner product then implies that $\pi^\vee \cong {}^c \pi$, which in turn implies that for all $\sigma \in \Aut(\CC)$, ${}^{c\sigma} \pi \cong {}^{\sigma c}\pi$, and therefore that the fixed field $\Q(\pi_f)$ of $\{\sigma \in \Aut(\CC): {}^{\sigma} \pi \cong \pi\}$ is CM.
\endproof
We can now formulate (and prove under Hypothesis \ref{algconj}) the appropriate higher-rank generalization of Weil's result that type $A$ Hecke characters of $F$, up to a finite-order twist, descend to $F_{cm}$.
\begin{prop-conj}\label{cmdescent}
Let $\pi$ be as in Hypothesis \ref{algconj}. Conjecturally, this should hold for all isobaric $C$- or $L$-algebraic automorphic representations of $\mr{GL}_n(\af)$. Then the infinity-type of $\pi$ descends to $F_{cm}$.  
\end{prop-conj}
\proof
Let $E$ denote the Galois closure (in $\CC$) of $\Q(\pi_f)$; $E$ is also a CM field. Extend each $\iota \colon F \into \CC$ to an embedding $\tilde{\iota} \colon \tilde{F} \into \CC$ of the Galois closure $\tilde{F}$. The image $\tilde{\iota}(\tilde{F}) \subset \CC$ does not depend on the extension. By the corollary, $\tilde{\iota}(\tilde{F})$ is linearly disjoint from $E$ over $\tilde{\iota}((\tilde{F})_{cm})$, and therefore we can find $\sigma \in \Aut(\CC/E)$ restricting (via $\tilde{\iota}$) to any element we like of $\Gal(\tilde{F}/(\tilde{F})_{cm})$; the collection of such $\sigma$ acts transitively on the set of embeddings $F \into \CC$ lying above a fixed $F \cap (\tilde{F})_{cm}= F_{cm} \into \CC$.\footnote{$G= \Gal(\tilde{F}/F_{cm})$ is generated by $H= \Gal(\tilde{F}/F)$ and $H'= \Gal(\tilde{F}/\tilde{F}_{cm})$, with $H'$ normal. The set $\Hom_{F_{cm}}(F, \CC)$ is permuted transitively by $G$, with $H$ acting trivially, so for any such embedding $x$, $H'x= H'Hx= Gx= \Hom_{F_{cm}}(F, \CC)$.} For any two such $\iota, \iota' \colon F \into \CC$ (related by $\sigma \in \Aut(\CC/E)$), we deduce, since ${}^{\sigma}\pi \cong \pi$ and hence ${}^\sigma M= M$, that $\mu_{\iota}= \mu_{\sigma^{-1} \iota}= \mu_{\iota'}$.   
\endproof
Any study of algebraic automorphic forms over non-CM fields will need to take this result into account.
\begin{rmk}\label{algconjrmk}
\begin{itemize}
\item If $\pi$ is regular, this yields an unconditional descent result for its infinity-type.
\item Let $F$ be a CM field and $\pi$ be a regular $L$ or $C$-algebraic cuspidal representation of $\mr{GL}_n(\af)$. Suppose that $\pi= \Ind_{L}^F(\pi_0)$ for some extension $L/F$. Then $L$ is CM. This modest consequence of the proposition suggests that in the study of regular automorphic representations/motives/Galois representations over CM fields, we will never have to grapple with the (less well-understood) situation over non-CM fields. In \S \ref{liemultfree} we discuss an abstract Galois-theoretic analogue.
\item If we know more about $\Q(\pi_f)$ than that it is CM (an extreme case: $\Q(\pi_f)=\Q$), the proof of the proposition yields a correspondingly stronger result. For the related formalism, see \S \ref{formal}.
\end{itemize}
\end{rmk}
\begin{rmk}\label{cmdescentextrapolation}
Let us also note that Proposition \ref{cmdescent} is a `seed' result (under Hypothesis \ref{algconj}); if we moreover assume functoriality (in a form that requires functorial transfers to be strong transfers-- i.e. of archimedean $L$-packets-- at infinity), then it implies:
\begin{itemize}
\item Let $G$ be a connected reductive group over a totally imaginary field $F$, and let $\pi$ be a cuspidal tempered $L$-algebraic automorphic representation of $G(\af)$. Then the infinity-type of $\pi$ descends to $F_{cm}$.
\item Let $G$ and $\pi$ be as above. Let $\tG \supset G$ be a connected reductive $F$-group obtained by enlarging the center of $G$ to a torus $\tZ$ (see \S \ref{liftingwalg}). Then the central character $\omega_\pi \colon Z_G(F) \bs Z_G(\af) \to \CC^\times$ extends to an $L$-algebraic Hecke character of $\tZ(\af)$. (Compare Proposition $13.3.1$.)
\end{itemize}
\end{rmk}
\section{W-algebraic representations}\label{walgsection}
The current section generalizes a second aspect of \cite{weil:characters}, discussing a higher-rank analogue of the type $A$, but not necessarily $A_0$, condition, and begins to motivate its arithmetic significance. As always, let $G$ be a connected reductive group over a number field $F$.\index{s}{$G$} We continue with the infinity-type notation of \S \ref{cmdescentsection}.

For many interesting questions about algebraicity of automorphic representations, and especially the interaction of algebraicity and functoriality, the framework of $C$ and $L$-algebraic representations does not suffice. Motivated initially by Weil's study of type-$A$ Hecke characters, we make the following definition:
\begin{defn}
Let $\pi$ be an automorphic representation of $G(\af)$. We say that $\pi$ is $W$-algebraic if for all $v \vert \infty$, $\mu_{\iota_v}$ and $\nu_{\iota_v}$ in fact lie in $\frac{1}{2}X^\bullet(T)$.\index{t}{$W$-algebraic automorphic representation} 
\end{defn}
\begin{eg}\label{mixedhecke}
For $\mr{GL}_1$, a unitary $W$-algebraic representation is precisely a unitary Hecke character of type $A$ in the sense of Weil (\cite{weil:characters}). Weil's type $A$ characters also include arbitrary twists $|\cdot|^r$ for $r \in \Q$, since these also yield $L$-series with algebraic coefficients. The $L$- and $C$-algebraic representations are the Hecke characters of type $A_0$. 
\end{eg}
\begin{eg}\label{mixed}
Consider a Hilbert modular form $f$ on $\mr{GL}_2{/F}$ with (classical) weights $\{k_\tau\}_{\tau \colon F \into \RR}$, where the positive integers $k_\tau$ are not all congruent modulo $2$; these are called `mixed-parity'.\index{t}{mixed-parity Hilbert modular form} $f$ then gives rise to an automorphic representation (in the unitary normalization, say) that is $W$-algebraic but neither $L$- nor $C$-algebraic. The previous example yields a special case: choose a quadratic CM (totally imaginary) extension $L/F$, and let $\psi$ be a unitary Hecke character of $L$ that is type $A$ but not type $A_0$. Then the automorphic induction $\Ind_L^F \psi$ yields an example on $\mr{GL}_2{/F}$. We will elaborate on the case of mixed-parity Hilbert modular forms in \S \ref{HMFs}.
\end{eg}
We are led to ask (compare \S \ref{cmdescentsection} and Conjectures 3.1.5 and 3.1.6 of \cite{buzzard-gee:alg}):
\begin{question}\label{Walg}
\begin{itemize}
\item Let $\pi$ be a $W$-algebraic automorphic representation on $G$. Does the $G(\mbf{A}_{F, f})$-module $\pi_f$ have a model over $\Qb$ (or $\Q^{cm}$)?. Note that we do not seek a model over a fixed number field. Alternatively, is $\pi_v$ defined over $\Qb$ for almost all finite places $v$? By analogy with the terminology of \cite{buzzard-gee:alg}, we might call the latter condition $W$-arithmetic, and ask whether some condition (which will not quite be $W$-algebraicity!) on infinity-types characterizes $W$-arithmeticity.\index{t}{$W$-arithmetic}
\item Similarly, we can ask the CM descent question of \S \ref{cmdescentsection} for $W$-algebraic representations. Example \ref{mixinginoneplace} below will be a cautionary tale: if non-induced cuspidal $\Pi$ with infinity-types as in that example exist, then there is no evidence that they would satisfy the analogue of Clozel's algebraicity conjecture. If they do not exist, or if they miraculously still satisfy the conclusion of Clozel's theorem, then we could confidently extend the CM descent conjectures to the $W$-algebraic case.
\end{itemize}
\end{question}

For tori, the most optimistic conjecture holds:
\begin{lemma}
For arbitrary $F$-tori, $W$-algebraic implies $W$-arithmetic.
\end{lemma}
\proof
After squaring, this follows from the corresponding result for $L$-algebraic/$L$-arithmetic (\S $4$ of \cite{buzzard-gee:alg}).
\endproof
Continuing with Example \ref{mixedhecke}, let us emphasize that the $W$-algebraic condition does not capture all automorphic representations with algebraic Satake parameters; we nevertheless want to make a case for isolating this condition, rather than allowing arbitrary rational parameters $\mu_v, \nu_v \in X^\bullet(T)_\Q$, which would be the na\"{i}ve analogue of type $A$. First, note that if a unitary $\pi$ is tempered at infinity with real parameters $\mu_v, \nu_v \in X^\bullet(T)_{\RR}$, then these parameters in fact lie in $\frac{1}{2}X^\bullet(T)$. In particular, the Ramanujan conjecture\index{t}{Ramanujan conjecture (archimedean)} implies that all cuspidal unitary $\pi$ on $\mr{GL}_n/F$ with real infinity-type are in fact $W$-algebraic.  It is easy to construct non-cuspidal automorphic representations with rational Satake parameters that do not twist to $W$-algebraic representations, but the point is that all such examples will be degenerate, so $W$-algebraicity is the condition of basic importance. A more interesting question, concerning the difference between $W$, $L$, and $C$-algebraicity, is the following:
\begin{eg}\label{mixinginoneplace}
Let $G= \mr{GL}_4/\Q$ (for example; there are obvious analogues for any totally real field). Let $F/\Q$ be real quadratic, and let $\pi$ as in Example \ref{mixed} be a mixed parity (unitary) Hilbert modular representation. It cannot be isomorphic to its $\Gal(F/\Q)$-conjugate, so $\Pi= \Ind_F^{\Q}(\pi)$ is cuspidal automorphic on $\mr{GL}_4/\Q$, and at infinity its Langlands parameter (restricted to $\CC^\times$) looks like
\[
z \mapsto 
\begin{pmatrix}
(z/ \bar{z})^{\frac{k_1-1}{2}} & 0 & 0 & 0 \\
0 & (z/ \bar{z})^{\frac{1-k_1}{2}} & 0 & 0 \\
0 & 0 & (z/ \bar{z})^{\frac{k_2-1}{2}} & 0 \\
0 & 0 & 0 & (z/ \bar{z})^{\frac{1-k_2}{2}}, 
\end{pmatrix}
\]
where $k_1$ and $k_2$ are the classical weights of $\pi$ at the two infinite places of $F$. This exhibits the `parity-mixing' within a single infinite place, which implies the following:
\begin{lemma}
For any functorial transfer $^{L}\mr{GL}_4 \xrightarrow{r} {}^{L}\mr{GL}_{N}$ arising from an irreducible representation $r$, not equal to any power of the determinant, of $\mr{GL}_4(\CC)$, the transfer $\mr{Lift}_{r}(\Pi)$ cannot be $L$-algebraic.
\end{lemma}
\proof
We have stated the lemma globally, and therefore conjecturally, but of course the analogous archimedean statement (which is well-defined since the local transfer of $\Pi_\infty$ via $r$ is known to exist) is all we are really interested in. The proof is elementary highest-weight theory.
\endproof
We will see (in \S \ref{HMFs}) that mixed-parity Hilbert modular representations are $W$-arithmetic, using the fact that they have $L$-algebraic functorial transfers; by contrast, if a representation $\Pi$ with the infinity-type of the above example is $W$-arithmetic, it may be harder to establish. Of course, in the above example, $\Pi$ is automorphically induced, so its $W$-arithmeticity is an immediate consequence of that of $\pi$. We are led to ask whether there exist non-automorphically induced examples of such $\Pi$, or more generally of cuspidal, non-induced $\Pi$ on $\mr{GL}_n/F$ that exhibit the `parity-mixing' within a single infinite place. The trace formula does not easily yield them, since such a $\Pi$ is not the transfer from a classical group of a form that is discrete series at infinity (in these cases there is a parity constraint on regular elliptic parameters).
\end{eg}

One further motivation for considering $W$-algebraic representations comes from studying the fibers of functorial lifts, and their algebraicity properties. For instance, given two mixed-parity Hilbert modular forms $\pi_1$ and $\pi_2$, with the (classical) weight of $\pi_{1, v}$ and $\pi_{2, v}$ having the same parity for each $v \vert \infty$, the tensor product $\Pi= \pi_1 \boxtimes \pi_2$ is $L$-algebraic. The $\pi_i$ are not themselves twists of $L$-algebraic automorphic representations, so $\Pi$ cannot be expressed as a tensor product of $L$-algebraic representations (this is proven in Corollary \ref{2,2fibers}). This is, however, the `farthest' from $L$-algebraic that the $\pi_i$'s can be:
\begin{prop}\label{tensordescent}
Let $\Pi$ be an $L$-algebraic cuspidal automorphic representation of $\mr{GL}_{nn'}(\af)$ for some number field $F$. Assume that $\Pi= \pi \boxtimes \pi'$ is in the image of $\mr{GL}_n \times \mr{GL}_{n'} \xrightarrow{\boxtimes} \mr{GL}_{nn'}$ for cuspidal automorphic representations $\pi$ and $\pi'$ of $\mr{GL}_n(\af)$ and $\mr{GL}_{n'}(\af)$.\footnote{To be precise, we want $\Pi$ to be a weak lift that is also a strong lift at archimedean places.} If $F$ is totally real or CM, then there are quasi-tempered (i.e., tempered up to a twist) $W$-algebraic automorphic representations $\pi$ of $\mr{GL}_n(\af)$ and $\pi'$ of $\mr{GL}_{n'}(\af)$ such that $\Pi= \pi \boxtimes \pi'$.  
\end{prop} 
\proof
The first part of the argument is similar to, and will make use of, Lemme $4.9$ of \cite{clozel:alg} (Clozel's `archimedean purity lemma'\index{t}{archimedean purity lemma of Clozel}), which shows that $\Pi_{\infty}$ itself is quasi-tempered. First assume $F$ is totally imaginary. Let $w_{\Pi}$ (the `motivic weight') be the integer such that $\Pi |\cdot|^{-w_{\Pi}/2}$ is unitary, and write $z^{\mu_v}\bar{z}^{\nu_v}$ and $z^{\mu'_v} \bar{z}^{\nu'_v}$ for the $L$-parameters at $v \vert \infty$ of $\pi$ and $\pi'$. Temperedness of $\Pi$ implies that for all $i, j= 1, \ldots, n$, 
\[
\mr{Re}(\mu_{v, i}+ \mu'_{v, j}+ \nu_{v, i}+ \nu'_{v, j})= w_{\Pi}.
\]
Fixing $j$ and varying $i$, we find that $\mr{Re}(\mu_{v, i}+ \nu_{v, i})$ is independent of $i$-- call it $w_{\pi_v}$. Similarly, for $w_{\pi'_v}= w_{\Pi}- w_{\pi_v}$, we have $\mr{Re}(\mu'_{v, i}+ \nu'_{v, i})= w_{\pi'_v}$ for all $i$. For each of $\pi$ and $\pi'$ there is a unique real number $r, r'$ such that each of $\pi |\cdot|^{-r/2}$, $\pi' |\cdot|^{-r'/2}$ is unitary. Thus, $z^{\mu_v- r/2} \bar{z}^{\nu_v-r/2}$ (likewise for $\pi'_v$) is a generic (by cuspidality) unitary parameter.\footnote{On $\mr{GL}_n$, $L$-packets are singletons, and we can identify which parameters correspond to generic representations.} This implies (see the proof of Lemme $4.9$ of \cite{clozel:alg}) that these parameters are sums of the parameters of \textit{unitary} characters and $2 \times 2$ complementary series blocks (what Clozel denotes $J(\chi, 1)$ and $J(\chi, \alpha, 1)$, where $\alpha \in (0, 1/2)$ is the complementary series parameter). Since $\mr{Re}(\mu_{v, i}+ \nu_{v, i})$ is independent of $i$, we can immediately rule out any complementary series factors, and we deduce that each $z^{\mu_v-r/2} \bar{z}^{\nu_v-r/2}$ is a direct sum of (parameters of) unitary characters, and that $w_{\pi_v}=r$ is independent of $v$ (and $w_{\pi'_v}= r'$), and $r+r'= w_{\Pi} \in \Z$. We may replace $\pi$ by $\pi |\cdot|^{\frac{r'-r}{2}}$ and $\pi'$ by $\pi' |\cdot|^{\frac{r-r'}{2}}$, so that still $\pi \boxtimes \pi' = \Pi$ but now each has half-integral `motivic weight' (namely, $\frac{1}{2}w_{\Pi} \in \frac{1}{2} \Z$). In particular, $\mr{Re}(\mu_v), \mr{Re}(\nu_v), \mr{Re}(\mu'_v)$, and $\mr{Re}(\nu'_v)$ all now consist of half-integers (it is easily seen that all of these half-integers are moreover congruent modulo $\Z$). Now, $\mr{Im}(\mu_{v, i}+ \mu'_{v, j})=0$ for all $i, j$, so there exists a $t_v \in \RR$ such that 
\[
\mr{Im}(\mu_{v, i})= \mr{Im}(\nu_{v, i})= -\mr{Im}(\mu'_{v, i})= - \mr{Im}(\nu'_{v, i})= t_v
\]
for all $i$. To conclude the proof, note that if there exists a Hecke character of $F$ with infinity-componenets $|z|_{\CC}^{i t_v}$ for all $v \vert \infty$, then we can twist $\pi$ and $\pi'$ to arrange that both be $W$-algebraic (in fact, either both $L$-algebraic or both $C$-algebraic). Looking at central characters, we are handed a Hecke character with infinity-type
\[
z \mapsto z^{\det(\mr{Re}(\mu_v))+ int_v} \bar{z}^{\det(\mr{Re}(\mu_v))+ int_v},
\]
and by Corollary \ref{HCstructurethm}, the desired Hecke character exists if and only if the set of half-integers $\det(\mr{Re}(\mu_v)):= \sum_i \mr{Re}(\mu_{v, i})$ is the infinity-type of a type $A$ Hecke character of $F$, i.e. descends to $F_{cm}$. Of course, when $F$ itself is CM, this is no obstruction, and the proof is complete. 
%Note also that if the fibers of the automorphic tensor product are (as for $\mr{GL}_2 \times \mr{GL}%_2 \to \mr{GL}_4$) just given by twisting, then this descent condition is actually equivalent to the %existence of the desired $W$-algebraic descents $\pi$ and $\pi'$. These aren't the fibers! Are %they if you restrict to the fiber over cuspidal $\Pi$???

Now assume $F$ is totally real. By a global base-change argument (as in Clozel's Lemme $4.9$), we may deduce temperedness of $\pi$ and $\pi'$ from the totally imaginary (or even CM) case, and to determine whether there are $W$-algebraic descents $\pi$ and $\pi'$, as usual it suffices to look at the archimedean $L$-parameters restricted to $W_{\overline{F}_v}$. The purity constraint  (see Corollary \ref{HCstructurethm}) on the half-integers $\det(\mr{Re}(\mu_v))$ forces these to be independent of $v$, and we can argue as in the CM case to finish the proof. 
\endproof
In the non-CM case, we can describe the infinity-types of possible Hecke characters; this gives rise to an explicit relation between the $\det(\mr{Re}(\mu_v))$ and $t_v$ that needs to be satisfied if $\pi$ and $\pi'$ are to have $W$-algebraic twists.
%Does CM descent for $\Pi$ help???
%\begin{question}
%Work out what you expect for multiple tensor products.
%\end{question}
Conversely, note that the tensor product transfer (only known to exist locally, of course) $\mr{GL}_n \times \mr{GL}_m \xrightarrow{\boxtimes} \mr{GL}_{nm}$ clearly preserves $L$- or $W$-algebraicity (but not $C$-algebraicity).

In \S \ref{liftingwalg}, we take up in detail the (algebraicity of the) fibers of the functorial transfer ${}^L\tG \to {}^LG$, where $G \subset \tG$ with torus quotient.

\section{Further examples: the Hilbert modular case and $\mr{GL}_2 \times \mr{GL}_2 \xrightarrow{\boxtimes} \mr{GL}_4$}\label{HMFs}
In this section we will elaborate on the example of mixed-parity Hilbert modular forms: we discuss $W$-arithmeticity in this context and make some initial forays on the Galois side. First, however, we recall known results on the automorphic tensor product $\mr{GL}_2 \times \mr{GL}_2 \xrightarrow{\boxtimes} \mr{GL}_4$ and provide a description of its fibers.

\subsection{General results on the $\mr{GL}_2 \times \mr{GL}_2 \xrightarrow{\boxtimes} \mr{GL}_4$ functorial transfer}\label{22tensor}
Ramakrishnan (\cite{ramakrishnan:tensor}) has proven the existence of the automorphic tensor product transfer in the case
$\mr{GL}_2 \times \mr{GL}_2 \xrightarrow{\boxtimes} \mr{GL}_4$.  Moreover, he establishes a cuspidality criterion:
\begin{thm}[Ramakrishnan, {\cite[Theorem M]{ramakrishnan:tensor}}]
Let $\pi_1$ and $\pi_2$ be cuspidal automorphic representations on $\mr{GL}_2/F$. Then $\pi_1 \boxtimes \pi_2$ is automorphic. The cuspidality criterion divides into two cases:
\begin{itemize}
\item
If neither $\pi_i$ is an automorphic induction, $\pi_1 \boxtimes \pi_2$ is cuspidal if and only if $\pi_1$ is not equivalent to $\pi_2 \otimes \chi$ for some Hecke character $\chi$ of $F$. 
\item If $\pi_1= \Ind_L^F(\psi)$ for a Hecke character $\psi$ of a quadratic extension $L/F$, then $\pi_1 \boxtimes \pi_2$ is cuspidal if and only if $\BC_L(\pi_2)$ is not isomorphic to its own twist by $\psi^{\theta} \psi^{-1}$, where $\theta$ is the non-trivial automorphism of $L/F$.
\end{itemize}
\end{thm}
Here is the Galois-theoretic heuristic for the first part of the cuspidality criterion. Assume $V_1$ and $V_2$ are two-dimensional, irreducible, non-dihedral representations. Non-dihedral implies that each $\Sym^2 V_i$ is irreducible. Suppose that $V_1 \otimes V_2 \cong W_1 \oplus W_2$. Taking exterior squares, we get
\[
(\Sym^2 V_1 \otimes \det V_2) \oplus (\det V_1 \otimes \Sym^2 V_2) \cong \wedge^2 W_1 \oplus (W_1 \otimes W_2) \oplus \wedge^2 W_2.
 \]
We may therefore assume $\dim W_1=3$, $\dim W_2=1$ (so rename $W_2$ as a character $\psi$), and thus
\[
V_1 \otimes V_2 (\psi^{-1}) \cong W_1(\psi^{-1}) \oplus 1,
\]
which implies that $V_1$ and $V_2$ are twist-equivalent. 

We now describe the fibers of this tensor product lift, beginning with a Galois heuristic. Suppose we have four irreducible $2$-dimensional representations (over $\Qlb$, say) satisfying
\[
V_1 \otimes V_2 \cong W_1 \otimes W_2.
\]
Taking exterior squares as above, we get
\[
(\Sym^2 V_1 \otimes \det V_2) \oplus (\det V_1 \otimes \Sym^2 V_2) \cong (\Sym^2 W_1 \otimes \det W_2) \oplus (\det W_1 \otimes \Sym^2 W_2)
\]
Assume $V_1$ is non-dihedral, so $\Sym^2 V_1$ is irreducible. Then we may assume
\[
\Sym^2 V_1 \otimes \det V_2 \cong \Sym^2 W_1 \otimes \det W_2
\]
and
\[
\det V_1 \otimes \Sym^2 V_2 \cong \det W_1 \otimes \Sym^2 W_2.
\]
Comparing determinants in this and the initial isomorphism, we find $\det V_1 \det V_2= \det W_1 \det W_2$, and so (for $i=1, 2$)
\[
\Ad^0(V_i) \cong \Ad^0(W_i).
\]
$V_i$ and $W_1$ therefore give rise to isomorphic two-dimensional projective representations, and so, by the (not very) long exact sequence in continuous cohomology associated to 
\[
1 \to \Qlb^\times \to \mr{GL}_2(\Qlb) \to \mathrm{PGL}_2(\Qlb) \to 1,\footnote{Note that the surjection admits a topological section.}
\]
$V_i$ and $W_i$ are twist-equivalent: we find a pair of characters $\psi_i$ such that $V_1 \cong W_1(\psi_1)$ and $V_2 \cong W_2(\psi_2)$. If $V_2$ is also non-induced, then $\psi_2= \psi_1^{-1}$, else $V_1 \otimes V_2$, and hence one of $V_1$ or $V_2$, is an induction.

Because the functorial transfers $\mr{GL}_2 \xrightarrow{\Sym^2} \mr{GL}_3$ (\cite{gelbart-jacquet:sym2}) and $\mr{GL}_4 \xrightarrow{\wedge^2} \mr{GL}_6$ (\cite{kim:wedge2}) are known, this argument works automorphically until the last step. To conclude a purely automorphic argument, we have to know the fibers of $\Ad^0$. These were determined by Ramakrishnan as a consequence of his construction of the tensor product lift (and his cuspidality criterion):
\begin{thm}[Ramakrishnan, \cite{ramakrishnan:tensor} Theorem 4.1.2]
Let $\pi$ and $\pi'$ be unitary cuspidal automorphic representations on $\mr{GL}_2/F$ such that $\Ad^0 \pi \cong \Ad^0 \pi'$. Then there exists a Hecke character $\chi$ of $F$ such that $\pi \cong \pi' \otimes \chi$.
\end{thm}
It would be remiss not to mention the special nature of this fiber description: by \cite{langlands-labesse}, it is essentially multiplicity one for $\mr{SL}_2$. In any case, from this and the preceding calculation, we conclude:
\begin{cor}\label{2,2fibers}
Let $\pi_1, \pi_2, \pi_1', \pi_2'$ be unitary cuspidal automorphic representations on $\mr{GL}_2/F$ satisfying $\pi_1 \boxtimes \pi_2 \cong \pi_1' \boxtimes \pi_2'$ (it suffices to assume this almost everywhere locally). Assume that $\pi_1 \boxtimes \pi_2$ is cuspidal (equivalently, satisfies Ramakrishnan's criterion). Up to reordering, we have $\Ad^0 \pi_i \cong \Ad^0 \pi_i'$ for $i= 1, 2$, and there are Hecke characters $\chi_i$ of $F$ such that $\pi_i \cong \pi_i' \otimes \chi_i$ ($i=1, 2$). Moreover, $\chi_1= \chi_2^{-1}$, so the fibers of the lift are just twists by Hecke characters.\index{t}{fibers of a functorial transfer!$\mr{GL}_2 \times \mr{GL}_2$ tensor product}
\end{cor}
\proof
If at least one, say $\pi_1$, of the $\pi_i$ is non-dihedral, then it only remains to check $\chi_1= \chi_2^{-1}$. Since $\pi_1 \boxtimes \pi_2 \cong (\pi_1 \boxtimes \pi_2)\otimes (\chi_1 \chi_2)$, we deduce from a theorem of Lapid-Rogawski-- which we elaborate on in Lemma \ref{cyclicautind} below--  that either $\chi_1 \chi_2=1$, or $\pi_1 \boxtimes \pi_2$ is an automorphic induction from the (quadratic or quartic) cyclic extension of $F$ cut out by $\chi_1 \chi_2$. In the latter case, there is a cyclic base change $L/F$ such that $\BC_L (\pi_1 \boxtimes \pi_2)= \BC_L(\pi_1) \boxtimes \BC_L(\pi_2)$ is non-cuspidal. If both $\pi_1$ and $\pi_2$ are non-dihedral, or if $\pi_2= \Ind_K^F(\psi_2)$ with $K$ a quadratic extension of $F$ not contained in $L$, then the cuspidality criterion implies  there exists a Hecke character $\psi$ of $L$ such that $\BC_L(\pi_1) \otimes \psi \cong \BC_L(\pi_2)$. Since $\pi_1$ is non-dihedral, $\psi$ is invariant under $\Gal(L/F)$; for cyclic extensions, there is no obstruction to descending invariant Hecke characters, so $\psi$ descends to a character of $F$ that we also write as $\psi$, i.e. $\BC_L(\pi_1 \otimes \psi) \cong \BC_L(\pi_2)$. By cyclic descent (also Theorem $B$ of \cite{lapid-rogawski:descent}), $\pi_1$ and $\pi_2$ are twist-equivalent, contradicting cuspidality of $\pi_1 \boxtimes \pi_2$. We conclude that $\chi_1 \chi_2=1$. 

We resume the above argument when $\pi_2$ is dihedral with $K \subset L$. If $L=K$, then $\chi_1 \chi_2$ is the quadratic character $\chi$ of $\Gal(K/F)$, so $\pi_2= \pi_2' \otimes (\chi_1^{-1} \chi)$, and thus $\pi_2= \pi_2' \otimes \chi_1^{-1}$ since $\pi_2 \otimes \chi = \pi_2$. If $L/F$ is quartic, then the result of Lapid-Rogawski still implies $\pi_1 \boxtimes \pi_2= \Ind_L^F(\psi)$, where $\psi$ is now a Hecke character of $L$. Base-changing to $L$ and comparing isobaric constituents, we again contradict the assumption that $\pi_1$ is non-dihedral.

Finally, we treat the case where both $\pi_1$ and $\pi_2$ are dihedral. Base-change implies that both $\pi_1'$ and $\pi_2'$ are dihedral as well, and moreover that we may assume that $\pi_i$ and $\pi_i'$ are both induced from the same field $L_i$, say by characters $\psi_i$ and $\psi_i'$. Write $\sigma_i$ for the non-trivial element of $\Gal(L_i/F)$. Then
\[
\Ind_{L_1}^F\left(\psi_1 \otimes \BC_{L_1}( \Ind_{L_2}^F(\psi_2)) \right)= \Ind_{L_1}^F\left(\psi'_1 \otimes \BC_{L_1}( \Ind_{L_2}^F(\psi'_2)) \right),
\]
and up to replacing $\psi_1$ by $\psi_1^{\sigma_1}$, we may assume that on $\mr{GL}_2/L_1$ we have
\[
\frac{\psi_1}{\psi'_1} \otimes \BC_{L_1}(\Ind_{L_2}^F(\psi_2))= \BC_{L_1}(\Ind_{L_2}^F(\psi'_2)).
\]
If $L_1 \neq L_2$, then we find
\[
\Ind_{L_1 L_2}^{L_1}(\frac{\psi_1}{\psi'_1}\cdot \psi_2)= \Ind_{L_1 L_2}^{L_1}(\psi'_2),
\]
and possibly replacing $\psi_2$ by its conjugate $\psi_2^{\sigma_2}$, we may assume that as Hecke characters of $L_1 L_2$, $\psi_1/\psi'_1= \psi'_2/\psi_2$. Let $\alpha$ be the Hecke character of $L_1$ given by $\psi_1/\psi'_1$. This equality shows $\alpha$ is $\sigma_1$-invariant, and therefore descends to a Hecke character of $F$. Consequently, $\pi_1= \pi_1' \otimes \alpha$, and $\pi_2= \pi'_2 \otimes \alpha^{-1}$ as automorphic representations on $\mr{GL}_2/F$. The case $L_1=L_2$ is similar, but easier, and we omit the details.
\endproof
The above corollary made use of a characterization of the image of cyclic automorphic induction (in the quadratic and quartic case). This characterization in the case of non-prime-degree does not follow from the results of \cite{arthur-clozel:basechange}, but it is now known thanks to work of Lapid-Rogawski.\index{t}{image of a functorial transfer!cyclic automorphic induction} The following is stated without proof in \cite{lapid-rogawski:descent}, as an easy application of their Statement B\footnote{The proof of which was conditional on versions of the fundamental lemma that are now known.}. We fill in the details:
\begin{prop}\label{cyclicautind}[Lapid-Rogawski]
Let $\Pi$ be a cuspidal automorphic representation on $\mr{GL}_n/F$, and let $L/F$ be a cyclic extension. Then $\Pi$ is automorphically induced from $L$ if (and only if) $\Pi \cong \Pi \otimes \omega$ for a Hecke character $\omega$ of $F$ that cuts out the extension $L/F$.
\end{prop}
\proof
Lapid-Rogawski establish the following, which implies the Proposition:\footnote{Which is in turn used to prove the more refined parts $(b)$ and $(c)$ of Statement B of \cite{lapid-rogawski:descent}.}
\begin{thm}[Statement B, part $(a)$ of \cite{lapid-rogawski:descent}]\label{l-r}
Let $E/F$ be a cyclic extension of number fields, with $\sigma$ a generator of $\Gal(E/F)$, and let $\omega$ be a Hecke character of $E$. Denote by $\omega_F$ its restriction to $\mbf{A}_F^\times \subset \mbf{A}_E^\times$. Suppose $\Pi$ is a cuspidal automorphic representation on $\mr{GL}_n/E$ satisfying $\Pi^{\sigma} \cong \Pi \otimes \omega$. Let $K/F$ be the extension (necessarily of order dividing $n$) of $F$ cut out by $\omega_F$, and let $L=KE$. Then $E \cap K=F$, and $[K:F]$ divides $n$.
\end{thm}
Suppose $\omega$ cuts out $L/F$, cyclic of order $m$, and suppose $m$ is divisible by a prime $\ell$ (and $m \neq \ell$). We see that $\Pi \cong \Pi \otimes \omega^{m/l}$ as well, and letting $E/F$ be the (cyclic degree $\ell$) extension cut out by $\omega^{m/l}$, the case of prime degree implies that $\Pi = \Ind_E^F(\pi)$ for some cuspidal representation $\pi$ on $\mr{GL}_{n}/E$. Moreover,
\[
\Ind_E^F(\pi) \cong \Ind_E^F(\pi \otimes (\omega \circ N_{E/F})),
\]
so the description of the fibers in the prime case implies $\pi$ and $\pi \otimes (\omega \circ N_{E/F})$ are Galois-conjugate: writing $\tau$ for a generator of $\Gal{L/F}$ (or for its restriction to $E$), there exists an integer $i$ such that
\[
\pi^{\tau^i} \cong \pi \otimes (\omega \circ N_{E/F}).
\]
Part (a) of Theorem \ref{l-r} implies that the restriction of $\omega \circ N_{E/F}$ to $\mbf{A}_F^\times/F^\times$, which is just $\omega^\ell$, cuts out an extension $K/F$ that is linearly disjoint from $E/F$,\footnote{The point is that $\Pi$ is cuspidal. For instance, in the simple case $\ell^2 \nmid m$, there's no need to appeal to the Lapid-Rogawski result: $\ell$ divides $\frac{m}{\ell}i$, so $\ell \vert i$, and then $\pi \cong \pi \otimes (\omega \circ N_{E/F})$.} so in particular $(\omega \circ N_{E/F})^\ell$ is a non-trivial Hecke character of $E$.
But now, since $\tau^\ell$ is trivial on $E$, we can iterate the previous conjugation relation to deduce 
\[
\pi \cong \pi \otimes (\omega \circ N_{E/F})^\ell.
\]
It now follows from induction on the degree $[L:F]$ that $\pi \cong \Ind_L^E(\pi_0)$ for some cuspidal $\pi_0$ on $\mr{GL}_{\frac{n}{[L:F]}}/L$, so we're done.

\endproof
\subsection{The Hilbert modular case}
Now we take up the case of Hilbert modular forms, starting with an observation of Blasius-Rogawski (\cite{blasius-rogawski:motiveshmfs}): while a `mixed parity' Hilbert modular representation $\pi$ (Definition \ref{mixedparitydef}) does not itself twist to an $L$-algebraic representation, its base changes to every CM field will.\footnote{This observation of Blasius-Rogawski is the only time I have seen substantive use made of type $A$ but not $A_0$ Hecke characters.} Since we work up to twist, we may assume the parameters of $\pi$ at places $v \vert \infty$ have the form (restricted to $\overline{F}_v^\times \subset W_{F_v}$ and implicitly invoking $\iota$ as above)
\[
z \mapsto 
\begin{pmatrix}
(z/\bar{z})^{\frac{k_v-1}{2}} & 0 \\
0 & (z/\bar{z})^{\frac{1- k_v}{2}}
\end{pmatrix}.
\]
\begin{defn}\label{mixedparitydef}
In the above notation, $\pi$ is \textit{mixed parity} if as $v$ varies over all $v \vert \infty$, $k_v$ takes both even and odd values.\index{t}{mixed-parity Hilbert modular representation}
\end{defn}
For a quadratic CM extension $L/F$, let $\psi \colon \mbf{A}_L^\times /L^\times \to \mathbf{S}^1$ be a Hecke character whose infinity type at $v$ is given by
\[
z \mapsto (z/|z|)^{k_v-1}
\]
for all $v$. Then $\BC_L(\pi) \otimes \psi$ is $L$-algebraic, with infinity type
\[
z \mapsto 
\begin{pmatrix}
(z/\bar{z})^{k_v-1} & 0 \\
0 & 1
\end{pmatrix}
\]
Given two such $\pi_1$ and $\pi_2$, both $W$- but not $L$-algebraic, and with $\Pi= \pi_1 \boxtimes \pi_2$ $L$-algebraic, Blasius-Rogawski (see Theorem 2.6.2 and Corollary 2.6.3 of \cite{blasius-rogawski:motiveshmfs}) can then associate an $\ell$-adic Galois representation $\rho_{\Pi, \ell}$ to $\Pi$ via the identity
\[
\BC_L(\Pi) \cong (\BC_L(\pi_1)(\psi_1)) \boxtimes (\BC_L(\pi_2)(\psi_2)) \otimes (\psi_1 \psi_2)^{-1},
\]
since all three tensor factors on the right-hand side are $L$-algebraic with associated Galois representations. (To get a Galois representation for $\Pi$ over $F$ itself, they vary $L$ and use the now-famous patching argument.) We pursue this a little farther, emphasizing the `$W$-algebraic' and non-de Rham/geometric aspects of this construction.
\begin{prop}\label{subtle}
Let $\Pi$ be the tensor product $\pi_1 \boxtimes \pi_2$ for $\pi_i$ as above corresponding to mixed parity Hilbert modular forms (but $\pi_1$ and $\pi_2$ having common weight-parity at each infinite place). Assume that $\Pi$ is cuspidal, and for simplicity assume that the $\pi_i$ are non-dihedral.
\begin{itemize}
\item[(i)] The $\ell$-adic Galois representation $\rho_{\Pi, \ell}$ is Lie irreducible. 
\item[(ii)] $\rho_{\Pi, \ell}$ is isomorphic to a tensor product $\rho_{1, \ell} \otimes \rho_{2, \ell}$ of two-dimensional continuous, almost everywhere unramified, representations of $\gal{F}$. No twist of either $\rho_{i, \ell}$ is de Rham (or even Hodge-Tate), and $\rho_{\Pi, \ell}$ is not a tensor product of geometric Galois representations (although it is after every CM base change, by construction).\footnote{Compare Proposition \ref{tensorformalism}.}
\item[(iii)] $\Ad^0(\rho_{i, \ell})$ is geometric, corresponding to the $L$-algebraic $\Ad^0(\pi_i)$ (suitably ordering the representations). In particular, a mixed-parity Hilbert modular representation $\pi$ gives rise to a geometric projective Galois representations $\bar{\rho}_{\pi, \ell} \colon \gal{F} \to \mr{PGL}_2(\Qlb)$; these are the representations predicted by Conjecture $3.2.2$ of \cite{buzzard-gee:alg} for irreducible $\mr{SL}_2(\af)$-constituents of $\pi|_{\mr{SL}_2(\af)}$.   
\item[(iv)] With the $\pi_i$ normalized to be $W$-algebraic, they are also $W$-arithmetic (see Question \ref{Walg}).
\end{itemize}
\end{prop}
\proof
To show $\rho_{\Pi, \ell}$ is a tensor product, observe that it is essentially self-dual, since 
\[
\Pi^\vee \cong \pi_1^\vee \boxtimes \pi_2^\vee \cong (\pi_1 \otimes \omega_{\pi_1}^{-1}) \boxtimes (\pi_2 \otimes \omega_{\pi_2}^{-1}) \cong \Pi \otimes (\omega_{\pi_1} \omega_{\pi_2})^{-1}.
\]
On the Galois side, we deduce that $\rho_{\Pi, \ell}^\vee \cong \rho_{\Pi, \ell} \otimes \mu_{\Pi}^{-1}$, where $\mu_{\Pi}$ is the (geometric) Galois character corresponding to the ($L$-algebraic) Hecke character $\omega_{\pi_1} \omega_{\pi_2}$. Since after one of these quadratic base changes $L/F$ $\rho_{\Pi, \ell}|_{\gal{L}}$ is irreducible and orthogonal (rather than symplectic), we see that it (as $\gal{F}$-representation) is orthogonal, i.e.
\[
\rho_{\Pi, \ell} \colon \gal{F} \to \mr{GO}_4(\Qlb).
\]
Writing $\mu$ for the multiplier character on $\mr{GO}_4$, and observing that $\det^2= \mu^4$ on this group, we find an exact sequence of algebraic groups
\[
1 \to \mathbb{G}_m \to \mr{GL}_2 \times \mr{GL}_2 \xrightarrow{\boxtimes} \mr{GO}_4 \xrightarrow{\det \cdot \mu^{-2}} \{\pm 1\} \to 1.
\]
The image of $\rho_{\Pi, \ell}$ is contained in the kernel $\mr{GSO}_4$ of $\det \cdot \mu^{-2}$: for this it suffices to know that $\rho_{\Pi, \ell}$ is isomorphic to a tensor product after two disjoint quadratic base changes, for then $(\det \cdot \mu^{-2}) \circ \rho_{\Pi, \ell}$ is a character of $\gal{F}$ that is trivial on a set of primes of density strictly greater than $1/2$, and therefore it is trivial. We can then apply Proposition \ref{tatelift} to deduce that $\rho_{\Pi, \ell}$ lifts across $\mr{GL}_2(\Qlb) \times \mr{GL}_2(\Qlb) \to \mr{GSO}_4(\Qlb)$, i.e. it is isomorphic to a tensor product $\rho_{1, \ell} \otimes \rho_{2, \ell}$. The same exact sequence implies that any expression of $\rho_{\Pi, \ell}$ as a tensor product has the form $(\rho_{1, \ell} \otimes \chi) \otimes (\rho_{2, \ell} \otimes \chi^{-1})$ for some continuous character $\chi \colon \gal{F} \to \Qlb^\times$. Lie-irreducibility also follows since the Zariski closure of the image of each $\rho_{i, \ell}$ contains $\mr{SL}_2$ (otherwise $\rho_{\Pi, \ell}$, and hence $\Pi$, is automorphically induced), and thus the Zariski closure of the image of $\rho_{\Pi, \ell}$ contains $\mr{SO}_4$, which acts Lie-irreducibly in its standard $4$-dimensional representation.

Now we show all twists $\rho_{i, \ell} \otimes \chi$ are almost everywhere unramified, but none are geometric (they fail to be de Rham). This follows from the Galois-theoretic arguments of the next section, but here we give a more automorphic argument. Note that both the de Rham and almost everywhere unramified conditions can be checked after a finite extension. By construction, for any $L=FK$ for $K/\Q$ imaginary quadratic, there exist ($L$-algebraic) cuspidal automorphic representations $\tau_i$ ($i=1, 2$) on $\mr{GL}_2/L$ with associated geometric Galois representations $\sigma_i \colon \gal{L} \to \mr{GL}_2(\Qlb)$ such that
\[
\rho_{1, \ell}|_{\gal{L}} \otimes \rho_{2, \ell}|_{\gal{L}} \cong \rho_{\Pi, \ell}|_{\gal{L}} \cong \sigma_1 \otimes \sigma_2.
\]
This implies there is a Galois character $\chi \colon \gal{L} \to \Qlb^\times$ such that
\begin{align*}
&\sigma_1(\chi) \cong \rho_{1, \ell}|_{\gal{L}}, \\
&\sigma_2(\chi^{-1}) \cong \rho_{2, \ell}|_{\gal{L}},
\end{align*}
and every abelian $\ell$-adic representation is unramified almost everywhere (easy), so each $\rho_{i, \ell}$ is almost everywhere unramified. We then have the equivalences 
\[
\text{$\rho_{i, \ell}$ is geometric $\iff$ $\rho_{i, \ell}|_{\gal{L}}$ is geometric $\iff$ $\chi$ is geometric.}\footnote{Moreover, $\rho_{i, \ell}$ is Hodge-Tate $\iff$ $\rho_{i, \ell}|_{\gal{L}}$ is Hodge-Tate $\iff$ $\chi$ is Hodge-Tate $\iff$ $\chi$ is geometric!}
\]
But if $\chi$ is geometric, then (Theorem \ref{FMLGL1}) there exists a type $A_0$ Hecke character $\chi_{\mbf{A}}$ of $L$ such that 
\begin{align*}
&\tau_1 \otimes \chi_{\mbf{A}} \sim_{w, \infty} \rho_{1, \ell} \\
&\tau_2 \otimes \chi_{\mbf{A}}^{-1} \sim_{w, \infty} \rho_{2, \ell},
\end{align*}
where $\sim_{w, \infty}$ is the notation of \S \ref{notation}. $\Gal(L/F)$-invariance and cyclic base change then imply there are cuspidal automorphic representations $\tilde{\tau}_i$ on $\mr{GL}_2/F$ lifting $\tau_1 \otimes \chi_{\mbf{A}}$ and $\tau_2 \otimes \chi_{\mbf{A}}^{-1}$, and therefore $\BC_L(\tilde{\tau}_1 \boxtimes \tilde{\tau}_2) \cong \BC_L(\Pi)$. This implies $\tilde{\tau}_1 \boxtimes \tilde{\tau}_2$ and $\Pi$ are twist-equivalent, and Corollary \ref{2,2fibers} ensures that (up to relabeling) $\pi_i$ and $\tilde{\tau}_i$ are twist-equivalent. But $\tilde{\tau}_i$ is $L$-algebraic since $\tau_i$ and $\chi_{\mbf{A}}$ are, yet it is easy to see (because $\pi$ is mixed-parity and $F$ is totally real) that no twist of $\pi_i$ can be $L$-algebraic.

That, suitably ordered, $\Ad^0(\rho_{i, \ell})$ is a geometric Galois representation whose local factors match those of the $L$-algebraic cuspidal representation $\Ad^0(\pi_i)$ follows easily from the earlier $\wedge^2$ argument (see \S \ref{22tensor})and the fact that our forms (representations) are all non-dihedral.

Part $(iv)$: Normalize $\pi_i$ to be $W$-algebraic. Then $\pi_i \boxtimes \pi_i= \Sym^2(\pi_i) \boxplus \omega_{\pi_i}$ is $L$-algebraic, and by \cite{blasius-rogawski:motiveshmfs} it has Satake parameters in $\Qb$ at unramified primes. Writing $\{\alpha_{i, v}, \beta_{i, v}\}$ for the parameters of $\pi_{i,v}$ (when unramified), we conclude that $\alpha_{i, v}^2, \beta_{i, v}^2 \in \Qb$, hence $\alpha_{i, v}, \beta_{i, v} \in \Qb$. This implies $\pi_{i, v}$ has a model over $\Qb$ for all such $v$,\footnote{Note that the fields of definition of $\pi_v$ and its Satake parameters may be different; but if one is algebraic, then both are.} answering the weaker version of Question \ref{Walg}. In this case the stronger version (that $\pi_{i, f}$ is defined over $\Qb$) follows as well, by a check (omitted) at the ramified primes.
\endproof
Now we prove a more refined result in the case $\pi_1= \pi_2$:
\begin{prop}\label{mixedparitygalois}
Let $\pi$ be a unitary, cuspidal, non-induced, mixed-parity Hilbert modular representation, so that $\Pi= \pi \boxtimes \pi$ has an associated Galois representation, the sum of Galois representations associated to the $L$-algebraic representations $\Sym^2(\pi)$ and $\omega_\pi$. As before, there are $\ell$-adic representations $\rho_{i, \ell}$ such that $\rho_{\Pi, \ell}= \rho_{1, \ell} \otimes \rho_{2, \ell}$.
\begin{itemize}
\item We can normalize the $\rho_{i, \ell}$ so that $\rho_{i, \ell}(fr_v)$ has characteristic polynomial in $\Qb[X]$ for all unramified $v$. Moreover, its eigenvalues in fact lie in $\Q^{cm}$, and are pure of some weight, which we may take to be zero.
\item For quadratic CM extensions $L= KF$ with $K/\Q$ a quadratic imaginary field in which $\ell$ is inert, there is a $\rho_{\ell} \colon \gal{F} \to \mr{GL}_2(\Qlb)$ for which $\mr{Sym}^2(\rho_{\ell}) \sim_{w, \infty} \mr{Sym}^2(\pi)$. Later on (Corollary \ref{mixedparitymismatch}), we will see this is not possible over a totally real field. 
\item Nevertheless, we can choose $\rho_\ell$ and a finite-order character $\chi$ such that $\Sym^2(\rho_\ell) \colon \gal{F} \to \mr{GO}_3(\Qlb)$ is, viewing $\mr{GO}_3(\Qlb)$ as the dual group $(\mathbb{G}_m \times \mr{SL}_2)^\vee(\Qlb)$, associated to an automorphic representation $(\omega_\pi \chi, \pi_0)$ of $(\mathbb{G}_m \times \mr{SL}_2) (\af)$, as in $3.2.2$ of \cite{buzzard-gee:alg}, with $\pi_0$ any irreducible constituent of $\pi|_{\mr{SL}_2(\af)}$.
\end{itemize}
\end{prop}
\proof
If we take $\pi_1= \pi_2= \pi$ in the unitary normalization (this will ensure $\pi \boxtimes \pi$ has `motivic weight zero,' with an implicit-- and provable-- application of Ramanujan at infinity) and decompose the Galois representation $\rho_{\Pi, \ell}$ associated to $\Pi= \pi \boxtimes \pi$ as $\rho_1 \otimes \rho_2$, then the now-familiar $\wedge^2$ argument shows that $\Ad^0(\rho_1) \cong \Ad^0(\rho_2)$. We can therefore write $\rho_{\Pi, \ell} \cong \rho \otimes (\rho \otimes \chi)$ for some Galois character $\chi$.  This $\chi$ need not be a square, but we can normalize $\rho$ so that $\chi$ is finite order (by Lemma \ref{Galoisroots}). Now, the Ramanujan conjecture is known for $\Pi$\footnote{Of course, $\Pi$ is not cuspidal, so literally this is for $\Sym^2(\pi)$ and $\omega_\pi$.}, and of course for $\chi$, so $\rho \otimes \rho$ is pure of weight zero (at unramified primes, say). Let $\alpha_v$ and $\beta_v$ be frobenius eigenvalues of $\rho(fr_v)$ at an unramified prime $v$, so that $\alpha_v^2$ and $\beta_v^2$ (being eigenvalues of $\Sym^2(\rho)$) have absolute value one under any embedding $\Qlb \into \CC$. The same then holds for $\alpha_v$ and $\beta_v$, so $\rho$ itself is pure of weight zero. Therefore the frobenius eigenvalues of $\rho$ lie in $\Q^{cm}$.

We conclude by showing that, after certain CM base-changes, we can take the finite-order character $\chi$ to be trivial. Restricting to $L=KF$ as above, we can find a type A Hecke character $\psi$ of $L$ such that $\pi \otimes \psi$ is $L$-algebraic, and then letting $(\widehat{\psi^2})_\ell$ be the geometric $\ell$-adic representation associated to $\psi^2$, we claim that, possibly replacing $\psi$ by $\psi^{-1}$, there exists a Galois character $\lambda$ of $L$ such that 
\[
\rho|_{\gal{L}} \otimes \lambda \cong \rho|_{\gal{L}} \otimes (\chi (\widehat{\psi^2})_\ell \lambda^{-1}).
\]
Over $L$, $\pi \cdot \psi^{-1}$ has associated Galois representation $r$, and then $r \cdot (\widehat{\psi^2})_\ell$ corresponds to $\pi \cdot \psi$. Then there exists a $\lambda$ such that either $r= \rho \cdot \lambda$ and $r \cdot (\widehat{\psi^2})_\ell= \rho(\chi \lambda^{-1})$, or $r= \rho \cdot (\chi \lambda^{-1})$ and $r \cdot (\widehat{\psi^2})_\ell= \rho \cdot \lambda$. Either way, plugging one equation into the other we get the claim.

$\rho$ is not an induction, so $\chi|_{\gal{L}}= \lambda^2 (\widehat{\psi^2})_\ell^{-1}$. In particular, $\chi|_{\gal{L}}$ is a square if and only if $(\widehat{\psi^2})_\ell$ admits a Galois-theoretic square-root. This was discussed in Lemma \ref{galsqrt} above: it need not always be the case (for instance, if $\ell$ is split in $L/\Q$), but it is if $\ell$ is inert in $K$. Choosing such an $L=KF$, we can now find $\chi_L \colon \gal{L} \to \Qlb$ such that
\[
\rho_{\Pi, \ell}|_{\gal{L}} \cong (\rho|_{\gal{L}} \otimes \chi_L)^{\otimes 2}.  
\]
Here $\chi_L= \lambda \hat{\psi}$ for $\hat{\psi}$ a square-root of $(\widehat{\psi^2})_\ell$, and we're done. 

For the final point, note that the inclusion $\mathbb{G}_m \times \mr{SO}_3 \to \mr{GO}_3$ is an isomorphism, with inverse map $g \mapsto (\frac{\det}{\mu}(g), \frac{\det}{\mu}(g)^{-1} g)$. In particular, $\Sym^2 \colon \mr{GL}_2 \to \mr{GO}_3$ becomes $g \mapsto (\det(g), \Ad^0(g)$ in these coordinates. The $\mr{SL}_2(\af)$-constituents of $\pi$ have Langlands parameters corresponding to the projectivization of those of $\pi$, and since $\Ad^0 \colon \mr{GL}_2 \to \mr{SO}_3 \cong \mr{PGL}_2$ is just the quotient map $\mr{GL}_2 \onto \mr{PGL}_2$, the claim follows.
\endproof
\begin{rmk}
In particular, the Proposition tells us that there are pure $\ell$-adic Galois representations no twist of which are geometric. They do not, however, live in compatible systems.
\end{rmk}
These propositions would have no content if the only examples of mixed parity Hilbert modular forms were the inductions described in Example \ref{mixed}; we end this section by showing 
\begin{lemma}\label{limmult}
Over any totally real field, non-CM mixed parity Hilbert modular forms exist.
\end{lemma}
\proof
Considering the (semi-simple, connected) $\Q$-group $G= \Res_{F/\Q}(\mr{SL}_2/F)$, we can apply Theorem 1B from \S 4 of Clozel's paper \cite{clozel:limmult}. We fix a discrete series representation $\pi_{\infty}$ of $G(\RR) \cong \prod_{\tau \colon F \into \RR} \mr{SL}_2(\RR)$ corresponding to the desired mixed-parity infinity-type. Then fix an auxiliary supercuspidal level $K_{p_0} \subset G(\Q_{p_0})$ for some prime $p_0$. At some other finite prime $p$, let $\pi_p$ be a fixed (twist of) Steinberg representation. Finally, let $S$ be a set of finite primes (disjoint from $p_0, p$) at which we will let the level go to infinity in the manner of Clozel's paper (written $K_S \to 1$). Letting $K$ be any fixed compact open subgroup away from $S, p, p_0, \infty$, Clozel proves that 
\[
\liminf_{K_S \to 1} \left[ \mathrm{vol}(K_S) \cdot \mathrm{mult}\left( \pi_{\infty} \otimes \pi_p, L^2_{cusp}(G(\Q) \slash G(\af))^{K_{p_0} \times K_S \times K}\right) \right] >0.
\]
(He also determines an explicit but non-optimal constant.) This suffices for the Lemma: it implies the existence of (infinitely many) cuspidal automorphic representations on $G$ that are Steinberg at $p$, and therefore not automorphically induced, and have the desired mixed parity at infinity.
\endproof 
\begin{rmk}
In \cite{shin:plancherel}, Shin derives, as a corollary of very general existence results for automorphic representations, an exact limit multiplicity formula for \textit{cohomological} Hilbert modular forms. Unfortunately, this excludes precisely the mixed-parity forms we are interested in. It seems, however, that his methods should extend to cover all discrete series infinity types.
\end{rmk}
\subsection{A couple of questions}
In Proposition \ref{subtle}, the key point was an understanding of which type $A$ Hecke characters-- or their Galois analogues-- could exist. I would like to raise here a couple questions in higher rank about the existence of automorphic forms with certain infinity-types.\index{t}{infinity-types of automorphic representations: questions about}
\begin{question}\label{existencequestions}
Are there automorphic representations $\pi$ on $\mr{GL}_2/F$ ($F$ totally real) such that at two infinite places $v_1$ and $v_2$,
\begin{itemize}
\item $\pi_{v_1}$ is discrete series and $\pi_{v_2}$ is limit of discrete series (excluding the easily-constructed dihedral cases, as in Example \ref{mixed})? Clozel's result does not apply to this case;
\item $\pi_{v_1}$ is discrete series (or limit of discrete series), and $\pi_{v_2}$ is principal series (looks like a Maass form)?
\item as in the previous item, but with algebraic infinity-type (so a `Plancherel density'-type result will not suffice)?
\end{itemize} 
\end{question}

\section{Galois lifting: Hilbert modular case}\label{GaloisHMFs}
\subsection{Outline}\label{outline}
We continue to elaborate on the examples of the previous section, now turning to some preliminary cases of a problem raised by Conrad in \cite{conrad:dualGW}. Recall from the introduction that he addresses lifting problems 
\[
\xymatrix{
& H'(\Qlb) \ar[d] \\
\gal{F} \ar@{-->}[ru]^{\tilde{\rho}} \ar[r]^{\rho} & H(\Qlb) \\
},
\]
where $H' \onto H$ is a surjection of linear algebraic groups with central kernel, and the key remaining question is:\index{s}{$H'$}\index{s}{$H$}\index{t}{Galois lifting problem}
\begin{question}
Suppose that the kernel of $\tilde{H} \onto H$ is a torus. If $\rho$ is geometric, when does there exist a geometric lift $\tilde{\rho}$?
\end{question}
Here is an outline of the coming sections:
\begin{itemize}
\item First we address, in the regular case, with $F$ totally real,\footnote{Modulo the difference between potential automorphy and automorphy!} Conrad's question for lifting across $\mr{GL}_2(\Qlb) \to \mr{PGL}_2(\Qlb)$, finding that all examples of $\rho$ not having geometric lifts are accounted for by mixed-parity Hilbert modular forms (a converse to Propositions \ref{subtle} and \ref{mixedparitygalois}).
\item To give a higher-rank example in which potential automorphy theorems still allow us to link the automorphic and Galois sides, we then carry out in \S \ref{spineg} an analogous discussion for $\mr{GSpin}_{2n+1}(\Qlb) \to \mr{SO}_{2n+1}(\Qlb)$ (\S \ref{l-adichodge} introduces the necessary background about Sen operators and `labeled Hodge-Tate-Sen weights' in $\ell$-adic Hodge theory; in the current section, we use more elementary terminology, although some manipulations are justified by the general theory).
\item Before proceeding to develop the general framework (\S \ref{generalGalois}), we take a detour (\S \ref{liftingwalg}) to discuss the automorphic analogue. This is in some sense more intuitive, and it motivates the Galois-theoretic solution. Note that in general we have neither constructions of automorphic Galois representations nor (potential) automorphy theorems, so we cannot unconditionally bridge the automorphic-Galois chasm. Proving a `modular lifting theorem' in this context is an important remaining question (see Question \ref{modularlifting}).
\end{itemize}
\subsection{$\mr{GL}_2(\Qlb) \to \mr{PGL}_2(\Qlb)$}
In this subsection we prove:
\begin{thm}\label{nodRlift}
Let $F$ be totally real. Suppose $\rho \colon \gal{F} \to \mr{PGL}_2(\Qlb)$ is a geometric representation with no geometric lift to $\mr{GL}_2(\Qlb)$. Then 
\begin{itemize}
\item There exists a lift $\tilde{\rho} \colon \gal{F} \to \mr{GL}_2(\Qlb)$ such that for all CM $L/F$, $\tilde{\rho}|_{\gal{L}}$ is the twist of a geometric Galois representation (in particular, $\rho|_{\gal{L}}$ has a geometric lift, which ought to be automorphic). 
\end{itemize}
Now assume moreover that $\Ad(\rho) \colon \gal{F} \to \mr{SO}_3(\Qlb) \subset \mr{GL}_3(\Qlb)$ satisfies the hypotheses of (the potential automorphy theorem) Corollary $4.5.2$ of \cite{blggt:potaut}\footnote{Essentially: there is a lift $\tilde{\rho}$ of $\rho$ such that $(\tilde{\rho} \mod \ell)|_{\gal{F(\zeta_\ell)}}$ is irreducible and $\tilde{\rho}$ is regular (if true for one lift, this is true for all), and $\ell$ is sufficiently large.} Then:
\begin{itemize}
\item After a totally real base-change $F'/F$, $\Ad(\rho)$ is automorphic, and more precisely
\item we may normalize $\tilde{\rho}$ such that $\Sym^2(\tilde{\rho})|_{\gal{F'}}$ is a geometric Galois representation corresponding to $\Sym^2(\pi) \otimes \chi$ for some mixed parity Hilbert modular representation $\pi$ on $\mr{GL}_2/F$ and *non-trivial* finite order character $\chi$. Furthermore, $\tilde{\rho}$ is totally odd. 
\end{itemize}
\end{thm}
\begin{rmk}
If $\Ad^0(\rho)$ is not regular, there are two possibilities: either it fails to be regular everywhere, in which case Fontaine-Mazur conjecture that $\rho$ (up to twist) has finite image. Such $\rho$ always have finite-image (hence geometric) lifts by Tate's theorem (Proposition \ref{tatelift}). If $\rho$ is totally odd, and $\Ad^0(\rho)$ exhibits a mixture of regular and irregular behavior, it should be related to the existence of a mixed-parity Hilbert modular form that is limit of discrete series at some infinite places, and genuine discrete series at others. I know of no such examples not arising from Hecke characters. The question of whether there exist non-dihedral $\rho$ that are even at some and odd at other infinite places is related to Question \ref{existencequestions}.
\end{rmk}
We now make our first use of non-Hodge-Tate Galois representations; recall that in \S \ref{l-adichodge} we have described the theory of (non-integral) Hodge-Tate-Sen weights. In fact, a reader unfamiliar with this general theory will have no difficulty following the arguments of the current section.
\begin{lemma}\label{halfZ}
For any number field $F$ and any geometric $\rho \colon \gal{F} \to \mr{PGL}_2(\Qlb)$, there is a lift $\tilde{\rho} \colon \gal{F} \to \mr{GL}_2(\Qlb)$ of $\rho$ whose Hodge-Tate-Sen (HTS) weights at all $v \vert \ell$ belong to $\frac{1}{2} \Z$. In this case, $\Sym^2(\tilde{\rho})$ and $\det(\tilde{\rho})$ are geometric.
\end{lemma}
\proof
Proposition \ref{tatelift} ensures that we can find some lift $\tilde{\rho}$. By assumption, $\Ad^0(\tilde{\rho})$ is geometric, so 
\[
\tilde{\rho} \otimes \tilde{\rho} \otimes \det(\tilde{\rho})^{-1} \cong \Ad^0(\tilde{\rho}) \oplus \mathbf{1}
\]
is geometric. We ask whether $\det(\tilde{\rho})$ is a square. There at least exists a Galois character $\psi \colon \gal{F} \to \Qlb^\times$ such that $\det(\tilde{\rho})= \psi^2 \chi$ with $\chi$ finite order (Lemma \ref{Galoisroots}). Thus
\[
[\tilde{\rho}\otimes \psi^{-1}]^{\otimes 2} \cong (\Ad^0(\tilde{\rho})\otimes \chi) \oplus \chi
\]
is also geometric, and in particular Hodge-Tate. Replacing $\tilde{\rho}$ with the new lift $\tilde{\rho} \otimes \psi^{-1}$ we conclude that $\Sym^2(\tilde{\rho})$ and $\det{\tilde{\rho}}$ are Hodge-Tate. The theorem of Wintenberger and Conrad (Theorem \ref{conradliftingP}) then proves they must be de Rham, since the original projective $\rho$ was. Any lift $\tilde{\rho}$ is almost everywhere unramified by Lemma 5.3 of \cite{conrad:dualGW}, so the proof is complete. 
\endproof
We now proceed in the same spirit as in the construction of $W$-algebraic forms via type $A$ Hecke characters. The following lemma is an \textit{ad hoc} version of Theorem \ref{fullmonodromylift}; it may help in navigating that rather formal proof:
\begin{lemma}\label{galtwist}
Consider $\tilde{\rho} \colon \gal{F} \to \mr{GL}_2(\Qlb)$ as produced by Lemma \ref{halfZ}.  Let $L/\Q$ be a quadratic CM extension. Then for all such $L$, the restriction of $\tilde{\rho}$ to $\gal{L}$ is the twist of a geometric Galois representation.
\end{lemma}
\proof
Recall that our fixed embeddings $\Qb \xhookrightarrow{\iota_\infty} \CC$ and $\Qb \xhookrightarrow{\iota_\ell} \Qlb$ yield a map 
\[
\iota_{\infty, \ell}^* \colon \bigcup_{v \vert \ell} \Hom_{\Ql}(F_v, \Qlb) \to \Hom_{\Q}(F, \CC).
\]
$\Sym^2(\tilde{\rho})$ is de Rham, so for each place $v \vert \ell$ of $F$ and each $\tau \colon F_v \into \Qlb$, it has $\tau$-labeled Hodge-Tate weights\footnote{By duality, or by the Hodge-Tate-Sen theory.} 
\[
\mathrm{HT}_{\tau}\left(\Sym^2(\tilde{\rho}|_{\gal{F_v}})\right)= \{A_{\tau}, \frac{A_\tau+B_\tau}{2}, B_{\tau}\},
\]
where the $A_\tau$ and $B_\tau$ are integers, necessarily having the same parity. We can distinguish the subset 
\[
\{A_\tau, B_\tau\} \subset \mathrm{HT}_{\tau}\left(\Sym^2(\tilde{\rho}|_{\gal{F_v}})\right)
\]
(trivially in the non-regular case, easily otherwise). It follows that $\tilde{\rho}|_{\gal{F_v}}$ is Hodge-Tate (hence de Rham) if and only if all $A_\tau$ and $B_\tau$ are even. Consider the map $\iota_{\infty, \ell}^*$ for $L$ as well as $F$; the two versions are compatible, so conjugate archimedean embeddings $\iota, \bar{\iota} \colon L \into \CC$ pull back to embeddings $L \into \Qlb$ that lie above a common $\tau \colon F \subset F_v \into \Qlb$. To each $\tau$ we then unambiguously associate the parity $\epsilon_\tau \in \{0, 1\}$ of $A_\tau$ (and $B_\tau$), and at the unique archimedean place $v_{\infty, \tau}$ of $L$ above $\iota_{\infty, \ell}^*(\tau)$ we define a character
\begin{align*}
\psi_{v_{\infty, \tau}} \colon &L_{v_{\infty, \tau}}^\times \to \CC^\times \\
&z \mapsto \iota(z)^{\frac{k_{\tau}}{2}} \bar{\iota}(z)^{-\frac{k_{\tau}}{2}},
\end{align*}
where $k_\tau$ is any integer with parity $\epsilon_\tau$.

$L$ is a CM field, so Weil (\cite{weil:characters}) tells us that there exists a type $A$ Hecke character $\psi$ of $L$ with components at archimedean places given by these $\psi_{v_{\infty, \tau}}$. We wish to twist $\tilde{\rho}|_{\gal{L}}$ by a Galois realization of $\psi$, shifting all of its HTS weights to be integers (recall this can be checked by looking at $\Sym^2(\tilde{\rho}|_{\gal{L}})$), but care is needed: there is no Galois representation canonically associated to a type $A$ (not $A_0$) Hecke character. Instead, we associate (canonically, having specified $\iota_\infty$ and $\iota_\ell$) an $\ell$-adic representation $(\widehat{\psi^2})_{\ell}$ to the type $A_0$ Hecke character $\psi^2$, and then we non-canonically extract, up to a finite-order character, a (Galois-theoretic, non-geometric) square root $\hat{\psi}$. Then $\tilde{\rho}|_{\gal{L}} \otimes \hat{\psi}$ is Hodge-Tate, since the analogues of the $A_\tau$ and $B_\tau$, but now for $\Sym^2(\tilde{\rho}|_{\gal{L}} \otimes \hat{\psi})= \Sym^2(\tilde{\rho}|_{\gal{L}}) \otimes (\widehat{\psi^2})_\ell$, are all even. Thus $\tilde{\rho}|_{\gal{L}} \otimes \hat{\psi}$ is geometric (again by Theorem \ref{conradliftingP}).
\endproof
This completes the proof of the first part of Theorem \ref{nodRlift}. We need deeper tools to make further progress. Let $\tilde{\rho}$ be the lift as above, with both $\Sym^2(\tilde{\rho})$ and $\det(\tilde{\rho})$ geometric. $\Sym^2(\tilde{\rho})$ factors through $\mr{GO}_3(\Qlb)$ with totally even multiplier character $\det(\tilde{\rho})^2$. We now make further assumptions on $\tilde{\rho}$ (or $\rho$) in order to apply a potential automorphy theorem:
\begin{hypothesis}\label{potauthyp}
\begin{enumerate}
\item Assume $\ell > 7$.
\item Assume $\Sym^2(\tilde{\rho})$ is `potentially diagonalizable' and regular at all $v \vert \ell$. For instance, we require that for each $v \vert \ell$ it be crystalline with $\tau$-labeled Hodge-Tate weights distinct and falling in the Fontaine-Laffaille range for all $\tau \colon F_v \into \Qlb$. 
\item Assume that the reduction mod $\ell$ of $\Sym^2(\tilde{\rho})|_{\gal{F(\zeta_{\ell})}}$ is irreducible. (We have taken $\ell \neq 2$, so this is equivalent to the analogous assumption for $\tilde{\rho}$.)
\end{enumerate}
\end{hypothesis}
\index{t}{potential automorphy theorem}Then we can can apply Theorem $C$ of \cite{blggt:potaut} to deduce that after base-change to some totally real field $F'$, $\Sym^2(\tilde{\rho})|_{\gal{F'}}$ is automorphic, corresponding to a RAESDC automorphic representation $\Pi$ on $\mr{GL}_3/F'$. Let us abusively denote by $\det(\tilde{\rho})$ the Hecke character (of type $A_0$) corresponding to $\det(\tilde{\rho})$. Then $\Pi \otimes \det(\tilde{\rho})^{-1}$ is self-dual with trivial central character, and (by, for instance, \cite{arthur:classical}-- details, in greater generality, appear in Lemma \ref{arthur}-- although in this case the result is older) there exists a cuspidal $\pi$ on $\mr{GL}_2/F'$ such that $\Ad^0(\pi) \cong \Pi \otimes \det(\tilde{\rho})^{-1}$. We may assume that $\omega_\pi$ has finite order (applying Lemma \ref{Heckelift} and Proposition \ref{llgeneralization}, since $F'$ is totally real), because the descent is originally to $\mr{SL}_2/F'$, and we then have control over the choice of extension to $\mr{GL}_2/F'$.\footnote{We only sketch this here, because these matters will be discussed in greater detail and generality in \S \ref{liftingwalg}.}
\begin{lemma}
Such a $\pi$ on $\mr{GL}_2/F'$ satisfying $\Ad^0(\pi) \cong \Pi \otimes \det(\tilde{\rho})^{-1}$ is necessarily a mixed parity Hilbert modular representation. 
\end{lemma}
\proof
It is clear from the archimedean L-parameters that $\pi$ is $W$-algebraic. If it were $L$-algebraic, there would be a corresponding Galois representation $\rho_{\pi}$\footnote{By \cite{blasius-rogawski:motiveshmfs}.} with $\Ad^0(\rho_\pi) \cong \Ad^0(\tilde{\rho})$, and thus $\rho_{\pi}$ and $\tilde{\rho}|_{\gal{F'}}$ would be twist-equivalent. Lemma \ref{galtotreal} then implies $\tilde{\rho}$ has all integral or all half-integral Hodge-Tate weights, which contradicts the assumption that $\rho$ has no geometric lift (note that both $F'$ and $F$ are totally real). Similarly, if $\pi$ were $C$-algebraic, then there would be a geometric representation corresponding to $\pi |\cdot|^{1/2}$, and we can use the same argument. We conclude that $\pi$ must be mixed parity.
\endproof
\begin{cor}\label{mixedparitymismatch}
For any $\rho \colon \gal{F} \to \mr{PGL}_2(\Qlb)$ as in Theorem \ref{nodRlift}, and satisfying the auxiliary potential automorphy Hypothesis \ref{potauthyp}, there exists a lift $\tilde{\rho} \colon \gal{F} \to \mr{GL}_2(\Qlb)$ and, after a totally real base change $F'/F$, a mixed-parity Hilbert modular representation $\pi$ on $\mr{GL}_2/F'$, such that 
\[
\Sym^2(\tilde{\rho})|_{\gal{F'}} \sim_{w, \infty} \Sym^2(\pi) \otimes \chi
\]
for some finite order Hecke character $\chi$ of $F'$. Moreover, any lift $\tilde{\rho}$ is totally odd, and the character $\chi$ is necessarily non-trivial. In contrast, restricting to $L'=F' K$ where $K/\Q$ is a quadratic imaginary extension in which $\ell$ is inert, we can find a lift $\tilde{\rho}$ such that $\Sym^2(\tilde{\rho})$ corresponds to $\BC_{L'/F'}(\Sym^2(\pi))$.

Conversely, starting with a mixed-parity $\pi$, we can produce such a $\chi$ and $\tilde{\rho}$.
\end{cor}
\proof
In the notation of the previous Lemma, we have 
\[
\Sym^2(\pi) \otimes \frac{\det(\tilde{\rho})}{\omega_\pi} \cong \Pi \sim_{w, \infty} \Sym^2(\tilde{\rho})|_{\gal{F'}}. 
\]
The character $\det(\tilde{\rho})/\omega_\pi$ is $L$-algebraic, hence equals $\chi |\cdot|^m$ for a finite-order character $\chi$ and an integer $m$.  Replacing $\pi$ by $\pi \otimes |\cdot|^{m/2}$, which is still mixed-parity, we are done except for the last part.

Assume instead that we can take $\chi=1$. In particular, $\omega_\pi \sim_{w, \infty} \det(\tilde{\rho})$, and thus $\det(\tilde{\rho})(c_v)= \omega_{\pi_v}(-1)$. By the description of discrete series representations of $\mr{GL}_{2}(\RR)$ and purity of $\omega_\pi$, we know that the Langlands parameter $\phi_v \colon W_{\RR} \to \mr{GL}_{2}(\CC)$ takes the form
\begin{align*}
&z \mapsto 
\begin{pmatrix}
(z/\bar{z})^{\frac{k_v-1}{2}} & 0 \\
0 & (z/\bar{z})^{\frac{1-k_v}{2}}
\end{pmatrix} |z|_{\CC}^{-w/2}\\
&j \mapsto
\begin{pmatrix}
0 & (-1)^{k_v-1} \\
1 & 0
\end{pmatrix}.
\end{align*}
Hence $\omega_{\pi_v}(-1)= \det(\phi_v)(j)= (-1)^{k_v}$. Since $\Ad^0(\tilde{\rho})|_{\gal{F'}} \sim_{w, \infty} \Ad^0(\pi)$ is regular algebraic self-dual cuspidal, Proposition $A$ of \cite{taylor:sign} shows that $\Ad^0(\tilde{\rho})$ is odd, and thus for all $v \vert \infty$, the eigenvalues of $\Ad^0(\tilde{\rho})(c_v)$ are $\{-1,1,-1\}$. We conclude that $\det(\tilde{\rho})(c_v)= -1$ for all $v \vert \infty$, proving the oddness claim for $\tilde{\rho}$, and contradicting, when we assume $\chi=1$, the fact that the $k_v$ are both odd and even.

The remaining claims, constructing from $\pi$ such a $\tilde{\rho}$ and $\chi$, were established in Proposition \ref{mixedparitygalois}.
\endproof
This concludes the proof of Theorem \ref{nodRlift}.

\section{Spin examples}\label{spineg}
To address Conrad's question in general, we will have to re-cast the arguments of \S \ref{GaloisHMFs} in terms of root data. In the meantime we give another pleasantly concrete example, but one that relies on root-theoretic manipulation, in the hopes of easing the transition to the general case; moreover, specializing to discrete series/regular examples, we can also still exploit known results about automorphic Galois representations. The previous example (lifting across the surjection $\mr{GL}_2(\Qlb) \to \mr{PGL}_2(\Qlb)$) is really just the first case of a family of spin examples,
\[
\xymatrix{
& \mr{GSpin}_{2n+1}(\Qlb) \ar[d] \\
\gal{F} \ar@{-->}[ru]^-{\tilde{\rho}} \ar[r]_-{\rho} & \mr{SO}_{2n+1}(\Qlb),}
\]
where $\rho$ is geometric but has no geometric lift $\tilde{\rho}$. $F$ will be a totally real field, and we will again see that over CM fields, geometric lifts do in fact exist. The case $n=1$ will amount to the contents of \S \ref{GaloisHMFs}.

We first recall the basic setup for (odd) $\mr{Spin}$ groups.\index{s}{$\mr{Spin}$}\index{s}{$\mr{GSpin}$} We will also make heavy use of this notation, and its obvious analogues for even $\mr{Spin}$ groups, later on in some of our examples of `motivic lifting' (see \S \ref{hyperlift} and \S \ref{conclusion}).
\[
\xymatrix @R=1.5pc {
& 1 \ar[d] & 1 \ar[d] & & \\
1 \ar[r] & \{\pm 1\} \ar[r]  \ar[d] & \mr{Spin}_{2n+1} \ar[r] \ar[d]& \mr{SO}_{2n+1} \ar[r] \ar@{=}[d]& 1 \\
1 \ar[r] & \mathbb{G}_m \ar[r] \ar[d]^{z \mapsto z^2}& \mr{GSpin}_{2n+1} \ar[r] \ar[d]^{\nu} & \mr{SO}_{2n+1} \ar[r] & 1 \\
& \mathbb{G}_m \ar@{=}[r] \ar[d] & \mathbb{G}_m \ar[d]& & \\
& 1 & 1 && 
}
\]
where $\nu$ is the Clifford norm. We can then identify
\[
\xymatrix{
\mr{GSpin}_{2n+1} \ar[d]^{\nu}  & \frac{\mathbb{G}_m \times \mr{Spin}_{2n+1}}{\{1, (-1, c)\}} \ar[l]^{\sim} \ar[d]^{(x,g) \mapsto x^2} \\
\mathbb{G}_m \ar@{=}[r] & \mathbb{G}_m}
\]
where $c$ is the non-trivial central element of $\mr{Spin}_{2n+1}$. This yields an identification of Lie algebras
\[
\mf{gspin}_{2n+1} \xleftarrow{\sim} \mf{gl}_1 \times \mf{so}_{2n+1},
\]
where $\mf{gl}_1$ is identified with the center of $\mf{gspin}_{2n+1}$. Alternatively, $\mr{GSpin}_{2n+1}$ is the dual group to $\mr{GSp}_{2n}$, and it will be convenient to keep both normalizations of (based) root data for $\mr{GSpin}_{2n+1}$-- one from the spin description, one from the dual description-- in mind. Let $(X, \Delta, X^\vee, \Delta^\vee)$ be the based root datum for $\mr{GSp}_{2n}$ (say defined with respect to $J_{2n}= 
\left( \begin{smallmatrix}
0 & 1_n \\
-1_n & 0 
\end{smallmatrix} \right)$),
with its diagonal maximal torus $T$ and the Borel $B \supset T$ for which $e_i-e_{i+1}$ and $2e_n-e_0$ form a set of positive simple roots, with
\begin{align*}
&e_i \colon diag(t_1, \ldots, t_n, \nu t_1^{-1}, \ldots, \nu t_n^{-1}) \mapsto t_i \\
&e_0 \colon diag(t_1, \ldots, t_n, \nu t_1^{-1}, \ldots, \nu t_n^{-1}) \mapsto \nu.
\end{align*}
Then for $\mr{GSpin}_{2n+1}$, we have the based root datum (with $X^\vee$ the character group)\index{t}{root data for $\mr{Spin}$ and $\mr{GSpin}$ groups}
\begin{align*}
&X^\vee= \bigoplus_{i=0}^n \Z e_i^{*} \\
&\Delta^\vee= \{\alpha_i^\vee= e_i^*-e_{i+1}^*\}_{i=1}^{n-1} \cup \{\alpha_n^\vee= e_n^*\} \\
&X= \bigoplus_{i=0}^n \Z e_i \\
&\Delta= \{\alpha_i= e_i-e_{i+1}\}_{i=1}^{n-1} \cup \{\alpha_n= 2e_n- e_0\},
\end{align*}
with $X$ and $X^\vee$ in the duality $\langle e_i, e_j^* \rangle= \delta_{ij}$. 

\label{spinnotation}Alternatively, since $\mr{Spin}_{2n+1}$ is the connected, simple, simply-connected group with Lie algebra $\mf{so}_{2n+1}$, we can write its character group as the weight lattice of $\mf{so}_{2n+1}$, i.e., as the submodule of $\oplus_{i=1}^n \Q \chi_i$ spanned by $\chi_1, \ldots, \chi_n, \frac{\chi_1+ \ldots+ \chi_n}{2}$, and its co-character lattice as those $\sum_{i=1}^n c_i \lambda_i$ such that $c_i \in \Z$ and $\sum c_i \in 2\Z$. The duality is $\langle \chi_i, \lambda_j \rangle= \delta_{ij}$. Letting $\chi_0$ and $\lambda_0$ generate $X^{\bullet}(\mr{GL}_1)$ and $X_{\bullet}(\mr{GL}_1)$, respectively, the description of $\mr{GSpin}_{2n+1}$ as $\frac{\mr{GL}_1 \times \mr{Spin}_{2n+1}}{\{1, (-1, c)\}}$ then leads to another description of the root datum as
\begin{align*}
&Y^\bullet = \bigoplus_{i=1}^n \Z \chi_i \oplus \Z (\chi_0+ \frac{\chi_1 + \ldots + \chi_n}{2}) \subset \bigoplus_{i=0}^n \Q \chi_i \\
&\Delta^\bullet = \{\chi_i-\chi_{i+1}\}_{i=1}^{n-1} \cup \{\chi_n\}\\
&Y_\bullet = \bigoplus_{i=1}^n \Z(\lambda_i + \frac{\lambda_0}{2}) \oplus \Z \lambda_0\\
&\Delta_\bullet= \{\lambda_i-\lambda_{i+1} \}_{i=1}^{n-1} \cup \{2 \lambda_n\}.
\end{align*}
We summarize by comparing the two descriptions:
\begin{lemma}
There is an isomorphism of based root data
\[
(X^\vee, \Delta^\vee, X, \Delta) \cong (Y^\bullet, \Delta^\bullet, Y_\bullet, \Delta_\bullet)
\]
given by
\begin{align*}
&e_i^* \mapsto \chi_i \quad \text{for $i=1, \ldots, n$;}\\
&e_0^* \mapsto \chi_0 - \frac{\chi_1 + \ldots +\chi_n}{2}; \\
&e_i \mapsto \lambda_i+ \frac{\lambda_0}{2} \quad \text{for $i=1, \ldots, n-1$;}\\
&e_0 \mapsto \lambda_0.
\end{align*}
\end{lemma}
For example, the center of $\mr{GSpin}_{2n+1}$ is generated by the co-character $e_0 \leftrightarrow \lambda_0$. The Clifford norm is given by the character $2\chi_0 \leftrightarrow 2e_0^*+ \sum_{1}^n e_i^*$. 
%The projection to $\mr{SO}_{2n+1}$ corresponds to including $\oplus_1^n \Z e_i$ and modding out by $\Z e_0^*$, $e_0^*$ being the co-character that generates the center $Z_{\mr{GSpin}_{2n+1}}$. Finally, the Clifford norm $\nu \in X(\mr{GSpin}_{2n+1})$ extends $2 e_0 \in X$ (see the above diagram).

Returning to our representation $\rho$, recall that to each $v \vert \ell$ and $\iota \colon \Qlb \into \CC_{F_v}$ we can associate the Sen operator 
\[
\Theta_{\rho, \iota}:= \Theta_{\rho|_{\gal{F_v}}, \iota} \in \mf{so}_{2n+1} \otimes_{\Qlb, \iota} \CC_{F_v}.
\]
(The place $v$ is implicit in $\iota$.) Since $\rho$ is Hodge-Tate, we can identify $\Theta_{\rho, \iota}$ up to conjugation with a diagonal element
\[
\left(\begin{smallmatrix}
m_1 & &&&&& \\
& \ddots &&&&& \\
& & m_n &&&& \\
& & & 0 & &&& \\
& & & & -m_n & & \\
& & &  & & \ddots & \\
& & & & & & -m_1
\end{smallmatrix}\right),
\]
where $m_j= m_j(\iota)$ are all integers.\footnote{We choose the `anti-diagonal' orthogonal pairing, so that $\mr{SO}_{2n+1}$ has a maximal torus consisting of diagonal matrices.}
\begin{prop}\label{spinlift0}
Let $\rho \colon \gal{F} \to \mr{SO}_{2n+1}(\Qlb)$ be a geometric representation as above, with $F$ totally real. Then $\rho$ lifts to a geometric, $\mr{GSpin}_{2n+1}$-valued representation if and only if the parity of $\sum_j m_j(\iota)$ is independent of $v \vert \ell$ and $\iota \colon \Qlb \into \CC_{F_v}$.
\end{prop}
\begin{rmk}
By itself, this result is formal. Later we will see how such $\rho$ arise, at least when $\rho$ is regular; but note that this Proposition makes no such assumption.
\end{rmk}
\proof
By Proposition \ref{tatelift}, some lift $\tilde{\rho}$, necessarily unramified almost everywhere, exists. All possible lifts differ from this $\tilde{\rho}$ by $e_0 \circ \chi$ for some $\chi \colon \gal{F} \to \Qlb^\times$. As before, we can write the multiplier character $\nu(\tilde{\rho}) \colon \gal{F} \to \Qlb^\times$ as $\chi^2 \chi_0$, where $\chi_0$ has finite order, and then replace $\tilde{\rho}$ by a twist having $\nu(\tilde{\rho})$ of finite-order (Lemma \ref{twist} will explain this procedure in general). Then, since finite-order characters have Hodge-Tate weights zero, functoriality of the Sen operator implies that for all $\iota \colon \Qlb \into \CC_{F_v}$, $\Theta_{\tilde{\rho}, \iota}$ maps to $(0, \Theta_{\rho, \iota})$ in $(\mf{gl}_1)_{\iota} \oplus (\mf{so}_{2n+1})_{\iota}$.\footnote{Writing $\mf{g}_\iota$ as a short-hand for $\mf{g} \otimes_{\Qlb, \iota} \CC_{F_v}$.}

Now consider the spin representation $\mr{GSpin}_{2n+1} \xrightarrow{r_{spin}} \mr{GL}_{2^n}$. In the above root datum notation this corresponds to the highest weight $-e_0^* \leftrightarrow -\chi_0 + \frac{\chi_1 + \ldots + \chi_n}{2}$.  The image of $\Theta_{\tilde{\rho}, \iota}$ is then a semi-simple element with one eigenvalue $\frac{m_1+ \ldots +m_n}{2}$, and all eigenvalues congruent to this (half-integer) modulo $\Z$. In particular, $r_{spin} \circ \tilde{\rho}|_{\gal{F_v}}$ is Hodge-Tate, and thus de Rham, if and only if for all $\iota \colon \Qlb \into \CC_{F_v}$, $m_1(\iota)+ \ldots + m_n(\iota)$ is even. If all (for all $v \vert \ell$ and all $\iota$) of these sums are odd, then we twist by a character $\chi$ with all Hodge-Tate-Sen weights $1/2$ (see Lemma \ref{galtotreal}) to get a new $\tilde{\rho}$, now geometric, lifting $\rho$. On the other hand, Lemma \ref{galtotreal} shows that for $F$ totally real we cannot twist $\tilde{\rho}$ in a similar way if some $\sum_j m_j(\iota)$ are even and others are odd.
\endproof
\begin{cor}
Let $\rho \colon \gal{F} \to \mr{SO}_{2n+1}(\Qlb)$ be as in Proposition \ref{spinlift0}. Let $L/F$ be a CM extension. Then $\rho|_{\gal{L}}$ has a geometric lift. 
\end{cor}
\proof
With the framework of the previous proof, this follows by the same argument as Lemma \ref{galtwist}.
\endproof
If we make additional assumptions so that potential automorphy theorems apply, we can of course say more. The next two lemmas merely cash in on some very deep recent results.
\begin{lemma}\label{cashin1}
Assume that $\rho \colon \gal{F} \to \mr{SO}_{2n+1}(\Qlb) \subset \mr{GL}_{2n+1}(\Qlb)$ as in the previous proposition moreover satisfies:
\begin{itemize}
\item For all $v \vert \ell$, $\rho|_{\gal{F_v}}$ is regular and potentially diagonalizable;
\item $\bar{\rho}|_{\gal{F(\zeta_\ell)}}$ is irreducible.
\item $\ell \geq 2(2n+2)$.
\end{itemize}
Then after some totally real base change $F'/F$, there exists a regular $L$-algebraic self-dual cuspidal automorphic representation $\pi$ of $\mr{GL}_{2n+1}(\mathbf{A}_{F'})$ such that $\pi \sim_{w, \infty} \rho$. 
\end{lemma}
\proof
This is immediate from Theorem $C$ ($=$Corollary 4.5.2) of \cite{blggt:potaut}, since $\rho$ is automatically totally odd self-dual.\footnote{Note that those authors always work with $C$-algebraic automorphic representations, so the statements of their theorems always have an extra, but easily unraveled, twist.}
\endproof
\begin{lemma}\label{arthur}
Let $\rho \sim_{w, \infty} \pi$ be as in Lemma \ref{cashin1}. Then $\pi$ descends to a cuspidal automorphic representation on $\mr{Sp}_{2n}(\afp)$, in a way compatible with unramified and archimedean local L-parameters.
\end{lemma}
\proof
This follows from Arthur's classification of automorphic representations of classical groups (\cite{arthur:classical}). Namely, $\rho^\vee \cong \rho$ implies that $\pi^\vee \cong \pi$, and $\det(\rho)=1$ implies that $\omega_{\pi}=1$. $\pi$ is cuspidal, so the associated formal $A$-parameter $\phi$ (see \S $1.4$ of \cite{arthur:classical}) is simple generic and so comes from a unique simple twisted endoscopic datum $G_\phi$ as in Theorem $1.4.1$ of \cite{arthur:classical}; since $2n+1$ is odd, the parameter $\phi$ therefore factors through either $\mr{SO}_{2n+1}(\CC)$ or $\mr{O}_{2n+1}(\CC)$, but the latter case is ruled out since $\omega_\pi=1$. It follows that $G_\phi= \mr{Sp}_{2n}/F$.\footnote{One would like to say that since $\rho$ lands in $\mr{SO}_{2n+1}(\Qlb)$, $\Sym^2(\rho)$ contains the trivial representation and therefore $L(s, \pi, \Sym^2)$ has a pole at $s=1$; this would allow us to descend $\pi$ to $\mr{Sp}_{2n}(\af)$ using \cite{cogdell-kim-ps-shahidi}. Unfortunately, nothing is known a priori about $L(s, V)$ where $\Sym^2(\rho)= \mbf{1} \oplus V$; in particular, we can't say immediately it is non-vanishing at $s=1$, so this argument doesn't seem to work. The argument we have given allows us to deduce that $L(s, \pi, \Sym^2)$ has a simple pole at $s=1$, and therefore $L(s, V)$ is non-vanishing at $s=1$.} The local statement follows from \cite[Theorem 1.4.2]{arthur:classical}.
\endproof
\begin{prop}\label{sp}
Continuing with the assumptions (and conclusions) of the previous two lemmas, $\rho|_{\gal{F'}}$ has a geometric lift to $\mr{GSpin}_{2n+1}(\Qlb)$ if and only if $\pi$ admits an $L$-algebraic extension to $\mr{GSp}_{2n}(\af)$. In the other direction, we have, as in Lemma \ref{limmult}, an automorphic construction of infinitely many such $\rho$, whose local behavior\footnote{i.e. inertial type} we can additionally specify at any finite number of places, that do not admit a geometric lift.
\end{prop}
\proof
(Some details of this argument are omitted, deferring to more general arguments in the next section.) At each $v \vert \infty$, we can write the $L$-parameter of $\pi_v$ as 
\begin{align*}
&\phi_v \colon W_{F_v} \to \mr{SO}_{2n+1}(\CC) \times W_{F_v} \\
&\phi_v|_{W_{\overline{F}_v}} \colon z \mapsto
\left( \begin{smallmatrix}
z^{m_1}\bar{z}^{l_1} & &&&&& \\
& \ddots &&&&& \\
& & z^{m_n}\bar{z}^{l_n} &&&& \\
& & & 1 & &&& \\
& & & & z^{-m_n}\bar{z}^{-l_n} & & \\
& & &  & & \ddots & \\
& & & & & & z^{-m_1}\bar{z}^{-l_1}.
\end{smallmatrix}\right),
\end{align*}
with all $m_i, l_i \in \Z$. Since the transfer of $\pi$ to $\mr{GL}_{2n+1}(\af)$ is cuspidal, Clozel's archimedean purity theorem (see Lemme $4.9$ of \cite{clozel:alg}) implies that $m_i+l_i=0$ for all $i$ (and all $v$). To extend $\pi$ in an $L$ or $W$-algebraic fashion, the key point is to construct an appropriate extension of the central character. Any lift $\tilde{\phi}_v$ to $\mr{GSpin}_{2n+1}(\CC)$ of the parameter $\phi_v$ must, on $W_{\overline{F}_v}$, take the form $z \mapsto z^{\tilde{\mu}_v}\bar{z}^{\tilde{\nu}_v}$, where
\begin{align*}
&\tilde{\mu}_v= \sum_{i=1}^n m_{v,i}(\lambda_i+ \frac{\lambda_0}{2}) + \mu_{v, 0} \lambda_0 \\
&\tilde{\nu}_v=-\sum_{i=1}^n m_{v, i}(\lambda_i + \frac{\lambda_0}{2}) + \nu_{v, 0} \lambda_0.
\end{align*}
The central character of this extension is (by pairing with $2\chi_0$, the Clifford norm; this procedure for computing central characters is explained in general in \cite{langlands:archllc})
\[
\omega_v \colon z \mapsto z^{\sum_{i=1}^n m_{v, i}+ 2\mu_{v, 0}} \bar{z}^{-\sum_{i=1}^n m_{v, i}+ 2 \nu_{v, 0}}.
\]
If we can choose an extension of $\pi$ to an automorphic representation $\tpi$ of $\mr{GSp}_{2n}(\afp)$ with finite-order central character, then this calculation shows that $\tilde{\phi}_v$ (the local $L$-parameter for $\tpi_v$) is always `$W$-algebraic,'\footnote{Using this term abusively for the obvious local analogue.} and it is `$L$-algebraic' if and only if $\sum_{i=1}^n m_{v, i}$ is even. The result would follow easily, twisting our given $\tpi$ by $|\cdot|^{1/2}$ if all $\sum_i m_{v, i}$ are odd (just as in the Galois case). That we can find such an extension $\tpi$ with finite-order central character will be proven in much more generality in Proposition \ref{totrealalglift}.

The construction of geometric $\rho \colon \gal{F} \to \mr{SO}_{2n+1}(\Qlb)$ having no geometric lift to $\mr{GSpin}_{2n+1}(\Qlb)$ (and with specified local behavior) follows as in Lemma \ref{limmult}, applying Clozel's limit multiplicity formula to $G= \mr{Sp}_{2n}$ over $F$, transferring these forms to $\mr{GL}_{2n+1}$ (via \cite{arthur:classical}), and then invoking the Paris Book Project, or, more precisely, Remark $7.6$ of \cite{shin:galoisreps}.
\endproof
We have therefore generalized some of the results of \S \ref{GaloisHMFs} (the case $n=1$).

%\chapter{Galois and automorphic lifting}\label{2}
\chapter{Galois and automorphic lifting}\label{2}
This chapter begins (\S \ref{liftingwalg}) by addressing the natural automorphic analogue of Conrad's lifting question. It is much easier to see what should be true in this setting, and the proofs are simpler as well. Equipped with the intuition coming from the automorphic case, we address the original Galois-theoretic question in \S \ref{generalGalois}. In \S \ref{descentcomparison}, we combine the results of \S \ref{liftingwalg} and \ref{generalGalois} to compare, \textit{modulo the Fontaine-Mazur-Langlands conjecture}, descent problems for certain automorphic representations and Galois representations. The closing section, \S \ref{monodromy}, is of a rather different nature, assembling a few results about the images of compatible systems of $\ell$-adic Galois representations. The results of this section continue the attempt to compare aspects of the automorphic and Galois-theoretic formalisms.

\section{Lifting $W$-algebraic representations}\label{liftingwalg}
In this section, we make some simplifying assumptions: $G$ will be a connected semi-simple split group over $F$\index{s}{$G$}, and $\tZ$\index{s}{$\tZ$} will be an $F$-split \textit{torus} containing $Z_G$. Then let $\tG= (\tZ \times G)/Z_G$\index{s}{$\tG$}, as before, with maximal torus $\tT= (\tZ \times T)/Z_G$\index{s}{$\tT$} and center $\tZ$. In each case $Z_G$ is embedded anti-diagonally. In fact, the cases of greatest interest are when $G$ is simply-connected (so $G^\vee$ is adjoint), such as $\mr{SL}_n$ or $\mr{Sp}_{2n}$, but we do not assume this. The assumptions that $G$ is semi-simple and that $\tZ$ is a torus are not essential, and in fact in \S \ref{generalGalois} we will work more generally. 

Now let $\pi$ be a (unitary) cuspidal automorphic representation of $G(\mbf{A}_F)$ that is $W$-algebraic. We are interested in the problem of lifting $\pi$ to a $W$-algebraic automorphic representation $\tpi$ of $\tG(\mbf{A}_F)$; when $\pi$ is moreover $L$-algebraic, we similarly ask whether an $L$-algebraic lift exists. By `lift,' we simply mean the most naive thing: the restriction $\tpi|_{G(\af)}$ contains $\pi$. Corollary \ref{it'salift} justifies this convention, showing that under the $L$-morphism ${}^L \tG \to {}^L G$, the $L$-packet of $\pi$ is a functorial transfer of the $L$-packet of $\tpi$. 
\subsection{Notation and central character calculation}\label{coordinates}
In order to do computations in terms of root data, we choose coordinates (using invariant factors) such that
\[
X^\bullet(\tZ)= \bigoplus_{i=1}^r \Z w_i \oplus \bigoplus_{j=1}^s \Z w_j',
\]
and the kernel of $X^\bullet(\tZ) \to X^\bullet(Z_G)$ is
\[
\bigoplus_{i=1}^r d_i \Z w_i \oplus \bigoplus_{j=1}^s \Z w_j'.
\]
This is $X^\bullet(S)$ for the torus $S= \tZ/Z_{G}$. We then write $X^\bullet(Z_G)= \oplus_{1}^r (\Z/d_i\Z) \bar{w}_i$. To relate parameters for $\tG$ to those for $G$, we use the Cartesian diagram
\[
\xymatrix{
X^\bullet(\tT) \ar[r] \ar[d] & X^\bullet(\tZ) \ar[d] \\
X^\bullet(T) \ar[r] & X^\bullet(Z_G),
}
\]
representing elements of $X^\bullet(\tT)$ as pairs $(\chi_T, \chi_{\tilde{Z}})$ that are congruent in $X^\bullet(Z_G)$. Extending scalars to $\Q$ (or any characteristic zero field), $X^\bullet(\tT)_\Q \xrightarrow{\sim} X^\bullet(T)_\Q \oplus X^\bullet(\tZ)_\Q$. We will usually (out of sad necessity) take $F$ to be either CM or totally real. In the former case (or whenever an archimedean place is complex) we compute local central characters (see page $21$ ff. of \cite{langlands:archllc}) as follows:

\index{t}{central character of an automorphic representation-- how to compute}Suppose $v$ is complex, so there is an isomorphism $\iota_v \colon F_v \xrightarrow{\sim} \CC$. Let 
\[
\rec_{v}(\pi_v) \colon z \mapsto \iota_v(z)^{\mu_v} \bar{\iota}(z)^{\nu_v} \in T^\vee(\CC)\footnote{From now on, we write $z, \bar{z}$ in place of $\iota_v(z), \bar{\iota}_v(z)$.}
\]
be the associated Langlands parameter, with $\mu_v, \nu_v \in X^\bullet(T)_\CC$ with $\mu_v-\nu_v \in X^\bullet(T)$. Write $\sum_{i=1}^r [\mu_{v}- \nu_v]_i \bar{w}_i$ for the projection of $\mu_v- \nu_v$ to $X^\bullet(Z_G)$. Choose any $\mu_{v, i}, \nu_{v, i} \in \CC$ with $\mu_{v, i}- \nu_{v, i}$ an integer projecting to $[\mu_v- \nu_v]_i \in \Z/d_i \Z$. We can then regard $\tilde{\mu}_v= (\mu_v, \sum_1^r \mu_{v, i} w_i)$ and $\tilde{\nu}_v=(\nu_v, \sum_1^r \nu_{v, i} w_i)$ as elements of $X^\bullet(\tT)_\CC$ parametrizing an extension of the local $L$-parameter of $\pi_v$ to a parameter for $\tG(F_v)$. Identifying, via our chosen basis, $\tZ(F_v)$ with $(\CC^\times)^{r+s}$, the central character $\omega_{\tpi_v}$ of this lift is simply
\[
\omega_{\tpi_v} \colon (z_1, \ldots, z_r, z_1', \ldots, z_s') \mapsto \prod_{i=1}^r z_i^{\mu_{v, i}- \nu_{v, i}} |z_i|_{\CC}^{\nu_{v, i}}.
\]
Restricting to $Z_G(F_v) \cong \mu_{d_1} \times \cdots \times \mu_{d_r} \subset (\CC^\times)^r$, we find that the central character $\omega_{\pi_v}$ is given by
\[
\omega_{\pi_v} \colon (\zeta_1, \ldots, \zeta_r) \mapsto \prod_{i=1}^r \zeta_i^{\mu_{v, i}- \nu_{v, i}},
\]
where each $\zeta_i$ is a $d_i^{th}$ root of unity. Clearly this is independent of the choice of lifts $\tmu_v, \tnu_v$.
\begin{defn}\label{cmheckelift}
If $F$ is a CM field, and $\pi$ is a (unitary) automorphic representation of $G(\af)$ with archimedean parameters $\mu_v, \nu_v$ as above at each $v \vert \infty$, then we define $\tomega= \tomega(\pi)$ to be any choice of Hecke character of $\tZ$ lifting $\omega_\pi$ and such that \index{s}{$\tomega$}
\[
\tomega_v(z_1, \ldots, z_s')= \prod_{i=1}^r (z_i/|z_i|)^{\mu_{v, i}-\nu_{v, i}}.
\]
This is a unitary, type $A$ Hecke character whose existence is assured by Lemma \ref{Heckelift}.
\end{defn}

When $v$ is real, the lifting process is more complicated (see \cite{langlands:archllc}), and the central character computation depends very much on where the $L$-parameter of $\pi_v$ sends an element of $W_\RR- W_\CC$. We can avoid this, however:
\begin{defn}\label{totrealheckelift}
Let $\pi$ be as in the previous definition, but now suppose $F$ is totally real. Then we define $\tomega$ to be any choice of finite-order Hecke character of $\tZ$ extending $\omega_{\pi}$ (existence again by Lemma \ref{Heckelift}).
\end{defn}

\subsection{Generalities on lifting from $G(\af)$ to $\tG(\af)$}
We now try to find an automorphic representation $\tpi \subset L^2_{cusp}(\tG(F) \bs \tG(\af), \tomega)$\footnote{That is, the space of measurable functions on $\tG(F) \bs \tG(\af)$ with (unitary) central character $\tomega$, and square-integrable modulo $\tZ(F_{\infty})^0$-- or, equivalently, modulo $A_{\tG}(\RR)^0$, where $A_{\tG}$ is the maximal $\Q$-split central torus in $\Res_{F/\Q}(\tG)$.} lifting $\pi$, where $\tomega$ is any (unitary) lift of the central character $\omega_\pi$. We say that we are in the Grunwald-Wang special case if one of the pairs $(F, d_i)$ is in the usual Grunwald-Wang special case\index{t}{Grunwald-Wang special case}. Let $H_{\af}= G(\af) \tZ(\af)$ (a normal subgroup of $\tG(\af)$), and let $H_F= H_{\af} \cap \tG(F)$.
\begin{lemma}
We are in the Grunwald-Wang special case if and only if $H_F$ strictly contains $G(F) \tZ(F)$.
\end{lemma}
\proof
Recall the characters $d_iw_i \in X^\bullet(S)= X^\bullet(\tG)$ (since $G$ is semi-simple). These induce an isomorphism
\[
H_F/ G(F) \tZ(F) \xrightarrow{\prod d_i w_i} \prod_{i=1}^r \left( F^\times \cap (\af^\times)^{d_i}\right)/(F^\times)^{d_i}.
\]
\endproof
There are various ways to show forms on $G(\af)$ extend to $\tG(\af)$; the argument we use here is modeled on one of Flicker for the case $(\tG, G)= (\mr{GSp}_{2n}, \mr{Sp}_{2n})$ (see Proposition $2.4.3$ of \cite{flicker:low-rank}).\index{t}{image of a functorial transfer!Tate's lifting problem}
\begin{prop}\label{llgeneralization}
Let $\pi$ be a cuspidal automorphic representation of $G(\af)$, with central character $\omega_\pi$. If we are not in the Grunwald-Wang special case, then for any extension $\tomega$ of $\omega_\pi$ to a Hecke character of $\tZ$, there exists a cuspidal automorphic representation $\tpi$ of $\tG(\af)$ lifting $\pi$, and having central character $\tomega$. If we are in the Grunwald-Wang special case, then for at least one extension of $\omega_\pi$ to 
\[
\prod_{i=1}^r C_F[d_i] \supset \prod_{i=1}^r \mu_{d_i}(F) \bs \mu_{d_i}(\af),
\]
and for any extension of this character to a Hecke character $\tomega$ of $\tZ$, there exists a cuspidal $\tpi$ lifting $\pi$ with central character $\tomega$.\footnote{This odd modification in the special case is presumably an artifact of my clumsy proof.}
\end{prop}
\begin{rmk}
 In particular, in all cases, for all Hecke characters $\tomega$ extending $\omega_\pi$, there exists a Hecke character $\tomega'$ having the same infinity-type, and a cuspidal representation $\tpi$ of $\tG(\af)$ lifting $\pi$, with central character $\tomega'$.
\end{rmk}
\proof
Choose an extension $\tomega \colon \tZ(F) \bs \tZ(\af) \to \mr{S}^1$ of $\omega_\pi$. By extending functions along $\tZ(\af)$ via $\tomega$, we embed the space $V_\pi$ of $\pi$ into the space of cusp-forms on $G(F)\tZ(F) \bs H_{\af}$, and thereby obtain an extension of $\pi$ to a representation $\pi_{\tomega}$ of $H_{\af}$. We construct an intertwining map
\[
\Ind_{H_{\af}}^{\tG(\af)} (\pi_{\tomega}) \xrightarrow{U} L^2_{cusp}(\tG(F) \bs \tG(\af), \tomega),
\] 
where the induction consists of functions $F \colon \tG(\af) \to V_{\pi}$ with the usual left-$H_{\af}$-equivariance, and the requirement of compact support modulo $H_{\af}$. Write $\delta_1 \colon V_{\pi} \to \CC$ for evaluation at $1$ of the cusp forms in $V_{\pi}$, and set
\[
L(F)= \sum_{G(F) \tZ(F) \bs \tG(F) \ni u} \delta_1(F_u),
\]
where we write $F_u$ for the value at $u$ of $F \in \Ind (\pi_{\tomega})$. This sum is in fact finite: $\tG(F)$ is discrete in $\tG(\af)$, the fibers of $G(F) \tZ(F) \bs \tG(F) \to H_{\af} \bs \tG(\af)$ are finite, and we have assumed $F$ has compact support modulo $H_{\af}$. It is also well-defined because  for $\gamma \in G(F) \tZ(F)$, 
\[
\delta_1(F_{\gamma u})= \delta_1(\gamma \cdot F_u)= F_u(\gamma)= F_u(1)= \delta_1(F_u).
\]
The map $U$ is then given by 
\[
U(F) \colon g \mapsto L(I(g)F),
\]
where $I(\cdot)$ denotes the $\tG(\af)$-action on the induction. $U$ is clearly $\tG(\af)$-equivariant and has output $U(F)$ which is left-$\tG(F)$-equivariant (by its construction as an average) and has central character $\tomega$; $U(F)$ is a cusp form by a simple calculation using the fact that the unipotent radical of any parabolic of $\tG$ is in fact contained in $G$.

To see that $U \neq 0$, take a non-zero form $f \in V_{\pi}$; we may assume (by translating) $f(1) \neq 0$. Then define $F \in \Ind (\pi_{\tomega})$ by 
\[
F_h \colon h' \mapsto f(h' h)
\]
if $h \in H_{\af}$, and $F_g=0$ for $g \notin H_{\af}$. Then 
\[
U(F)(1)= \sum_u \delta_1(F_u)= \sum_u
\begin{cases}
 f(u) & \text{if $u \in H_{\af}$;}\\
 0 & \text{otherwise.}
\end{cases}
\]
So, if we are not in the special case, we just get $U(F)(1)= f(1) \neq 0$. If we are in the special case, then for each pair $(F, d_i)$ in the special case there is an element $\alpha_i \in F^\times$ which is everywhere locally a $d_i^{th}$-power-- say $\alpha_i= \beta_i^{d_i}$, but not globally a $d_i^{th}$-power. By abuse of notation, we also write $\alpha_i$ for an element of $\tG(F)$ such that $(d_iw_i)(\alpha_i)$ is this element of $F^\times$ and the other characters $d_j w_j$ ($j \neq i$) and $w_j'$ of $\tG$ are trivial on $\alpha_i$.\footnote{We use the assumption that $G$ is split over $F$.} Regarding $\beta_i$ as an element of the $i^{th}$ componenet of $\tZ(\af)$, we can therefore write $\alpha_i= \beta_i \cdot \gamma_i$, where each $\gamma_i$ lies in $G(\af)$. If only one pair $(F, d_i)$ is in the special case, then the above expression for $U(F)(1)$ becomes $f(1)+ f(\alpha_i)= f(1)+ \tomega(\beta_i)f(\gamma_i)$. This time we normalize $f$ so that $f(\gamma_i) \neq 0$ (rather than $f(1) \neq 0$), and so necessarily either $f(1)+\tomega(\beta_i) f(\gamma_i) \neq 0$, or this expression is non-zero for any $\tomega$ extending the *other* extension of $\omega_\pi$ to a character of $C_F[d_i]$. Thus we can choose at least one initial extension to $C_F[d_i]$, and thereafter the argument proceeds as in the non-exceptional case. If multiple $(F, d_i)$ are in the special case, the same argument applies: arrange some $f(\gamma_i) \neq 0$, and if $U(F)(1)=0$, change the extension of $\omega_\pi$ to $C_F[d_i]$ just in this $i^{th}$ component.

Finally, given that $U \neq 0$, we take $\tpi$ to be the image of any irreducible constituent of $\Ind_{H_{\af}}^{\tG(\af)} (\pi_{\tomega})$ that survives under $U$, and we claim this is the desired extension of $\pi_{\tomega}$. The isomorphism classes of $H_{\af}$-representations appearing in the restriction to $H_{\af}$ of the full induction are precisely the conjugates $\pi_{\tomega}^g$, for $g \in \tG(\af)$, and $\tG(\af)$-stability of $\tpi$ implies that all of these in fact appear in $\tpi|_{H_{\af}}$. In particular, this latter restriction contains $\pi_{\tomega}$.
\endproof
We need to supplement this by understanding to what extent $\pi$ is a `functorial transfer' of $\tpi$, with respect to the $L$-homomorphism ${}^L \tG \onto {}^LG$. Somewhat more generally, let $G \xrightarrow{\eta} \tG$ be a morphism of connected reductive groups over $F_v$, with abelian kernel and cokernel; here we either take $v$ to be archimedean, or finite such that $G$ and $\tG$ are unramified over $F_v$. By assumption, there is a dual $L$-homomorphism ${}^L \eta \colon {}^L \tG \to {}^L G$. In both the archimedean and unramified cases, the $L$-parameters $\phi_v$ and $\tilde{\phi}_v$, as well as the corresponding $L$-packets $\Pi_{\phi_v}$ and $\Pi_{\tilde{\phi}_v}$, of $\pi_v$ and $\tpi_v$ have been defined (as will be explained in Corollary \ref{it'salift}). 
\begin{cor}\label{it'salift}
Continue in the setting of the proposition. \index{t}{weak transfer of automorphic representations}Then $\pi$ is a weak transfer of $\tpi$, which is also a strong transfer at archimedean places. More precisely, for all places $v$ of $F$ that are either archimedean or such that $G$ is unramified at $v$, the restriction $\tpi_v|_{G(F_v)}$ is a finite direct sum of elements of $\Pi_{\phi_v}$, and $\tilde{\phi}_v$ reduces to $\phi_v$ (i.e. ${}^L \eta \circ \tilde{\phi}_v= \phi_v$ up to $G^\vee(\CC)$-conjugacy). \index{t}{$L$-packet!archimedean}\index{t}{$L$-packet!unramified}
\end{cor}
\proof
For infinite places $v$, the assertion follows from desideratum $(iv)$ on page $30$ of \cite{langlands:archllc}. Since the construction of the correspondence is inductive, the verification necessarily stretches out through \S 3 of that paper; for the first (and most important) case of discrete series, see page $43$. 

Now we treat the unramified case. Suppose $\tpi_v$ is unramified, i.e. there exists a hyperspecial maximal compact subgroup $\tilde{K}_v$ of $\tG(F_v)$ such that $\tpi_v^{\tilde{K}_v} \neq 0$. Then for some unramified character $\tilde{\chi}_v$ of $\tT(F_v)$, $\tpi_v$ is a sub-quotient of $I_{\tilde{B}(F_v)}^{\tG(F_v)} \tilde{\chi}_v$. The natural map $G(F_v)/B(F_v) \to \tG(F_v)/\tilde{B}(F_v)$ is an isomorphism, so restriction of functions gives an isomorphism of $G(F_v)$ representations
\[
\left( I_{\tilde{B}(F_v)}^{\tG(F_v)} \tilde{\chi}_v \right)|_{G(F_v)} \xrightarrow{\sim} I_{B(F_v)}^{G(F_v)}(\tilde{\chi}_v|_{B(F_v)}).
\]
Write $\chi_v$ for $\tilde{\chi}_v|_{B(F_v)}$. It is an unramified character of $T(F_v)$, and under the identification of unramified characters of $T(F_v)$ with $\Hom(X_{\bullet}(T)^{\gal{F_v}}, \CC^\times)$ (and the analogue for $\tilde{T}$), and thus with the space of unramified $L$-parameters, $\chi_v$ corresponds to the parameter $\phi_v= {}^L\eta \circ \tilde{\phi_v}$. By definition of unramified $L$-packets (see \cite[\S 10.4]{borel:L}), to see that any constituent $\pi_v$ of $\tilde{\pi}_v|_{G(F_v)}$ lies in the packet $\Pi_{\phi_v}$, we need only check that $\pi_v$ has invariants under some hyperspecial maximal compact subgroup of $G(F_v)$. To see this, let $u$ be a non-zero vector in $\tpi_v^{\tilde{K}_v}$. Decomposing $\tpi_v|_{G(F_v)}= \oplus_{1}^M \pi_i^{\oplus m_i}$ (with the $\pi_i$ distinct isomorphism classes; in fact, the multiplicities $m_i$ are all equal), we see that $u$ lies in one of the isotypic componenets $\pi_i^{\oplus m_i}$, since the induction $I(\chi_v)$ can only have one $K_v:= G(F_v) \cap \tilde{K}_v$- invariant line. Fixing a decomposition $u= \sum_{j=1}^{m_i} u_j$ with $u_j$ in the $j^{th}$ copy of $\pi_i$, we find that the $u_j$ are themselves also $K_v$-invariant, hence $\pi_i$ contains a $K_v$-invariant vector (this in turn implies the multiplicities $m_i$ must equal $1$). Since $\tG(F_v)$ acts transitively on the isomorphism classes $\pi_i$, each $\pi_i$, and in particular our given $\pi_v$, contains some vector of the form $\tpi_v(g)u$, which is $gK_v g^{-1}$-invariant; $gK_v g^{-1}$ is also a hyperspecial maximal compact subgroup, so we're done.
\endproof
\begin{rmk}
According to expected properties of local $L$-packets, $\pi$ should be a strong transfer of $\tpi$ (see \S $10.3$ of \cite{borel:L}), but of course this statement is meaningless until the local Langlands correspondence is known for $G$ and $\tG$. 
\end{rmk}
Eventually, we will also want to understand the ambiguity in the choice of $\tpi$.\index{t}{fiber of a functorial transfer!Tate's lifting problem} Strictly speaking, we only need this in one direction of Proposition \ref{totrealalglift} below, but as a general problem, its importance is basic. We will need to assume something from local representation theory, which has also previously been conjectured, more generally for quasi-split $G$ and $\tG$, by Adler and Prasad (\cite{adler-prasad}. 

\begin{conj}\label{locmultone}
Let $v$ be any place $F$, and let $G$ and $\tG$ be as above. Then for any irreducible smooth representation $\tpi_v$ of $\tG(F_v)$, the restriction $\tpi_v|_{G(F_v)}$ decomposes with multiplicity one.
\end{conj}
The motivation for this conjecture is the uniqueness of Whittaker models for quasi-split groups, and indeed the conjecture holds for generic $\tpi_v$. The analogous statement for $v$ archimedean is straightforward in many cases: for instance, for a simple $G$, it is easy because $\RR^\times/(\RR^\times)^n$ is cyclic for all $n$.  Moreover, the structure of archimedean L-packets may be well-enough understood to deduce the conjecture in general for archimedean $v$, but we do not pursue this here. For $v$ non-archimedean, we have verified it in the unramified case (see the proof of Lemma \ref{it'salift}), and in general the conjecture is known for the following pairs $(\tG, G)$: $(\mr{GL}_n, \mr{SL}_n)$ (from the theory of, possibly degenerate, Whittaker models; see \cite{langlands-labesse} for the $n=2$ case), and any pair $(GU(V), U(V))$ where $(V, \langle, \rangle)$ is a vector space over $F_v$ equipped with a non-degenerate symmetric or skew-symmetric form $\langle, \rangle$, and $GU(V)$, respectively $U(V)$, denote the similitude and isometry groups of the pairing (see Theorem $1.4$ of \cite{adler-prasad}). So, whatever the status of the general conjecture, the discussion below applies to some interesting cases.
\begin{lemma}\label{fibers}
Assume Conjecture \ref{locmultone}. Suppose $\tpi$ and $\tpi'$ are two lifts of $\pi$ (i.e. their restrictions contain $\pi$) with the same central character $\tomega$. 
\begin{itemize}
\item There exist continuous idele characters $\alpha_i \colon \af^\times \to \CC^\times$, for $i=1, \ldots, r$, such that, in the notation of \S \ref{coordinates},
\[
\tpi \cong \tpi' \cdot \left( \prod_{i=1}^r \alpha_i \circ d_i w_i \right),
\]
and each $\alpha_i^{d_i}$ factors as a genuine Hecke character $\alpha_i^{d_i} \colon \af^\times/F^\times \to \CC^\times$.
\item Since Hecke characters satisfy purity, so do the characters $\alpha_i$. 
\end{itemize}
\end{lemma}
\begin{rmk}
The point of the Lemma is that the idele characters $\alpha_i$ \textit{need not be Hecke characters}, i.e. need not factor through $C_F$. We are particularly interested in the constraint on the infinity type of $\alpha_i$, which by the lemma must be just as rigid as the constraints for Hecke characters. This will be applied in Proposition \ref{totrealalglift}.
\end{rmk}
\proof
As before, write $H_{\af}= G(\af) \tZ(\af)$, and write $H_v= G(F_v) \tZ(F_v)$ for its local analogue. By hypothesis, both restrictions $\tpi|_{H_{\af}}$ and $\tpi'|_{H_{\af}}$ contain $\pi \boxtimes \tomega$, so this holds everywhere locally as well. Lemma $2.4$ of \cite{gelbart-knapp:l-indistinguishability} (applying Conjecture \ref{locmultone}) implies that for all finite $v$, there exists a smooth character $\alpha_v \colon \tG(F_v)/H_v \to \CC^\times$ such that $\tpi_v \cong \tpi'_v \cdot \alpha_v$. By duality (for finite abelian groups), we can extend this to a character $\bar{\alpha}_v$:
\[
\xymatrix{
1 \ar[r] & \tG(F_v)/H_v \ar[d]_{\alpha_v} \ar[r] & \prod_{i=1}^r F_v^\times/(F_v^\times)^{d_i} \ar[ld]_{\bar{\alpha}_v}\\
& \CC^\times & \\
}
\]
We express $\alpha_v$ in coordinates as 
\[
\alpha_v= \left(\alpha_{i, v} \circ (d_i w_i)\right)_{i=1}^r
\]
for smooth characters $\alpha_{i,v}$ of $F_v^\times$ trivial on $(F_v^\times)^{d_i}$. Globally, we then have 
\[
\tpi \cong \tpi' \cdot \left( \prod_{i=1}^r (\otimes_v' \alpha_{i, v}) \circ d_i w_i \right).
\]
Here the $\alpha_i:= \otimes_v' \alpha_{i, v}$ are continuous\footnote{Almost all $\alpha_{i, v}$ are unramified, since the same holds for $\tpi_v$ and $\tpi'_v$.} characters $\af^\times \to \CC^\times$, but they need not be Hecke characters. Taking central characters, however, we see that each $\alpha_i^{d_i}$ \textit{is} in fact a Hecke character.
\endproof
\begin{rmk}
Under Conjecture \ref{locmultone}, one should be able to refine the arguments of this section to produce a multiplicity formula for cuspidal automorphic representations of $G(\af)$ in terms of those of $\tG(\af)$, as in Lemma $6.2$ of \cite{langlands-labesse}). One can also axiomatize the passage between local Langlands conjectures for $G$ and for $\tG$; this would include verifying compatibility of the conjectural formulae for the sizes of $L$-packets (a template, in the case of $(\mr{GSp}, \mr{Sp})$, is given in the paper of \cite{gan-takeda:llcsp(4)}; their arguments will clearly apply much more generally).
\end{rmk}

\subsection{Algebraicity of lifts: the ideal case}
Now we return to questions of algebraicity, handling the CM and totally real cases in turn. In each case, we carefully choose an extension of the central character of $\pi$ (as described in Definitions \ref{cmheckelift} and \ref{totrealheckelift}), and then find a $\tpi$ as in Proposition \ref{llgeneralization} with that central character, and whose restriction to $G(\af)$ contains $\pi$. This will yield enough information about the archimedean $L$-parameters of $\tpi$ to deduce the desired algebraicity statements. To give a clean argument with broad conceptual scope, we will assume that $\pi$ is tempered at infinity\footnote{No assumption at finite places.}. As remarked before, this is not a serious assumption for algebraic representations: it is satisfied for forms having cuspidal transfer to some $\mr{GL}_N$. For non-tempered forms, analogous results for the discrete spectrum can be deduced from Arthur's conjectures.
\begin{prop}\label{cmextension}
Let $F$ be CM, with $\pi, \tpi$ as above, and with $\tomega$ the (type $A$) extension of $\omega_{\pi}$ described in Definition \ref{cmheckelift}. Assume $\pi_\infty$ is tempered.
\begin{enumerate}
\item If $\pi$ is $L$-algebraic, then $\tpi$ is $L$-algebraic. In particular, there exists an $L$-algebraic lift of $\pi$.
\item If $\pi$ is $W$-algebraic, then $\tpi$ is $W$-algebraic.
\end{enumerate}
\end{prop}
\proof
We use the notation of \S \ref{coordinates}. We may of course take all $\mu_{v, i}$ and $\nu_{v, i}$ to be integers. By the choice of central character $\tomega$, the archimedean $L$-parameter for $\tpi_v$ corresponds to 
\[
\tmu_v=(\mu_v, \sum \frac{\mu_{v, i}- \nu_{v, i}}{2} w_i) \in X^\bullet(T)_\CC \oplus X^\bullet(\tZ)_\CC,
\] 
and 
\[
\tnu_v=(\nu_v, \sum \frac{\nu_{v, i}- \mu_{v, i}}{2} w_i) \in X^\bullet(T)_\CC \oplus X^\bullet(\tZ)_\CC.
\] 
We write this parameter as an obviously integral term plus a defect:
\[
\tmu_v= (\mu_v, \sum \mu_{v, i} w_i)+ (0, -\sum \frac{\mu_{v, i}+ \nu_{v, i}}{2} w_i),
\]
and likewise for $\tnu_v$. Note that this lies in $X(\tT)$ if and only if for all $i=1, \ldots, r$, we have $\frac{\mu_{v, i}+ \nu_{v, i}}{2} \in d_i \Z$. The discussion is so far general; if we now assume $\pi_v$ is tempered, then for all $\lambda \in X_\bullet(T)$, the character
\[
z \mapsto z^{\langle \mu_v, \lambda \rangle} \bar{z}^{\langle \nu_v, \lambda \rangle}
\]
is unitary, i.e.
\[
Re(\langle \mu_v + \nu_v, X_\bullet(T) \rangle)=0.
\]
$W$-algebraic representations of course have real infinitesimal character, so $\mu_v= -\nu_v$, and therefore in the initial choice of lift we may assume $\mu_{v, i}=- \nu_{v, i}$ for all $i$. Then obviously $\tmu_v$ is $L$-algebraic.

If $\pi$ is only $W$-algebraic, then we have to check that $2 \tmu_v= (2 \mu_v, \sum(\mu_{v, i}- \nu_{v, i}) w_i)$ lies in $X^\bullet(\tT)$. This element of $X^\bullet(T) \oplus X^\bullet(\tZ)$ represents an element of $X^\bullet(\tT)$ if and only if $2 \mu_v$ maps to $\sum_i [\mu_v-\nu_v]_i \bar{w}_i$ in $X^\bullet(Z_G)$. The latter is also the image of $\mu_v- \nu_v$, so it is equivalent to ask that $\mu_v+\nu_v \in X^\bullet(T)$ map to zero in $X^\bullet(Z_G)$. As above, temperedness of $\pi_v$ guarantees this, so we are done.
\endproof
\begin{rmk}\label{arthurgripe1}
This result would immediately extend to totally imaginary fields if we knew that $\omega_{\pi}$ extended to a type $A_0$ Hecke character of $\tZ(\af)$. As noted in Remark \ref{cmdescentextrapolation}, this would follow from the (`CM descent') conjectures of \S \ref{cmdescentsection}. 
\end{rmk}
So we turn to the case of totally real $F$. I am grateful to Brian Conrad for urging a coordinate-free formulation of the obstruction in part $2$ of the proposition.
\begin{prop}\label{totrealalglift}
Now suppose $F$ is totally real, with $\pi, \tomega, \tpi$ as before ($\tomega$ finite-order as in Definition \ref{totrealheckelift}). Alternatively, let $F$ be arbitrary, but assume that $\omega_\pi$ admits a finite-order extension $\tomega$. Continue to assume $\pi_\infty$ is tempered.
\begin{enumerate}
\item If $\pi$ is $L$-algebraic, then it admits a $W$-algebraic lift $\tpi$.
\item Assume Conjecture \ref{locmultone} for the `only if' direction of this statement. Assume $F$ is totally real, and consider this $W$-algebraic lift $\tpi$. Then the images of $\mu_v$ and $\nu_v$ under $X^\bullet(T) \to X^\bullet(Z_G)$ lie in $X^\bullet(Z_G)[2]$, and $\pi$ admits an $L$-algebraic lift if and only if these images are independent of $v \vert \infty$. 
\end{enumerate}
\end{prop}
\proof
We continue with the parameter notation of the previous proof. Let $\pi$ be $L$-algebraic. First we check that since $\tomega$ can be chosen finite-order (automatic in the totally real case, but an additional assumption at complex places) $\mu_v- \nu_v$ maps to zero in $X^\bullet(Z_G)$. If $v$ is imaginary, then $\tomega$ cannot be finite order unless $\omega_{\pi_v}$ is trivial, hence $\mu_v$ and $\nu_v$ themselves are trivial in $X^\bullet(Z_G)$. If $v$ is real, then in $T^\vee \subset G^\vee$ we have the relation
\[
\phi(j)z^{\mu_v} \bar{z}^{\nu_v} \phi(j)^{-1}= z^{\nu_v} \bar{z}^{\mu_v},
\]
writing $\phi$ for the $L$-parameter and $j \in W_{F_v}- W_{\overline{F}_v}$. Now, $\phi(j) \in N_{G^\vee}(T^\vee)$ represents an element $w$ of the Weyl group of $G^\vee$, and this yields the relations $w \mu_v= \nu_v$ and $w \nu_v= \mu_v$. Thus $\mu_v-\nu_v= \mu_v-w\mu_v$, which lies in the root lattice $Q$ of $G$. [This holds for $w\chi-\chi$ for any $\chi \in X^\bullet(T)$: reduce to the case of simple reflections, where it is clear from the defining formula.] Restricting characters to $Z_G$ factors through a perfect duality $Z_G \times X^\bullet(T)/Q \to \Q/\Z$, so $\mu_v-w\mu_v \in Q$ has trivial image in $X^\bullet(Z_G)$.

Any lift with finite-order central character has parameters $\tmu_v= (\mu_v, 0)$ and $\tnu_v= (\nu_v, 0)$; this pair yields a well-defined representation of $W_{\overline{F}_v}$, i.e. $\tmu_v- \tnu_v \in X^\bullet(\tT)$, and its projections to $T$ and $\tZ$ are what they have to be, so they are the only possible parameters; that this extends (possibly non-uniquely) to an $L$-parameter on the whole of $W_{F_v}$ follows from general theory (Langlands Lemma), but we do not need the details of this extension. Now we use the assumption that $\pi_v$ is tempered: as in the previous proposition, we see that $\mu_v= -\nu_v$. Hence $2 \mu_v$ maps to zero in $X^\bullet(Z_G)$, and so $2(\mu_v, 0) \in X^\bullet(\tT)$, i.e. $(\mu_v, 0)$ is $W$-algebraic.

For $F$ totally real, we now show the second part of the proposition. By Lemma \ref{fibers}, if a second lift $\tpi'$ is $L$-algebraic, then 
\[
\tpi' \cong \tpi \cdot \left( \prod_{i=1}^r \alpha_i \circ d_i w_i \right),
\]
where each $\alpha_i^{d_i}$ is a Hecke character of our totally real field, and therefore (since $\tpi$ and $\tpi'$ both have rational infinitesimal character) takes the form $\chi_i \cdot |\cdot|^{d_i y_i}$ for a finite-order character $\chi_i$ and a rational number $y_i$. The infinitesimal character of $\tpi'_v$ (for $v \vert \infty$) then corresponds to $(\mu_v, \sum y_i d_i w_i) \in X^\bullet(\tilde{T})_{\Q}$, which is integral if and only if $y_i \in \frac{1}{d_i} \Z$ and $y_i d_i$ is the image of $\mu_v$ in $\Z/d_i \Z$ for each $i$. The claim follows easily.
\endproof
\begin{cor}\label{counterexamples}
Let $F$ be totally real and $\pi$ be $L$-algebraic on the $F$-group $G$.
\begin{enumerate}
\item Suppose $\pi_\infty$ is tempered. If all $d_i$ are odd, then $\tpi$ is $L$-algebraic. In particular, if the simple factors of $G$ are all type $A_{2n}$ (i.e. $\mr{SL}_{2n+1}$), or of type $E_6$, $E_8$, $F_4$, $G_2$,\footnote{$A_{2n}$ and $E_6$ are the interesting examples here, since in the other cases the adjoint group is simply-connected. The index $[X^\bullet(T):Q]$ for the simply-connected simple group of type $E_6$ is $\Z/3\Z$.} then any $L$-algebraic $\pi$ has an $L$-algebraic lift. 
\item Assume $F \neq \Q$ (still totally real). For $F$-groups $G$ that are simple simply-connected (and split) of type $B_{n}$, $C_n$, $D_{2n}$, and $E_7$, there exist, assuming Conjecture \ref{locmultone} for the pair $(G, \tilde{G})$, $L$-algebraic $\pi$ on $G$, tempered at infinity, that admit no $L$-algebraic lift to $\tG$. 
\end{enumerate}
\end{cor}
\proof
The first part is immediate from the proof of the proposition: for $d_i$ odd, if the image of $2\mu_v$ in $\Z/d_i\Z$ is trivial, then so of course is the image of $\mu_v$. For the second part, we use existence results for automorphic forms that are discrete series at infinity as in Lemma \ref{limmult} and Proposition \ref{sp}. For a semi-simple, simply-connected group, $\rho$ lies in the weight lattice, so discrete series representations will be $L$-algebraic. Choose a quotient $G \onto G' \onto G_{ad}$, with $T'$ the induced maximal torus of $G'$, such that $X^\bullet(T')/Q \cong X^\bullet(Z_G)[2]$ under the isomorphism $X^\bullet(T)/Q \xrightarrow{\sim} X^\bullet(Z_G)$.Then we can find a discrete series representation of $G(F_{\infty})$ such that at two infinite places $v_1$ and $v_2$, the parameters have the form (on $\CC^\times$)
\[
z \mapsto z^{\rho + \mu_{v_i}} \bar{z}^{-\rho-\mu_{v_i}},
\]
where the $\mu_{v_i}$ are distinct in $X^\bullet(T')/Q$.\footnote{We have shifted what was previously denoted $\mu_v$ by $\rho$ in this case because of the description of discrete series $L$-parameters as arising from elements of $\rho + X^\bullet(T)$. Note that $\rho$ mapst to $X^\bullet(Z_G)[2]$ since $2 \rho$ is in the root lattice $Q$.} This ensures that the $\rho+\mu_{v_i}$ have distinct images in $X^\bullet(Z_G)[2]$. It remains therefore to check that the split real forms of these groups, except for $A_{2n-1}$, actually admit discrete series. Type $C_n$ was treated previously. For type $B_n$, the (real) split group is $\mr{Spin}(n, n+1)$, which always admits discrete series (real orthogonal groups $\mr{SO}(p, q)$ have discrete series whenever $pq$ is even). Likewise, for type $D_{n}$, the split form $\mr{Spin}(n, n)$ has discrete series if and only if $n$ is even.  The split real form of $E_7$ also admits discrete series: its maximal compact subgroup is $\mr{SU}(8)$, which contains a compact Cartan isomorphic to $(\mr{S}^1)^7$ (see the table in Appendix $4$ of \cite{knapp:beyond}). 
\endproof
\begin{rmk}\label{arthurgripe2}
\begin{itemize}
\item Arthur constructs in \cite{arthur:LF} a conjectural candidate for the automorphic Langlands group $\mc{L}_F$,\index{s}{$\mc{L}_F$}\index{t}{automorphic Langlands group} also giving an `automorphic' construction of a candidate for the motivic Galois group $\mc{G}_F$. (In \S \ref{motivatedprelims} we will discuss motivic Galois groups in detail, explaining both conjectural and unconditional constructions.) The ability to extend $L$-algebraic representations of $G$ to $L$-algebraic representations of $\tG$ is essential to Arthur's conjectural construction (see \cite[\S 6]{arthur:LF}) of a morphism $\mc{L}_F \to \mc{G}_F$. Roughly, the factor of the motivic Galois group $\mc{G}_F$ associated to $\pi$ (or rather its almost-everywhere system of Hecke eigenvalues; this construction is restricted to suitably `primitive' $\pi$) is defined as a sub-group $\mc{G}_{\pi}$ of $\tG^\vee \times \mc{T}_F$, with $\mc{T}_F$ the Taniyama group (see \S \ref{taniyama} for background on $\mc{T}_F$)\index{s}{$\mc{T}_F$}\index{t}{Taniyama group}; then $\mc{L}_F$ and $\mc{G}_F$ are constructed as suitable fiber products over varying pairs $(G, \pi)$.\footnote{For non-split $G$, one works with the semi-direct product $\tG^\vee \rtimes \mc{T}_F$, with $\mc{T}_F$ acting via its projection to the global Weil group $W_F$.} The natural motivic Galois representation corresponding to `the motive of $\pi$' is the projection to $G^\vee$, and implicit in the construction is the fact that this lifts to $\tG^\vee$. Corollary \ref{counterexamples} then gives concrete examples of automorphic representations (including ones that ought to be `primitive,' arranging the local behavior suitably) for which the Tannakian formalism cannot behave in this way; that is, while Arthur's construction is consistent with our discussion in the CM case (Proposition \ref{cmextension}), it must be modified for totally real fields $F$. 
\item If our lifting results were generalized to quasi-split groups, we could presumably include $D_n$ for $n$ odd in the second part of the Proposition, since the non-split quasi-split form has signature $(n+1, n-1)$, for which the associated orthogonal group admits discrete series. Similarly, we could include type $A_{2n-1}$ by using the quasi-split unitary group $\mr{SU}(n, n)$ instead of $\mr{SL}_{2n}$.
\end{itemize}
\end{rmk}

\section{Galois lifting: the general case}\label{generalGalois}
Now we discuss a general framework for Conrad's lifting problem. We consider lifting problems of the form
\[
\xymatrix{
& H'(\Qlb) \ar[d] \\
\gal{F} \ar@{-->}[ru]^-{\tilde{\rho}} \ar[r]_-{\rho} & H(\Qlb),}
\]
for any surjection $H' \onto H$, with central torus kernel, of linear algebraic groups over $\Qlb$. In this picture, we want to understand when a geometric $\rho$ does or does not admit a geometric lift $\tilde{\rho}$. We will see that the general case reduces to the case where $H'$ and $H$ are connected reductive, and so to reap the psychological benefits of the work of \S \ref{liftingwalg}, we will begin by considering the case $H'= \tG^\vee$, $H= G^\vee$, where $G \subset \tG$ is an inclusion of connected reductive\footnote{Unlike in \S \ref{liftingwalg}, whose results were largely intended to motivate the results of this section, we no longer require $G$ to be semi-simple.} split $F$-groups, constructed by extending the center $Z_G$ of $G$ to a central torus $\tZ$; that is\index{s}{$G$}\index{s}{$\tG$}\index{s}{$\tZ$}\index{s}{$S$}
\[
\tG= (G \times \tZ)/Z_G.
\] 
The quotient $\tZ/Z_G$ is a torus $S$, and we get an exact sequence
\[
1 \to S^\vee \to \tilde{G}^\vee \to  G^\vee \to 1,
\]
and a canonical isomorphism of Lie algebras $\tilde{\mf{g}}^\vee \cong \mf{g}^\vee \oplus \mf{s}^\vee$. We caution the reader that, even if $H'$ and $H$ are assumed connected reductive, this is \textit{not} the most general situation: namely, the center $\tZ$ need not be connected in order to have $\tG^\vee \onto G^\vee$ a surjection with central torus kernel.\footnote{Consider for instance $G= \mr{Spin}_{2n} \subset \tG= \mr{GSpin}_{2n}$, which dualizes to a surjection $\mr{GSO}_{2n} \onto \mr{PSO}_{2n}$.} Nevertheless, at least for $F$ totally imaginary, we will reduce the general analysis to this case. To start, we elaborate on some of the associated group theory, slightly recasting the notation of the previous section to take into account the fact that $G$ may have positive-dimensional center. The exact sequence
\[
1 \to Z_G \to \tZ \to S \to 1
\]
gives rise, by applying $X^\bullet(\cdot)$ and then $\Hom(\cdot, \Z)$, to an exact sequence
\[
1 \to \Hom(X^\bullet(Z_G), \Z) \to X_{\bullet}(\tZ) \to X_{\bullet}(S) \to \Ext^1(X^\bullet(Z_G), \Z) \to 1.
\]
Using invariant factors, it is convenient to fix a basis of $X^\bullet(\tZ)$ such that the inclusion $X^\bullet(S) \subset X^\bullet(\tZ)$ is in coordinates
\[
\xymatrix{
X^\bullet(S) \ar[r] \ar@{=}[d] & X^\bullet(\tZ) \ar[r] \ar@{=}[d] & X^\bullet(Z_G) \ar@{=}[d] \\
\bigoplus_{i=1}^r \Z d_i w_i \ar[r] & \bigoplus_{i=1}^{r+s} \Z w_i \ar[r] & \bigoplus_{i=1}^r \Z/d_i \Z \bar{w}_i \oplus \bigoplus_{j=1}^s \Z \bar{w}_{r+j} \\
}
\]
Let $w_i^*$ denote the dual basis for $X^\bullet(\tZ^\vee)$, and let $\nu_i^*$ denote the basis of $X^\bullet(S^\vee)$ dual to $d_i w_i$; in particular, $w_i^*$ maps to $d_i \nu_i^*$ under $X^\bullet(\tZ^\vee) \to X^\bullet(S^\vee)$. 

We work with the maximal torus $\tT= (T \times \tZ)/Z_G$ of $\tG$, and deduce from its definition an exact sequence
\[
1 \to X^\bullet(\tT) \to X^\bullet(T) \oplus X^\bullet(\tZ) \to X^\bullet(Z_G) \to 1,
\]
and thus a (crucial) exact sequence \label{defectsequence}
\begin{equation}\label{defect}
1 \to \Hom(X^\bullet(Z_G), \Z) \to X^\bullet(T^\vee) \oplus X^\bullet(\tZ^\vee) \to X^\bullet(\tT^\vee) \to \Ext^1(X^\bullet(Z_G), \Z) \to 1.
\end{equation}
Note also that there is a canonical isomorphism $X^\bullet(\tZ^\vee) \xrightarrow{\sim} X^\bullet(\tG^\vee)$, which follows from exactness of the row in the diagram:
\[
\xymatrix{
& & 1 \ar[d] & & \\
& & S^\vee \ar[d] & & \\
1 \ar[r] & G^\vee_{sc} \ar[r] \ar[dr] & \tG^\vee \ar[r] \ar[d] & \tZ^\vee \ar[r] & 1\\
& & G^\vee \ar[d] & & \\
& & 1 & & \\
}
\]
where $G^\vee_{sc}$ denotes the simply-connected cover of the derived group of $G^\vee$.\index{s}{$G^\vee_{sc}$}
Recall that we index Sen operators $\Theta_{\rho, \iota}$ by embeddings $\iota \colon \Qlb \into \CC_{F_v}$, where $v$ is a place above $\ell$. Given a geometric $\rho \colon \gal{F} \to G^\vee(\Qlb)$, Conrad's result (see Theorem \ref{conradliftingP}, and Remark \ref{aeunr}) implies that a lift $\tilde{\rho}$ is geometric if and only if it is Hodge-Tate at all places above $\ell$; we will use this from now on without comment. That is, we need to arrange $\tilde{\rho}$ such that each Sen operator $\Theta_{\tilde{\rho}, \iota} \in \mr{Lie}(\tG^\vee)_{\iota}$ is conjugate to an element of $\mr{Lie}(\tT^\vee)_{\iota}$ that pairs integrally with all of $X^\bullet(\tT^\vee)$\footnote{This condition is independent of the way $\Theta_{\tilde{\rho}, \iota}$ is conjugated into $\mr{Lie}(\tT^\vee)_\iota$, since: $(1)$ we already know that $\Theta_{\tilde{\rho}, \iota}$ pairs integrally with the roots, which all lie in $X^\bullet(T^\vee)$; $(2)$ the ambiguity in conjugating into $\mr{Lie}(\tT^\vee)_{\iota}$ is an element of the Weyl group; $(3)$ for any weight $\lambda \in X^\bullet(\tT^\vee)$ and $w$ in the Weyl group, $w \lambda-\lambda$ lies in the root lattice. Compare the proof of Proposition \ref{totrealalglift}.} under the natural map
\[
X^\bullet(\tT^\vee) \xrightarrow{\mr{Lie}} \Hom( \mr{Lie}(\tT^\vee), \Qlb).
\] 
Here is our starting-point:  
\begin{lemma}
There exists some lift $\tilde{\rho}$ of $\rho$. Any other lift is of the form 
\[
\tilde{\rho}(\sum_{i=1}^r \nu_i \circ \chi_i) \colon g \mapsto \tilde{\rho}(g) \cdot \prod_{i=1}^r (\nu_i \circ \chi_i)(g),
\]
where the $\nu_i= d_i w_i$ range over the above basis of $ X_\bullet(S^\vee)$, and each $\chi_i \colon \gal{F} \to \Qlb^\times$ is a continuous character.
\end{lemma}
\proof
A lift exists by Proposition \ref{tatelift}. Although continuous cohomology does not in general have good $\delta$-functorial properties, short exact sequences do give (not very) long exact sequences on $H^0$ and $H^1$, so any lift has the form $\tilde{\rho}(\chi)$ for some $\chi \colon \gal{F} \to S^\vee(\Qlb)$. We compose with the dual characters $\nu_i^*$ in $X^\bullet(S^\vee)$ to put $\chi$ in the promised form.
\endproof
The following lemma is the general substitute for choosing lifts with finite-order Clifford norm in the spin examples; this result is also implicit in the proof of \ref{tatelift}, but a little warm-up with our notation is perhaps helpful.
\begin{lemma}\label{twist}
Let $\rho \colon \gal{F} \to G^\vee(\Qlb)$ be a geometric representation. Then there exists a lift $\tilde{\rho} \colon \gal{F} \to \tilde{G}^\vee(\Qlb)$ such that, for all $v \vert \ell$ and all $\iota \colon \Qlb \into \CC_{F_v}$, the Sen operator $\Theta_{\tilde{\rho}, \iota}$ pairs integrally with all of $X^\bullet(\tZ^\vee) \cong X^\bullet(\tG^\vee)$.
\end{lemma}
\proof
It suffices to find a lift whose composition with all elements of $X^\bullet(\tZ^\vee)$ is Hodge-Tate. We use the bases of the various character groups specified above. In particular, composing an initial lift $\tilde{\rho}$ with the various $w_i^* \in X^\bullet(\tG^\vee)$, $i=1, \ldots, r$, we can write
\[
w_i^* \circ \tilde{\rho}= \chi_i^{d_i} \chi_{i, 0} \colon \gal{F} \to \Qlb^\times,
\]
where the $\chi_i$ and $\chi_{i, 0}$ are Galois characters with $\chi_{i, 0}$ finite-order. Then we consider the new lift
\[
\tilde{\rho}'= \tilde{\rho}( \sum_{i=1}^r (d_i w_i) \circ \chi_i^{-1}),
\]
which has the advantage that $w_i^* \circ \tilde{\rho}'= (w_i^* \circ \tilde{\rho})\cdot \chi_i^{-d_i}$ is finite-order for all $i=1, \ldots, r$. Moreover, for the characters $w_{r+j}^*$, $j=1, \ldots, s$, namely, the sub-module $\Hom(X^\bullet(Z_G), \Z) \subset X^\bullet(\tZ^\vee)$, the compositions $w_{r+j}^* \circ \tilde{\rho}'$ are all geometric, since $\rho$ is. Therefore $\alpha \circ \tilde{\rho}'$ is geometric for all $\alpha \in X^\bullet(\tZ^\vee)$, as desired.
\endproof
Returning to equation \eqref{defect} on page \pageref{defectsequence}, we see that the obstruction to geometric lifts then comes from $\Ext^1(X^\bullet(Z_G)_{tor}, \Z)$. For any weight $\lambda \in X^\bullet(\tT^\vee)$, there is a positive integer $d \in \Z$ such that $d \lambda \in X^\bullet(T^\vee) \oplus X^\bullet(\tZ^\vee)$, so with a $\tilde{\rho}$ as produced by Lemma \ref{twist}, $\Theta_{\tilde{\rho}, \iota}$ pairs integrally with $d \lambda$ for all $\iota$. We obtain a well-defined class, independent of the choice of $\tilde{\rho}$ as constructed in the proof of Lemma \ref{twist} (namely, with the compositions $w_i^* \circ \tilde{\rho}$ finite-order for $i = 1, \ldots, r$),
\[
\langle \lambda, \Theta_{\tilde{\rho}, \iota} \rangle \in \Q /\Z,
\]
which can also be interpreted as the common value modulo $\Z$ of the eigenvalues of $\Theta_{r_\lambda \circ \tilde{\rho}, \iota}= \mr{Lie}(r_\lambda) \circ \Theta_{\tilde{\rho}, \iota}$; here, and in what follows, we denote by $r_\lambda$ the irreducible representation of $\tG^\vee$ associated to the (dominant) weight $\lambda \in X^\bullet(\tT^\vee)$. This pairing factors through a map
\[
\Ext^1(X^\bullet(Z_G), \Z) \to \Q/\Z,
\]
and since for any finite abelian group $A$, the long exact sequence associated to $0 \to \Z \to \Q \to \Q/\Z \to 0$ yields an isomorphism
\[
\Hom(A, \Q/\Z) \xrightarrow{\sim} \Ext^1(A, \Z),
\] 
we can make the following definition:
\begin{defn}\label{torsionclass} Let $\theta_{\rho, \iota}$\index{s}{$\theta_{\rho, \iota}$} be the element of $X^\bullet(Z_G)_{tor}$ canonically corresponding to the above map $\Ext^1(X^\bullet(Z_G), \Z) \to \Q/\Z$.
\end{defn}

To make further progress, we need to assume that $\rho$ satisfies certain Hodge-Tate weight symmetries. 
\begin{hypothesis}\label{HTsym}
Let $H$ be a linear algebraic group and $\rho \colon \gal{F} \to H(\Qlb)$ a geometric Galois representation with connected reductive algebraic monodromy group $H_{\rho}=(\overline{\rho(\gal{F})})^{Zar}$. We formulate the following Hodge-Tate symmetry hypothesis for such a $\rho$: 
\begin{itemize}
\item Let $r$ be any irreducible algebraic representation $r \colon H_{\rho} \to \mr{GL}(V_r)$. Then:
\begin{enumerate}
\item For $\tau \colon F \into \Qlb$, the set $\mr{HT}_\tau(r \circ \rho)$ depends only on $\tau_0= \tau|_{F_{cm}}$.
\item Writing $\mr{HT}_{\tau_0}$ for this set common to all $\tau$ above $\tau_0$, there exists an integer $w$ such that\footnote{Interpreted in the obvious way.} 
\[
\mr{HT}_{\tau_0 \circ c}(r \circ \rho)= w- \mr{HT}_{\tau_0}(r \circ \rho),
\] 
for $c$ the unique complex conjugation on $F_{cm}$.
\end{enumerate}
\end{itemize}
For $\rho$ whose algebraic monodromy group is reductive but not necessarily connected, the corresponding hypothesis is simply that some finite restriction (with connected monodromy group) of $\rho$ satisfies the above.
\end{hypothesis}
In particular, we note for later use that the lowest weight in $\mr{HT}_{\tau_0 \circ c}(r \circ \rho)$ is $w$ minus the highest weight in $\mr{HT}_{\tau_0}(r \circ \rho)$. We will see that to establish lifting results for geometric representations $\rho \colon \gal{F} \to G^\vee(\Qlb)$, we will only need to know Hypothesis \ref{HTsym} for some easily identifiable finite collection of compositions $r_\lambda \circ \rho$, $\lambda \in X^\bullet(T^\vee)$, but we do not make that explicit here. Most important, Hypothesis \ref{HTsym} should in fact be no additional restriction on $\rho$, because of the following conjecture, which would follow from various versions of the conjectural Fontaine-Mazur-Langlands correspondence and Tate conjecture (Conjecture \ref{FML}):
\begin{conj}\label{galoisdescent}
Let $F$ be a number field, and let $\rho \colon \gal{F} \to \mr{GL}_N(\Qlb)$ be an irreducible geometric Galois representation. Then:
\begin{enumerate}
\item For $\tau \colon F \into \Qlb$, the set $\mr{HT}_\tau(\rho)$ depends only on $\tau_0= \tau|_{F_{cm}}$. (This will still hold if $\rho$ is geometric but reducible.) 
\item Writing $\mr{HT}_{\tau_0}$ for this set common to all $\tau$ above $\tau_0$, there exists an integer $w$ such that
\[
\mr{HT}_{\tau_0 \circ c}(\rho)= w- \mr{HT}_{\tau_0}(\rho),
\] 
for $c$ the unique complex conjugation on $F_{cm}$.
\end{enumerate}
\end{conj}
Unfortunately, for an abstract Galois representation, this conjecture will be extremely difficult to establish. The next lemma explains it in the automorphic case; for a motivic variant, see Corollary \ref{motivichodgesymmetry}.\footnote{But note for now that in the most basic motivic cases, where the Galois representation is given by $H^j(X_{\overline{F}}, \Qlb)$ for some smooth projective variety $X/F$, the Hodge-Tate symmetries are immediate (even when this representation is reducible) from the $\ell$-adic comparison isomorphism of \cite{faltings:crystalline}; if $j$ is even, the symmetries similarly hold for primitive cohomology. What is not obvious is that if this Galois representation decomposes, that the irreducible factors all satisfy the conjecture.}
\begin{lemma}
Suppose $\rho \colon \gal{F} \to \mr{GL}_N(\Qlb)$ is an irreducible geometric representation. If $\rho$ is automorphic in the sense of Conjecture \ref{FML}, corresponding to a cuspidal automorphic representation $\pi$ of $\mr{GL}_N(\af)$, and if we assume that Proposition \ref{cmdescent} is unconditional for $\pi$ (i.e., admit Hypothesis \ref{algconj}), then Conjecture \ref{galoisdescent} holds for $\rho$.
\end{lemma}
\proof
This is immediate from the passage between infinity-types and Hodge-Tate weights, Conjecture \ref{cmdescent}, and Clozel's archimedean purity lemma (which was proven as part of Proposition \ref{tensordescent}). 
\endproof

We can now understand when geometric lifts ought to exist; the proof proceeds by reduction to the following key case:
\begin{prop}\label{fullmonodromylift}
Let $F$ be a totally imaginary field, and let $\rho \colon \gal{F} \to G^\vee(\Qlb)$ be a geometric representation with algebraic monodromy group equal to the whole of $G^\vee$. Assume $\rho$ satisfies Hypothesis \ref{HTsym}. Then $\rho$ admits a geometric lift $\tilde{\rho} \colon \gal{F} \to \tG^\vee(\Qlb)$.
\end{prop}
\proof
Choose a lift $\tilde{\rho}$ as supplied by Lemma \ref{twist}. Recall that our weight-bookkeeping is done within the diagram
\[
\xymatrix{
 & 0 & & \\
X^\bullet(\tZ^\vee) \ar[r] & X^\bullet(S^\vee) \ar[u] \ar[r] & \Ext^1(X^\bullet(Z_G), \Z) \ar[r] & 0 \\
 & X^\bullet(\tT^\vee) \ar[u] & & \\
 & X^\bullet(T^\vee) \ar[u] & & \\
 & 0 \ar[u] & & \\
}
\]
We have the elements $\nu_i^* \in X^\bullet(S^\vee)$ dual to $d_i w_i \in X^\bullet(S)$. Their images in $\Ext^1(X^\bullet(Z_G), \Z)$ form a basis. Let $\lambda_i$ be a (dominant weight) lift to $X^\bullet(\tT^\vee)$ that also satisfies $\langle \lambda_i, d_j w_j \rangle= \delta_{i j}$. Note that the value $\langle \Theta_{\tilde{\rho}, \iota}, \lambda_i \rangle \in \Q/\Z$ does not depend on the choice of lift of $\nu_i^*$, and it clearly lies in $\frac{1}{d_i}\Z/\Z$, so we write it in the form $\frac{k_{\iota, i}}{d_i}+ \Z$ for an integer $k_{\iota, i}$. By considering the geometric representation $r_{d_i \lambda_i}\circ \tilde{\rho}$, we find that $\Theta_{r_{d_i \lambda_i} \circ \tilde{\rho}, \iota}$ has eigenvalues that depend only on $\tau_\iota$ (by Lemma \ref{HTS}; see that lemma for the notation $\tau_\iota$ as well), and thus we can write $k_{\tau, i}$ in place of $k_{\iota, i}$. (Equivalently, we can work with the elements $\theta_{\rho, \tau} \in X^\bullet(Z_G)_{tor}$.) These classes $\frac{k_{\tau, i}}{d_i} \mod \Z$ serve as both highest and lowest $\tau$-labeled Hodge-Tate weights (modulo $\Z$) for $r_{\lambda_i} \circ \tilde{\rho}$; we deduce that, modulo $d_i \Z$, the highest and lowest $\tau$-labeled weights of $r_{d_i \lambda_i} \circ \tilde{\rho}$ are both congruent to $k_{\tau, i}$\footnote{Note that all eigenvalues of $\Theta_{r_{\lambda_i} \circ \tilde{\rho}}$ are congruent modulo $\Z$; this does not imply that all elements of $\mr{HT}_{\tau}(r_{d_i \lambda_i} \circ \tilde{\rho})$ are congruent modulo $d_i \Z$, but this congruence does hold for the highest and lowest weights.}. 

Now, the geometric Galois representations $(r_{d_i \lambda_i}) \circ \tilde{\rho}$ for $i=1, \ldots, r$ are irreducible, because we have assumed that $G^\vee$ is the monodromy group of $\rho$. Applying Hypothesis \ref{HTsym}, we deduce that $k_\tau$ depends only on $\tau_0= \tau|_{F_{cm}}$, along with the symmetry relation 
\[
k_{\tau_0, i}+ k_{\tau_0 \circ c, i} \equiv w_i \mod d_i
\]
for some integer $w_i$ (for all $\tau_0 \colon F_{cm} \into \Qlb$). 

This relation allows us, by Lemma \ref{galchar}, to find Galois characters $\hat{\psi}_i \colon \gal{F} \to \Qlb^\times$ with $\mr{HT}_\tau(\hat{\psi}_i) \in \frac{k_{\tau, i}}{d_i} + \Z$ for all $\tau$. We then form the twist $\tilde{\rho}':= \tilde{\rho}(\sum_i (d_i w_i) \circ \hat{\psi}^{-1}_i)$; recall that $\langle d_i w_i, \lambda_j \rangle= \delta_{ij}$. This new lift is then geometric:
\[
\langle \Theta_{\tilde{\rho}', \iota}, \lambda_i \rangle= \langle \Theta_{\tilde{\rho}, \iota}, \lambda_i \rangle+ \langle \Theta_{\sum (d_j w_j) \circ \hat{\psi}_j^{-1}, \iota}, \lambda_i \rangle \equiv \frac{k_{\tau_\iota, i}}{d_i} - \frac{k_{\tau_\iota, i}}{d_i} \equiv 0 \mod \Z.
\]
(Recall it suffices to check $\tilde{\rho}'$ is Hodge-Tate, by Conrad's result, quoted here as Theorem \ref{conradliftingP}.)
\endproof
\begin{cor}\label{fullmonodromytotreal}
Let $\rho \colon \gal{F} \to G^\vee(\Qlb)$ be geometric with algebraic monodromy group $G^\vee$; maintain the notation of the previous proof, but now suppose $F$ is totally real. Then $\rho$ has a geometric lift if and only if for varying $\iota$, the elements $\theta_{\rho, \iota}$ (see Definition \ref{torsionclass}) in $X^\bullet(Z_G)_{tor}$ are independent of $\iota$.
\end{cor}
\begin{rmk}
More concretely, to determine whether $\rho$ has a geometric lift, apply the following criterion:
\begin{enumerate}
\item Ignore any $i$ for which $d_i$ is odd; these do not obstruct geometric lifting;
\item Then for fixed $i$ the integers $k_{\tau, i} \mod d_i$ are all $\frac{w_i}{2}$ translated by a two-torsion class in $\Z/d_i \Z$;
\item $\rho$ has a geometric lift if and only if each of these classes $k_{\tau, i}$ (or, equivalently, the associated two-torsion class) is independent of $\tau$.
\end{enumerate}
\end{rmk}
\proof
By the previous proof and Lemma \ref{galtotreal}, $\rho$ has a geometric lift if and only if the classes $k_{\tau, i} \mod d_i$ (fixed $i$, varying $\tau$) are independent of $\tau$. The weight-symmetry relation becomes $2k_{\tau, i} \equiv w_i \mod d_i$, and the claim follows easily.
\endproof
Over totally imaginary fields, we can now reduce the general lifting problem to the special case of full monodromy:
\begin{thm}\label{anymonodromylift}
Let $F$ be totally imaginary, and let $\pi \colon H' \onto H$ be any surjection of linear algebraic groups with central torus kernel. Suppose $\rho \colon \gal{F} \to H(\Qlb)$ is a geometric representation, with arbitrary image, satisfying Hypothesis \ref{HTsym}. Then $\rho$ admits a geometric lift $\tilde{\rho} \colon \gal{F} \to H'(\Qlb)$. 
\end{thm}
\proof
We may assume $H'$ and $H$ are reductive since $\pi$ induces an isomorphism on unipotent parts (we take this observation from Conrad, who has exploited this reduction in the arguments of \cite{conrad:dualGW}).

We next show that the theorem holds if $H_{\rho}:= \overline{\rho(\gal{F})}^{Zar} \subset H$ is connected. In that case, let $H'_{\rho}$ be the preimage in $H'$ of $H_{\rho}$. Then $H'_{\rho} \to H_{\rho}$ is a surjection of connected reductive groups with central torus kernel, and we may write $H_{\rho}= G^\vee$, $H'_{\rho}= \tG^\vee$ where $G \subset \tG$ is an inclusion of connected reductive groups of the form $\tG= (G \times \tZ)/Z_G$ for some inclusion of $Z_G$ into a multiplicative group $\tZ$. If $\tZ$ is not connected, Proposition \ref{fullmonodromylift} does not immediately apply, so we embed $\tZ$ into a torus $\widetilde{\mathbf{Z}}$, with corresponding inclusions $G \subset \tG \subset \widetilde{\mathbf{G}}$. Dually, $\widetilde{\mathbf{G}}^\vee \onto G^\vee$ is a quotient to which we can apply Proposition \ref{fullmonodromylift}, and then projecting from $\widetilde{\mathbf{G}}^\vee$ to $\tG^\vee$, we obtain a geometric lift $\tilde{\rho} \colon \gal{F} \to H'_{\rho}(\Qlb) \subset H'(\Qlb)$ of $\rho$.

For arbitrary $\rho$ (i.e. $H_{\rho}$ not necessarily connected), let $F'/F$ be a finite extension such that $\overline{\rho(\gal{F'})}^{Zar}= H_{\rho}^0$ is connected. Over $F'$, we can therefore find a geometric lift $\tilde{\rho}_{F'} \colon \gal{F'} \to H'_{\rho}(\Qlb) \subset H'(\Qlb)$, letting $H'_{\rho}$ as before denote the preimage of $H_{\rho}$ in $H'$. Thus, for \textit{any} lift $\tilde{\rho}_0 \colon \gal{F'} \to H'(\Qlb)$ of $\rho|_{\gal{F'}}$, there exists a character $\hat{\psi} \colon \gal{F'} \to S^\vee(\Qlb)$ such that $\tilde{\rho}_0 \cdot \hat{\psi}$ is geometric. In particular, letting $\tilde{\rho} \colon \gal{F} \to H'(\Qlb)$ be any lift over $F$ itself (with rational Hodge-Tate-Sen weights), there is a $\hat{\psi}_{F'} \colon \gal{F'} \to S^\vee(\Qlb)$ such that $\tilde{\rho}|_{\gal{F'}} \cdot \hat{\psi}_{F'}$ is geometric. But by Corollary \ref{galchar}, there is a Galois character $\hat{\psi} \colon \gal{F} \to S^\vee(\Qlb)$ whose labeled Hodge-Tate-Sen weights descend those of $\hat{\psi}_{F'}$. It follows immediately that $\tilde{\rho} \cdot \hat{\psi} \colon \gal{F} \to H'(\Qlb)$ is a geometric lift of $\rho$.
\endproof
\begin{rmk}\label{liftsupp}
\begin{enumerate}
\item The theorem also lets us make explicit precisely which sets of labeled Hodge-Tate weights can be achieved in a geometric lift $\tilde{\rho}$. We will exploit this in \S \ref{hyperlift}.
\item We will not treat the general totally real case here, since a somewhat different approach seems more convenient in that case. For a succinct, coordinate-free treatment of the totally real case, and another perspective on the arguments of this section, see \cite{stp:parities}, where the following general result is obtained: for any $H' \onto H$ as in our lifting setup, we write $H^0= G^\vee$ and $(H')^0= \tG^\vee$, where $\tG= (G \times \tZ)/Z_G$; here $\tZ$ is not necessarily connected. We can define in this generality elements
\[
\theta_{\rho, \tau} \in \mr{coker} \left( X^{\bullet}(\tZ)_{tor} \to X^{\bullet}(Z_G)_{tor} \right),
\]
and $\rho$ admits a geometric lift to $H'$ if and only if the $\theta_{\rho, \tau}$ are independent of $\tau$. The argument described in \cite{stp:parities} has the disadvantage of not making explicit the parity obstruction (assuming `Hodge symmetry') found in Proposition \ref{fullmonodromytotreal}; for this reason we have retained the two different expositions of the totally real case.
\end{enumerate}
\end{rmk}
The method of proof of Proposition \ref{fullmonodromylift} and Theorem \ref{anymonodromylift} also implies the following local result, which emerged from a conversation with Brian Conrad. A proof of a more difficult result, where $F_v$ is replaced by a $p$-adic field with algebraically closed residue field, will appear in \S $3.6$ of \cite{chai-conrad-oort:cm}. 
\begin{cor}\label{HTlift}
Let $H' \onto H$ be a central torus quotient, and let $\rho \colon \gal{F_v} \to H(\Qlb)$ be a Hodge-Tate representation of $\gal{F_v}$, for $F_v/\Q_\ell$ finite. Then there exists a Hodge-Tate lift $\tilde{\rho} \colon \gal{F_v} \to H'(\Qlb)$.
\end{cor} 
\proof
We sketch a proof, replacing appeal to Hypothesis \ref{HTsym} and the existence of certain Hecke characters by the simpler observation that local class field theory lets us find the necessary twisting characters $\hat{\psi} \colon \gal{F_v} \to \Qlb^\times$ by hand. That is, $\mc{O}_{F_v}$ sits in $I_{F_v}^{ab}$ in finite index, and on the former we can define the character
\[
x \mapsto \prod_{\tau \colon K \into \Qlb} \tau(x)^{k_\tau}
\]
for any integers $k_\tau$; then up to a finite-order character we take a $d^{th}$ root for an integer $d$ to build (having extended to all of $\gal{F_v}$) characters with any prescribed set of rational $\tau$-labeled Hodge-Tate weights ($\tau$ running over all $F_v \into \Qlb$). This observation suffices (invoking Lemma \ref{HTS} in the full generality of Hodge-Tate, rather than merely de Rham, representations) for the previous arguments to carry through.
\endproof
This Corollary combined with \cite[Proposition 6.5]{conrad:dualGW} implies the following stronger result:
\begin{cor}
Let $H' \onto H$ be a central torus quotient, and let $\rho \colon \gal{F_v} \to H(\Qlb)$ be a representation of $\gal{F_v}$, for $F_v/\Q_\ell$ finite, satisfying a basic $p$-adic Hodge theory property\footnote{Namely: crystalline, semi-stable, de Rham, or Hodge-Tate.} $\mbf{P}$. Then there exists a lift $\tilde{\rho} \colon \gal{F_v} \to H'(\Qlb)$ also satisfying $\mbf{P}$.
\end{cor}

\section{Applications: comparing the automorphic and Galois formalisms}\label{descentcomparison}
In \S \ref{generalGalois}, we took the Fontaine-Mazur-Langlands conjecture (or rather its weakened form Hypothesis \ref{HTsym}) relating geometric Galois representations $\gal{F} \to \mr{GL}_n(\Qlb)$ to $L$-algebraic automorphic representations of $\mr{GL}_n(\af)$ as input to establish some of our lifting results. Now we want to apply these lifting results to give some evidence for the relationship between automorphic forms and Galois representations on groups other than $\mr{GL}_n$. We will touch on the Buzzard-Gee conjecture, certain cases of the converse problem, and some general thoughts about comparing descent problems on the ($\ell$-adic) Galois and automorphic sides. We first digress to discuss what `automorphy' of an ${}^L G(\Qlb)$-valued representation even means.
\subsection{Notions of automorphy}\label{automorphynotions}
As usual we have fixed $\iota_\infty \colon \Qb \into \CC$ and $\iota_\ell \colon \Qb \into \Qlb$. It is sometimes more convenient simply to fix an isomorphism $\iota_{\ell, \infty} \colon \CC \to \Qlb$, and to regard $\Qb$ as the subfield of algebraic numbers in $\CC$. $G$ is a connected reductive $F$-group, and we take a $\Qb$-form of the $L$-group ${}^L G$. For an automorphic representation $\pi$ of $G(\af)$ and $\rho \colon \gal{F} \to {}^L G(\Qlb)$, always assumed continuous and composing with $\gal{F}$-projection to the identity, there are (at least) four relations between such $\rho$ and $\pi$ that might be helpful, of which only the first two can be stated unconditionally. First we describe analogues that restrict to either the automorphic or Galois side, borrowing some ideas from \cite{lapid:sl_n} (and, inevitably, \cite{langlands-labesse}). For automorphic representations $\pi$ and $\pi'$ of $G(\af)$, three notions of similarity are:\index{t}{equivalence of automorphic representations and $\ell$-adic representations}
\begin{itemize}
\item $\pi \sim_w \pi'$ if almost everywhere locally, the (unramified, say) $L$-parameters are $G^\vee(\CC)$-conjugate.\index{s}{$\sim_w$}
\item $\pi \sim_{w, \infty} \pi'$ if almost everywhere locally, and at infinity, the $L$-parameters are conjugate. This also makes sense unconditionally, since archimedean local Langlands is known.\index{s}{$\sim_{w, \infty}$}
\item $\pi \sim_{ew} \pi'$ if everywhere locally, the $L$-parameters are $G^\vee(\CC)$-conjugate; this only makes sense if one knows the local Langlands conjecture for $G$.\index{s}{$\sim_{ew}$}
\item $\pi \sim_{s} \pi'$ is the most fanciful condition: if the conjectural automorphic Langlands group $\mc{L}_F$ exists, so $\pi$ and $\pi'$ give rise to representations $\mc{L}_F \to {}^L G(\CC)$, then this condition requires these representations to be globally $G^\vee(\CC)$-conjugate.\index{s}{$\sim_s$}
\end{itemize}
A fifth notion would compare $L$-parameters in a particular finite-dimensional representation of $G^\vee$.

We can unconditionally make the same sort of comparisons between $\ell$-adic Galois representations, writing $\rho \sim_w \rho'$, $\rho \sim_{w, \infty} \rho'$, $\rho \sim_{ew} \rho'$, and $\rho \sim_s \rho'$. By equivalence `at infinity' here, we mean that at real places the actions of complex conjugation are conjugate, and at places above $\ell$, the associated Sen operators (i.e. labeled Hodge-Tate data) are conjugate. Since it is conjectured, but totally out of reach, that frobenius elements act semi-simply in a geometric Galois representation, we should only compare `frobenius semi-simplifications' in these definitions of local equivalence. For some nice examples, Lapid's paper (\cite{lapid:sl_n}) studies the difference between $\sim_w, \sim_{ew}$, and $\sim_s$ for certain Artin representations (and, when possible, the corresponding comparison on the automorphic side).

We then have corresponding ways to relate an $\ell$-adic $\rho$ and an automorphic $\pi$:
\begin{itemize}
\item  $\rho \sim_w \pi$
\item $\rho \sim_{w, \infty} \pi$ 
\item  $\rho \sim_{ew} \pi$ 
\item  $\rho \sim_s \pi$
\end{itemize}
First, write $\rho \sim_w \pi$ if for almost all unramified $v$ (for $\rho$ and $\pi$), $\rho|_{W_{F_v}}^{ss}$ is $G^\vee(\Qlb)$-conjugate to $\rec_v(\pi_v) \colon W_{F_v} \to G^\vee(\CC) \rtimes \gal{F}$; implicit is the assumption that the local parameter lands in $G^\vee(\Qb) \rtimes \gal{F}$, so that we can apply $\iota_\ell \circ \iota_{\infty}^{-1}$. Note that this definition does not distinguish between $\pi$ and other elements of its (conjectural) global $L$-packet $L(\pi)$. The other relations are straightforward modifications (for compatibility with complex conjugation, we take the condition in Conjecture $3.2.1$ of \cite{buzzard-gee:alg}), except for $\rho \sim_s \pi$. Writing $\mc{G}_{F, E}(\sigma)$ for the motivic Galois group for motives\footnote{Either for absolutely Hodge cycles, for motivated cycles, or, assuming the standard conjectures, for homological cycles. Again, we will deal more precisely with motivic Galois groups in \S \ref{motivatedprelims}.}over $F$ with $E$-coefficients, using a Betti realization via $\sigma \colon F \into \CC$, one might hope that after fixing $E \into \CC$ as well, one would obtain a map of pro-reductive groups over $\CC$, $\mc{L}_F \to (\mc{G}_{F, E})(\sigma) \otimes_E {\CC}$ (compare Remark \ref{arthurgripe2}). If $\pi$ corresponds to a representation $\rec(\pi)$ of $\mc{L}_F$, and $\rho$ arises from (completing at some finite place of $E$) a representation $\rho_E$ of $\mc{G}_{F, E}(\sigma)$, it makes sense to ask whether $\rec(\pi)$ factors through $\mc{G}_{F, E}(\sigma)(\CC)$, and whether the resulting representation is globally $G^\vee(\CC)$-conjugate to (the complexification via $E \into \CC$ of) $\rho_E$. Of course, any discussion of the automorphic Langlands group and its relation with the motivic Galois group is pure speculation; but these heuristics do provide context for the basic problems raised in Question \ref{modularlifting} and Conjecture \ref{motivicliftconjecture}, as well as the work of \S \ref{hyperlift}.

Although we don't actually require it, it is helpful to keep in mind a basic lemma of Steinberg, which implies that $\sim_w$ can be checked by checking in all finite-dimensional representations:
\begin{lemma}\label{Steinberg}
Let $x$ and $y$ be two semi-simple elements of a connected reductive group $G^\vee$ over an algebraically closed field of characteristic zero. If $x$ and $y$ are conjugate in every (irreducible) representation of $G^\vee$, or even merely have the same trace, then they are in fact conjugate in $G^\vee$.
\end{lemma}
\proof
For semi-simple groups, this is Corollary $3$ (to Theorem $2$) in Chapter $3$ of \cite{steinberg:conjugacyclasses}; the proof extends to the reductive case (and even more generally, see Proposition $6.7$ of \cite{borel:L}). The key point is that the characters of finite-dimensional representations of $G^\vee$ restrict to a basis of the ring of Weyl-invariant regular functions on a maximal torus; these in turn separate conjugacy classes in the torus.
\endproof
\subsection{Automorphy of projective representations} 
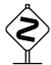
\hspace*{\fill}

\dbend \\

Throughout this section, \textit{we assume the Fontaine-Mazur-Langlands conjecture on automorphy of geometric} ($\mr{GL}_N$-valued) \textit{Galois representations}. It suffices to take a version that matches unramified (almost everywhere) and Hodge-theoretic parameters; to be precise, assume Part 3 of Conjecture \ref{FML}, and note that this includes the requirement that cuspidality is equivalent to irreducibility under the automorphic-Galois correspondence. We will show how our lifting results-- both automorphic and Galois-theoretic-- give rise to a `Fontaine-Mazur-Langlands'-type correspondence between algebraic automorphic representations of $\mr{SL}_N(\af)$ and $\mr{PGL}_N(\Qlb)$-valued geometric $\gal{F}$-representations.\footnote{For some unconditional results in this direction, see \cite[\S 4]{stp:parities}.} The starting point is the following consequence of the results of \S \ref{generalGalois}:
\begin{cor}
Let $F$ be totally imaginary. Then any geometric $\rho \colon \gal{F} \to \mr{PGL}_n(\Qlb)$ is weakly automorphic, i.e. there exists an $L$-algebraic automorphic representation $\pi$ of $\mr{SL}_n(\af)$ such that $\rho \sim_{w, \infty} \pi$ (or $\rho \sim_{ew} \pi$, if we assume a form of Fontaine-Langlands-Mazur that matches local factors everywhere).
\end{cor}
\proof
We have seen that $\rho$ lifts to a geometric $\tilde{\rho} \colon \gal{F} \to \mr{GL}_n(\Qlb)$, which by assumption is automorphic, corresponding to some $\tpi$ on $\mr{GL}_n/F$. The irreducible constituents of $\tpi|_{\mr{SL}_n(\af)}$ form a global $L$-packet whose local unramified parameters correspond to those of $\rho$. 
\endproof
We now give descent arguments that extend this automorphy result to $F$ totally real. First we need a couple of elementary lemmas.
\begin{lemma}\label{cft}
Let $L/F$ be a cyclic, degree $d$, extension of number fields, with $\sigma$ a generator of $\Gal(L/F)$. Let $\chi$ be a Hecke character of $L$, and let $\delta$ be any Hecke character whose restriction to $C_F \subset C_L$ is $\delta= \delta_{L/F}$, a fixed order $d$ character that cuts out the extension $L/F$. Assume that $\chi^{1+\sigma + \ldots +\sigma^{d-1}}=1$. Then for a unique integer $i=0, \ldots, d-1$, $\chi \delta^i$ is of the form $\psi^{\sigma-1}$ for a Hecke character $\psi$ of $L$. 
\end{lemma}
\proof
We may assume $\chi$ is unitary. Write $C^D$ as usual for the Pontryagin dual of a locally compact abelian group $C$. We have the following exact sequences:
\[
\xymatrix{
& 1 & & \\
& \Gal(L/F) \ar[u] & & \\
1 \ar[r] & C_F \ar[u] \ar[r] & C_L \ar[r]^{\sigma-1} & C_L \\
& C_L \ar[u]^{N_{L/F}} & &,
}
\]
dualizing to
\[
\xymatrix{
& 1 \ar[d] & & \\
& \Gal(L/F)^D \ar[d] & & \\
1 \ar[r] & C_F^D \ar[d]^{N_{L/F}} \ar[l] & C_L^D \ar[l]_{res} & C_L^D \ar[l]_{\sigma-1} \\
& C_L^D  & &.
}
\]
By assumption, $N_{L/F} \circ res(\chi)=1$, so $res(\chi)= \delta^{-i} \in \Gal(L/F)^D$ for some integer $i$, unique modulo $d$. Then $res(\chi \delta^i)=1$, and we are done by exactness of the horizontal diagram.
\endproof
We need a special case of an $\ell$-adic analogue of the remark after Statement $A$ of \cite{lapid-rogawski:descent};\footnote{In the proof of Corollary \ref{slautomorphy} below, we could replace appeal to this lemma by simply citing Statement $B$ of \cite{lapid-rogawski:descent}.} that remark is in turn the (much easier) analogue,  for complex representations of the Weil group, of the main result of their paper. We first record the simple case that we need, and then out of independent interest we prove a general $\ell$-adic analogue of the Lapid-Rogawski result. 
\begin{lemma}\label{simpleLR}
Let $L/F$ be a quadratic CM extension of a totally real field $F$, with $\sigma \in \gal{F}$ generating $\Gal(L/F)$ Suppose $\hat{\psi} \colon \gal{L} \to \Qlb^\times$ is a Galois character such that $\hat{\psi}^{1-\sigma}$ is geometric. Then there exists a unitary, type $A$ Hecke character $\psi$ of $L$ such that $\psi^{1-\sigma}$ is the type $A_0$ Hecke character associated to $\hat{\psi}^{1-\sigma}$.
\end{lemma}
\proof
Write $(\hat{\psi}^{1-\sigma})_{\mbf{A}}$ for the Hecke character associated to $\hat{\psi}^{1-\sigma}$. By the previous lemma, it suffices to check that $(\hat{\psi}^{1-\sigma})_{\mbf{A}}$ is trivial on $C_F \subset C_L$. If a finite place $v$ of $L$ is split over a place $v_F$ of $F$, and unramified for $\hat{\psi}$, then for a uniformizer $\varpi_v$ of $F_{v_F}$ (embedded into the $L_v$ and $L_{\sigma v}$ components of $\mbf{A}_L^\times$), 
\[
(\hat{\psi}^{1-\sigma})_{\mbf{A}}(\varpi_v, \varpi_v)= \hat{\psi}^{1-\sigma}(fr_v) \hat{\psi}^{1-\sigma}(fr_{\sigma v})= \frac{\hat{\psi}(fr_v)}{\hat{\psi}(fr_{\sigma v})} \cdot \frac{\hat{\psi}(fr_{\sigma v})}{\hat{\psi}(fr_v)}= 1.
\]
Similarly for $v$ inert, $(\hat{\psi}^{1-\sigma})_{\mbf{A}}(\varpi_v)= \frac{\hat{\psi}(fr_v)}{\hat{\psi}(\sigma fr_v \sigma^{-1})}=1$. The Hecke character $(\hat{\psi}^{1-\sigma})_{\mbf{A}}|_{C_F}$ is therefore trivial. To see that we may choose $\psi$ to be unitary and type $A$, we invoke Corollary \ref{HCstructurethm}: decomposing $\psi$ as in that result, both the $|\cdot|^w$ and `Maass' components descend to the totally real subfield $F$, so dividing out by them yields a new $\psi$, now unitary type $A$, and with $\psi^{1-\sigma}$ unchanged.
\endproof
\begin{lemma}\label{pain}
Let $L/F$ be cyclic of degree $d$, with $\sigma \in \gal{F}$ restricting to a generator of $\Gal(L/F)$. Suppose $\rho \colon \gal{L} \to \mr{GL}_n(\Qlb)$ is an irreducible continuous representation satisfying $\rho^\sigma \cong \rho \cdot \chi$ for some character $\chi \colon \gal{L} \to \Qlb^\times$. Suppose further that $\chi$ is geometric (in particular, this holds if $\rho$ is geometric), necessarily of weight zero, so we may regard it as a character $\chi_{\mbf{A}}= \prod_{w \in |L|} \chi_w \colon C_L \to \CC^\times$. Then the (finite-order) restriction of $\chi_{\mbf{A}}$ to $C_F \subset C_L$ cannot factor through a non-trivial character of $\Gal(L/F)$. 
\end{lemma}
\proof
Iterating the relation $\rho^\sigma \cong \rho \cdot \chi$, we obtain $\rho \cong \rho \cdot \chi^{1+\sigma+ \ldots+ \sigma^{d-1}}$, so that $\chi^{1+\sigma+\ldots+\sigma^{d-1}}$ is finite-order. Each of the characters $\chi^{\sigma^i}$ has some common weight $w$, so $d\cdot w=0$, and thus $w=0$. Moreover, writing as usual $p_{\iota_w}$ for the algebraic parameter of $\chi_w$ (with respect to a choice $\iota_w \colon L \into \CC$ representing the place $w$), we have $\sum p_{\iota_w}=0$ as $\iota_w$ ranges over a $\Gal(L/F)$-orbit of such embeddings. If $L$ has a real embedding, then $\chi$ has finite-order, and the passage from $\chi \colon \gal{L}^{ab} \to \Qlb^\times$ to $\chi_{\mbf{A}} \colon C_L \to \CC^\times$ is simply via the reciprocity map $C_L \onto \gal{L}^{ab}$.\footnote{Throughout this argument we implicitly use our fixed embeddings $\iota_\ell, \iota_\infty$, but omit any reference to them for notational simplicity.} In particular, $\chi(x)= \chi_{\mbf{A}}(x')$ for any representative $x'$ in $C_L$ of the image of $x$ in $\gal{L}^{ab}$. If on the other hand $L$ is totally imaginary, then continuing to write $x'= (x'_w)_{w \in |L|}$, we have
\[
\chi(x)= \prod_{w \nmid \ell \infty} \chi_w(x'_w) \prod_{w \vert \ell} \left( \chi_w(x'_w) \prod_{\tau \colon L_w \into \Qlb} \tau(x'_w)^{p_{\iota_{\infty, \ell}^*(\tau)}} \right).
\]
If we further assume that the representative $x'$ can be chosen in $C_F \subset C_L$, with elements $x_v' \in F_v^\times$ giving rise to all $x_w'$ for $w \vert v$, then we can rewrite 
\[
\prod_{w \vert \ell} \prod_{\tau \colon L_w \into \Qlb} \tau(x'_w)^{p_{\iota_{\infty, \ell}^*(\tau)}}
\]
as
\[
\prod_{v \vert \ell} \prod_{\tau \colon F_v \into \Qlb} \prod_{w \vert v, \tilde{\tau} \vert \tau} \tau(x'_v)= \prod_{v \vert \ell} \prod_{\tau \colon F_v \into \Qlb} \tau(x'_v)^{\sum_{w \vert v, \tilde{\tau}\vert \tau} p_{\iota_{\infty, \ell}^*(\tilde{\tau})}}= 1.
\]
The last equality follows since $\chi^{1+\sigma+\ldots+ \sigma^{d-1}}$ is finite-order, and we are summing $p_{\iota_{\infty, \ell}^*(\tilde{\tau})}$ over a full $\Gal(L/F)$-orbit. A similar argument shows that 
\[
\prod_{w \vert \infty} \chi_w(x'_w)=1.
\]
We conclude that in this case ($x' \in C_F$), $\chi(x)$ can be computed simply as $\chi_{\mbf{A}}(x')$.

Let $V$ be the space on which $\rho$ acts. Then the isomorphism $\rho^\sigma \cong \rho \cdot \chi$ yields an operator $A \in \Aut(V)$ satisfying $A \rho^\sigma= \rho \cdot \chi A$. Fix $g \in \gal{L}$, and for any $x \in \gal{L}$ we compute
\begin{align*}
\tr(\rho(g)A)&=\tr(\rho^\sigma(x) \rho(g) A \rho^\sigma(x)^{-1})\\
&= \tr(\rho^\sigma(x) \rho(g) \rho(x^{-1}) \chi(x^{-1}) A)\\
&= \chi(x^{-1}) \tr(\rho(^\sigma{x} g x^{-1}) A).
\end{align*}
(Here $^\sigma x= \sigma x \sigma^{-1}$.) So, if we can find an $x \in \gal{L}$ such that $^\sigma x g x^{-1}=g$ and $\chi(x^{-1}) \neq 1$, then we will have $\tr(\rho(g) A)=0$. Doing this for all $g \in \gal{L}$, we see that by Schur's Lemma $\rho$ cannot be irreducible, else $A=0$. Now, $y= g^{-1} \sigma \in \gal{F}$ satisfies $^\sigma y g y^{-1}=g$, so $x= y^d \in \gal{L}$ does as well. It suffices to show that if $\chi|_{C_F}$ cuts out the extension $L/F$, then $\chi(x) \neq 1$. In fact, the image of $x$ in $\gal{L}^{ab}$ lies in the image of the transfer $\mr{Ver} \colon \gal{F}^{ab} \to \gal{L}^{ab}$. Explicitly,
\[
\mr{Ver}(y)= \prod_{i=0}^{d-1} \sigma^i (g^{-1} \sigma) \phi(\sigma^i g^{-1} \sigma)^{-1},
\]
where $\phi \colon \gal{F} \to \{\sigma^i\}_{i=0, \ldots, d-1}$ records the representative of the $\gal{L}$-coset of an element of $\gal{F}$. It is then easily seen\footnote{Note that $\phi(\sigma^{d-1} g^{-1} \sigma)=1$, while otherwise $\phi(\sigma^{i}g^{-1} \sigma)= \sigma^{i+1}$.} that 
\[
\mr{Ver}(y)= (g^{-1} \sigma)^d= x,
\]
so by class field theory $x \in \gal{L}^{ab}$ is represented by an element $x'$ of $C_F \subset C_L$ under the reciprocity map $\rec_L$. This element is a generator of $C_F/N_{L/F} C_L$, since $y$ lifts a generator of $\Gal(L/F)$, and thus $\chi(x)= \chi_{\mbf{A}}(x') \neq 1$ if $\chi_{\mbf{A}}|_{C_F}$ factors through a non-trivial character of $\Gal(L/F)$. 
\endproof
Finally we can (conditionally) prove automorphy of geometric projective representations over totally real fields.
\begin{cor}\label{slautomorphy}
Let $F$ be totally real. Continue to assume Part 3 of Conjecture \ref{FML} (Fontaine-Mazur-Langlands). Then for any geometric $\rho \colon \gal{F} \to \mr{PGL}_n(\Qlb)$, there exists an $L$-algebraic $\pi$ on $\mr{SL}_n/F$ such that $\rho \sim_w \pi$. 
\end{cor}
\proof
We will first treat the case of $\rho$ having irreducible lifts to $\mr{GL}_n(\Qlb)$. Choose a lift $\tilde{\rho}$ with finite-order determinant, a CM quadratic extension $L/F$, and, by Theorem \ref{fullmonodromylift}, a Galois character $\hat{\psi} \colon \gal{L} \to \Qlb^\times$ such that $\tilde{\rho}|_{\gal{L}}\cdot \hat{\psi}^{-1}$ is geometric. Let $\tilde{\pi}$ be the cuspidal automorphic representation of $\mr{GL}_n(\mbf{A}_L)$ corresponding to this geometric twist. Write $\sigma$ for the nontrivial element of $\Gal(L/F)$, so that $\tpi^\sigma \cong \tpi \cdot \chi$ where $\chi$ is the Hecke character corresponding to the geometric Galois character $\hat{\psi}^{1-\sigma}$. Appealing either to Lemma \ref{simpleLR} or Lemma \ref{pain}, we can write $\chi= \psi^{1-\sigma}$ for a unitary type $A$ Hecke character $\psi$. Then $\tpi \cdot \psi$ is $\sigma$-invariant, and by cyclic (prime degree) descent, there is a cuspidal representation $\pi$ of $\mr{GL}_n(\af)$ whose base-change is $\tpi \cdot \psi$. We want to compare the projectivization of the unramified parameters of $\pi$ with the unramified restrictions $\rho|_{\gal{F_v}}$. 

To do so, we repeat the argument but instead with infinitely many (disjoint) quadratic CM extensions $L_i/F$, showing that in all cases the descent $\pi$ to $\mr{GL}_n(\af)$ gives an $L$-packet for $\mr{SL}_n(\af)$ that is independent of the field $L_i$. $\tilde{\rho}$ is still a fixed lift with finite-order determinant, and we can write, for each $\tau \colon F \into \Qlb$,  $\mr{HT}_{\tau}(\tilde{\rho}) \in \frac{k_\tau}{n}+ \Z$ for some integer $k_\tau$, which we fix (rather than just its congruence class mod $n$). For some integer $w$, we have the purity relation $2 k_\tau \equiv w \mod n$, as follows, for instance, from (geometric) liftability after a CM base-change;\footnote{After such a base-change $L/F$, the character $\hat{\psi}$ twisting $\tilde{\rho}$ to a geometric representation will have Hodge-Tate-Sen weights congruent to $\frac{k_\tau}{n} \in \Q/\Z$ at both embeddings $L \into \Qlb$ above $\tau$; the integer $w$ is then the weight of the Hecke character associated to $\hat{\psi}^n$.} as with $k_\tau$, we fix an actual integer $w$, not just the congruence class. Now, for each such $\tau$, let $\iota \colon F \into \CC$ be the archimedean embedding associated via $\iota_\infty, \iota_\ell$ (elsewhere denoted $\iota^*_{\infty, \ell}(\tau)$). For each $L_i$, fix an embedding $\tau(i) \colon L_i \into \Qlb$ extending $\tau$, so that the other extension is $\tau(i) \circ c$. Likewise write $\iota(i)$ and $\iota(i)\circ c= \overline{\iota(i)}$ for the corresponding complex embeddings. We can then construct Galois characters $\hat{\psi}_i \colon \gal{L_i} \to \Qlb^\times$ such that 
\begin{align*}
\mr{HT}_{\tau(i)}(\hat{\psi}_i)&= \frac{k_{\tau}}{n} \\ 
\mr{HT}_{\tau(i) \circ c}(\hat{\psi}_i)&= \frac{w-k_{\tau}}{n},
\end{align*} 
such that $\tilde{\rho}|_{\gal{L_i}} \cdot \hat{\psi}_i^{-1}$ is geometric, corresponding to an $L$-algebraic cuspidal $\tpi_i$ on $\mr{GL}_n/L_i$. As before, we find a Hecke character $\psi_i$ of $L_i$ such that $\psi_i^{1-\sigma_i}$ is the type $A_0$ Hecke character associated to $\hat{\psi}_i^{1-\sigma_i}$; here we write $\sigma_i$ for the non-trivial element of $\Gal(L_i/F)$, but of course all the $\sigma_i$ are just induced by complex conjugation. Again, for all $i$ we find cuspidal automorphic representations $\pi_i$ of $\mr{GL}_n(\af)$ such that $\BC_{L_i/F}(\pi_i)= \tpi_i \cdot \psi_i$. Restricting to composites $L_i L_j$, we have the comparison
\[
\tilde{\rho}|_{L_i L_j} \hat{\psi}_i^{-1}|_{L_i L_j} \cdot \left( \frac{\hat{\psi}_i}{\hat{\psi}_j} |_{L_i L_j}\right)= \tilde{\rho}|_{L_i L_j} \hat{\psi}_j^{-1}|_{L_i L_j},
\]
and thus
\[
\BC_{L_i L_j/L_i}(\tpi_i \psi_i) \cdot \BC_{L_i L_j} \left( \frac{\psi_j}{\psi_i} \right)\cdot \BC_{L_i L_j}\left( \frac{\hat{\psi}_i}{\hat{\psi}_j} \right)= \BC_{L_i L_j/ L_j}( \tpi_j \psi_j),\footnote{Here $\frac{\hat{\psi}_i}{\hat{\psi}_j}$ restricted to $\gal{L_i L_j}$ is geometric, so we abusively write this for the associated Hecke character as well.}
\]
so finally
\[
\BC_{L_i L_j}(\pi_i) \cdot \BC_{L_i L_j} \left( \frac{\psi_j}{\psi_i} \cdot \frac{\hat{\psi}_i}{\hat{\psi}_j} \right)= \BC_{L_i L_j}(\pi_j).
\]

If the character $\frac{\psi_j}{\psi_i} \cdot \frac{\hat{\psi}_i}{\hat{\psi}_j}$ is finite-order-- in the next paragraph, we check that we may assume this-- it cuts out a cyclic extension $L'/L_i L_j$, and we have $\BC_{L'}(\pi_i)= \BC_{L'}(\pi_j)$. $L'/F$ is solvable, however, so the characterization of the fibers of solvable base-change in \cite{rajan:solvable} implies that $\pi_i$ and $\pi_j$ are twist-equivalent, hence that $\pi_i|_{\mr{SL}_n(\af)}$ and $\pi_j|_{\mr{SL}_n(\af)}$ define the same $L$-packet of $\mr{SL}_n(\af)$. Let us denote by $\pi_0$ any representative of this global $L$-packet. Now consider places $v$ of $F$ that are split in a given $L_i/F$. The semi-simple part $\rho(fr_v)^{ss}$ is equal (in $\mr{PGL}_n(\Qlb)$) to $(\tilde{\rho}\hat{\psi}^{-1}(fr_w))^{ss}$ for any $w \vert v$, and this is conjugate in $\mr{GL}_n(\Qlb)$ to $\iota_{\ell, \infty}\left(\rec_w(\tpi_w)(fr_w)\right)$, whose projectivization lies in the same $\mr{PGL}_n(\Qlb)$-conjugacy class as $\iota_{\ell, \infty}\left(rec_v(\pi_{0,v})(fr_v)\right)$. This verifies that for all such $v$, $\rho(fr_v)^{ss}$ is $\mr{PGL}_n(\Qlb)$-conjugate to $\iota_{\ell, \infty}\left(rec_v(\pi_{0,v})(fr_v)\right)$. Varying $L_i/F$, and remembering that $\pi_0$ is independent of this variation, we get the same result for all $v$ split in any single quadratic CM extension $L_i/F$ (we have to throw out a finite number of such $L_i$ to ensure our representations remain cuspidal/irreducible), we conclude that $\rho \sim_w \pi$.

To finish the proof, we must check that $\frac{\psi_j}{\psi_i} \cdot \frac{\hat{\psi}_i}{\hat{\psi}_j}$ may indeed be assumed finite-order. First, recall that each $\psi_i$ may be taken unitary and type $A$; in this case, the infinity-type is determined by the relation $\psi_i^{1-\sigma_i}= \hat{\psi}_i^{1-\sigma_i}$. Explicitly (using the above notation for the various embeddings), $\hat{\psi_i}^n$ corresponds to a Hecke character of $L_i$ with infinity-type (where we abusively denote $\iota(i)(z)$ by simply $z$)
\[
\rec_{\iota(i)}(\hat{\psi_i}^n) \colon z \mapsto z^{k_\tau} \bar{z}^{w-k_\tau},
\]
so
\[
\rec_{\iota(i)}(\psi_i^{1-\sigma_i}) \colon z \mapsto z^{\frac{2k_\tau-w}{n}} \bar{z}^{\frac{w-2k_\tau}{n}}.
\]
(Recall that $2k_\tau \equiv w \mod n$.) We then have, under our assumptions,
\[
\rec_{\iota(i)}(\psi_i) \colon z \mapsto \left( \frac{z}{|z|} \right)^{\frac{2k_\tau-w}{n}}.
\]
Of course, $\rec_{\iota(i) \circ c}$ is the same but with $\frac{w-2k_\tau}{n}$ in the exponent. To make the parameter comparison after restriction to a composite $L_i L_j$, we use the following notation for embeddings of $L_i L_j$ into $\Qlb$ and $\CC$, lying above the given $\tau$ and $\iota$:
\begin{align*}
&\text{$\tau_1$ extends $\tau(i)$ and $\tau(j)$,}\\
&\text{$\tau_2$ extends $\tau(i)$ and $\tau(j) \circ c$,}\\
&\text{$\iota_1$ extends $\iota(i)$ and $\iota(j)$,}\\
&\text{$\iota_2$ extends $\iota(i)$ and $\iota(j) \circ c$}.
\end{align*}
We then have the conjugate embeddings $\tau_1 \circ c$, etc. Computing the $\tau_k$-labeled weights of $\frac{\hat{\psi}_i}{\hat{\psi}_j}$, and translating them to the infinity-type at the place corresponding to $\iota_k$, with $\iota_k$ as the chosen embedding $L_i L_j \into \CC$, we then find
\begin{align*}
\rec_{\iota_1}\left(\frac{\hat{\psi}_i}{\hat{\psi}_j}\right)& \colon z \mapsto 1 \\
\rec_{\iota_2}\left(\frac{\hat{\psi}_i}{\hat{\psi}_j}\right)& \colon z \mapsto \left( z/\bar{z} \right)^{\frac{2k_\tau-w}{n}},
\end{align*}
whereas
\begin{align*}
\rec_{\iota_1}\left( \frac{\psi_j}{\psi_i} \right)& \colon z \mapsto 1\\
\rec_{\iota_2}\left( \frac{\psi_j}{\psi_i} \right)& \colon z \mapsto \left( z/\bar{z} \right)^{\frac{w-2k_\tau}{n}}.
\end{align*}
We conclude that $\frac{\psi_j}{\psi_i} \cdot \frac{\hat{\psi}_i}{\hat{\psi}_j}$ is, with our normalization of the $\psi_i$, in fact finite-order, and the proposition follows.

Finally, we quickly treat the case of general $\rho$, having reducible lifts. If $\tilde{\rho}$ as above decomposes $\tilde{\rho}= \oplus_{i=1}^m \tilde{\rho}_i$, say with $\tilde{\rho}_i$ of dimension $n_i$, geometricity of $\rho$ implies that over any CM $L/F$ the same Galois character $\hat{\psi}$ twists $\tilde{\rho}_i$, for all $i$, to a geometric representation. We can therefore use, for all $i$, the same Hecke character $\psi$ such that $\psi^{1-\sigma}= \hat{\psi}^{1-\sigma}$. As above, we invoke automorphy of $\tilde{\rho}_i \hat{\psi}^{-1}$, and, twisting by $\psi$, descend to a cuspidal automorphic representation $\Pi_i$ of $\mr{GL}_{n_i}/F$. The same local check (for $v$ split in $L/F$) as above, but now crucially relying on the fact that $\psi$ and $\hat{\psi}$ were independent of $i$, shows $\rho(fr_v)$ is $\mr{PGL}_n(\Qlb)$-conjugate to $\iota \left( \rec_v(\boxplus_{i=1}^m \Pi_{i, v})(fr_v) \right)$.
\endproof
By a similar argument, we can `construct' the Galois representations (assuming of course the $\mr{GL}_N$ correspondence) associated to (tempered) $L$-algebraic $\pi$ on $\mr{SL}_n/F$ for $F$ CM (or imaginary, assuming Conjecture \ref{cmdescent}) or totally real. By Proposition \ref{cmextension}, we are reduced to the case of $F$ totally real.
\begin{prop}\label{BGeg}
Continue to assume Fontaine-Mazur-Langlands. Let $F$ be a totally real field, and let $\pi$ be an $L$-algebraic cuspidal automorphic representation of $\mr{SL}_n(\af)$. Assume that $\pi_\infty$ is tempered. Then there exists a (not necessarily unique) projective representation $\rho \colon \gal{F} \to \mr{PGL}_n(\Qlb)$ satisfying $\rho \sim_w \pi$.\end{prop}
\proof
By Proposition \ref{totrealalglift}, there exists a $W$-algebraic cuspidal (and tempered at $\infty$) $\tpi$ on $\mr{GL}_n/F$ lifting $\pi$. For all but finitely many quadratic CM $L/F$, we can find a type $A$ Hecke character $\psi$ such that $\BC_{L/F}(\tpi) \cdot \psi$ is $L$-algebraic and cuspidal on $\mr{GL}_n/L$, hence corresponds to an irreducible geometric representation $\tilde{\rho}_L \colon \gal{L} \to \mr{GL}_n(\Qlb)$. Conjugating, we find $\tilde{\rho}_L^\sigma \equiv \tilde{\rho}_L \cdot \widehat{(\psi^{\sigma-1})}$, the twist being by the geometric character associated to the type $A_0$ Hecke character $\psi^{\sigma-1}$. We wish to write $\hat{\chi}= \widehat{(\psi^{\sigma-1})}$ in the form $\hat{\psi}^{\sigma-1}$ for some Galois character $\hat{\psi} \colon \gal{L} \to \Qlb^\times$ (reversing the process in Corollary \ref{slautomorphy}). Once this is managed, we have $(\tilde{\rho}_L \cdot \hat{\psi}^{-1})$ is $\sigma$-invariant, hence descends to $\rho \colon \gal{F} \to \mr{PGL}_n(\Qlb)$. Arguing as in Corollary \ref{slautomorphy}, we find that this projective descent is independent of $L/F$, and, again as in that proof, by varying $L/F$ we obtain the compatibility $\rho \sim_w \pi$.

To construct $\hat{\psi}$, we take as first approximation a Galois character $\hat{\psi}$ such that $\hat{\psi}^2$ equals $\widehat{(\psi^2)}$ up to a finite-order character $\chi_0$; recall that $\psi^2$ is type $A_0$, so we can attach the Galois character $\widehat{(\psi^2)}$. Then 
\[
(\hat{\psi}^{\sigma-1})^2= (\hat{\psi}^2)^{\sigma-1}= (\widehat{(\psi^2)} \chi_0)^{\sigma-1}= \hat{\chi}^2 \chi_0^{\sigma-1},
\]
and consequently $\hat{\chi}$ agrees with $\hat{\psi}^{\sigma-1}$ up to a finite-order character. Twisting $\tilde{\rho}_L$, we may therefore assume that $\hat{\chi}$ has finite-order. It still satisfies $\hat{\chi}^{1+\sigma}=1$, so invoking Lemma \ref{cft} and (a simple case of) Lemma \ref{pain} we find a finite-order Hecke character,\footnote{Writing $\hat{\chi}= \psi^{\sigma-1}$ where the infinity-components of $\psi$ have the form $z^p \bar{z}^q$, we see $z^{p-q}\bar{z}^{q-p}$ is the corresponding component of $\hat{\chi}$, hence that $p=q$ at each infinite place. Twisting $\psi$ by the base-change of a character of the totally real field $F$, we can then assume it is finite-order.} which may therefore be directly regarded as a Galois character, that casts $\hat{\chi}$ in the desired form.
\endproof
Here is another example comparing the Tannakian formalisms:
\begin{prop}\label{tensorformalism}
Continue to assume Fontaine-Langlands-Mazur and, in the totally imaginary but non-CM case, Conjecture \ref{cmdescent}. Let $\Pi$ be a cuspidal $L$-algebraic representation of $\mr{GL}_n(\af)$, and suppose that $\rho_\Pi \cong \rho_1 \otimes \rho_2$, where $\rho_i \colon \gal{F} \to \mr{GL}_{n_i}(\Qlb)$. Then there exist cuspidal automorphic representations $\pi_i$ of $\mr{GL}_{n_i}(\af)$ such that $\Pi= \pi_1 \boxtimes \pi_2$.
\end{prop}
\begin{rmk}
As the examples in \S \ref{HMFs} show, sometimes the $\pi_i$ cannot be taken $L$-algebraic. Nevertheless, by Proposition \ref{tensordescent}, they can always be taken $W$-algebraic.
\end{rmk}
\proof
First suppose $F$ is totally imaginary. For all $v \vert \ell$, the fact that $\rho_1 |_{\gal{F_v}} \otimes \rho_2 |_{\gal{F_v}}$ is de Rham implies that locally these $\gal{F_v}$-representations are twists of de Rham representations. To see this, we apply Corollary \ref{HTlift} (to find a Hodge-Tate lift) and Theorem \ref{conradliftingP} (to find a de Rham lift, given that a Hodge-Tate lift exists) to the lifting problem (with central torus kernel) $\mr{GL}_{n_1} \times \mr{GL}_{n_2} \xrightarrow{\boxtimes} G_{n_1, n_2} \subset \mr{GL}_{n_1 n_2}$, where $G_{n_1, n_2}$ denotes the image of the tensor product map. In particular, the projectivizations of the $\rho_i|_{\gal{F_v}}$ are de Rham, so the global projective representations $\rho_i \colon \gal{F} \to \mr{PGL}_{n_i}(\Qlb)$ are geometric. We know that these geometric projective representations have geometric lifts, and we may therefore assume our original $\rho_i$ were in fact geometric. They then correspond to $L$-algebraic $\pi_i$, and we have $\Pi= \pi_1 \boxtimes \pi_2$.

For $F$ totally real, we perform a descent similar to previous arguments. Restricting to CM $L/F$, we find a Galois character $\hat{\psi}$ such that $\rho_1 \cdot \hat{\psi}^{-1}$ and $\rho_2 \cdot \hat{\psi}$ are geometric, corresponding to $L$-algebraic cuspidal $\pi_i$ on $\mr{GL}_{n_i}/L$. Writing $\hat{\psi}^{1-\sigma}= \psi^{1-\sigma}$ for a Hecke character $\psi$, we find $\pi_1 \cdot \psi$ and $\pi_2 \cdot \psi^{-1}$ are $\sigma$-invariant, so descend to cuspidal representations $\bar{\pi}_i$ of $\mr{GL}_{n_i}/F$. Since $\BC_{L/F}(\Pi)= \BC_{L/F}(\bar{\pi}_1 \boxtimes \bar{\pi}_2)$, we deduce that $\Pi$ and $\bar{\pi}_1 \boxtimes \bar{\pi}_2$ are twist-equivalent, from which the result follows. (The same argument applies to $F$ that are neither totally real nor totally imaginary: just replace the restriction to CM extensions $L/F$ with restrictions to totally imaginary quadratic extensions $L/F$.)
\endproof
These examples (and, for instance, Corollary \ref{mixedparitymismatch}) motivate a comparison of the images of $r \colon {}^L H \to {}^L G$ on the automorphic and Galois sides, when $r$ is an $L$-morphism with central kernel. The most optimistic expectation (for $H$ and $G$ quasi-split) is that if $\ker(r)$ is a central torus, then the two descent problems for ($L$-algebraic) $\Pi$ and (geometric) $\rho_\Pi$ are equivalent; whereas if $\ker(r)$ is disconnected, there is an obstruction to the comparison, that nevertheless can be killed after a finite base-change. If $G$ is not $\mr{GL}_n/F$, then one will have to decide whether weak equivalence ($\pi \sim_w \rho$) suffices to connect the descent problems, or whether some stronger link (the mysterious $\pi \sim_s \rho$) must be postulated.
%\begin{conj}
%Let $H/F$ and $G/F$ be connected reductive quasi-split $F$-groups, and let $r \colon {}^L H \to {}^L G$ be an $L$-homomorphism with central kernel. Suppose that $\Pi$ is an $L$-algebraic automorphic representation of $G(\af)$ corresponding to a geometric representation $\rho \colon \gal{F} \to {}^L G(\Qlb)$. 
%\begin{itemize}
%\item If $\ker(r)$ is a central torus, then $\Pi$ descends to $H$ if and only if $\rho$ factors through ${}^L H(\Qlb)$.
%\item If $\ker(r)$ is disconnected, then there is an obstruction to comparing descent of $\Pi$ and $\rho$; but if one descends, then after a finite base-change, so will the other.
%\end{itemize}
%\end{conj}

\section{Monodromy of abstract Galois representations}\label{monodromy}
In this section we discuss some general results about monodromy of $\ell$-adic Galois representations. Much of the richness of this subject comes from its blending of two kinds of representation theories, that of finite groups, and that of connected reductive algebraic groups. We will see (Proposition \ref{lietensorart}) that the basic lifting result (Proposition \ref{tatelift}) allows us to some extent to understand how these two representation theories interact. In \S \ref{liemultfree} we develop more refined results in the `Lie-multiplicity-free' case (see Definition \ref{liemultfreedefn}); this situation encapsulates the essential difficulties of independence-of-$\ell$ questions, such questions being trivial for Artin representations.
\subsection{A general decomposition}
The following result is a simple variant of a result of Katz (\cite[Proposition 1]{katz:duke}), which he proves for lisse sheaves on affine curves over finite fields. We can replace Katz's appeal to the Lefschetz affine theorem by Proposition \ref{tatelift}. Recall that a Galois representation is \textit{Lie irreducible} if it is irreducible after restriction to every finite-index subgroup (i.e., the connected component or Lie algebra of its algebraic monodromy group acts irreducibly).\index{t}{Lie irreducible}
\begin{prop}\label{lietensorart}
Let $F$ be any number field, and let $\rho \colon \gal{F} \to \mr{GL}_{\Qlb}(V)$ be an irreducible representation of dimension $n$. Then either $\rho$ is induced, or there exists $d \vert n$, a Lie irreducible representation $\tau$ of dimension $n/d$, and an Artin representation $\omega$ of dimension $d$ such that $\rho \cong \tau \otimes \omega$. Consequently, any (irreducible) $\rho$ can be written in the form
\[
\rho \cong \Ind_{L}^F (\tau \otimes \omega)
\]  
for some finite $L/F$ and irreducible representations $\tau$ and $\omega$ of $\gal{L}$, with $\tau$ Lie-irreducible and $\omega$ Artin.
\end{prop}
\proof
Let $\mc{G}$ denote the algebraic monodromy group of $\rho$, with $\mc{G}^0$ the connected component of the identity. Abusively writing $\rho$ for the representation $\mc{G} \into \mr{GL}(V)$, we may assume $\rho|_{\mc{G}^0}$ is isotypic (else $\rho$ is induced, and we are done). If $\rho|_{\mc{G}^0}$ is irreducible, then $\rho$ itself is Lie-irreducible, so again we are done. Therefore, we may assume that $\rho|_{\mc{G}^0} \cong \tau_0^{\oplus d}$ for some $d \geq 2$, with $\tau_0$ an irreducible representation of $\mc{G}^0$, and consequently a Lie-irreducible representation of $\gal{L}$ for any $L/F$ sufficiently large that $\rho(\gal{L}) \subset \mc{G}^0(\Qlb)$. Since the irreducible $\gal{L}$-representation $\tau_0$ is $\gal{F}$-invariant, it extends to a projective representation of $\gal{F}$. By the basic lifting result (Proposition \ref{tatelift}), this projective representation lifts to an honest $\gal{F}$-representation $\tau_1$, so for some character $\chi \colon \gal{L} \to \Qlb^\times$, 
\[
\tau_1^{\oplus d}|_{\gal{L}} \cong \rho|_{\gal{L}} \otimes \chi.
\]
The character $\alpha:= \det(\rho)/\det(\tau_1^{\oplus d})$ of $\gal{F}$ has $\gal{L}$-restriction equal to $\chi^{-n}$, and over $F$ itself we can find characters $\alpha_1, \alpha_0 \colon \gal{F} \to \Qlb^\times$, with $\alpha_0$ finite-order, such that $\alpha= \alpha_1^n \alpha_0$. Then $(\chi \alpha_1|_{\gal{L}})^n= (\alpha^{-1} \alpha_1^n)|_{\gal{L}}= \alpha_0^{-1}|_{\gal{L}}$, and replacing $\tau_1$ by $\tau_1 \otimes \alpha_1$, and $L$ by a finite extension trivializing $\alpha_0$, we find a Lie-irreducible representation $\tau$ of $\gal{F}$ and a finite extension $L$ of $F$ such that $\tau^{\oplus d}|_{\gal{L}} \cong \rho|_{\gal{L}}$. The $\gal{F}$-representation
\[
\omega:= \Hom_{\gal{L}}(\tau, \rho).
\]
is therefore a $d$-dimensional Artin representation,\footnote{As $\gal{L}$-representation, $\omega \cong \Hom_{\gal{L}}(\tau|_{\gal{L}}, \tau|_{\gal{L}}^{\oplus d}) \cong \Qlb^{\oplus d}$.} and the natural map $\tau \otimes \omega \to \rho$ (i.e. $v \otimes \phi \mapsto \phi(v)$) is an isomorphism of $\gal{F}$-representations.
\endproof
\begin{cor}\label{lieisotypic}
Let $\rho \colon \gal{F} \to \mr{GL}_{\Qlb}(V)$ be a semi-simple representation (not necessarily irreducible), and suppose that $\rho$ is Lie-isotypic, i.e. for all $F'/F$ sufficiently large, $\rho|_{\gal{F'}}$ is isotypic. Then there exists a Lie-irreducible representation $\tau$ and an Artin representation (possibly reducible) $\omega$, both of $\gal{F}$, such that $\rho \cong \tau \otimes \omega$.
\end{cor}
\proof
Decompose $\rho$ into irreducible $\gal{F}$-representations as $\oplus_1^r \rho_i$. Each $\rho_i$ is Lie-isotypic: there exists $L/F$ and integers $m_i$ such that $\rho_i|_{\gal{L}} \cong \tau_0^{\oplus m_i}$ for all $i$, where $\tau_0$ is a Lie-irreducible representation independent of $i$. By the argument of the previous proposition, after possibly enlarging $L$ we find a $\gal{F}$-representation $\tau$ whose restriction to $L$ is isomorphic to $\tau_0$, and then there are Artin representations $\omega_i$ of $\gal{F}$ such that $\tau \otimes \omega_i \cong \rho_i$. Consequently, 
\[
\rho \cong \tau \otimes (\bigoplus_1^r \omega_i).
\]
\endproof
\begin{rmk}
\begin{itemize}
\item In general, the field $L$ in Proposition \ref{lietensorart} is not unique, even up to $\gal{F}$-conjugacy. Examples of such non-uniqueness should not arise in the Lie-multiplicity free case (see \S \ref{liemultfree}), but in the Artin case there are easy examples arising from the representation theory of finite groups. Consider, for instance, the quaternion group $Q_8= \{\pm 1, \pm i, \pm j, \pm k\}$. The (unique) irreducible two-dimensional representation of $Q_8$ can be written in the form $\Ind_{\langle x \rangle}^{Q_8}(\varepsilon_x)$, where $\langle x \rangle$ denotes one of the subgroups generated by $i$, $j$, or $k$, and $\varepsilon_x$ is a generator of the character group of $\langle x \rangle$. None of these subgroups is conjugate to any of the others.\footnote{Although these groups are not conjugate, they are related by (outer) automorphisms of $Q_8$, but applying outer automorphisms in this fashion will not in general preserve an irreducible induced character: consider the principal series of $\mr{GL}_2(\mbb{F}_p)$ and the outer automorphism of $A_1$.} It would be interesting to achieve a more systematic understanding of these ambiguities.
\item These structure theorems for Galois representations should have an analogue on the automorphic side. In fact, Tate has shown ($2.2.3$ of \cite{tate:background}) the analogue of Proposition \ref{lietensorart} for representations of the Weil group $W_F$, where it takes the particularly simple form that any irreducible, non-induced $\rho \colon W_F \to \mr{GL}_n(\CC)$ is isomorphic to $\omega \otimes \chi$, where $\omega$ is an Artin representation and $\chi \colon W_F \to \CC^\times$ is a character. A basic question is whether we should expect Proposition \ref{lietensorart} to hold for `representations of $\mc{L}_F$.' If we had the formalism of $\mc{L}_F$, then to carry out the argument of the proposition with complex representations of $\mc{L}_F$ in place of $\ell$-adic representations of $\gal{F}$ requires two ingredients:
\begin{itemize}
\item That a homomorphism $\mc{L}_F \to \mr{PGL}_n(\CC)$ lifts to a homomorphism $\mc{L}_F \to \mr{GL}_n(\CC)$; but this is `implied' by Proposition \ref{llgeneralization}. I should mention in this respect the theorem of Labesse (\cite{labesse:lifting}), which establishes the analogue for lifting homomorphisms $W_F \to {}^L G$ across surjections ${}^L \tG \to {}^L G$ with central torus kernel. 
\item That the analogue of the character $\alpha$ in the proof of Proposition \ref{lietensorart} can be written as a finite-order twist $\alpha_0$ of the $n^{th}$ power of a character $\alpha_1$. This is \textit{not} automatic, as it is for $\ell$-adic characters, but in this case we can exploit Lemma \ref{Heckeroots}, which applies to the $\alpha$ of the Proposition, since there the restriction to $L$ is (continuing with the notation of the Proposition) $\chi^{-n}$.
\end{itemize}
This discussion motivates the following conjecture, whose formulation of course requires assuming deep cases of functoriality:
\begin{conj}
Let $\pi$ be a cuspidal automorphic representation of $\mr{GL}_n(\af)$; assume $\pi$ is not automorphically induced from any non-trivial extension $L/F$. Then there exist cuspidal automorphic representations $\tau$ and $\omega$ of, respectively, $\mr{GL}_d(\af)$ and $\mr{GL}_{n/d}(\af)$ such that $\pi= \tau \boxtimes \omega$, with the following properties:
\begin{itemize}
\item for all finite extensions $L/F$, the base-change $\mr{BC}_{L/F}(\tau)$ remains cuspidal;
\item for some finite extension $L/F$, $\BC_{L/F}(\omega)$ is isomorphic to the isobaric sum of $n/d$ copies of the trivial representation.
\end{itemize} 
\end{conj}
\end{itemize}
\end{rmk}
\subsection{Lie-multiplicity-free representations}\label{liemultfree}
In this section, we focus on the cases antithetical to that of Artin representations, putting ourselves in the following situation. Let $F$ be any number field, and recall from definition \ref{wcsystem} the definition of a (weakly) compatible system of $\lambda$-adic representations of $\gal{F}$, with coefficients in a number field $E$. Let 
\[
\rho_{\lambda} \colon \gal{F} \to \mr{GL}_n(\Elb)
\] 
be such a (semi-simple, continuous) compatible system. Write $V_\lambda$ for the space on which $\rho_\ell$ acts. 
\begin{defn}\label{liemultfreedefn}\index{t}{Lie multiplicity free} 
We say that $V_\lambda$ is Lie-multiplicity free if after any finite restriction $L/F$, $V_\lambda|_{\gal{L}}$ is multiplicity-free. Equivalently, 
\[
\varinjlim_{L/F}\End_{\Elb[\gal{L}]}(V_{\lambda})
\] 
is commutative. We will often abbreviate `Lie-multiplicity-free' to `LMF.'\index{t}{LMF}
\end{defn}

Cases to keep in mind are Hodge-Tate regular (see Definition \ref{regular} and the preceding discussion of \S \ref{formal}) $V_{\lambda}$, or $V_{\lambda}$ of the form $H^1(A_{\overline{F}}, \Q_\ell)$ where $A/F$ is an abelian variety with $\End^0(A_{\overline{F}})$ a commutative $\Q$-algebra (by Faltings's proof of the Tate conjecture). Elementary representation theory yields:
\begin{lemma}\label{elementaryLMF}
\begin{enumerate}
\item Suppose $V_{\lambda}$ is irreducible. Then $V_{\lambda}$ is LMF if and only if it can be written 
\[
V_{\lambda} \cong \Ind_{L(\lambda)}^F (W_{\lambda}), 
\]
where we write $L(\lambda)$\index{s}{$L(\lambda)$} to show the \textit{a priori} dependence on $\lambda$ if $V_{\lambda}$ belongs to a compatible system, and where $W_{\lambda}$ is a Lie-irreducible $\Elb$-representation of $\gal{L(\lambda)}$, all of whose $\Gamma_F$-conjugates remain distinct after any finite restriction. 
\item Let $W_{\lambda}$ be an irreducible representation of $\gal{L}$, and assume that $V_{\lambda}= \Ind_{L}^F (W_{\lambda})$ is LMF. Then $V_{\lambda}$ is irreducible.
\end{enumerate}
\end{lemma}
\proof
For $(1)$, restrict to a finite-index subgroup of $\gal{F}$ over which $V_{\lambda}$ decomposes into a direct sum of Lie-irreducible representations; take one such factor, and consider its stabilizer in $\gal{F}$-- $V_{\lambda}$ is then induced from this subgroup. For $(2)$, Mackey theory implies we need to check that $W_{\lambda}|_{g\gal{L}g^{-1} \cap \gal{L}}$ and $(gW_{\lambda})|_{g\gal{L}g^{-1} \cap \gal{L}}$ are disjoint for all $g \in \gal{F}- \gal{L}$. These two representations occur as distinct factors in the $V_{\lambda}|_{g\gal{L}g^{-1} \cap \gal{L}}$, so they are disjoint since $V_{\lambda}$ is LMF.
\endproof
For general (possibly reducible) LMF representations, there is a decomposition into a sum of terms as in the lemma. If $V_{\lambda}$ belongs to a compatible system, we expect that the number of such factors should be independent of $\lambda$; this is an extremely difficult problem (unlike the corresponding question for Artin representations). Let us indicate the difficulties through an example. 
\begin{eg}
Suppose $f$ is a holomorphic cuspidal Hecke eigenform on the upper half-plane of some weight $k \geq 2$ and level $N$ and nebentypus $\epsilon$. We normalize $f$ so that its $q$-expansion at the cusp $\infty$, $f= \sum a_n(f) q^n$, has leading coefficient $a_1(f)=1$. Then work of Eichler-Shimura-Deligne (see \cite{deligne:galmod}) yields a compatible system
\[
\rho_{f, \lambda} \colon \gal{\Q} \to \mr{GL}_2(\Elb)
\]
of $\lambda$-adic representations (here $E$ is the number field generated by the $a_n(f)$) characterized by the property that for all $p \nmid N$, the characteristic polynomial of $\rho_{f, \lambda}(fr_p)$ is equal to
\[
X^2-a_p(f)X+p^{k-1}\epsilon(p).
\]
Moreover, the eigenvalues of $\rho_{f, \lambda}(fr_p)$ are $p$-Weil numbers of weight $k-1$, i.e. have absolute value $p^{\frac{k-1}{2}}$ in all complex embeddings, and $\rho_{f, \lambda}$ is Hodge-Tate of weights $0, 1-k$. Since $f$ is cuspidal, we expect these Galois representations to be irreducible. This can be proven \textit{because the Fontaine-Mazur-Langlands conjecture is known for the group} $\mr{GL}_1$. Suppose $\rho_{f, \lambda} \cong \chi_1 \oplus \chi_2$. Then (up to re-ordering $\chi_1$, $\chi_2$) Theorem \ref{FMLGL1} implies\footnote{This special case of Theorem \ref{FMLGL1} is notably easier than the general case. It suffices to show that a character $\chi \colon \gal{\Q} \to \Qlb^\times$ with all labeled Hodge-Tate weights equal to zero is finite-order; this follows from global class field theory and the corresponding statement for characters $\gal{\Q_{\ell}} \to \Qlb^{\times}$. This latter statement is a slight improvement of a fundamental theorem of Tate: see \cite[\S III.A-3]{serre:ladic}.} that $\chi_1$ must be a finite-order character and $\chi_2$ must be the product of a finite-order character and the $(1-k)^{th}$ power of the cyclotomic character. This contradicts the fact that the eigenvalues of $\rho_{f, \lambda}(fr_p)$ are $p$-Weil numbers of weight $k-1$, so we have proven the asserted irreducibility. 
\end{eg}

Therefore in this section we pursue a much more modest goal: restricting to the case of irreducible compatible systems, we will be able to say something about independence of $\lambda$ of the fields $L(\lambda)$ of Lemma \ref{elementaryLMF}. 

Our basic strategy is that the places $v$ for which $\tr(\rho_{\lambda}(fr_v)$ equals zero should detect the field $L(\lambda)$. This is in marked contrast to the case of Artin representations: any irreducible (non-trivial) representation of a finite group has elements acting with trace zero. Our main tool will be the following (slight weakening of a) theorem of Rajan:
\begin{thm}[Theorem $3$ of \cite{rajan:sm1}]\label{rajan}
Let $E$ be a mixed characteristic non-archimedean local field, and let $H/E$ be an algebraic group. Let $X$ be a subscheme of $H$ (over $E$), stable under the adjoint action of $H$. Suppose $\rho \colon \gal{F} \to H(E)$ is a Galois representation, unramified almost everywhere, and let $C= X(E) \cap \rho(\gal{F})$. Denote by $H_{\rho} \subset H$ the algebraic monodromy group $\overline{\rho(\gal{F})}^{Zar}$, and let $\Phi= H_{\rho}/H_{\rho}^0$ denote its group of connected components. For $\phi \in \Phi$, we write $H^{\phi}$ for the corresponding component. 
\begin{itemize}
\item Let $\Psi= \{\phi \in \Phi \vert H^\phi \subset X\}$. Then the density of the set of places $v$ of $F$ with $\rho(fr_v) \in C$ is precisely $|\Psi|/|\Phi|$.
\end{itemize}
\end{thm}
Rajan applies this to prove\footnote{In addition to the main result of his paper, a beautiful `strong multiplicity one' theorem for $\ell$-adic representations.} (Theorem $4$ of \cite{rajan:sm1}) that an irreducible, but \textit{Lie reducible}, representation necessarily has a positive density of frobenii acting with trace zero; note that this also follows immediately from \v{C}ebotarev and Proposition \ref{lietensorart}, which is a more robust version of Rajan's result (basically combining his argument with Proposition \ref{tatelift}). Our next two results establish a converse, also extending Corollaire $2$ to Proposition $15$ of \cite{serre:cebotarev} to its natural level of generality: that result handles the case of connected monodromy groups.
\begin{prop}\label{tracezero}
Let $ \rl \colon \gal{F} \to \mr{GL}_n(\Elb)$ be a continuous, semi-simple, LMF representation. Decompose $V_{\lambda}$ as above, so
\[
V_{\lambda} \cong \bigoplus_{i=1}^{r_{\lambda}} \Ind_{L(\lambda)_i}^F (W_{\lambda, i})
\]
for Lie-irreducible representations $W_{\lambda, i}$ of $\gal{L(\lambda)_i}$. Then:
\begin{enumerate}
\item Up to a density zero set of places,
\[
\{v \in |F|: \tr(\rl(fr_v))=0\}= \{v: fr_v \notin \bigcup_i \bigcup_{\sigma \in S_{\lambda, i}} \sigma \Gamma_{L(\lambda)_i} \sigma^{-1}\},
\]
where $S_{\lambda, i}$ is a set of representatives of $\Gamma_F/\Gamma_{L(\lambda)_i}$.
\item Further assume that $\rl$ belongs to a compatible system $\{\rl\}$ of $\lambda$-adic representations of $\gal{F}$ (although in contrast to definition \ref{wcsystem}, we need not assume here that the $\rl$ are geometric). Then up to a set of density zero, the set of places of $F$ which have a split factor in $L(\lambda)_i$ for some $i$ is independent of $\lambda$. If we further assume that all $V_\lambda$ are absolutely irreducible $(r_\lambda=1)$ and all $L(\lambda)/F$ are Galois, then $L(\lambda)$ is independent of $\lambda$.\footnote{In the non-Galois case, see Exercise $6$ in \cite{cassels-frohlich}!}
\end{enumerate}
\end{prop}
\proof
For the first part of the Proposition, we ignore the underlying `coefficient' number field $E$ and just view $\rho_{\lambda}$ as valued in $\mr{GL}_n(E)$ for some sufficiently large finite extension $E$ of $\Q_{\ell}$. The ``$\supseteq$'' direction follows from the usual formula for the trace of an induced representation. To establish the reverse inclusion, let us consider, for each non-empty subset $I \subset \coprod_i S_{\lambda, i}$,  the set $\mc{X}_I$ of places $v$ such that $\tr(\rl(fr_v))=0$, and $fr_v \in \sigma\Gamma_{L(\lambda)_{i(\sigma)}} \sigma^{-1}$ if and only if $\sigma \in I$ [Notation: if $\sigma \in I$, then $\sigma \in S_{\lambda, i}$ for a unique $i=: i(\sigma)$]. Also set
\[
\Gamma_I= \bigcap_{\sigma \in I} \sigma \Gamma_{L(\lambda)_{i(\sigma)}} \sigma^{-1},
\]
and let $\mc{G}_I$, resp. $\mc{G}$, denote the algebraic monodromy group of $\rl|_{\Gamma_I}$, resp. $\rl$. To establish the ``$\subseteq$'' (up to density zero) direction, we must show that $\mc{X}_I$ has density zero for every non-empty $I$. For all $I$, $\mc{G}_I$ contains $\mc{G}^0$, the identity component of $\mc{G}$. We apply Rajan's Theorem (\ref{rajan} above) to $\rl|_{\Gamma_I}$: if $\mc{X}_I$ has positive density, then there is a full connected component $T \mc{G}^0 \subset \mc{G}_I$ on which the trace vanishes (here $T$ is some coset representative for the component). Representing endomorphisms of $V_{\lambda}$ in block-matrix form corresponding to the decomposition 
\[
V_{\lambda}|_{\rl^{-1}(\mc{G}^0)}= \bigoplus_{i, \sigma} \sigma W_{\lambda, i},
\] 
we have
\[
\tr \left[T \cdot \left( 
\begin{array}{cccc}
* & 0 & \cdots & 0 \\
0 & * & \cdots & 0 \\
\vdots & \vdots & \ddots & 0 \\
0 & 0 & \cdots & *
\end{array} \right) \right]=0,
\]
where the $*$'s represent \textit{arbitrary} elements of the $\End(\sigma W_{\lambda, i})$. For this, we use the fact that these constituents are (absolutely) irreducible and distinct: either apply Wedderburn theory to the semi-simple $\Elb$-algebra $\Elb[[\rl(\rl^{-1}(\mc{G}^0))]]$, or apply Schur's lemma (using absolute irreducibility) and an algebra version of Goursat's lemma (using multiplicity-freeness). Then the above matrix equation implies the automorphism $T$, written in block-matrix form, has all zeros along the (block)-diagonal. But $I$ is non-empty, so $T \in \mc{G}_I$ preserves at least one subspace $\sigma W_{\lambda, i}$, and so its block-diagonal entry corresponding to that sub-space must be non-zero (invertible). This contradiction forces all components of $\mc{G}_I$ to have non-zero trace (generically), for all non-empty $I$, and thus the ``$\subseteq$'' direction is established. 

Part $2$ now follows from independence of $\lambda$ of $\tr(\rl(fr_v))$ and \v{C}ebotarev, noting that
\[
\{v: fr_v \in \bigcup_i \bigcup_{\sigma \in S_{\lambda, i}} \Gamma_{\sigma(L(\lambda)_i)}\}
\]
is the set of places $v$ of $F$ that have at least one split factor in some $L(\lambda)_i$. 
\endproof
We make a few remarks about the limitations of this method:
\begin{rmk}
\begin{enumerate}
\item Without any information about how the various $V_{\lambda}$ decompose into irreducible sub-representations, this result yields frustratingly little, since it is easy to find disjoint collections of number fields $\{L_i\}_{i \in I}$ and $\{L'_i\}_{i \in I'}$ such that the union of primes with a split factor (or even split) in the various $L_i$ equals the union of those with a split factor in the various $L'_i$. For, the simplest example, take $L_1= \Q$, $L'_1= \Q(i)$, $L'_2= \Q (\sqrt{2})$, $L'_3= \Q (\sqrt{-2})$.
\item Even when the $V_{\lambda}$ are irreducible, and even assuming that one field $L(\lambda_0)$ is Galois, care must be taken when the other inducing fields $L(\lambda)$ are not (known to be) Galois. We see that $L(\lambda_0)$ is contained in $L(\lambda)$ for all $\lambda$, but equality does not follow, as the following example from finite group theory shows. We want an inclusion $H \leq K \lhd G$ of groups with $K$ normal in $G$, $H$ a proper subgroup of $K$, and $\cup_{g \in G} gHg^{-1}=K$. Taking $K \lhd S_4$ to be the copy of the Klein four-group given by the $(2, 2)$-cycles, and $H$ to be the subgroup generated by one of these permutations, meets the requirements. Note that such examples require $K$ to have non-trivial outer automorphism group: if $G$-conjugation acts by $K$-inner automorphisms on $K$, then $\cup_G gHg^{-1}= \cup_K kH^k{-1}=K$ implies $H=K$, since no finite group is the union of conjugates of a proper subgroup.
\item The group-theoretic counterexample of the previous item should not arise in practice: if we assume that $W_{\lambda_0}$ can be put in a compatible-system, say with $\lambda$-adic realization $W_{\lambda, 0}$, then conjecturally $W_{\lambda, 0}$ will be Lie irreducible as well, and then the isomorphism
\[
\Ind_{L(\lambda)}^F(W_{\lambda}) \cong \Ind_{L(\lambda_0)}^F(W_{\lambda, 0})
\]
implies, by Mackey theory, that there is a non-zero $\gal{L(\lambda_0)}$-morphism
\[
\Ind^{L(\lambda_0)}_{s(L(\lambda))}(sW_{\lambda}) \onto W_{\lambda, 0}
\]
for some $s \in \gal{F}$. By part $(2)$ of Lemma \ref{elementaryLMF}, this induction is irreducible, so this map is an isomorphism. But $W_{\lambda, 0}$ is (conjecturally) Lie-irreducible, so $L(\lambda_0)=L(\lambda)$.
\end{enumerate}
\end{rmk}
In any case, the following corollary is the promised converse to Rajan's result:
\begin{cor}\label{LI}
An irreducible, LMF representations $\rho_{\ell} \colon \gal{F} \to \mr{GL}_n(\Qlb)$ is Lie irreducible precisely when the set of $v$ with $\tr(\rho_{\ell}(fr_v))=0$ has density-zero. In particular, in a compatible system of irreducible, LMF representations, Lie-irreducibility is independent of $\lambda$. 
\end{cor}
\begin{rmk}
If we know general automorphic base-change for $\mr{GL}_n$, we can formulate the conditions `Lie irreducible' and `Lie multiplicity free' on the automorphic side. This result then suggests how to tell whether a `LMF' automorphic representation is automorphically induced. Finding an intrinsic characterization of the image of automorphic induction, even conjecturally, is a mystery (in contrast to its close cousin base-change), so it may come as a surprise that there should be such a simple condition at the level of Satake parameters, for this broad class of LMF representations.   
\end{rmk}

I originally developed Proposition \ref{tracezero} to prove that a regular compatible system of representations of $\gal{F}$ for $F$ a CM field, if induced, is necessarily induced from a CM field. See Remark \ref{algconjrmk} for an automorphic analogue. Here is a partial result; it is another application of the 'Hodge-theory with coefficients' in \S \ref{formal}. For the most natural formulation of the result, recall the definition of purity of a weakly compatible system (Definition \ref{purewcs}).
\begin{lemma}
Let $F$ be a CM field, and let $\mc{R}= \{\rho_{\lambda} \colon \gal{F} \to \mr{GL}_n(\overline{E_{\lambda}})\}$ be a weakly compatible system of $\lambda$-adic representations with coefficients in a number field $E$. Assume that the $\rho_{\lambda}$ are (almost all) Hodge-Tate regular, and pure. Finally, assume that there is a single number field $L$ such that for $\lambda$ above a set of rational primes $\ell$ of density one, 
\[
\rho_{\lambda} \cong \Ind_L^F(r_{\lambda}),
\]
for some $\overline{E_{\lambda}}$-representation $r_{\lambda}$ of $\gal{L}$. Then $L$ is CM.
\end{lemma}
\proof
We may assume $E/\Q$ is Galois. In yet another variant of the theme of \S \ref{cmdescentsection}, purity implies we may take the number field $E$ to be CM: simply observe that for any choice $c$ of complex conjugation in $\Gal(E/\Q)$, the characteristic polynomials $Q_v(X)$ of $\rho_{\lambda}(fr_v)$ satisfy
\[
{}^c Q_v(X)= X^n Q_v(q_v^w/X)/Q_v(0),
\]
and thus for any two complex conjugations $c, c'$, ${}^{cc'}\mc{R} \cong \mc{R}$. It follows that for a density one set of $v$, $Q_v(X)$ has coefficients in $E_{cm}$. The \v{C}ebotarev density theorem then yields a well-defined weakly compatible system (in the sense of Definition \ref{wcsystem}) of $\lambda$-adic representations with coefficients in $E_{cm}$. (This argument is taken from \cite[Lemma 1.1, 1.2]{patrikis-taylor:irr}.\footnote{Our definition of compatible system is somewhat weaker than that given in \cite{patrikis-taylor:irr}; for our purposes, we can just ignore the part of the proof of \cite[Lemma 1.1]{patrikis-taylor:irr} that computes Hodge numbers.}). Thus we take $E$ to be CM. By regularity, we may also assume that there is a CM extension $E'/E$ such that for all finite-index subgroups $H$ of $\gal{F}$, and all primes $\lambda$ of $E'$, all sub-representations $r \subset \rho_{\lambda}|_{H}$ are actually defined over $E'_{\lambda}$: this is an elementary argument (see \cite[Lemma 5.3.1(3)]{blggt:potaut}, with the CM refinement of \cite[Lemma 1.4]{patrikis-taylor:irr}), the key idea being that regularity gives us an abundant supply of elements in the image of $\rho_{\lambda}$ with distinct eigenvalues, and that $\rho_{\lambda}$ can then be defined over the extension of $E$ generated by these eigenvalues. Thus, after enlarging $E$, we may assume all $\rho_{\lambda}$ and $r_{\lambda}$ act on $E_{\lambda}$-vector spaces. 

For simplicity enlarge $E$ to contain $L_{cm}$, the maximal CM subfield of $L$, and take its Galois closure-- the result remains CM. If $L \neq L_{cm}$, then we can find a (positive-density) set of $\ell$ (unramified in $L$ and for the system $\rl$) which are split in $E$ (with, say, $\lambda \vert \ell$) but not in (the Galois closure of $L$, hence) $L$. Consider a non-split prime $w \vert \ell$ of $L$, above a place $v$ of $F$, and the restriction
\[
r_{\lambda}|_{\gal{L_w}} \colon \gal{L_w} \to \mr{GL}_n(E_{\lambda})= \mr{GL}_n(\Q_\ell);
\]
By Lemma \ref{dRInd}, $\Ddr(\rho_{\lambda}|_{\gal{F_v}})$ is the image under the forgetful functor (from filtered $L_w$-vector spaces to filtered $F_v$-vector spaces) of $\Ddr(r_{\lambda}|_{\gal{L_w}})$. Since $L_w$ does not embed in $E_{\lambda}= \Ql$, we can invoke Corollary \ref{notregular} to show $\rho_{\lambda}$ is not regular, a contradiction. Therefore $L= L_{cm}$.
\endproof
Combining this with Proposition \ref{tracezero}, since regular clearly implies LMF,\footnote{This is the one case in which the LMF condition is provably independent of $\ell$.} we deduce: 
\begin{cor}\label{regularoverCMfields}
Let $F$ be a CM field, and let $\mc{R}= \{\rho_{\lambda}\}_{\lambda}$ be an absolutely irreducible, pure, Hodge-Tate regular weakly compatible system of representations of $\gal{F}$. Suppose that when we write 
\[
\rho_{\lambda} \cong \Ind_{L(\lambda)}^F(W_{\lambda}),
\]
where $W_{\lambda}$ is Lie irreducible, the extensions $L(\lambda)/F$ are Galois. Then the field $L= L(\lambda)$ is independent of $\lambda$, and $L$ is itself CM.
\end{cor}
\begin{rmk}
One way of interpreting this result is that to study regular motives, compatible systems, or algebraic automorphic representations over CM fields, we will never have to leave the comfort of CM fields; essentially all progress in the study of automorphic Galois representations is currently restricted to this context. This is in particular the case for Galois representations occurring as irreducible sub-quotients of the cohomology of a Shimura variety.
\end{rmk}

%\chapter{Motivic lifting}
\chapter{Motivic lifting}\label{motivicliftingpart}
As noted in the introduction (see Question \ref{modularlifting}), the results of Part \ref{2} raise more questions than they resolve. In this chapter, we discuss some cases of the motivic analogue of Conrad's lifting question; this is also the natural framework for the problem of finding \textit{compatible} lifts of a compatible system of Galois representations. Ideally, we would be able to work in Grothendieck's category of pure motives, defined using the relation of homological equivalence on algebraic cycles (see \S \ref{homologicalmotives} for a brief review). This category can only be proven to have the desirable categorical properties-- namely, equivalence to the category of representations of some pro-reductive (`motivic Galois') group-- if we assume Grothendieck's Standard Conjectures on algebraic cycles, which are far out of reach. (For the basic formalism (as relevant for the theory of motives) of algebraic cycles and precise statement of the Standard Conjectures, see \cite{kleiman:algcycles}.) We therefore need an unconditional variant of the motivic Galois formalism, and we adopt Andr\'{e}'s approach, using his theory of motivated cycles (\cite{andre:motivated}). In \S \ref{motivatedprelims} we review Andr\'{e}'s theory, prove some supplementary results needed for the application to motivic lifting, and then treat the motivic lifting problem in the potentially abelian case. In \S \ref{hyperlift} we prove an arithmetic refinement of Andr\'{e}'s work (\cite{andre:hyperkaehler}) on, roughly speaking, the motivated theory of hyperk\"{a}hler varieties. This provides a motivic analogue of Theorem \ref{anymonodromylift} in many non-abelian examples. Finally, in \S \ref{conclusion}, we speculate on a generalized Kuga-Satake construction, of which the results of \S \ref{hyperlift} are the `classical' case; we then prove this for $H^2$ of an abelian variety, generalizing known results for abelian surfaces.

\section{Motivated cycles: generalities}\label{motivatedprelims}
\subsection{Lifting Hodge structures}\label{liftinghodge}
In trying to produce a motivic analogue either of Wintenberger's or of my lifting theorem, one is naturally led to try to lift all cohomological, rather than merely the $\ell$-adic, realizations of a `motive.' The easiest such lifting problem is for real Hodge-structures, which are parametrized by representations of the Deligne torus $\bS= \Res_{\CC/\RR}(\mathbb{G}_m)$. Recall that $X^\bullet(\bS)= \Z \alpha_1 \oplus \Z \alpha_2$, where $\alpha_1$ and $\alpha_2$ are the first and second projections in the isomorphism
\begin{align*}
\bS(\CC)= (\CC \otimes_{\RR} \CC)^{\times} &\xrightarrow{\sim} \CC^{\times} \times \CC^{\times} \\
z_1 \otimes z_2 & \mapsto (z_1 z_2, \overline{z_1} z_2).
\end{align*}
Here $\bS(\RR) \subset \bS(\CC)$ is $\CC^{\times}= (\CC \otimes_{\RR} \RR)^{\times}$, and the $\Gal(\CC/\RR)$-action, induced by $z_1 \otimes z_2 \mapsto z_1 \otimes \overline{z_2}$, is given by $c \colon (w, z) \mapsto (\bar{z}, \bar{w})$. In particular, 
\[
(c \cdot \alpha_1)(w, z)= c(\alpha_1(\bar{z}, \bar{w})=z,
\] 
i.e. $(c \cdot \alpha_1)= \alpha_2$, and similarly $(c \cdot \alpha_2)= \alpha_1$. The group of characters over $\RR$, denoted $X^\bullet_{\RR}(\bS)$, is then $\Z(\alpha_1+\alpha_2)$, where $\alpha_1+\alpha_2=N$ is the norm, on $\RR$-points satisfying $N(z)= z \bar{z}$.

Let $H'_{0} \to H_0$ be a surjection of linear algebraic groups over $\RR$ with kernel equal to a central torus $Z_0$. We are interested in the lifting problem for algebraic representations over $\RR$:
\[
\xymatrix{
& H'_0 \ar[d]^{\pi} \\
\bS \ar@{-->}[ur]^{\tilde{h}} \ar[r]_{h} & H_0. \\
}
\]
Any such $h$ lands in some (typically non-split) maximal torus $T_0$, and any lift will land in the preimage $\pi^{-1}(T_0)=: T'_0$, which is a maximal torus of $H'_0$. We are therefore reduced to studying the dual diagram of free $\Z$-modules with $\Gal(\CC/\RR)= \gal{\RR}$-action:
\[
\xymatrix{
& 0 \\
& X^\bullet(Z_0) \ar[u] \\
& X^\bullet(T'_0) \ar@{-->}[dl] \ar[u] \\
X^\bullet(\bS) & X^\bullet(T_0) \ar[l] \ar[u] \\
& 0 \ar[u] \\
}
\]
As short-hand, we denote the vertically-aligned character groups, from bottom to top, by $Y$, $Y'$, and $L$, so we in fact are studying the sequence
\[
0 \to \Hom_{\gal{\RR}}(L, X^\bullet(\bS)) \to \Hom_{\gal{\RR}}(Y', X^\bullet(\bS)) \to \Hom_{\gal{\RR}}(Y, X^\bullet(\bS)) \to \Ext^1_{\gal{\RR}}(L, X^\bullet(\bS)) \to \ldots
\]	 
Any real torus is isomorphic to a product of copies of $\mathbb{G}_m$, $\bS$, and $\bS^1= \ker(N \colon \bS \to \mathbb{G}_m)$. The character group $X^\bullet(\bS^1)$ is $X^{\bullet}(\bS)/\Z(\alpha_1 + \alpha_2)$. A generator is the image of $\alpha_1- \alpha_2$, which on $\RR$-points is simply the character $z \mapsto z/\bar{z}$ of the (analytic) unit circle $\mr{S}^1$. Complex conjugation acts as $-1$ on $X^\bullet(\bS^1)$. We can therefore completely address the lifting problem by understanding morphisms and extensions between these three basic $\Z[\gal{\RR}]$-modules. The case of immediate interest to us will be when $Z_0$ is split, so $L$ is just some number of copies of $\Z$ with trivial $\gal{\RR}$-action. Now,
\[
\Ext^1_{\gal{\RR}}(\Z,  X^\bullet(\bS))= H^1(\gal{\RR}, X^\bullet(\bS))= 0,
\]
since $\ker(1+c)= \im(c-1)= \Z(\alpha_1-\alpha_2)$. Therefore any $h \colon \bS \to H_0$ lifts, and the ambiguity in lifting is a collection of elements (of order equal to the rank of $Z_0$) of $\Hom_{\gal{\RR}}(\Z, X^\bullet(\bS)) \cong \Z(\alpha_1+ \alpha_2)$.\footnote{Note that if we dealt with representations of $\bS^1$ we would find an obstruction in $H^1(\gal{\RR}, X^\bullet(\bS^1))= \Z/2\Z$.}
\begin{eg}
In \S \ref{KSreview} we will consider the following setup: $V_{\RR}$ will be an orthogonal space with signature $(m-2, 2)$,\footnote{This part of the discussion applies to any signature $(p, q)$ with at least one of $p$ or $q$ even, so that $\mr{SO}(V_{\RR})$ has a compact maximal torus.}. Write $m=2n$ or $m=2n+1$. The lifting problem will be
\[
\xymatrix{
& \mr{GSpin}(V_{\RR}) \ar[d]^{\pi} \\
\bS \ar@{-->}[ur]^{\tilde{h}} \ar[r]_{h} & \mr{SO}(V_{\RR}), \\
}
\]
where $h$ lands in a maximal anisotropic torus $T_0 \cong (\bS^1)^n$. We can write
%The preimage $T'_0$ is the $\RR$-torus corresponding to the element of 
%\[
%\Ext^1_{\gal{\RR}}(\Z, T_0) \cong H^1(\gal{\RR}, X^\bullet(\bS^1))^n \cong (\Z/2\Z)^n
%\]
%that is non-trivial in each component.  
\[
X^\bullet(T'_0)= \left( \bigoplus_{i=1}^n \Z \chi_i \right) \oplus \Z(\chi_0 + \frac{\sum{\chi_i}}{2}),\footnote{The characters $\chi_i$ are conjugate in $\mr{SO}(V_{\CC})$ to the characters denoted $\chi_i$ in \S \ref{spineg} (see page \pageref{spinnotation}). The torus $T_0$ is built out of copies of $\mr{SO}(2)$ embedded in $\mr{SO}(m-2, 2)$, and these are just the usual characters 
\[
\left( \begin{matrix}
\cos(\theta) & \sin(\theta) \\
-\sin(\theta) & \cos(\theta)
\end{matrix} \right) \mapsto e^{i \theta}.
\]}
\]
where conjugation acts by $-1$ on each $\chi_i$, $i=1, \ldots, n$, and trivially on $\chi_0$. Here $\oplus_{i=1}^n \Z \chi_i$ is the submodule $X^\bullet(T_0)$. A morphism $\bS \to T_0 \subset \mr{SO}(V_{\RR})$ is given in coordinates by $\chi_i \mapsto m_i(\alpha_1-\alpha_2)$, for some integers $m_i$. A lift to a morphism $\bS \to \mr{GSpin}(V_{\RR})$ then amounts to an extension 
\[
\chi_0 + \frac{\sum_{i=1}^n \chi_i}{2} \mapsto \frac{\epsilon_0}{2}(\alpha_1+ \alpha_2)+ \frac{\sum_{i=1}^n m_i}{2}(\alpha_1-\alpha_2),
\]
where $\epsilon_0$ is any integer having the same parity as $\sum m_i$. The Clifford norm $N$ is given by the character $2 \chi_0$, so the composition of such a lift with the Clifford norm is $\epsilon_0(\alpha_1+\alpha_2)$. In the $K3$ (or hyperk\"{a}hler) examples to be considered in the next section, $m_1=1$ and all other $m_i=0$, so we find there is a unique lift 
\[
\xymatrix{
& \mr{GSpin}(V_{\RR}) \ar[d]^{\pi} \\
\bS \ar@{-->}[ur]^{\tilde{h}} \ar[r]_{h} & \mr{SO}(V_{\RR}), \\
}
\]
where $N \circ \tilde{h} \colon \bS \to \mathbb{G}_m$ is any \textit{odd} power of the usual norm $\bS \to \mathbb{G}_m$. The Kuga-Satake theory takes the norm itself, which then gives rise to weight $1$ Hodge structures.
\end{eg}
\subsection{Motives for homological equivalence and the Tannakian formalism}\label{homologicalmotives}
Our aim in this sub-section is to review Grothendieck's construction of the category $\mc{M}^{hom}_F$ of pure motives for homological equivalence over a field $F$. As a preliminary, we provide some background on the theory of neutral Tannakian categories. We can then describe the output of the Standard Conjectures, taking for simplicity $F$ to be an abstract (i.e., not embedded) field of characteristic zero, small enough to be embedded in $\CC$: $\mc{M}^{hom}_F$ is (conjecturally) a graded, semi-simple, $\Q$-linear neutral Tannakian category. With the exception of \S \ref{AVlifting}, in the rest of this paper we will not work directly with $\mc{M}_F^{hom}$, so much of this discussion serves only to orient the reader. For more background, the reader should consult \cite{andre:introduction}, especially Chapter 4, or \cite{scholl:motives}.

\subsubsection{Neutral Tannakian categories}
We begin with an overview of the theory of neutral Tannakian categories; \cite{deligne-milne} is a very readable and thorough introduction to which we refer the reader for details (the original source is \cite{saavedrarivano}). Note that we will always work with neutral Tannakian categories, which simplifies the theory considerably; for deeper aspects in the non-neutral case, see \cite{deligne:tannakian}. Let $E$ be a field. The prototypical neutral Tannakian category is the category $\mr{Vec}_{E}$ of finite-dimensional vector spaces over $E$. $\mr{Vec}_E$ is an abelian category with a notion of tensor product (namely, the usual tensor product of $E$-vector spaces) satisfying certain natural requirements: associativity, commutativity, and existence of a unit ($E$ regarded as an $E$-vector space) for the tensor product. Moreover, every object of $\mr{Vec}_E$ has a dual (more generally, internal $\Hom$'s exist), and the endomorphisms of the unit object are just the field $E$ itself. The main theorem of the theory in this context is the trivial observation that $\mr{Vec}_E$ is equivalent to the category of finite-dimensional representations (over $E$) of the trivial group.

The general definition merely formalizes the notions in the previous paragraph; we will quickly make the necessary definitions, and then give a number of examples (Example \ref{tannakianeg} below).
\begin{defn}
Let $\mc{C}$ be a category, and let $\otimes \colon \mc{C} \times \mc{C} \to \mc{C}$ be a functor satisfying the usual axioms of a symmetric monoidal category, namely:
\begin{enumerate}
\item (Associativity constraint) There is a functorial isomorphism
\[
A_{X, Y, Z} \colon X \otimes (Y \otimes Z) \to (X \otimes Y) \otimes Z
\]
satisfying the \textit{pentagon axiom} (\cite[1.0.1]{deligne-milne}).
\item (Commutativity constraint) There is a functorial isomorphism
\[
C_{X, Y} \colon X \otimes Y \to Y \otimes X
\]
such that $C_{Y, X} \circ C_{X, Y}= \id_{X \otimes Y}$, and that is compatible with the associativity constraint in the sense of the \textit{hexagon axiom} (\cite[1.0.2]{deligne-milne}).
\item (Unit object) There is a pair $(\mbf{1}, e)$ consisting of an object $\mbf{1}$ and an isomorphism $e \colon \mbf{1} \to \mbf{1} \otimes \mbf{1}$ such that the functor $\mc{C} \to \mc{C}$ given by $X \mapsto \mbf{1} \otimes X$ is an equivalence of categories.
\end{enumerate}
We call such a $(\mc{C}, \otimes)$, equipped with its associativity and commutativity constraints (but omitted from the notation), a \textit{tensor category}, for short.\index{t}{tensor category}

If $(\mc{C}, \otimes)$ and $(\mc{C}', \otimes')$ are tensor categories (whose associativity and commutativity constraints we will write as $A, C$ and $A', C'$, respectively; unit elements will be $(\mbf{1}, e)$ and $(\mbf{1}', e')$), then a \textit{tensor functor} from $(\mc{C}, \otimes)$ to $(\mc{C}', \otimes')$ is a pair $(F, k)$ consisting of a functor $F \colon \mc{C} \to \mc{C}$ and a functorial isomorphism $k_{X, Y} \colon F(X) \otimes F(Y) \to F(X \otimes Y)$ satisfying the following three compatibilities:
\begin{enumerate}
\item For all objects $X, Y, Z$ of $\mc{C}$, the diagram
\[
\xymatrix{
FX \otimes (FY \otimes FZ) \ar[d]^{A'} \ar[r]^{\id \otimes k} & FX \otimes F(Y \otimes Z) \ar[r]^k & F(X \otimes (Y \otimes Z)) \ar[d]^{F(A)} \\
(FX \otimes FY) \otimes FZ \ar[r]^{k \otimes \id} & F(X \otimes Y) \otimes FZ \ar[r]^{k} & F((X \otimes Y) \otimes Z)
}
\]
is commutative.
\item For all objects $X, Y$ of $\mc{C}$, the diagram
\[
\xymatrix{
FX \otimes FY \ar[r]^k \ar[d]^{C'} & F(X \otimes Y) \ar[d]^{F(C)} \\
FY \otimes FX \ar[r]^k & F(Y \otimes X)
}
\]
is commutative.
\item If $(\mbf{1}, e)$ is a unit object of $\mc{C}$, then $(F(\mbf{1}), F(e))$ is a unit object of $\mc{C}'$.
\end{enumerate}
\end{defn}
We have spelled out the precise conditions for $(F, k)$ to be a tensor functor because of its relevance for the structure of $\mc{M}^{hom}_F$: the K\"{u}nneth Standard Conjecture essentially `corrects' the fact that the Betti realization
\[
H_B \colon \mc{M}^{hom}_F \to \mr{Vec}_{\Q}
\]
with the natural K\"{u}nneth isomorphism $k_{X, Y} H_B(X) \otimes H_B(Y) \xrightarrow{\sim} H_B(X \times Y)$ does \textit{not} intertwine the (obvious) commutativity constraints on $\mc{M}^{hom}_F$ and $\mr{Vec}_{\Q}$. See the discussion surrounding Conjecture \ref{kunnethSC} for details.

For any unit object $(\mbf{1}, e)$, we obtain isomorphisms $l_X \colon \mbf{1} \otimes X \xrightarrow{\sim} X$ and $r_X \colon X \otimes \mbf{1} \xrightarrow{\sim} X$. 

For any additive tensor category $(\mc{C}, \otimes)$ (for which we require $\otimes$ to be bi-additive), and any unit object $(\mbf{1}, e)$, $\End_{\mc{C}}(\mbf{1})= R$ is a commutative ring, unique up to unique isomorphism, and $\mc{C}$ is naturally an $R$-linear category. 
\begin{defn}
Let $(\mc{C}, \otimes)$ be a tensor category. It is \textit{rigid} if for every object $X$ there is a `dual' object $X^\vee$ along with evaluation and co-evaluation morphisms
\begin{align*}
X^\vee \otimes X \xrightarrow{\mr{ev}_X} \mbf{1} \\
\mbf{1} \xrightarrow{\mr{coev}_X} X \otimes X^\vee
\end{align*}
such that the composites (we suppress the unit isomorphisms and the associators)
\begin{align*}
X \xrightarrow{\mr{coev}_X \otimes \id_X} X \otimes X^\vee \otimes X \xrightarrow{\id_X \otimes \mr{ev}_X} X \\
X^\vee \xrightarrow{\id_{X^\vee} \otimes \mr{coev}_X} X^\vee \otimes X \otimes X^\vee \xrightarrow{\mr{ev}_X \otimes \id_{X^\vee}} X^\vee
\end{align*}
are $\id_X$ and $\id_{X^\vee}$, respectively.
\index{t}{rigid tensor category}
\end{defn}
Note that in the absence of a commutativity constraint, $X^\vee$ is what would be called a `right dual.'

Finally, we come to the main definition:
\begin{defn}
Let $E$ be a field. A \textit{neutral Tannakian category} over $E$ is a rigid abelian tensor category $(\mc{C}, \otimes)$ such that $\End(\mbf{1})= E$, and for which there exists a faithful, exact, $E$-linear tensor functor $\omega \colon \mc{C} \to \mr{Vec}_E$.\index{t}{neutral Tannakian category} Such an $\omega$ is called a \textit{fiber functor}.\index{t}{fiber functor} 
\end{defn}
The main theorem of the theory of neutral Tannakian categories just says that every Tannakian category is equivalent to the category of representations of some affine group scheme:
\begin{thm}[{\cite[Theorem 2.11]{deligne-milne}}]\label{maintannakianthm}
Let $(\mc{C}, \otimes)$ be a neutral Tannakian category over $E$, equipped with a fiber functor $\omega \colon \mc{C} \to \mr{Vec}_E$. Then the functor of $E$-algebras given by tensor-automorphisms of $\omega$ (see \cite[1.9, 1.11]{deligne:tannakian}) is represented by an affine group scheme $G$, and the functor $\mc{C} \to \Rep_{E}(G)$ defined by $\omega$ is an equivalence of tensor categories.\footnote{That is, it is a tensor functor that is also an equivalence of categories; this ensures (\cite[Proposition 1.11]{deligne-milne}) that there exists an inverse of $F$ that is also a tensor functor.}
\end{thm}
Informally, we might call $G$ the `Galois group' of $(\mc{C}, \otimes)$, whence the later terminology `motivic Galois group' in the (conjecturally Tannakian) case of $\mc{M}_F^{hom}$.

At last, some examples:
\begin{eg}\label{tannakianeg} 
\begin{enumerate}
\item The following example is crucial for the theory of motives. The category $\mr{Vec}_E^{\Z/2}$ of $\Z/2\Z$-graded vector spaces over $E$, with the usual graded tensor product, and with commutativity constraint given by the Koszul sign rule, i.e. $C_{V, W} \colon V \otimes W \xrightarrow{\sim} W \otimes V$ given on homogeneous tensors by $v \otimes w \mapsto (-1)^{\deg(v) \cdot \deg(w)} w \otimes v$, is a rigid $E$-linear abelian tensor category. Note, however, that the forgetful functor $\mr{Vec}_E^{\Z/2\Z} \to \mr{Vec}_E$ (with the obvious isomorphisms $k_{X, Y}$) is not a tensor functor, because it fails to intertwine the commutativity constraints. Indeed, $\mr{Vec}_E^{\Z/2\Z}$ is \textit{not} a Tannakian category. Any rigid tensor category has an intrinsic notion of rank for every object $X$, given by the element of $\End(\mbf{1})$ arising as the composition
\[
\mbf{1} \xrightarrow{\mr{coev}_X} X \otimes X^\vee \xrightarrow{C_{X, X^\vee}} X^\vee \otimes X \xrightarrow{\mr{ev}_X} \mbf{1}.
\]
For a $\Z/2\Z$-graded vector space $V= V_0 \oplus V_1$, it is easy to check that the rank in $\mr{Vec}_E^{\Z/2\Z}$ is $\dim V_0-\dim V_1$. But this notion of rank is preserved by any tensor functor, so an obviously necessary condition to admit a fiber functor is that all objects have non-negative rank. (In fact, a deep result of Deligne gives a converse: \cite[7.1 Th\'{e}or\`{e}me]{deligne:tannakian}.)
\item The category $\Q-\mr{HS}^{pol}$ of pure polarizable $\Q$-Hodge structures is a neutral Tannakian category over $\Q$, with fiber functor just given by taking the underlying $\Q$-vector space of a Hodge structure (here we take the commutativity constraint on Hodge structures to be the same as the usual constraint on vector spaces). This is a useful `toy model' for the theory of pure motives over $\CC$. In particular, Theorem \ref{maintannakianthm} identifies $\Q-\mr{HS}^{pol}$ with the category of representations of some affine group scheme $\mr{MT}$ over $\Q$; $\mr{MT}$ is the `universal Mumford-Tate group.' We use this example to explain how properties of a Tannakian category reflect those of its Galois group (see \cite[\S 2 `Properties of $G$ and $\Rep(G)$']{deligne-milne}). From its Hodge-theoretic description, $\mr{MT}$ is obviously connected. This is equivalent (\cite[Corollary 2.22]{deligne-milne}) to the condition that for every non-trivial representation $X$ of $\mr{MT}$, the strictly full subcategory of $\Q-\mr{HS}^{pol}$ whose objects are isomorphic to sub-quotients of some $X^{\oplus n}$ is not stable under $\otimes$. But if $X$ is a non-trivial Hodge structure, which we may assume pure of weight zero (otherwise the result is clear by comparing weights), then some $H^{p, -p}(X)$ is non-zero, for some $p \neq 0$. Then $H^{mp, -mp}(X^{\otimes m}) \neq 0$, which for $m$ large enough cannot be the case for any sub-quotient of any $X^{\oplus n}$. In a similar spirit, $\mr{MT}$ is pro-reductive, since (\cite[Proposition 2.23]{deligne-milne}), since the category $\Q-\mr{HS}^{pol}$ is semi-simple (by polarizability). 
\end{enumerate}
\end{eg}
\subsubsection{Homological motives}\label{hommotivessection}
We begin by sketching the construction of Grothendieck's category of (pure) homological motives over a field $F$. References for more details are \cite[\S 1]{kleiman:algcycles} and \cite[\S 1]{scholl:motives}. Recall that for simplicity we take $F$ to be a field that can be embedded in $\CC$ (or even a specified subfield of $\CC$)-- this does not affect the construction, but does affect whether the category of motives is (conjecturally) \textit{neutral} Tannakian.  Recall that for a fixed Weil cohomology theory $H^*$ (see \cite[\S 1.2]{kleiman:algcycles}), homological equivalence defines an \textit{adequate equivalence relation} on algebraic cycles on smooth projective varieties; for example, since $F\subset \CC$, we will always have at our disposal Betti cohomology, $H_B^*$.\index{s}{$H_B^*$} For a smooth projective $X/F$, we let $A^*_{hom}(X)$ denote the graded, by codimension, $\Q$-algebra of ($\Q$-linear combinations of) algebraic cycles for homological equivalence.\index{s}{$A^*_{hom}(X)$} `Adequate' ensures that the intersection product is well-defined as a linear map $A^{r}(X) \otimes A^s(X) \to A^{r+s}(X)$. For any two smooth projective $F$-varieties $X$ and $Y$, with $X$ connected of (equi-)dimension $d$, we can then define the space of degree $r$ correspondences\index{t}{algebraic correspondence}\index{s}{$C^*_{hom}(X, Y)$}
\[
C^r_{hom}(X, Y)= A^{d+r}_{hom}(X \times Y).
\]
We obtain a composition of correspondences
\[
C^r_{hom}(X, Y) \otimes C^s_{hom}(Y, Z) \to C^{r+s}_{hom}(X, Z)
\]
given by 
\[
f \otimes g = g \circ f \mapsto p_{13, *}\left( p_{12}^* f \cdot p_{23}^* g \right),
\]
where the $p_{ij}$ are the projections from $X \times Y \times Z$ to the products of two factors. In particular, $C^0_{hom}(X, X)$ is a $\Q$-algebra.
\begin{defn}\label{hommotivesdef}
The category $\mc{M}_F^{hom}$ of motives over $F$ for homological equivalence has as objects triples $(X, p, m)$, where $X$ is a smooth projective variety over $F$, $m$ is an integer, and $p$ is an idempotent correspondence in $C^0_{hom}(X, X)$. Morphisms in $\mc{M}_F^{hom}$ are defined by\index{t}{motive for homological equivalence}
\[
\Hom_{\mc{M}^{hom}_F}( (X, p, m), (Y, q, n))= q C^{n-m}_{hom}(X, Y) p.
\]
\end{defn}
$\mc{M}^{hom}_F$ is an additive, $\Q$-linear, pseudo-abelian category (\cite[Theorem 1.6]{scholl:motives}). There is a bi-additive functor $\otimes \colon \mc{M}_F^{hom} \times \mc{M}_F^{hom} \to \mc{M}_F^{hom}$ given on objects by
\[
(X, p, m) \otimes (Y, q, n)= (X \times Y, p \otimes q, m+n)
\]
(the fiber product is over $F$; we do not here specify $\otimes$ on morphisms). This satisfies associativity and commutativity constraints induced by the tautological isomorphisms $(X \times Y) \times Z \cong X \times (Y \times Z)$ and $C \colon X \times Y \cong Y \times X$. A unit object for $\otimes$ is given by $(\Spec F, \id, 0)$, and we can define dual objects by (again for simplicity take $X$ of equi-dimension $d$)
\[
(X, p, m)^\vee= (X, {}^t p, d-m);
\] 

In sum, these structures make $(\mc{M}^{hom}_F, \otimes)$ into a rigid $\Q$-linear tensor category with $\End_{\mc{M}_F^{hom}}(\mbf{1})= \Q$. Grothendieck's Standard Conjectures anticipate that moreover $\mc{M}_F^{hom}$ should be a semi-simple neutral Tannakian category. One must prove that
\begin{itemize}
\item $\mc{M}^{hom}_F$ is abelian;
\item the abelian category $\mc{M}^{hom}_F$ is moreover semi-simple;
\item $\mc{M}^{hom}_F$ possesses a fiber functor to $\mr{Vec}_{\Q}$.
\end{itemize}
We begin by discussing the question of a fiber functor. The answer is simple: $(\mc{M}_F^{hom}, \otimes)$, with the commutativity constraint defined above, does \textit{not} have a fiber functor! For simplicity, let us take $F \subset \CC$, so that there is a Betti realization 
\[
H_B \colon \mc{M}_F^{hom} \to \mr{Vec}_{\Q}.
\]
$H_B$ is not a tensor functor, however: the diagram
\[
\xymatrix{
H_B(X) \otimes H_B(Y) \ar[d]^{C'} \ar[r]^{k}_{\sim} & H_B(X \times Y) \ar[d]^{H_B(C)} \\
H_B(Y) \otimes H_B(X) \ar[r]^k_{\sim} & H_B(Y \times X),
}
\]
given by applying the na\"{i}ve commutativity constraints and the K\"{u}nneth isomorphism, only commutes up to sign, since cup-product is anti-commutative. More precisely, $H_B$ gives a tensor functor valued in $\mr{Vec}_{\Q}^{\Z/2\Z}$, and since objects of $\mc{M}_F^{hom}$ can then have negative rank, it cannot be neutral Tannakian with the commutativity constraint induced by $X \times Y \xrightarrow{\sim} Y \times X$. The \textit{K\"{u}nneth Standard Conjecture} \index{t}{K\"{u}nneth Standard Conjecture} would give $\mc{M}_F^{hom}$ a grading (corresponding to cohomological degree), allowing for a modified commutativity constraint, for which $H_B \colon \mc{M}_F^{hom} \to \mr{Vec}_{\Q}$ \textit{is} a tensor functor. More precisely, for any Weil cohomology theory $H^*$, the isomorphism (again taking $d= \dim X$)
\[
H^{2d}(X \times X)(d) \cong \bigoplus_i H^{2d-i}(X) \otimes H^{i}(X)(d) \cong \bigoplus_i H^i(X)^\vee \otimes H^i(X) \cong \bigoplus_i \End(H^i(X))
\] 
gives a cohomology class $\pi^i_X$\index{s}{$\pi^i_X$} corresponding to the composition $H^*(X) \onto H^i(X) \into H^*(X)$.
\begin{conj}\label{kunnethSC}
For all $i= 0, 1, \ldots, 2d$, the cohomology class $\pi_X^i$ is algebraic, i.e. lies in the image of the cycle class map $A_{hom}^d(X \times X) \to H^{2d}(X \times X)(d)$.
\end{conj}
The modified commutativity constraint on $\mc{M}_F^{hom}$ would then be given as follows: if the original constraint is $C \colon M \otimes N \xrightarrow{\sim} N \otimes M$, we can decompose $C = \oplus_{r, s} C^{r, s}$ with
\[
C^{r, s} \colon \pi^rM \otimes \pi^s N \xrightarrow{\sim} \pi^s N \otimes \pi^r M,
\]
and then the correct commutativity constraint is $C'= \oplus_{r, s} (-1)^{rs} C^{r, s}$.\index{t}{modified commutativity constraint for $\mc{M}_F^{hom}$}

The K\"{u}nneth Standard Conjecture is in fact implied by a stronger conjecture, the Lefschetz Standard Conjecture, which would also show that the primitive decomposition of cohomology' makes sense in the category $\mc{M}^{hom}_F$. To state this, let $X$ be a smooth projective variety over $F$, and fix an ample line bundle $\eta$ on $X$, giving rise to the Lefschetz operator 
\[
L=L_{\eta, H^*} \colon H^i(X)(r) \to H^{i+2}(X)(r+1),  
\]
and the hard Lefschetz isomorphisms
\[
L^{d-i} \colon H^i(X)(r) \xrightarrow{\sim} H^{2d-i}(X)(d-i+r)
\]
for all $i \leq d$. As always with taking cup-product with the class of an algebraic cycle, these isomorphisms are given by algebraic correspondences.
\begin{conj}[Lefschetz Standard Conjecture]
For all $i \leq d$, the inverse of $L^{d-i}$ is given by algebraic correspondences.
\end{conj}
Note that \textit{a priori} both the K\"{u}nneth and Lefschetz Standard Conjectures depend on the choice of Weil cohomology. For the classical cohomologies (Betti, \'{e}tale, de Rham) that are related by comparison isomorphisms (respecting cycle class maps), these conjectures do not depend on the choice of $H^*$.

Now we move to the question of whether $\mc{M}_F^{hom}$ is abelian. Here there is a marvelous theorem of Jannsen:
\begin{thm}[ {\cite[Theorem 1]{jannsen:numerical}} ]
For an adequate equivalence relation $\sim$, the category of motives for $\sim$-equivalence is abelian semi-simple if and only if $\sim$ is numerical equivalence. In particular, $\mc{M}^{hom}_F$ is abelian if and only if numerical and homological equivalence coincide.
\end{thm}
Indeed, long before Jannsen proved his theorem, Grothendieck conjectured:
\begin{conj}
Homological and numerical equivalence coincide.
\end{conj}
We should also remark that semi-simplicity is closely related to the so-called Hodge Standard Conjecture (not to be confused with the Hodge conjecture), which in characteristic zero is a simple consequence of the Hodge index theorem (see \cite[\S3]{kleiman:algcycles}).

In summary, under the Standard Conjectures, we can equip $\mc{M}_F^{hom}$ with a commutativity constraint for which it is a semi-simple neutral Tannakian category, hence:
\begin{conj}
Let $F$ be a subfield of $\CC$, for simplicity. There is a pro-reductive $\Q$-group $\mc{G}_F^{hom}$, the \textit{motivic Galois group} for pure motives over $F$, and an equivalence (given by any choice of fiber functor $\mc{M}_F^{hom} \to \mr{Vec}_{\Q}$; an example is $H_B^*$)
\[
\mc{M}_F^{hom} \cong \Rep_{\Q}(\mc{G}_F^{hom}).
\]
\end{conj}
\subsection{Motivated cycles}\label{motivatedcycles}
In the present Part \ref{motivicliftingpart} of these notes, most of our results can be stated solely in terms of abelian varieties, but both the strongest assertions and the proofs will require invoking some version of the motivic Galois formalism; in the absence of the Standard Conjectures, we use Andr\'{e}'s theory of motives for motivated cycles (\cite{andre:motivated}). In this subsection, we provide a brief review, elaborating on some points for later application. Although Andr\'{e}'s theory is developed over any field $F$, to simplify we will always take $F$ to be an abstract (i.e., not embedded) field of characteristic zero, small enough to be embedded in $\CC$, and which we will eventually specify to be a number field. 

In \cite{andre:motivated}, Andr\'{e} defines a $\Q$-linear category of motives for `motivated cycles'  whose construction mirrors the classical construction of (Grothendieck) motives for homological equivalence, but circumvents the standard conjectures by formally enlarging the group of `cycles' to include the Lefschetz involutions. Here is the precise definition of the space of motivated cycles:\index{t}{motivated cycles}
\begin{defn}\label{motivatedcycle}
Fix a reference Weil cohomology $H^\bullet$. A motivated cycle on $X$ with coefficients in $E$ is an element of $H^\bullet(X)$ of the form 
\[
\pr^{XY}_{X, *}(\alpha \cup \ast_L \beta),
\] 
where
\begin{itemize}
\item $Y$ is a smooth projective $F$-scheme, with polarization $\eta_{Y}$ giving rise to a `product' polarization $\eta_{X \times Y}= [X] \otimes \eta_Y+ \eta_X \otimes [Y]$, with corresponding Lefschetz involution $\ast_L$, on $X \times Y$;
\item $\alpha$ and $\beta$ are algebraic cycles with $E$-coefficients on $X \times Y$.
\end{itemize}
We denote by $A^\bullet_{mot}(X)_E$ the $E$-vector space of motivated cycles on $X$ with $E$-coefficients. $A^\bullet_{mot}(X)$ will always mean the case $E=\Q$.
\end{defn}
As in \S \ref{hommotivessection}, we can then define spaces $C^\bullet_{mot}(X, Y)_E$\index{s}{$C^\bullet_{mot}(X, Y)_E$}\index{t}{motivated correspondences} of motivated correspondences (see also \cite[D\'{e}finition 2]{andre:motivated}); for $X=Y$, this construction yields a graded $E$-algebra containing the Lefschetz involutions and the K\"{u}nneth projectors $\pi^i_X$. 

The first main result of this theory asserts that the spaces of motivated cycles in a precise sense do not depend on the Weil cohomology $H^\bullet$ used to define them:
\begin{thm}[Th\'{e}or\`{e}me $0.3$ of \cite{andre:motivated}]\label{foundations}    
For any smooth projective $F$-scheme $X$, let $A^\bullet_{mot}(X)$ be the graded $\Q$-algebra  of `motivated cycles,' constructed with respect to a fixed Weil cohomology $H^\bullet$. $A^\bullet_{mot}(X)$ contains the classes of algebraic cycles modulo homological equivalence, and there is a $\Q$-linear injection
\[
cl_H \colon A^\bullet_{mot}(X) \to H^{2\bullet}(X)
\]
extending the cycle class map for $H$. $A^\bullet_{mot}(X)$ has the following properties:
\begin{itemize}
%\item for every element $\xi$ of $A^\bullet_{mot}(X)$, there is a smooth projective $F$-scheme $Y$, algebraic cycles $\alpha$ and $\beta$ on $X \times Y$, and a polarization on $X \times Y$ induced by a `product' of polarizations of $X$ and $Y$, and with corresponding Lefschetz involution $\ast_L$, such that
%\[
%cl_H(\xi)= \pr^{XY}_{X, *} (cl_H(\alpha) \cup \ast_{L} cl_H(\beta)).
%\] 
\item $A^\bullet_{mot}(X)$ depends bifunctorially (push-forward and pull-back) on $X$, satisfying the usual projection formula.
\item (See \cite[\S 2.3]{andre:motivated}) Two Weil cohomologies related by comparison isomorphisms yield canonically and functorially isomorphic algebras $A^\bullet_{mot}(X)$ of motivated cycles. In particular, all the classical cohomology theories yield the same $A^{\bullet}_{mot}(X)$.
\end{itemize}
\end{thm}
From here, Andr\'{e} defines a category $\mc{M}_F$\index{s}{$\mc{M}_F$} of motives for motivated cycles\index{t}{motives for motivated cycles} exactly as in Definition \ref{hommotivesdef}, replacing algebraic correspondences with motivated correspondences.\footnote{Andr\'{e} defines more generally a category of motives `modeled on' a full sub-category $\mc{V}$ of the category of smooth projective $F$-schemes, with $\mc{V}$ assumed stable under products, disjoint union, and passage to connected components. This amounts to restricting the auxiliary varieties $Y$ permitted in Definition \ref{motivatedcycle}. We will always take $\mc{V}$ to be all smooth projective $F$-schemes.}  Andr\'{e} shows the endomorphism algebras of objects of $\mc{M}_F$ are semi-simple (finite-dimensional) $\Q$-algebras, from which it follows (see \cite{jannsen:numerical}) that $\mc{M}_F$ is an abelian semi-simple category. Since the K\"{u}nneth projectors $\pi_X^i$ are motivated cycles, we can define the modified commutativity constraint described in \S \ref{hommotivessection}, thereby making $\mc{M}_F$ into a neutral Tannakian category over $\Q$:
%\begin{defn}
%Let $\mc{M}_F$, the category of motives over $F$ for motivated cycles, be the category whose objects are triples $M= (X, q, n)$ consisting of a smooth projective $F$-scheme $X$, an idempotent motivated correspondence $q \in C^0_{mot}(X, X)$, and a `Tate twist' $n$, which is a choice of integer for each connected component of $X$. For short-hand (and intuition), we write $M= qH(X)(n)$. The morphisms in $\mc{M}_F$ are given by
%\[
%\Hom_{\mc{M}_F}(qH(X)(n), pH(Y)(m))= p C^{m-n}_{mot}(X, Y) q.
%\] 
%\end{defn}
\begin{thm}[Th\'{e}or\`{e}me $0.4$ of \cite{andre:motivated}]
$\mc{M}_F$ is a neutral Tannakian category over $\Q$. It is graded, semi-simple, and polarized. Every classical cohomology factors through $\mc{M}_F$.
\end{thm}
\begin{rmk}
We will often use the short-hand $H(X)$ for the object $(X, \mr{id}, 0)$, and $H^i(X)$ for $(X, \pi_X^i, 0)$ (recall the notation from Definition \ref{hommotivesdef}). It will be clear from context when $H(X)$ refers to an object of $\mc{M}_F$, and when it refers to the output of any particular Weil cohomology theory $H^*$.
\end{rmk}
In particular, this sets in motion the formalism of motivic Galois groups, and the theory becomes a very useful circumvention of the standard conjectures. Perhaps its most unsatisfactory feature-- present also in the theory of absolute Hodge cycles-- is that its `motives' are not known to give rise to compatible systems of $\ell$-adic representations. 

Now we specify fiber functors and no longer regard $F$ as an abstract field. For any $\sigma \colon F \into \CC$, $\mc{M}_F$ is Tannakian and neutralized by the $\sigma$-Betti fiber functor (denoted $H_{\sigma}$), so we obtain its Tannakian group, `the' motivic Galois group, $\mc{G}_{F}(\sigma)$. For the $\ell$-adic fiber functor $X \mapsto H_{et}(X_{\overline{F}}, \Ql)$, we denote by $\mc{G}_{F, \ell}$ the corresponding motivic group (over $\Ql$). It is more convenient to choose an embedding $\sigma \colon \overline{F} \into \CC$, since this allows us, via the comparison isomorphisms
\[
H_{et}(X_{\overline{F}}, \Ql) \xrightarrow[\sigma^*]{\sim} H_{et}(X_{\overline{F}} \otimes_{\overline{F}, \sigma} \CC, \Ql) \cong H_{\sigma}(X, \Q) \otimes_{\Q} \Ql,
\] 
to deduce an isomorphism $\mc{G}_{F}(\sigma) \otimes_{\Q} \Ql \cong \mc{G}_{F, \ell}$. Eventually, $F$ will simply be regarded as a subfield of $\CC$, with $\overline{F}$ its algebraic closure in $\CC$; in that case, we will omit the $\sigma$ from the notation. For now, however, we retain it.

As a variant on this formalism, to any object (or collection of objects) $M$ of $\mc{M}_F$, we can associate the smallest Tannakian subcategory $\langle M \rangle^\otimes$ generated by $M$ (fully faithfully embedded in $\mc{M}_F$), and we can then look at its Tannakian group, denoted $\mc{G}^M_F(\sigma)$.\index{s}{$\mc{G}^M_F$}\index{t}{motivic Galois group of an object $M$ of $\mc{M}_F$} A particularly important example comes from taking $M$ to range over (the objects of $\mc{M}_F$ associated to) all finite \'{e}tale $F$-schemes; this defines the sub-category of \textit{Artin motives}.\index{t}{Artin motive} The fully faithful inclusion of the subcategory $\mc{M}^{art}_F$ of Artin motives over $F$ induces a surjection
\[
\mc{G}_F(\sigma) \onto \mc{G}^{art}_F(\sigma) \cong \gal{F},
\]
and when combined with the base-change functor $\mc{M}_F \to \mc{M}_{\overline{F}}$, we obtain an exact sequence\footnote{See Proposition $6.23a,c$ of \cite{deligne-milne}. This is a corrected, TeXed version of the original article in \cite{DMOS}; it is available at http://www.jmilne.org/math/xnotes/index.html. Note that part $b$ is a modification of the (only conjectural) statement in the original article.} of pro-algebraic groups
\[
1 \to \mc{G}_{\overline{F}}(\sigma) \to \mc{G}_F(\sigma) \to \gal{F} \to 1.
\]
By Proposition $6.23d$ of \cite{deligne-milne}, for all $\ell$ there are continuous sections (homomorphisms) $\gal{F} \to \mc{G}_F(\Ql)$ (after unwinding everything, this is simply the statement that $\gal{F}$ acts on $\ell$-adic cohomology). 

We can similarly compare the `arithmetic' and `geometric' versions of the motivic Galois group of any object $M$ of $\mc{M}_F$: there is a commutative diagram, where the vertical morphisms are surjective:
\[
\xymatrix{
1 \ar[r] & \mc{G}^{M_{\overline{F}}}_{\overline{F}}(\sigma) \ar[r] & \mc{G}^{M}_F(\sigma) & & \\
1 \ar[r] & \mc{G}_{\overline{F}}(\sigma) \ar[u] \ar[r] & \mc{G}_F(\sigma) \ar[u] \ar[r] & \gal{F} \ar[r] & 1.
}
\]
\begin{lemma}\label{connected}
Assume that $\mc{G}^{M_{\overline{F}}}_{\overline{F}}(\sigma)$ is connected. Then there exists a finite extension $F'/F$ such that $\mc{G}^{M_{F'}}_{F'}(\sigma)$ is connected.
\end{lemma}
\proof
Fix a prime $\ell$, a section $s_{\ell} \colon \gal{F} \to \mc{G}_F(\Ql)$, and a finite extension $F'/F$ such that the Zariski closure of the $\ell$-adic representation $\rho_{M_{\ell}} \colon \gal{F'} \to \mr{GL}(M_{\ell})$ is connected. The image of $\mc{G}^{M_{F'}}_{F'} \otimes_{\Q} \Ql$ in $\mr{GL}(M_{\ell})$ is then equal to the product of the two connected groups $\overline{\rho_{M_{\ell}}(\gal{F'})}^{Zar}$ and $\left( \mc{G}^{M_{\overline{F}}}_{\overline{F}} \otimes_{\Q} \Ql \right)$, hence is itself connected.
\endproof
\begin{rmk}
If in $\langle M_{\overline{F}}\rangle^\otimes$ all Hodge cycles are motivated, $\mc{G}^{M_{\overline{F}}}_{\overline{F}}(\sigma)$ is connected.
\end{rmk}

The assertion that Hodge cycles are motivated is a (still fiercely difficult) weakened version of the Hodge conjecture. We now describe what is essentially the only understood case of this problem. Many of the results of \cite{andre:hyperkaehler} rest on earlier work of Andr\'{e} (Th\'{e}or\`{e}me $0.6.2$ of \cite{andre:motivated}) showing that on a complex abelian variety, all Hodge cycles are motivated (a variant of Deligne's result that Hodge cycles on a complex abelian variety are absolutely Hodge). One useful consequence is:
\begin{cor}[Andr\'{e}]\label{AVMT}
Let $A/\CC$ be an abelian variety, and let $M$ be the motive $H^1(A)$ (an object of $\mc{M}_{\CC}$). Then the motivic group $\mc{G}^M_{\CC}$ (for the Betti realization) is equal to the Mumford-Tate group $MT(A)$, and in particular is connected. 
\end{cor} 
\proof
\index{t}{Mumford-Tate group of a $\Q$-Hodge structure}\index{s}{$MT(A)$}$MT(A)$, recall, is the smallest $\Q$-sub-group of $\mr{GL}(H^1_B(A_{\CC}, \Q))$ whose $\RR$-points contain the image of the $\bS$-representation corresponding to the $\RR$-Hodge structure $H^1_B(A_{\CC}, \RR)$; by the general theory of Mumford-Tate groups, this is equal to
\begin{itemize}
\item the Tannakian group for the Tannakian category of $\Q$-Hodge structures generated by $M_B= H^1_B(A_{\CC}, \Q)$; and
\item the subgroup of $\mr{GL}(M_B)$ fixing exactly the Hodge tensors in every tensor construction $T^{m, n}(M_B):= (M_B)^{\otimes m} \otimes (M_B^\vee)^{\otimes n}$.
\end{itemize}
Similarly, the motivic group $\mc{G}^M_{\CC}$ is the subgroup of $\mr{GL}(M_B)$ fixing exactly the motivated cycles in every tensor construction $T^{m, n}(M_B)$. Of course (on any variety) all motivated cycles are Hodge cycles, so there is a quite general inclusion $MT(A) \subset \mc{G}^M_{\CC}$. Applying Andr\'{e}'s result to all powers of $A$, we can deduce the reverse inclusion: if $t \in T^{m, n}(M_B)$ is a Hodge cycle, then (by weight considerations) $m=n$, and viewing this tensor space (via K\"{u}nneth and polarization) inside $H^{2m}(A^{2m}, \Q)(m)$, we see that $t$ is motivated.
\endproof
The following point is implicit in the arguments of \cite{andre:hyperkaehler}, but we make it explicit in part to explain an important foundational point in the theory of motivated cycles.
\begin{cor}\label{MT}
Let $F$ be a subfield of $\CC$, with algebraic closure $\overline{F}$ in $\CC$, and let $A/F$ be an abelian variety. Let $M$ be the object of $\mc{M}_{\overline{F}}$ corresponding to $H^1(A_{\overline{F}})$. Then $\mc{G}^M_{\overline{F}}$ (for the $\overline{F} \subset \CC$ Betti realization) is connected, equal to $MT(A_{\CC})$.
\end{cor}
\proof
We deduce this from the previous result and the following general lemma, which is of course implicit in \cite{andre:motivated}:
\begin{lemma}\label{indsepclosed}
Let $L/K$ be an extension of algebraically closed fields, and let $M$ be an object of $\mc{M}_K$, with base-change $M|_{L}$ to $L$. Then via the canonical isomorphism $H_{et}(M, \Ql) \cong H_{et}(M|_L, \Ql)$, the $\ell$-adic motivic groups $\mc{G}^M_{K, \ell}$ and $\mc{G}^{M|_L}_{L, \ell}$ agree.
\end{lemma}
\proof
This follows from the fact ($2.5$ \textit{Scolie} of \cite{andre:motivated}) that the comparison isomorphism (for $M= H^*(X)$) identifies the spaces of motivated cycles $A^*(X) \xrightarrow{\sim} A^*(X_L)$. This follows from standard spreading out techniques, but we provide some details, since they are omitted from \cite{andre:motivated}. Writing $L= \underset{\lambda}{\varinjlim} K_{\lambda}$ as the directed limit of its finite-type $K$-sub-algebras allows all the data required to define a motivated cycle on $X_L$ (a smooth projective $Y/L$, algebraic cycles on $X_L \times Y$, the Lefschetz involution on cohomology of $X_L \times Y$) to be descended to some $S_{\lambda}= \Spec(K_{\lambda})$. The general machinery allows us to assume (enlarging $\lambda$) that we have $Y_{\lambda}/S_{\lambda}$ smooth projective (of course $X_{\lambda}= X \otimes_K K_{\lambda}$ is smooth projective), a relatively ample invertible sheaf $\eta_{\lambda}$ of `product-type' on $X_{\lambda} \times_{S_{\lambda}} Y_{\lambda}$, and various closed $S_{\lambda}$-subschemes 
\[
Z \into \mc{X}:= X_{\lambda} \times_{S_{\lambda}} Y_{\lambda}
\]
that are smooth over $S_{\lambda}$ and whose linear combinations define spread-out versions of the algebraic cycles on $X_L \times Y$. The purity theorem in this context\footnote{Th\'{e}or\`{e}me $3.7$ of Artin's Exp. XVI of \cite{sga4.3}} yields an isomorphism 
\[
H^{2j}_{Z}(\mc{X}, \Ql)(j) \xrightarrow{\sim} H^0(Z, \Ql),
\]
hence a cycle class $[Z]$ in $H^{2j}(\mc{X}, \Ql)(j)$. Regarding $\bar{x} \colon \Spec(L) \to S_{\lambda}$ as a geometric point over some scheme-theoretic point $x$, and letting $\bar{s} \colon \Spec(K) \to S_{\lambda}$ be a geometric point over a scheme-theoretic closed point $s$ lying in the closure of $x$, the cospecialization map $H^{2j}(\mc{X}_{\bar{s}}, \Ql) \xrightarrow{\sim} H^{2j}(\mc{X}_{\bar{x}}, \Ql)$ is an isomorphism ($\mc{X}/S_{\lambda}$ being smooth proper), and it carries the cycle class $[Z_{\bar{s}}]$ to the cycle class $[Z_{\bar{x}}]$, since these are both restrictions of $[Z]$, and the diagram
\[
\xymatrix{
& H^{2j}(\mc{X}_{\bar{s}}, \Ql) \ar[dd] \\
H^{2j}(\mc{X}, \Ql) \ar[ur] \ar[dr] & \\
& H^{2j}(\mc{X}_{\bar{x}}, \Ql)
}
\] 
commutes.
\endproof
\endproof
\subsection{Motives with coefficients}\label{motiveswithcoefficients}
We will need the flexibility of working with related categories of motives with coefficients. Let $E$ be a field of characteristic zero. Given any $E$-linear abelian (or in fact just additive, pseudo-abelian) category $\mc{M}$, and any finite extension $E'/E$, we can define the category $\mc{M}_{E'}$ of objects with coefficients in $E'$ as either of the following (this discussion is taken from \S $2.1$ of \cite{deligne:valeurs} and $3.11, 3.12$ of \cite{deligne-milne}):
\begin{enumerate}
\item The category of `$E'$-modules in $\mc{M}$,' whose objects are pairs $(M, \alpha)$ of an object $M$ of $\mc{M}$ and an embedding $\alpha \colon E' \to \End_{\mc{M}}(M)$, and whose morphisms are those commuting with these $E'$-structures.
\item The pseudo-abelian envelope of the category whose objects are formally obtained from those of $\mc{M}$ (writing $M_{E'}$ for the object in $\mc{M}_{E'}$ arising from $M$ in $\mc{M}$), and whose morphisms are 
\[
\Hom_{\mc{M}_{E'}}(M_{E'}, N_{E'})= \Hom_{\mc{M}}(M, N) \otimes_{E} E'.
\]
This construction is valid for infinite-dimensional $E'/E$.
\end{enumerate}
To pass from the first to the second description, let $(M, \alpha)$ be as in $(1)$, so that $\End_{\mc{M}_{E'}}(M_{E'})= \End_{\mc{M}}(M) \otimes_E E'$ contains, via $\alpha$, $E' \otimes_E E'$. This $E'$-algebra (via the left factor) is isomorphic to a product of fields, and there is a unique projection $e_{id} \colon E' \otimes_E E' \onto E'$ in which $x \otimes 1$ and $1 \otimes x$ both map to $x$. Then $e_{id}(M_{E'})$ is the object of $(2)$ corresponding to $(M, \alpha)$.

There is a functor $\mc{M} \to \mc{M}_{E'}$, which in the first language is $M \mapsto (M \otimes_{E} E', id_M \otimes id_{E'})$.\footnote{See $2.11$ of \cite{deligne-milne} for a precise description of $M \otimes_E E'$.} If $\mc{M}$ is semi-simple, then so is $\mc{M}_{E'}$. If $\mc{M}$ is a neutral Tannakian category over $E$ with fiber functor $\omega$, then we can make $\mc{M}_{E'}$ into a neutral Tannakian category over $E'$. Define an $E'$-valued fiber functor $\omega_{E'} \colon \mc{M} \to \mr{Vec}_{E'}$ by $\omega_{E'}(M)= \omega(M) \otimes_E E'$. There is a diagram commuting up to canonical isomorphism,
\[
\xymatrix{
\mc{M} \ar[r] \ar[rd]_{\omega_{E'}} & \mc{M}_{E'} \ar[d]^{\omega'_{E'}} \\
& \mr{Vec}_{E'}, 
}
\] 
where 
\[
\omega'_{E'}(M, \alpha):= \omega_{E'}(M) \underset{\alpha \otimes \id_{E'}, E' \otimes_{E} E'}{\otimes} E'.
\] 
This $\omega'_{E'}$ neutralizes $\mc{M}_{E'}$, and so we can define the associated Tannakian group $\mc{G}_{E'}= \Aut^{\otimes}(\omega'_{E'})$. We wish to compare $\mc{G}_{E'}$ with $\mc{G}= \Aut^{\otimes}(\omega)$. Note that $\mc{M}_{E'}$ is equivalent to $\Rep_{E'}(\mc{G})$ and $\mc{G} \otimes_{E} E'= \Aut^{\otimes}(\omega_{E'})$. Then the composition of functors
\[
\mc{M}_{E'} \xrightarrow[\sim]{\omega} \Rep_{E'}(\mc{G}) \xrightarrow[\sim]{F} \Rep_{E'}(\mc{G} \otimes_{E} E')
\]
given on $\mc{M}_{E'}$ by\footnote{$F$ just sends an object $V$ of $\Rep_{E'}(\mc{G})$, corresponding to an $E$-homomorphism $\mc{G} \to \Res_{E'/E}(\mr{GL}_{E'}(V))$, to $E' \otimes_{E' \otimes_E E'} V$.}
\[
(X, \alpha) \mapsto  (\omega(X), \omega(\alpha)) \mapsto  E' \underset{E' \otimes E', \omega(\alpha) \otimes 1}{\otimes} (\omega(X) \otimes E')
\]
is just $\omega'_{E'}$, whence a tensor-equivalence
\[
\Rep_{E'}(\mc{G}_{E'}) \xrightarrow{\sim} \Rep_{E'}(\mc{G} \otimes_{E} E').
\]

We apply these constructions to Andr\'{e}'s category of motives over $F$ for motivated cycles, to form the variant $\mc{M}_{F, E}$, motives over $F$ with coefficients in $E$, for $E$ any finite extension of $\Q$. We denote the corresponding motivic Galois group (for the $\sigma$-Betti realization) by $\mc{G}_{F, E}(\sigma)$; its $E$-linear (pro-algebraic) representations correspond to objects of $\mc{M}_{F, E}$, and it is naturally isomorphic to $\mc{G}_{F, \Q}(\sigma) \otimes_{\Q} E$. For any object $M$ of $\mc{M}_{F, E}$, we also have the corresponding motivic group $\mc{G}^M_{F, E}(\sigma)$.

As a more concrete variant, we can start with the ($\Q$-linear, semi-simple) isogeny category $AV^0_F$ of abelian varieties over $F$, and form the $E$-linear, semi-simple category $AV^0_{F, E}$ of isogeny abelian varieties over $F$ with complex multiplication by $E$. There is a (contravariant) functor $AV^0_{F, E} \to \mc{M}_{F, E}$. Faltings' theorem implies the following $E$-linear variant: fix an embedding $E \into \Qlb$, and consider an object $A$ of $AV^0_{F, E}$; then the natural map
\[
\End_{AV^0_{F, E}}(A) \otimes_E \Qlb \to \End_{\Qlb[\gal{F}]}(H^1(A_{\overline{F}}, \Ql) \otimes_E \Qlb)
\]
is an isomorphism.\footnote{Simply take the usual isomorphism with $\Q$ in place of $E$, restrict to those endomorphism commuting with $E \to \End_{AV^0_{F}}(A)$, and then project to the $E \into \Qlb$ component of the resulting $E \otimes_{\Q} \Qlb$-module.}

\subsection{Hodge symmetry in $\mc{M}_{F, E}$}\label{hodgesymmetry}
Our next goal is to show that objects of $\mc{M}_{F, E}$ satisfy (unconditionally) the Hodge-Tate weight symmetries, and the more fundamental symmetries in the rational de Rham realization, needed for our abstract Galois lifting results. I expect these results are well-known to experts, but they do not seem to have been explained in their natural degree of generality, and in a context in which they can be \textit{proven} unconditionally\footnote{Our discussion is an analogue of the standard discussion in the theory of motives of CM type-- see \S $7$ of \cite{serre:motivicgalois}. The Galois symmetries themselves have also been postulated in \S $5$ of \cite{blggt:potaut}, but since that paper is only concerned with Galois representations arising from essentially conjugate self-dual automorphic representations, it has no need to deal with this symmetry as a general principle.}.  So, let $M$ be an object of $\mc{M}_{F, E}$, of rank $r$. The de Rham realization $M_{dR}$ is a filtered $F \otimes_{\Q} E$-module, which is moreover free of rank $r$. As in \S \ref{formal}, we define the $\tau \colon F \into \overline{E}$-labeled weights of $M$ as follows: $\mr{HT}_{\tau}(M)$ is an $r$-tuple of integers $h$, with $h$ appearing with multiplicity
\[
\dim_{\overline{E}} \gr^h\left(M_{dR} \otimes_{F \otimes_{\Q} E, \tau \otimes 1} \overline{E} \right).
\] 
In talking of $\tau$-labeled weights, there is alway an ambient over-field, in this case $\overline{E}$, but we will at times want to change this to either $\CC$ or $\Qlb$; that will require \textit{fixing} an embedding $\iota \colon \overline{E} \into \Qlb$ or $\iota \colon \overline{E} \into \CC$. In either case, we can then speak of $\mr{HT}_{\iota \tau}(M)$, with $E$ embedded into $\CC$ or $\Qlb$ via $\iota$, and there is an obvious equality $\mr{HT}_{\tau}(M)= \mr{HT}_{\iota \tau}(M)$.

The essential point, which is the motivic analogue of Corollary \ref{automorphiccmcoefficients}, is  closely related to the fact that $\mc{G}_{F}$ splits\footnote{i.e., for any algebraic quotient $\mc{G}_F \onto H$ over $\Q$, the connected component $H^0$ has maximal torus that splits} over $\Q^{cm}$ (see $4.6(iii)$ of \cite{andre:motivated} for this assertion). For lack of reference for the proof, we give some details:
\begin{lemma}\label{betticmdescent}
Let $N$ be an object of $\mc{M}_{F, E}$. Then there exists an object $N_0$ of $\mc{M}_{F, E_{cm}}$ such that $N \cong N_0 \otimes_{E_{cm}} E$. Consequently, $\mr{HT}_{\tau}(N)$ depends only on the restriction of $\tau$ to $F_{cm}$.
\end{lemma}
\proof
By the formalism of Lemma \ref{hodgetheory}, the second claim follows from the first. By Proposition $3.3$ (the analogue of the Hodge index theorem) of \cite{andre:motivated}, $C^0_{mot}(X, X)$ is endowed with a positive-definite, $\Q$-valued symmetric form, which we call $\langle \cdot , \cdot \rangle$. For any sub-object $M \subset H(X)$ in $\mc{M}_F$, there follows a decomposition $H(X)= M \oplus M^{\perp}$, by positivity and the fact that $\mc{M}$ is a semi-simple abelian category. $\langle \cdot , \cdot \rangle$ therefore restricts to a positive form on $M$ itself, and in particular every simple object of $\mc{M}_F$ carries a positive definite form. Again by semi-simplicity, any such simple object $M$ has $\End(M)$ isomorphic to a division algebra $D$, on which we now have an involution $'$ (transpose with respect to $\langle \cdot , \cdot \rangle$) and a trace form $\tr_M \colon D \to \Q$ (given by $\phi \mapsto \tr(\phi|_{H_B(M)}$) such that $\tr_M(\phi \phi')>0$ for all non-zero $\phi \in D$. As in the proof of the Albert classification (see \cite[\S 21, Theorem 2]{mumford:abvars}), this implies that the center of $D$ is either a totally real or a CM field. Thus, $\End(M)$ splits over a CM field; in particular, for any number field $E$ and any factor $N$ of $M \otimes_{\Q} E$ (in $\mc{M}_{F, E}$), $N$ can in fact be realized as (the scalar extension of) an object of $\mc{M}_{F, E_{cm}}$.
\endproof
\begin{rmk}\label{hodgeQcoefficients}
When $E= \Q$, Lemma \ref{hodgetheory} shows that $\mr{HT}_{\tau}(N)$ is independent of $\tau$; this is essentially the assertion that for a smooth projective $X/F$, the Hodge numbers of $X$ do not depend on the choice of embedding $F \into \CC$. The next few results (culminating in Corollary \ref{motivichodgesymmetry}) all have corresponding strengthenings when $E= \Q$.
\end{rmk}

We first record the de Rham-Betti version of the desired symmetry:\index{t}{Hodge symmetry}
\begin{lemma}\label{hodgesym}
Let $N$ be an object of $\mc{M}_{F, E}$ lying in the $k$-component of the grading. Then for any choice of complex conjugation $c$ in $\Gal(\overline{E}/\Q)$, 
\[
\mr{HT}_{c \circ \tau}(N)= \{ k- h: h \in \mr{HT}_{\tau}(N) \}.
\]
\end{lemma}
\proof
Any such object is built by taking one of the form $M= H^k(X) \otimes E$, for $X/F$ a smooth projective variety and applying an idempotent $\alpha \in C^0_{mot}(X, X)_E$.\footnote{And a Tate twist, but the statement of the lemma is obviously invariant under Tate twists.} Fix an embedding $\iota \colon \overline{E} \into \CC$, so that $\iota \circ \tau \colon F \into \CC$. There is a functorial (Betti-de Rham) comparison isomorphism
\[
M_{dR} \otimes_{F, \iota \tau} \CC \cong M_{B, \iota \tau} \otimes_{\Q} \CC,
\]
which commutes with the action of $C^0_{mot}(X, X)_{E}$ on $M_{dR}$ and $M_{B, \iota \tau}$, in particular making this an isomorphism of free $E \otimes_{\Q} \CC$-modules. It induces an isomorphism on the corresponding gradeds (with $q=k-p$):
\[
\gr^p(M_{dR} \otimes_{F, \iota \tau} \CC) \cong H^{p, q}_{\iota \tau}(M) \cong H^{p, q}(X_{\iota \tau}) \otimes_{\Q} E,
\]
again $E \otimes_{\Q} \CC$-linear. Just for orientation amidst the formalism, this says that
\[
\mr{HT}_{\tau}(M_{dR})= \mr{HT}_{\iota \tau}(M_{dR})= \left\{p: H^{p,q}(X_{\iota \tau}) \neq 0, \text{counted with multiplicity $\dim_{\CC} H^{p, q}(X_{\iota \tau})$}\right\}.
\]
Now, the projection $X_{\iota \tau} \times_{\CC, c} \CC \to X_{\iota \tau}$ induces a transfer of structure isomorphism on rational Betti cohomology, and then an $E \otimes_{\Q} \CC$-linear isomorphism
\[
F_{\infty} \colon H^{p, q}_{\iota \tau}(M) \to H^{q, p}_{c \iota \tau}(M).
\]
Let $\alpha \in C^0_{mot}(X, X)_{E}$ be an $E$-linear motivated idempotent correspondence defining an object $N= \alpha M$ of $\mc{M}_{F, E}$. There is a commutative diagram of $E \otimes_{\Q} \CC$-linear morphisms
\[
\xymatrix{
H^{p,q}_{\iota \tau}(M) \ar[r]^{F_{\infty}} \ar[d]_{\alpha_{\iota \tau}} & H^{q,p}_{c \iota \tau}(M) \ar[d]^{\alpha_{c \iota \tau}}\\
H^{p, q}_{\iota \tau}(M) \ar[r]^{F_{\infty}} & H^{q, p}_{c \iota \tau}(M).
}
\] 
Consider the images of the two vertical maps. Since the horizontal maps are isomorphisms, and all the maps are $E \otimes_{\Q} \CC$-linear, we deduce an isomorphism
\[
\alpha_{\iota \tau} H^{p, q}_{\iota \tau}(M) \otimes_{E \otimes_{\Q} \CC, \iota} \CC \xrightarrow{\sim} \alpha_{c \iota \tau} H^{q,p}_{c \iota \tau}(M) \otimes_{E \otimes_{\Q} \CC, \iota} \CC.
\]
The left-hand side is isomorphic to
\[
\gr^p \left( (\alpha M)_{dR} \otimes_{ F \otimes_{\Q} E, \iota \tau \otimes \iota} \CC \right),
\]
and letting $c'$ be a complex conjugation in $\Gal(\overline{E}/\Q)$ such that $\iota c' \tau= c \iota \tau$, the right-hand side is isomorphic to
\[
\gr^{k-p} \left( (\alpha M)_{dR} \otimes_{F \otimes_{\Q} E, \iota c' \tau \otimes \iota} \CC \right).
\] 
It follows that
\[
\mr{HT}_{c' \tau}(\alpha M)= \left\{ k- h: h \in \mr{HT}_{\tau}(\alpha M) \right\}.
\]
But by the previous lemma, $\mr{HT}_{c' \tau}(\alpha M)= \mr{HT}_{c \tau}(M)$, so we are done.
\endproof
Next we observe (as in \S $2.4$ of \cite{andre:motivated}, but with a few more details) that `$p$-adic' comparison isomorphisms hold unconditionally in $\mc{M}_{F, E}$; this is one of the main reasons for working with motivated cycles rather than absolute Hodge cycles.
\begin{lemma}
Let $N$ be an object of $\mc{M}_{F, E}$. Then for all $v \vert \ell$, the free $E \otimes_{\Q} \Ql$-module $M_{\ell}$ is a de Rham representation of $\gal{F_v}$, and there is a functorial (with respect to morphisms in $\mc{M}_{F, E}$) isomorphism of filtered $F_v \otimes_{\Q} E \cong F_v \otimes_{\Ql} (\Ql \otimes_{\Q} E)$-modules
\[
F_v \otimes_F M_{dR} \cong \Ddr(M_{\ell}|_{\gal{F_v}}).
\]
\end{lemma}
\proof
It suffices to check for objects of the form $M= H(X) \otimes_{\Q} E$, for $X/F$ smooth projective. Faltings' theorem\footnote{Of which there are now several proofs, including the necessary corrections to the original \cite{faltings:crystalline}. A `simple' recent proof is \cite{beilinson:derivedderham}.} provides the comparison
\[
H^*_{dR}(X_{F_v}) \otimes_{\Q} E \cong \Ddr \left(H^*(X_{\overline{F_v}}, \Ql)\right) \otimes_{\Q} E ,
\]
which is an isomorphism of filtered $F_v \otimes_{\Q} E$-modules. Let $\alpha \in C^0_{mot}(X_1, X_2)_E$ be a motivated correspondence. Let $d_1= \dim(X_1)$. $\alpha$ has realizations $\alpha_{dR}$ and $\alpha_{\ell}$ in cohomology groups that are also compared by Faltings:\footnote{Note that $\alpha_{\ell}$ is fixed by $\gal{F_v}$}
\[
H^{2d_1}_{dR}((X_1 \times X_2)_{F_v})(d_1)) \otimes_{\Q} E \cong \Ddr \left( H^{2d_1}((X_1 \times X_2)_{\overline{F_v}}, \Ql)(d_1) \right) \otimes_{\Q} E,
\]
and the claim is that $\alpha_{dR}$ maps to $\alpha_{\ell}$. The comparison isomorphism is compatible with cycle class maps,\footnote{See part $(b)$ of the proof of Theorem $3.6$ of \cite{beilinson:derivedderham}.} and our cycle $\alpha$ is spanned by elements of the form $\pr^{X_1 \times X_2\times Y}_{X_1 \times X_2, *}(\beta \cup \ast \gamma)$, with $Y/F$ another smooth projective variety and $\beta$ and $\gamma$ algebraic cycles on $X_1 \times X_2 \times Y$. Applying the comparison isomorphism also to $X_1 \times X_2 \times Y$, and applying compatibility with cycle classes, cup-product, (hence) Lefschetz involution, and the projection $\pr^{X_1 \times X_2 \times Y}_{X^2, *}$\footnote{This last point because compatible with pull-back induced by morphisms of varieties and, again, Poincare duality.}, we can deduce that $\alpha_{dR}$ maps to $\alpha_{\ell}$.
\endproof
\begin{cor}\label{motivichodgesymmetry}
Let $M$ be an object of $\mc{M}_{F, E}$. Fix an embedding $\iota \colon \overline{E} \into \Qlb$, and use this to identify $E$ as a subfield of $\Qlb$. Then for all $\tau \colon F \into \Qlb$, $\mr{HT}_{\tau}(M_{\ell} \otimes_E \Qlb)$ depends only on $\tau_0= \tau|_{F_{cm}}$. If $M$ is moreover pure of some weight $k$, then
\[
\mr{HT}_{\tau_0 \circ c}(M_{\ell}\otimes_E \Qlb)= \left\{ k- h: h \in \mr{HT}_{\tau_0}(M_{\ell} \otimes_E \Qlb) \right\}.
\] 
\end{cor}
\proof
Write $\iota^{-1} \tau$ for the embedding $F \into \overline{E}$ induced by $\tau$ and $\iota$. By the previous lemma, $\mr{HT}_{\tau}(M_{\ell} \otimes_{E} \Qlb)= \mr{HT}_{\iota^{-1} \tau}(M)$. The latter depends only on $\tau|_{F_{cm}}$, and when $M$ is pure satisfies the required symmetry, by Lemmas \ref{betticmdescent} and \ref{hodgesym}.
\endproof
Finally, we can show that the Galois lifting results of \S \ref{generalGalois} apply to representations arising from objects of $\mc{M}_{F, E}$. In an idealized situation in which all Hodge cycles are motivated (or motivic groups over $\overline{F}$ are connected), this will always be the case; in general there is a complication, which I have not yet been able to avoid, arising from the fact that the proof of Theorem \ref{fullmonodromylift} requires an initial reduction to the case of connected monodromy group. Fix an embedding $E \into \Qlb$, inducing a place $\lambda$ of $E$. For short-hand, we denote by $\mc{G}_{F, \Qlb}$ the base-change to $\Qlb$ of the $E_{\lambda}$-group $\mc{G}_{F, \lambda}$ defined by the $\lambda$-adic \'{e}tale fiber functor on the category $\mc{M}_{F, E}$, and we make the analogous definition of $\mc{G}^M_{F, \Qlb}$ for objects $M$ of $\mc{M}_{F, E}$.
\begin{cor}\label{galoisliftmotiviccase}
Let $F$ be a totally imaginary field. Let $M$ be an object of $\mc{M}_{F, E}$ and $E \into \Qlb$ a fixed embedding, to which we associate the $\ell$-adic representation $M_{\ell} \otimes_E \Qlb$, which is a representation both of $\gal{F}$ and of the motivic group $\mc{G}^M_{F, \Qlb}$. Make the following assumption
\begin{itemize}
\item $\mc{G}^{M_{\overline{F}}}_{\overline{F}, \Qlb}$ is connected.
\end{itemize}
Suppose that $H' \to H$ is any central torus quotient of linear algebraic groups over $\Qlb$, and that there is a factorization 
\[
\mc{G}^M_{F, \Qlb} \into H \into \mr{GL}_{\Qlb}(M_{\ell} \otimes_E \Qlb),
\]
where of course the composition of these inclusions is the natural representation. Then there exists a geometric lift $\tilde{\rho}$
\[
\xymatrix{
& & H'(\Qlb) \ar[d] \\
\gal{F} \ar@{-->}[rru]^{\tilde{\rho}} \ar[r]_-\rho & \mc{G}^M_{F, \Qlb}(\Qlb) \ar[r] & H(\Qlb).
}
\]
Similarly, if $F$ is any number field, but $E$ is $\Q$, then without any assumption on the motivic group $\mc{G}^{M_{\overline{F}}}_{\overline{F}}$, such a $\tilde{\rho}$ exists.
\end{cor}
\begin{rmk}
The assumption that $\mc{G}^{M_{\overline{F}}}_{\overline{F}, \Qlb}$ is connected is expected always to hold. In particular, this should not be necessary for the conclusion of the corollary to hold.
\end{rmk}
\proof
Corollary \ref{motivichodgesymmetry} tells us that composing $\rho$ with any irreducible algebraic representation of $\mc{G}^M_{F, \Qlb}$ yields an $\ell$-adic representation satisfying the Hodge symmetries of Conjecture \ref{galoisdescent}. But the arguments of Proposition \ref{fullmonodromylift} and Theorem \ref{anymonodromylift} do not apply directly, because they assume these symmetries after composition with irreducible representations of the connected component $(\overline{\rho(\gal{F})}^{Zar})^0$.\footnote{If this group were, contrary to conjecture, not equal to $\mc{G}^{M_{\overline{F}}}_{\overline{F}, \Qlb}$, then it is not obvious in general how to relate the two putative Hodge symmetries.} The simplest (but imperfect) way around this is to assume $\mc{G}^{M_{\overline{F}}}_{\overline{F}, \Qlb}$ is connected, and then run through the arguments of Proposition \ref{fullmonodromylift} and Theorem \ref{anymonodromylift} with $\mc{G}^{M_{F'}}_{F', \Qlb}$, for $F'$ sufficiently large (see Lemma \ref{connected}), in place of $\overline{\rho(\gal{F'})}^{Zar}$: reduce to the connected case (replacing $F$ by $F'$) as in Theorem \ref{anymonodromylift}, and then note that the proof of Proposition \ref{fullmonodromylift} only requires identifying a connected reductive group containing the image of $\rho$ such that composition with irreducible representations of this group yields Galois representations with the desired symmetry (whether or not these Galois representations are themselves irreducible). 

The second assertion (when $E=\Q$) follows by the same argument: by Remark \ref{hodgeQcoefficients}, the $\tau$-labeled Hodge-Tate co-character $\mu_{\tau}$ of $\rho$ is independent of $\tau$. It can be interpreted as a co-character of $(\overline{\rho(\gal{F})}^{Zar})^0$, and composing with finite-dimensional representations of this group, we obtain the necessary Hodge-Tate symmetries to apply the argument of Proposition \ref{fullmonodromylift}.
\endproof
\begin{rmk} 
Similarly, there is a motivic version of Corollary \ref{fullmonodromytotreal}.
\end{rmk}

\subsection{Motivic lifting: the potentially CM case}\label{taniyama}
One case of the desired motivic lifting result is immediately accessible, when the corresponding Galois representations are potentially CM. We denote by $\mc{CM}_F$ the Tannakian category of motives (for motivated cycles) over $F$ generated by Artin motives and potentially CM abelian varieties. As before an embedding $\sigma \colon F \into \CC$ yields, through the corresponding Betti fiber functor, a Tannakian group $\mc{T}_F(\sigma)= \Aut^{\otimes}(H_{\sigma}|_{\mc{CM}_F})$.  For a number field $E$, we can also consider, as in \S \ref{motiveswithcoefficients}, the category $\mc{CM}_{F, E}$ of potentially CM motives over $F$ with coefficients in $E$, with its Tannakian group\index{s}{$\mc{T}_{F, E}$}\index{t}{Taniyama group} $\mc{T}_{F, E}(\sigma)$.\footnote{From now on, $\sigma$ will be implicit.} As with $\mc{G}_F$, there is a sequence
\[
1 \to \mc{T}_{\overline{F}, E} \to \mc{T}_{F, E} \to \gal{F} \to 1,
\]
where again the projection $\mc{T}_{F, E} \to \gal{F}$ has a continuous section $s_{\lambda}$ on $E_{\lambda}$-points. Deligne showed in \cite{DMOS} (Chapter IV: `\textit{Motifs et groupes de Taniyama}') that $\mc{T}_{\Q}$ is isomorphic to the Taniyama group constructed by Langlands as an explicit extension of $\gal{\Q}$ by the connected Serre group.\footnote{More precisely, equipped with the data of the projection to $\gal{\Q}$, a section on $\mbf{A}_{F, f}$-points, and the co-character $\mathbb{G}_{m, \CC} \to \mc{T}_{\CC}$ giving the Hodge filtration on objects of $\mc{CM}_{\CC}$, $\mc{T}_{\Q}$ is uniquely isomorphic to the Taniyama group, equipped with its corresponding structures.}

\begin{prop}\label{taniyamalift}
Let $F$ be totally imaginary. Let $H' \to H$ be a surjection of linear algebraic groups over a number field $E$ with central torus kernel. For any homomorphism $\rho \colon \mc{T}_{F, E} \to H$, there is a finite extension $E'/E$ such that $\rho$ lifts to a homomorphism
\[
\xymatrix{
& H'_{E'} \ar[d] \\
\mc{T}_{F, E'} \ar@{-->}[ur]^{\tilde{\rho}} \ar[r]^{\rho} & H_{E'}.
}
\]
\end{prop}
\proof
Fix a finite place $\lambda$ of $E$ and consider the $\lambda$-adic representation $\rho_{\lambda}=\rho \circ s_{\lambda} \colon \gal{F} \to H(E_{\lambda})$. Since $\mc{T}_{\overline{F}}$ is isomorphic to the connected Serre group, we can unconditionally apply Corollary \ref{galoisliftmotiviccase} to find, for some finite extension $E'_{\lambda}/E_{\lambda}$ a geometric lift $\tilde{\rho}_{\lambda} \colon \gal{F} \to H'(E'_{\lambda})$ of $\rho_{\lambda}$. Note that $\tilde{\rho}_{\lambda}$ is potentially abelian, since $\rho_{\lambda}$ is, and the kernel of $H' \to H$ is central. Proposition $\text{IV.D}.1$ of \cite{DMOS} implies that $\tilde{\rho}_{\lambda}$ arises from a homomorphism (of groups over $E'_{\lambda}$) $\mc{T}_{F} \otimes_{\Q} E'_{\lambda} \cong \mc{T}_{F, E} \otimes_E E'_{\lambda} \to H' \otimes_E E'_{\lambda}$. This in turn must be definable over some finite extension $E'/E$ ($E'$ is thus embedded in $E'_{\lambda}$), i.e. there is a homomorphism $\tilde{\rho} \colon \mc{T}_{F, E'} \to H' \otimes_E E'$ whose extension to $E'_{\lambda}$ gives rise to $\tilde{\rho}_{\lambda}$. We explain this technical point in Lemma \ref{spreadout} below. To check that $\tilde{\rho}$ is actually a lift of $\rho \otimes_E E'$, it suffices to observe:
\begin{itemize}
\item The restriction to $\mc{T}_{\overline{F}, E'}$ is a lift: the restrictions of $\tilde{\rho}$ and $\rho$ are simply the algebraic homomorphisms of the connected Serre group (tensored with $E'$) corresponding to the labeled Hodge-Tate weights of $\tilde{\rho}_{\lambda}$ and $\rho_{\lambda}$.
\item The $\lambda'$-adic lift, for the place $\lambda'$ of $E'$ induced by $E' \into E'_{\lambda}$, is a lift (by construction).
\item It suffices to check that $\tilde{\rho}$ is a lift on $E'_{\lambda'}$-points. First, it suffices to check that $\tilde{\rho} \otimes_{E'} E'_{\lambda'}$ lifts $\rho \otimes_E E'_{\lambda'}$. This in turn can be checked on $E'_{\lambda'}$-points: we are free to replace $\mc{T}_{F, E'} \otimes_{E'} E'_{\lambda'}$ by some finite-type quotient $\mc{T}$ in which the $E'_{\lambda'}$-points are Zariski-dense. Then the closed subscheme $\mc{T} \underset{H \times H}{\times} H \into \mc{T}$ (where the two maps to $H \times H$ are the diagonal $H \to H \times H$ and the product of the two maps $\mc{T} \to H \times H$ given by $\rho$ and $\tilde{\rho}$ followed by the quotient $H' \to H$) contains $\mc{T}(E'_{\lambda'})$, hence equals $\mc{T}$.
\end{itemize}
Then we are done, since $\mc{T}_{F, E'}(E'_{\lambda'})= \mc{T}_{\overline{F}, E'}(E'_{\lambda'}) \cdot s_{\lambda'}(\gal{F})$.
\endproof
Here is the promised lemma showing that $\tilde{\rho}$ may be defined over a finite extension $E'$ of $E$:
\begin{lemma}\label{spreadout}
Let $K/k$ be an extension of algebraically closed fields of characteristic zero. Let $T$ be a (not necessarily connected) reductive group over $k$, and let $\pi \colon H' \to H$ be a surjection (defined over $k$) of reductive $k$-groups. Suppose that $\rho \colon T \to H$ is a $k$-morphism, and that the scalar extension $\rho_K$  lifts to a $K$-morphism $\tilde{\rho} \colon T_K \to H'_K$. Then $\rho$ lifts to a $k$-morphism $T \to H'$.
\end{lemma}
\proof
This is a simple spreading-out argument. Namely, writing $K$ as the direct limit of its finite-type $k$-sub-algebras, we can descend the relation $\pi_K \circ \tilde{\rho}= \rho_K$ to some finitely-generated $k$-sub-algebra $R$ of $K$, in particular obtaining $\tilde{\rho}_R \colon T_R \to H'_R$ such that $\tilde{\rho}_R \otimes_R K= \tilde{\rho}$. We then choose a non-zero $k$-point $\alpha \colon R \to k$ (a morphism of $k$-algebras), and specializing the relation $\pi_R \circ \tilde{\rho}_R= \rho_R$ via $\alpha$, we obtain a lift $\tilde{\rho}_R \otimes_{R, \alpha} k$ of $\rho$ over $k$. 
%By a basic result of Vinberg (Proposition $10$ of \cite{spreadout:invariants}), any homomorphism $T_K \to H'_K$ is $H'(K)$-conjugate to a homomorphism defined over $k$: so, there exists $\tilde{\rho}_0 \colon T \to H'$ and $\tilde{h} \in H'(K)$ such that $\tilde{h} \tilde{\rho} \tilde{h}^{-1}= \tilde{\rho}_{0, K}$. Composing with $\pi_K$, and writing $h= \pi_K(\tilde{h})$, $\rho_0= \pi_K \circ \tilde{\rho}_{0, K}$, we have $h \rho_K h^{-1}= \rho_{0, K}$. We claim that there is an $h_0 \in H(k)$ such that $h_0 \rho h_0^{-1}= \rho_0$. Granted this, we can lift $h_0$ to an element $\tilde{h}_0$ of $H'(k)$, and then $\tilde{h}_0^{-1} \tilde{\rho}_0 \tilde{h}_0$ is a lift of $\rho$ defined over $k$.

%To see the claim, note that for all $\sigma \in \Gal(K/k)$, ${}^{\sigma}h \rho_K {}^{\sigma}h^{-1}= \rho_{0, K}$, and thus $c \colon \sigma \mapsto h^{-1} \cdot {}^{\sigma}h$ defines a cohomology class in $H^1(\Gal(K/k), \mr{Cent}_{H(K)}(\rho_K))$, which is of course trivial, since $K/k$ is an extension of separably closed fields. Consequently, for some element $z$ of this centralizer, $c(\sigma)= z^{-1} \cdot {}^{\sigma} z$, and the element $hz^{-1}$ is then $\Gal(K/k)$-invariant. Hence $h_0 = hz^{-1}$ lies in $H(k)$ and conjugates $\rho$ to $\rho_0$.
\endproof
\begin{rmk}
\begin{itemize}
\item We have stated Proposition \ref{taniyamalift} only for imaginary fields for simplicity, but of course there is a variant for general number fields taking into account Corollary \ref{fullmonodromytotreal}. When $F$ is totally real, there will be such $\rho$ that do not lift: simply take a type $A$ Hecke character $\psi$ of a quadratic CM extension $L/F$ as in Example \ref{mixed}. Then $\Ad^0(\Ind_L^F(\psi))$ will have a potentially abelian Galois $\ell$-adic realizations $\gal{F} \to \mr{SO}_3(\Qlb)$, which arises from a representation of $\mc{T}_F$ by Proposition $\text{IV.D.}1$ of \cite{DMOS}. As we have seen, these representations do not lift geometrically to $\mr{GSpin}_3$. 
\item The method of checking that $\tilde{\rho}$-- once we know it exists-- actually lifts $\rho$ will recur in \S \ref{hyperlift}; this division into a geometric argument (lifting $\rho|_{\mc{T}_{\overline{F}, E}}$ and a Galois-theoretic argument (lifting $\rho \circ s_{\lambda}$) seems to be the natural way to make arguments about motivic Galois groups over number fields (compare Lemma \ref{connected}).
\item On the automorphic side, representations of the global Weil group ought to parametrize `potentially abelian' automorphic representations. In that case, the analogous lifting result is a theorem of Labesse (\cite{labesse:lifting}).
\end{itemize}
\end{rmk}
Before proceeding to more elaborate examples, let us clarify here that, even though we study lifting through \textit{central} quotients, the nature of these problems is in fact highly non-abelian. That is, suppose we have a lift $\tilde{\rho} \colon \mc{G}_{F, E} \to H'_{E}$ of $\rho \colon \mc{G}_{F, E} \to H_E$ (we may assume $H'$ and $H$ are reductive; let us further suppose they are connected), and let $r'$ and $r$ be irreducible faithful representations of $H'_E$ and $H_E$. If the derived group of $H$ is not simply-connected, but its simply-connected cover $H_{sc}$ injects into $H'$, then typically $r' \circ \tilde{\rho}$ will not lie in the Tannakian subcategory of $\mc{M}_{F, E}$ generated by $r \circ \rho$ and all potentially CM motives, i.e. objects of $\mc{CM}_{F, E}$ (of course, it does in the example of Proposition \ref{taniyamalift}). We make this precise at the Galois-theoretic level:
\begin{lemma}\label{notCM}
Let $H' \to H$, $r'$, and $r$ be as above. Let $\rho \colon \gal{F} \to H(\Qlb)$ be a geometric Galois representation having a geometric lift $\tilde{\rho} \colon \gal{F} \to H'(\Qlb)$. Assume that \begin{itemize}
\item the algebraic monodromy group of $\rho$ is $H$ itself; 
\item the kernel of $H_{sc} \cap \overline{(\tilde{\rho}(\gal{F}))}^{Zar} \to H$ is non-zero.
\end{itemize}
Then $r' \circ \tilde{\rho}$ is not contained in the Tannakian sub-category of semi-simple geometric $\Qlb$-representations of $\gal{F}$ generated by $r \circ \rho$ and all potentially abelian geometric representations.
\end{lemma}
\proof
If $r' \circ \tilde{\rho}$ were contained in this category, then there would exist an irreducible potentially abelian representation $\tau$ of $\gal{F}$ and an injection $r' \circ \tilde{\rho} \into \tau \otimes (r_1 \circ \rho)$ for some irreducible algebraic representation $r_1$ of $H$. By Proposition \ref{lietensorart}, $\tau$ has the form $\Ind_L^F (\psi \cdot \omega)$ for some character $\psi$ and irreducible Artin representation $\omega$ of $\gal{L}$. By Frobenius reciprocity, there is a non-zero map $(r' \circ \tilde{\rho})|_{\gal{L}} \to \psi \cdot (r_1 \circ \rho)|_{\gal{L}} \otimes \omega$; both sides are irreducible,\footnote{For the right-hand side, see Proposition \ref{lietensorart}: the tensor product of a Lie-irreducible and an Artin representation is irreducible; of course, for the purposes of this lemma, we could just further restrict $L$.} so $\omega$ is one-dimensional, and absorbing $\omega$ into $\psi$ we may assume $r' \circ \tilde{\rho}|_{\gal{L}} \xrightarrow{\sim} \psi \cdot (r_1 \circ \rho)|_{\gal{L}}$. Comparing algebraic monodromy groups, we obtain a contradiction, since by assumption $H_{sc} \cap \overline{(\tilde{\rho}(\gal{F}))}^{Zar}$ cannot inject into $\mathbb{G}_m \times H$.
\endproof
In \S \ref{hyperlift}, we will see many examples of $\rho$ with full $\mr{SO}$ monodromy group; their lifts to $\mr{GSpin}$ will satisfy the conclusions of Lemma \ref{notCM}.
\section{Motivic lifting: the hyperk\"{a}hler case}\label{hyperlift} 
\subsection{Setup}\label{HKsetup}
The aim of this section is to produce a lifting not merely at the level of a single $\ell$-adic representation, but of actual motives, in a very special family of cases whose prototype is the primitive second cohomology of a $K3$ surface over $F$. So that the reader has some examples to keep in mind, we recall the following definitions:
\begin{defn}\label{HK}
Let $F$ be a subfield of $\CC$. A \textit{hyperk\"{a}hler variety} $X$ over $F$ is a geometrically connected and simply-connected smooth projective variety over $F$ of even dimension $2r$ such that $\Gamma(X, \Omega^2_X)$ is one-dimensional, generated by a differential form $\omega$ for which $\omega^r$ is non-vanishing at every point of $X$ (i.e., as a linear functional on the top wedge power of the tangent space).\index{t}{hyperk\"{a}hler variety} A \textit{K3 surface} over $F$ is a hyperk\"{a}hler variety $X/F$ of dimension 2.\index{t}{K3 surface} 
\end{defn}
\begin{eg}
The simplest example of a K3 surface is a smooth quartic hypersurface in $\mathbb{P}^3$. Higher-dimensional examples of hyperk\"{a}hlers are notoriously difficult to produce. One standard family is gotten by starting with any K3 surface $X$, and then for any integer $r \geq 1$ considering the Hilbert scheme $X^{[r]}$ of $r$ points on $X$ (more precisely, the moduli of closed sub-schemes of length $r$); each $X^{[r]}$ is a hyperk\"{a}hler variety, with of course $X^{[1]}=X$. (This construction is due to Beauville: see \cite[\S 6]{beauville:HK}.)
\end{eg}

More generally, we work in the axiomatized setup of \cite{andre:hyperkaehler}. By a \textit{polarized variety} over a subfield $F$ of $\CC$ we mean a pair $(X, \eta)$ consisting of variety $X/F$ and an $F$-rational ample line bundle $\eta$ on $X$.\index{t}{polarized variety} For a fixed $k \leq \dim X$, we will consider the `motive' $Prim^{2k}(X)(k)$ (omitting the $\eta$-dependence from the notation), as an object of $\mc{M}_F$. Cup-product with $\eta$ lets us endow the motive $H^{2k}(X)(k)$ with the quadratic form 
\[
\langle x, y \rangle_{\eta}= (-1)^{k} x \cup y \cup \eta^{\dim X- 2k} \in H^{2\dim X}(X)(\dim X) \cong \Q,
\]
and then we can define $Prim^{2k}(X)(k)$ as the orthogonal complement of $H^{2k-2}(X)(k-1) \cup \eta$.\index{t}{primitive cohomology, as a motive}\index{s}{$Prim^{2k}(X)(k)$}

The two Weil cohomologies we use are $\ell$-adic (with coefficients in $\Ql$ or some extension inside $\Qlb$) and Betti cohomology; we will occasionally take integral Betti cohomology, where we will always work modulo torsion, so that $Prim^{2k}(X_{\CC}, \Z)(k)$ is a sub-module of $H^{2k}(X_{\CC}, \Z)(k)/{torsion}$. In that case, the primitive lattice defines a polarized (by $\langle \cdot , \cdot \rangle_{\eta}$) integral Hodge structure of weight $0$; we write $h^{p,q}$ for the Hodge numbers. Andr\'{e} proves his theorems, which include versions of the Shafarevich and Tate conjectures, under the following axioms:\index{s}{$A_k$}\index{s}{$B_k$}\index{s}{$B_k^+$}
\begin{itemize}
\item[$A_k$:] $h^{1,-1}=1$, $h^{0,0}>0$, and $h^{p,q}=0$ if $|p-q|>2$.
\item[$B_k$:] There exists a smooth connected $F$-scheme $S$, a point $s \in S(F)$, and a smooth projective morphism $f \colon \mc{X} \to S$ such that:
\begin{itemize}
\item $X \cong \mc{X}_s$;
\item the Betti class $\eta_{B} \in H^2(X_{\CC}^{an}, \Z)(1)/{torsion}$ extends to a section of $R^2 f_{\CC *}^{an} \Z(1)/{torsion}$;
\item letting $\tilde{S}$ denote the universal cover of $S(\CC)$, and $D$ denote the period domain of Hodge structures on $V_{\Z}:= Prim^{2k}(X_{\CC}, \Z)(k)$ polarized by $\langle \cdot , \cdot \rangle_{\eta}$, we require that the image of $\tilde{S} \to D$ contain an open subset.
\end{itemize}
\item[$B_k^{+}$:] For each $t \in S(\CC)$, every Hodge class in $H^{2k}(\mc{X}_t, \Q)(k)$ is an algebraic class.
\end{itemize}
$A_k$ is essential to the method, which relies on studying the associated Kuga-Satake abelian variety, which exists for Hodge structures of this particular form. $B_k$ is a statement about deforming $X$ into a `big' family. $B_k^+$ is of course a case of the Hodge conjecture, which is always known when $k=1$ (the theorem of Lefschetz). But we provide these axioms merely for orientation. Of interest is the following collection of varieties for which they are known to hold:
\begin{prop}
The axioms $A_1$, $B_1$, and $B_1^+$ are satisfied by: abelian surfaces, polarized surfaces of general type with $h^{2,0}(X)=1$ and $\mc{K}_X \cdot \mc{K}_X=1$, and polarized hyperk\"{a}hler varieties with $b_2>3$ (in particular, $K3$ surfaces). Cubic fourfolds (polarized via $\mc{O}_{\mathbb{P}^5}(1)$) satisfy $A_2$, $B_2$, and $B_2^+$.
\end{prop}
This relies on the work of many people; see \S $2$ and \S $3$ of \cite{andre:hyperkaehler}. Finally, Andr\'{e} observes that these axioms are independent of the choice of embedding $F \into \CC$; note that for $A_k$ this is a special case of Remark \ref{hodgeQcoefficients}. 
\begin{rmk}
The proposition does not apply to hyperk\"{a}hler varieties with Betti number $3$, since their (projective) deformation theory is not sufficiently robust. In fact, it is believed that such hyperk\"{a}hlers do not exist. Regardless, in that case one can still say something about the motivic lifting problem, since the $\ell$-adic representation $H^2(X_{\overline{F}}, \Qlb)$ is potentially abelian. In principle, this should force the underlying motivic Galois representation to factor through $\mc{T}_F$, at which point we would apply Proposition \ref{taniyamalift}; but this argument would require checking that two $\mc{G}_F$-representations with isomorphic $\ell$-adic realizations are themselves isomorphic, i.e. an unknown case of the Tate conjecture. At least we have (as in Proposition \ref{taniyamalift}):
\begin{lemma}\label{b2=3}
Let $X/F$ be a smooth projective variety with $b_2(X)=3$. If necessary, replace $F$ with a quadratic extension trivializing $\det H^2(X_{\overline{F}}, \Ql)$. Then the $\ell$-adic realization $\rho_{\ell} \colon \mc{G}_F \to \mr{SO}(H^2(X_{\overline{F}}, \Qlb))$ has a lift $\tilde{\rho}_{\ell}$ to $\mr{GSpin}(H^2(X_{\overline{F}}, \Qlb))$ that arises as the $\ell$-adic realization of a representation of $\mc{T}_{F, E}$ for suitable $E$.
\end{lemma}
\end{rmk}
\begin{notation}
\begin{itemize}
\item From now on we view $F$ as a subfield of $\CC$, with $\overline{F}$ its algebraic closure in $\CC$. These embeddings will be used implicitly to define Betti realizations, motivic Galois groups (for Betti realizations), and \'{e}tale-Betti comparisons, for varieties (or motives) over extensions of $F$ inside $\overline{F}$. 
\item We write $V_{\Q}$ for $Prim^{2k}(X_{\CC}, \Q)(k)$, and for a prime $\ell$, we will write $V_{\ell}$ for the $\ell$-adic realization $Prim^{2k}(X_{\overline{F}}, \Q_{\ell})(k)$. As before, $Prim^{2k}(X)(k)$ will be used to indicate the underlying motive, i.e. object of $\mc{M}_F$.\index{s}{$V_{\Q}$}\index{s}{$V_{\ell}$}
\item For fields $E$ containing $\Q$, we will sometimes denote the extension of scalars $V_{\Q} \otimes_{\Q} E$ by $V_{E}$ (or similarly for $\Ql$ and $V_{\ell}$). This is the Betti realization of an object $Prim^{2k}(X)(k)_E$ of $\mc{M}_{F, E}$.\index{s}{$V_E$}
\item The same subscript conventions will hold for other motives we consider, especially the direct factors of $Prim^{2k}(X)(k)$ given by the algebraic cycles (to be denoted $Alg$) and its orthogonal complement, the transcendental lattice ($T$). These motives will be discussed in \S \ref{prelimred}. In the meantime, we recall that the rational Hodge structure $Prim^{2k}(X_{\CC}, \Q)(k)$ has a $\Q$-subspace spanned by its Hodge classes; the orthogonal complement $T_{\Q}$ of this subspace is the \textit{transcendental subspace} is again a polarized $\Q$-Hodge structure.\index{t}{transcendental lattice, or subspace, of the Hodge structure $Prim^{2k}(X_{\CC}, \Q)(k)$} 
\item For an object $M$ of $\mc{M}_F$, we will sometimes write $\rho^M$ for the associated motivic Galois representation.\index{s}{$\rho^{M}$, for a motive $M$: the associated motivic Galois representation} For the $\ell$-adic realization of this motivic Galois representation (i.e., the representation of $\gal{F}$ on $H_{\ell}(M)$), we write $\rho^M_{\ell}$.\index{s}{$\rho^M_{\ell}$}
\item
First assume $\dim V_{\Q}=m$ is odd. The group $\mr{GSpin}(V_{\Q})$ may not have a rationally-defined spin representation, but it does after some extension of scalars, and over any suitably large field $E$, we denote by $r_{spin} \colon \mr{GSpin}(V_{E}) \to \mr{GL}(W_E)$ this algebraic representation. If $m$ is even, we similarly denote by $W_E= W_{+, E} \oplus W_{-, E}$ the direct sum of the two half-spin representations. In contrast to $V_{E}= V_{\Q} \otimes_{\Q} E$, this $W_E$ is not necessarily an extension of scalars from an underlying $\Q$-space.\index{s}{$r_{spin}$}\index{s}{$W_E$}\index{s}{$W_{+, E}$}\index{s}{$W_{-, E}$} 
\end{itemize}
\end{notation}
Possibly after replacing $F$ by a quadratic extension, we may assume that the $\gal{F}$-representation $V_{\ell}$ is special orthogonal, i.e. $\rho_{\ell} \colon \gal{F} \to \mr{SO}(V_{\ell}) \subset \mr{O}(V_{\ell})$. Since for almost all finite places $v$, the frobenius $fr_v$ acts on $H^{2k-2}(X_{\overline{F}}, \Ql)$ and $H^{2k}(X_{\overline{F}}, \Ql)$ with eigenvalues that are independent of $\ell$ (\cite{deligne:weil1}), and trivially on the $\ell$-adic Chern class $\eta_{\ell}$, the eigenvalues of $fr_v$ on $Prim^{2k}(X_{\overline{F}}, \Ql)(k)$ are independent of $\ell$. Consequently, if $\det \rho_{\ell}$ is trivial for one $\ell$, then it is for all $\ell$. We will from now on assume, for technical simplicity, that this determinant condition is satisfied.

\subsection{The Kuga-Satake construction}\label{KSreview}
We review the classical Kuga-Satake construction and outline Andr\'{e}'s refinement, which implies a potential version (i.e., after replacing $F$ by a finite extension) of our motivic lifting result. See \S $4$ and \S $5$ of \cite{andre:hyperkaehler}. His approach is inspired by that of \cite{deligne:weilK3}, in which Deligne used the Kuga-Satake construction (in families) to reduce the Weil conjectures for $K3$ surfaces to the (previously known) case of abelian varieties. Let $V_{\Z}$\index{s}{$V_{\Z}$} be the polarized quadratic lattice $Prim^{2k}(X_{\CC}, \Z)(k)$ of the previous subsection. Write $m=2n$ or $m=2n+1$ for its rank. Basic Hodge theory implies that the pairing on $V_{\RR}$ is negative-definite on the sub-space $(H^{1, -1} \oplus H^{-1, 1})_{\RR}$, and positive-definite on $H^{0,0}_{\RR}$, hence has signature $(2, m-2)$. Recalling the discussion of \S \ref{liftinghodge}, this real Hodge structure yields a homomorphism $h \colon \bS \to \mr{SO}(V_{\RR})$ that lifts uniquely to a homomorphism $\tilde{h} \colon \bS \to \mr{GSpin}(V_{\RR})$ whose composition $N_{spin} \circ \tilde{h}$ with the Clifford norm $N_{spin}$ is the usual norm $\bS \to \mathbb{G}_{m, \RR}$. 

Let $C(V_{\Z})$ and $C^{+}(V_{\Z})$ denote the Clifford algebra and the even Clifford algebra associated to the quadratic space $V_{\Z}$.\index{s}{$C(V_{\Z})$}\index{s}{$C^+(V_{\Z})$} Let $L_{\Z}$ be a free left $C^+(V_{\Z})$-module of rank one.\index{s}{$L_{\Z}$} Denote by $C^+$ the ring $\End_{C^+(V_{\Z})}(L_{\Z})^{op}$.\index{s}{$C^+$} Because of the `op,' $C^+$ naturally acts on $L_{\Z}$ on the right, and we then correspondingly have $\End_{C^+}(L_{\Z}) \cong C^+(V_{\Z})$. The choice of a generator $x_0$ of $L_{\Z}$ induces a ring isomorphism
\begin{align*}
C^+(V_{\Z}) &\xrightarrow{\phi_{x_0}} C^+ \\
c \mapsto &(b x_0 \mapsto bc x_0).
\end{align*}
Via $\tilde{h}$ and the tautological representation $\mr{GSpin}(V_{\RR}) \into C^+(V_{\RR})^{\times}$, $L_{\RR}:= L_{\Z} \otimes_{\Z} \RR$ acquires a Hodge structure of type $(1,0), (0,1)$; the Hodge type is easily read off from the fact that $\End_{C^+}(L_{\Z}) \cong C^+(V_{\Z})$ is also an isomorphism of Hodge-structures,\footnote{Giving $C^+(V_{\Z})$ the Hodge structure induced from that on the tensor powers of $V_{\Z}$.} or by recalling the remarks of \S \ref{liftinghodge} and listing the weights of the spin representation. Concretely (as in the original construction of Satake \cite{satake:clifford} and Kuga-Satake \cite{kuga-satake}), we can choose an orthogonal basis $e_1, e_2$ of $(H^{1, -1} \oplus H^{-1, 1})_{\RR}$, normalized so that $\langle e_i, e_i \rangle =-1$ for $i=1, 2$. Then the automorphism of $C^+(V_{\RR})$ given by multiplication by $e_1 e_2$ is a complex structure, by the defining relations for the Clifford algebra. In fact, this integral Hodge structure is polarizable: this can be shown explicitly, or by a very soft argument (apply Proposition $2.11.b.iii'$ of \cite{deligne:weilK3}). We therefore obtain a complex abelian variety, the Kuga-Satake abelian variety $KS(X)=KS(X_{\CC}, \eta, k)$, associated to our original $V_{\Z}$.\index{s}{$KS(X)$}\index{t}{Kuga-Satake construction} The right-action of $C^+$ on $L_{\Z}$ commutes with the Hodge structure, so $C^+$ acts as endomorphisms of $KS(X)$. Generically, this will be the full endomorphism ring; we will have to be attentive later to how much bigger the endomorphism ring can be. 

One of the main technical ingredients in \cite{andre:hyperkaehler} is the following descent result, which uses rigidity properties of the Kuga-Satake construction:\index{s}{$A_{F'}$: Andr\'{e}'s descent of the classical Kuga-Satake abelian variety}
\begin{lemma}[Main Lemma $1.7.1$ of \cite{andre:hyperkaehler}]\label{andremainlemma}
Let $(X, \eta)$ be a polarized variety over a subfield $F$ of $\CC$, satisfying properties $A_k$ and $B_k$. Then there exists an abelian variety $A_{F'}$ over some *finite* extension $F'/F$ such that
\begin{itemize}
\item The base-change $A_{\CC}$ is the Kuga-Satake variety $KS(X)$;
\item there is a subalgebra $C^+$ of $\End(A_{F'})$ and an isomorphism of $\Z_{\ell}[\gal{F'}]$-algebras
\[
\End_{C^+}\left(H^1(A_{\overline{F}}, \Z_{\ell})\right) \cong C^+(Prim^{2k}(X_{\overline{F}}, \Z_{\ell})(k)).
\]
\end{itemize}
\end{lemma}
We subsequently write $L_{\ell}$ for the $\ell$-adic realization $H^1(A_{F'} \otimes \overline{F}, \Ql)$.\index{s}{$L_{\ell}$} The main result of the `motivated' theory of hyperk\"{a}hlers, (the foundation of Andr\'{e}'s results on the Tate and Shafarevich conjectures) is:
\begin{thm}[see Theorem $6.5.2$ of \cite{andre:hyperkaehler}]\label{motivatedHK}
For some finite extension $F'/F$, $Prim^{2k}(X_{F'})(k)$ is a direct factor (in $\mc{M}_F$) of $\underline{End} \left( H^1(A_{F'}) \right)$, and both $Prim^{2k}(X_{F'})(k)$ and $H^1(A_{F'})$ have connected motivic Galois group.
\end{thm}
Write $\rho^A \colon \mc{G}_{F'} \to \mr{GL}(L_{\Q})$ and $\rho^V \colon \mc{G}_F \to \mr{SO}(V_{\Q})$ for the motivic Galois (for the Betti realization) representations associated to $H^1(A_{F'})$ and $Prim^{2k}(X)(k)$.\index{s}{$\rho^A$}\index{s}{$\rho^V$} We will somewhat sloppily use the same notation $\rho^A$ and $\rho^V$ for various restrictions of these representations (eg, to $\mc{G}_{\overline{F}}$). Note that by Corollary \ref{MT}, the image of $\rho^A$ is the Mumford-Tate group of $KS(X)$, which is easily seen to be contained in $\mr{GSpin}(V_{\Q})= \{x \in C^+(V_{\Q})^{\times}: xV_{\Q}x^{-1} = V_{\Q} \}$. More precisely: 
\begin{cor}\label{potentialmotiviclift}
For some finite extension $F'/F$, $\rho^A$ factors through $\mr{GSpin}(V_{\Q})$ and lifts $\rho^V|_{\mc{G}_{F'}}$. A fortiori, $\rho^A_{\ell}$ lifts $\rho^V_{\ell}|_{\gal{F'}}$.
\end{cor}
\proof
Over $\CC$ (or $\overline{F}$), the image of $\rho^A$ is $MT(A_{\CC})$, hence contained in $\mr{GSpin}(V_{\Q})$. We can therefore compare the two maps 
\[
\mc{G}_{\CC}  \xrightarrow{(\rho^V,  \pi \circ \rho^A)} \mr{SO}(V_{\Q}) \times \mr{SO}(V_{\Q}) \to \mr{GL}\left(C^+(V_{\Q})\right) \times \mr{GL}\left(C^+(V_{\Q})\right),
\]
where $\pi$ denotes the projection $\mr{GSpin}(V_{\Q}) \to \mr{SO}(V_{\Q})$. Under the motivated isomorphism 
\[
C^+(Prim^{2k}(X_{\CC})(k)) \cong \underline{End}\left( H^1(A) \right)
\]
the adjoint action of $\rho^A$ on $\End(L_{\Q})$ agrees with the action of $\pi \circ \rho^A$ on $C^+(V_{\Q})$, so the two compositions $\mc{G}_{\CC} \to \mr{GL}\left(C^+(V_{\Q})\right)$ above coincide. Now, if $m$ is odd, $\mr{SO}(V_{\Q}) \to \mr{GL}(C^+(V_{\Q}))$ is injective, so $\rho^V= \pi \circ \rho^A$; that is, at least over $\CC$, $\rho^A$ lifts $\rho^V$. If $m$ is even, we deduce that the compositions 
\[
\mc{G}_{\CC}  \xrightarrow{(\rho^V, \pi \circ \rho^A)} \mr{SO}(V_{\Q}) \times \mr{SO}(V_{\Q}) \to \mr{GL}(\wedge^2 V_{\Q}) \times \mr{GL}(\wedge^2 V_{\Q})
\]
agree.\footnote{The filtration on $C^+(V_{\Q})$ given by the image of $V^{\otimes \leq 2i}$ is motivated, since $\mr{O}(V)$-stable, so $\rho^V$ and $\pi \circ \rho^A$ coincide on the $i=1$ graded piece, which is just $\wedge^2(V_{\Q})$.} The kernel of $\wedge^2 \colon \mr{SO}(V_{\Q}) \to \mr{GL}(\wedge^2 V_{\Q})$ is $\{\pm 1\}$ (and central), so we see that $\rho^V$ and $\pi \circ \rho^A$ agree up to twisting by some character $\chi \colon \mc{G}_{\CC} \to \{ \pm 1\} \subset \mr{SO}(V_{\Q})$. But clearly this character factors through $\mc{G}^{M}_{\CC}$, where $M= Prim^{2k}(X_{\CC})(k) \oplus H^1(A)$, and by Theorem \ref{motivatedHK}, $\mc{G}^M_{\CC}$ is connected. Therefore $\chi$ is trivial, and $\rho^A$ lifts $\rho^V$ as $\mc{G}_{\CC}$-representations. By Lemma \ref{indsepclosed}, the same holds for the corresponding $\mc{G}_{\overline{F}}$-representations, and then as in Lemma \ref{connected} the same holds over some finite extension $F'/F$.
\endproof
This result does not always hold with $F'=F$, and it is our task in the coming sections to achieve a motivic descent over $F$ itself. The Kuga-Satake variety is highly redundant, and it is technically convenient to work with a smaller `spin' abelian variety, many copies of which constitute the Kuga-Satake variety. Since the rational Clifford algebra $C^+_{\Q}$ may not be split, this requires a finite extension of scalars $E/\Q$, after which we can work in the isogeny category $AV^0_{F', E}$ of abelian varieties over $F'$ with $E$-coefficients. We take $F'$ as in the Corollary, and now to ease the notation, we write simply $A$ for $A_{F'}$.\index{s}{$A$: short-hand for the abelian variety $A_{F'}$} 
\begin{lemma}\label{spinAV}
There exists a number field $E$ and an abelian variety $B/F'$ with endomorphisms by $E$ such that there is a decomposition in $AV^0_{F', E}$:
\begin{align*}
&A \otimes_{\Q} E \cong B^{2^n} \quad \text{if $m=2n+1$;}\\
&A \otimes_{\Q} E \cong B^{2^{n-1}} \quad \text{if $m=2n$.}
\end{align*} 
The $\ell$-adic realization $H^1(B_{\overline{F}}, \Ql)$ is isomorphic to the composite $r_{spin} \circ \rho^A_{\ell}$ as $(E \otimes \Q_{\ell})[\gal{F'}]$-modules, where as before $\rho^A_{\ell}$ denotes the representation $\gal{F'} \to \mr{GSpin}(V_{\ell} \otimes E)$ obtained from $L_{\ell} \otimes E$, and $r_{spin}$ denotes either the spin ($m$ odd) or sum of half-spin ($m$ even) representations of $\mr{GSpin}(V_{\ell} \otimes E)$. When $m=2n$, $B$ decomposes in $AV^0_{F', E}$ as $B_{+} \times B_{-}$, corresponding to the two half-spin representations.\index{s}{$B$: `spin representation' abelian variety}\index{s}{$B_+$}\index{s}{$B_-$}
\end{lemma}
\proof
Choose a number field $E$ splitting $C^+(V_{\Q})$. Then, letting $W_E$ denote either the spin ($m$ odd) representation or the direct sum $W_{+, E} \oplus W_{-, E}$ of the two irreducible half-spin representations ($m$ even), $C^+(V_E)$ is isomorphic as $\mr{GSpin}(V_E)$-representations to $W_E^{2^n}$, or to $W_{E, +}^{2^{n-1}} \oplus W_{E, -}^{2^{n-1}}$. As $E$-algebra it is then isomorphic either to $\End(W_E) \cong M_{2^n}(E)$ or $\End(W_{+, E}) \oplus \End(W_{-, E}) \cong M_{2^{n-1}}(E)\oplus M_{2^{n-1}}(E)$.\footnote{Note that the field $E$ can be made explicit if we know the structure of the quadratic lattice $V_{\Z}$. See Example \ref{K3eg} for a case where $E=\Q$.} Using the orthogonal idempotents in $C^+_{E} \cong C^+(V_{E})$, we decompose the object $A_{F'} \otimes_{\Q} E$ of $AV^0_{F', E}$ into $2^n$ (when $m$ is odd) or $2^{n-1}$ (when $m$ is even) copies of an abelian variety $B/F'$ with complex multiplication by $E$.
\endproof
We state a simple case of the main result of this section; for the proof, see \S \ref{genericcase}. More general versions will be proven in stages, depending on the complexity of the transcendental lattice of $V_{\Q}$.
\begin{thm}\label{motiviclift}
Let $(X, \eta)$ be a polarized variety over a number field $F \subset \CC$ for which $Prim^{2k}(X_{\CC}, \Z)(k)$ satisfies axioms $A_k$, $B_k$, and $B_k^+$. Possibly enlarging $F$ by a quadratic extension, we may assume as above that for all $\ell$, $\det V_{\ell}=1$, so that the $\ell$-adic representation $\rho^V_{\ell}$ maps $\gal{F}$ to $\mr{SO}(V_{\ell})$. We make the following hypothesis on the monodromy:
\begin{itemize}
\item The transcendental lattice $T_{\Q}$ has $\End_{\Q-Hodge}(T_{\Q})= \Q$.
\item $\det(T_{\ell})=1$.\footnote{Unlike the first hypothesis, this determinant condition is a technicality; it can again be arranged after an (independent of $\ell$) quadratic extension. Of course, in the generic case in which $V_{\Q}= T_{\Q}$, it is no additional hypothesis.}
\end{itemize}
Then there exists: 
\begin{itemize}
\item a finite extension $E/\Q$
\item an object $B$ of $AV^0_{F, E}$;
\item an Artin motive $M$ in $\mc{M}_{F, E}$
\item for each prime $\lambda$ of $E$, a lift $\tilde{\rho}_{\lambda}$
\[
\xymatrix{
&   \mr{GSpin}(V_{\ell} \otimes E_{\lambda}) \ar[d] \ar[r]^{r_{spin}} & \mr{GL}_{E_{\lambda}}(W_{E_{\lambda}}) \\
\gal{F} \ar@{-->}[ru]^{\tilde{\rho}_{\lambda}} \ar[r] & \mr{SO}(V_{\ell}) \subset \mr{SO}(V_{\ell} \otimes E_{\lambda})
}
\]
\end{itemize}
such that the composite $r_{spin} \circ \tilde{\rho}_{\lambda}$ is isomorphic to some number of copies\footnote{Which can be made explicit.} of $M_{\lambda} \otimes \left(H^1(B_{\overline{F}}, \Ql) \underset{E \otimes \Ql}{\otimes} E_{\lambda} \right)$ as $E_{\lambda}$-linear $\gal{F}$-representations. In particular, the system of lifts $\{\tilde{\rho}_{\lambda}\}_{\lambda}$ is weakly compatible. It is even `motivic' in the sense of arising from a lift
\[
\xymatrix{
& \mr{GSpin}(V_E) \ar[d] \\
\mc{G}_{F, E} \ar@{-->}[ru]^{\tilde{\rho}} \ar[r]_{\rho} & \mr{SO}(V_E).
}
\]
\end{thm}
This theorem is an optimal (up to the $\mr{O}(V_{\ell}) \supset \mr{SO}(V_{\ell})$ distinction) arithmetic refinement of the (\textit{a priori} highly transcendental) Kuga-Satake construction, showing the precise sense in which it descends to the initial field of definition $F$.

\subsection{A simple case}\label{explicitvariant}
To achieve the refined descent of Theorem \ref{motiviclift}, the basic idea is to apply our Galois-theoretic lifting results (which apply over $F$) to deduce $\gal{F}$-invariance of the abelian variety $B$ that we know to exist over $F'$; Faltings's isogeny theorem (\cite{faltings:endlichkeit}) implies that this invariance is realized by actual isogenies; then we apply a generalization of a technique used by Ribet (\cite{ribet:Qcurves}; our generalization is Proposition \ref{generalribet}) to study elliptic curves over $\Qb$ that are isogenous to all of their $\gal{\Q}$-conjugates (so-called ``$\Q$-curves").\index{t}{$\Q$-curves} Ribet's technique applies to elliptic curves without complex multiplication, and we will have to keep track of monodromy groups enough to reduce the descent problem to one for an absolutely simple abelian variety. In some cases, a somewhat `softer' argument than Ribet's works-- we give an example in Lemma \ref{softdescent}-- but in addition to being satisfyingly explicit, the Ribet method seems to be more robust. 

But first we prove Theorem \ref{motiviclift} in the simplest case, when the Hodge structure $V_{\Q}$ is `generic,' $\dim V_{\Q}= 2n+1$ is odd, and the even Clifford algebra $C^+(V_{\Q})$ is split over $\Q$. Our working definition of `generic' will be that $V_{\Q}$ contains no trivial $\Q$-Hodge sub-structures (i.e., it is equal to its transcendental lattice), and that $\End_{\Q-Hodge}(V_{\Q})= \Q$.
\begin{eg}\label{K3eg}
If $X/F$ is a generic $K3$ surface, then these hypotheses are satisfied. The $K3$-lattice $H^2(X_{\CC}, \Z)$ is an even unimodular lattice of rank $22$ whose signature over $\RR$ is (by Hodge theory) $(3, 19)$. The classification of even unimodular lattices implies it is isomorphic (over $\Z$) to $(-E_8)^{\oplus 2} \oplus U^{\oplus 3}$, where $E_8$ and $U$ are the $E_8$-lattice and the hyperbolic plane lattice. Over $\Q$, the orthogonal complement of the ample class $\eta$ is isomorphic to $(-E_8)^{\oplus 2} \oplus U^{\oplus 2} \oplus \langle -q(\eta) \rangle$, where $q$ is the quadratic form and $\langle \alpha \rangle$ denotes the one-dimensional quadratic space with a generator whose square is $\alpha$. Since $Prim^2(X_{\CC}, \Q)$ is odd-dimensional, the basic structure theory of Clifford algebras (see Chapter $9$, especially Theorem $2.10$, of \cite{scharlau}) implies that $C^+(Prim^2(X_{\CC}, \Q)$ is a central simple algebra over $\Q$, whose Brauer class is simply twice the Brauer class of $E_8$ plus twice the (trivial) Brauer class of $U$. This is obviously the trivial class, so $C^+(V_{\Q})$ is in this case isomorphic to $M_{2^{10}}(\Q)$. More generally, this argument applies to $V_{\Q}= Prim^2(X_{\CC}, \Q)(1)$ for $X$ a hyperk\"{a}hler satisfying:
\begin{itemize}
\item $b_2(X)>3$;
\item $H^2(X_{\CC}, \Z)$ is even;
\item the number of copies of the $E_8$-lattice, i.e. the number $\frac{b_2(X)-6}{8}$, is even.
\end{itemize} 
\end{eg}
Now, our `generic' hypothesis implies that the Hodge group of $V_{\Q}$ is the full $\mr{SO}(V_{\Q})$ (more generally, see Zarhin's result, quoted as Proposition \ref{zarhin}, below); it follows without difficulty that the Mumford-Tate group $MT(A_{\CC})$ is the full $\mr{GSpin}(V_{\Q})$, and $\End^0(A_{\CC})= C^+_{\Q}$. By our second simplifying hypothesis, this Clifford algebra is isomorphic to a matrix algebra $M_{2^n}(\Q)$, and, writing $W_{\Q}$ for the spin representation of $\mr{GSpin}(V_{\Q})$, we have two isomorphisms of $\mr{GSpin}(V_{\Q})$-representations:
\begin{align*}
C^+(V_{\Q}) &\cong W_{\Q}^{2^n};\\
C^+(V_{\Q})_{ad} &\cong \End(W_{\Q}),
\end{align*}
where $\mr{GSpin}(V_{\Q})$ acts on the first $C^+(V_{\Q})$ by left-multiplication, and on $C^+(V_{\Q})_{ad}$ by conjugation (i.e., through the natural $\mr{SO}(V_{\Q})$-action).

We may certainly enlarge $F'$ to a finite extension over which all endomorphisms of $A_{\CC}$ are defined, and the complex multiplication by $C^+_{\Q}$ then gives an isogeny decomposition $A \sim B^{2^n}$, where $B/F'$ is an abelian variety with $\End^0(B)= \Q$. We can take the spin representation $W_{\Q}$ to be equal to $H^1(B_{\CC}, \Q)$ (and, extending scalars and invoking comparison isomorphisms, we get identifications with other cohomological realizations of $B$-- in particular, the $\ell$-adic realization $W_{\lambda}= W_{\Q} \otimes_{\Q} \Ql$).   
\proof[Proof of Theorem \ref{motiviclift} when $V_{\Q}$ is generic]\label{ribet}
We sketch the argument, with some details postponed until later sections-- the goal here is to review Ribet's method and outline the argument in a simple case. We will now apply the technique of \cite{ribet:Qcurves} in combination with our abstract Galois-theoretic lifting results. By Corollary \ref{galoisliftmotiviccase}, there exists a lift 
\[
\tilde{\rho}_{\ell} \colon \gal{F} \to \mr{GSpin}(V_{\Qlb})
\]
of $\rho^V_{\ell}$, and we can normalize this lift so that in the spin representation its labeled Hodge-Tate weights $\mr{HT}_{\tau}(r_{spin} \circ \tilde{\rho}_{\ell})$ are of `abelian variety'-type, i.e. $2^{n-1}$ zeroes and $2^{n-1}$ ones, for each $\tau \colon F \into \Qlb$.\footnote{For details, see Lemma \ref{spinlift}.} This normalization determines $\tilde{\rho}_{\ell}$ up to finite-order twist, and it implies that $\tilde{\rho}_{\ell}|_{\gal{F'}}$ is a finite-order twist of $\rho^A_{\ell}$, since they both lift $\rho^V_{\ell}|_{\gal{F'}}$ with the same Hodge-Tate data. We may therefore replace $F'$ by a finite extension and assume 
\[
\rho^A_{\ell} = \tilde{\rho}_{\ell}|_{\gal{F'}}.
\]
Since $\tilde{\rho}_{\ell}$ begins life over $F$, we see that $\rho^A_{\ell}$ is $\gal{F}$-conjugation invariant. The composition of $\rho^A_{\ell}$ with the Clifford representation is $2^n$ copies of the $\ell$-adic representation $\rho^B_{\ell}$ associated to $B$, so $\rho^B_{\ell}$ is also $\gal{F}$-invariant. By Faltings's theorem, for each $\sigma \in \gal{F}$, there exists an isogeny $\mu_{\sigma} \colon {}^{\sigma}B \to B$; we can and do arrange that $\mu_{\sigma}$ is defined first for a (finite) system of representatives $\sigma_i$ in $\gal{F}$ for $\Gal(F'/F)$, and then defined in general by $\mu_{\sigma_i h}= \mu_{\sigma_i}$ for all $h \in \gal{F'}$. The collection of $\mu_{\sigma}$ yields an obstruction class 
\[
[c_B] \in H^2(\gal{F}, \End^0(B)^{\times})= H^2(\gal{F}, \Q^\times),\footnote{The local constancy of the isogenies $\mu_{\sigma}$ implies this is a continuous cohomology class, with $\Q^\times$ equipped with the discrete topology.}
\]
given by 
\[
c_B(\sigma, \tau)= \mu_{\sigma} \circ {}^{\sigma}\mu_{\tau} \circ \mu_{\sigma \tau}^{-1}.
\]
That is, $c_B$ measures the failure of the diagram
\[
\xymatrix{
{}^{\sigma \tau}B \ar[r]^{\mu_{\sigma \tau}} \ar[d]_{{}^{\sigma}\mu_{\tau}} & B \\
{}^{\sigma}B \ar[ur]_{\mu_{\sigma}} & 
}
\]
to commute. Now, the class $c_B$ may be non-trivial, and the abelian variety $B$ may not descend (up to isogeny) to $F$.\footnote{That triviality of this class is equivalent to isogeny-descent is Theorem $8.2$ of \cite{ribet:Qcurves}.} Nevertheless, Tate's vanishing result $H^2(\gal{F}, \Qb^\times)=0$ tells us that there is a continous $1$-cochain $\alpha \colon \gal{F} \to \Qb^\times$ whose coboundary equals $c_B$, i.e. $c_B(\sigma, \tau)= \alpha(\sigma \tau) \alpha(\tau)^{-1} \alpha(\sigma)^{-1}$. By continuity, $\alpha$ is locally constant with respect to $\gal{F''}$ for some finite $F''/F'$, and it takes values in some finite extension $\Q(\alpha)$ of $\Q$. We now consider the restriction of scalars abelian variety $C:= \Res_{F''/F}(B)$.\footnote{Here $B$ is really the base-change to $F''$, but we omit this not to clutter the notation.} $C$ is an abelian variety over $F$ with endomorphism algebra $\mc{R}= \End^0(C)$, which as $\Q$-vector space is isomorphic to 
\[
\Hom_{F''}\left(\prod_{\sigma \in \Gal(F''/F)} {}^{\sigma}B, B\right)= \bigoplus_{\sigma} \Q \mu_{\sigma},
\]
where again we use the fact that $\End^0(B_{\overline{F}})= \Q$. Write $\lambda_{\sigma}$ for the element of $\mc{R}$ corresponding to $\mu_{\sigma}$ under this isomorphism. 
\begin{lemma}[Lemma $6.4$ of \cite{ribet:Qcurves}]
The algebra structure of $\mc{R}$ is given by $\lambda_{\sigma} \lambda_{\tau}= c_B(\sigma, \tau) \lambda_{\sigma \tau}$, so there is a $\Q$-algebra homomorphism $\alpha \colon \mc{R} \to \Q(\alpha)$ given by the $\Q$-linear extension of $\lambda_{\sigma} \mapsto \alpha(\sigma)$.
\end{lemma}
Since the isogeny category of abelian varieties over $F$ is a semi-simple abelian category, we can form the object 
\[
M= \Res_{F''/F}(B) \otimes_{\mc{R}, \alpha} \Q(\alpha).
\]
We regard $M$ as an object of $AV^0_{F, \Q(\alpha)}$, with $\Q(\alpha)$-rank (in the obvious sense) equal to the $\Q$-rank of $B$. Moreover, for any place $\lambda \vert \ell$ of $\Q(\alpha)$, the $\lambda$-adic realization 
\[
\rho^M_{\lambda} \colon \gal{F} \to \mr{GL}_{\Q(\alpha)_{\lambda}}(M_{\lambda})
\] 
has projectivization isomorphic to the canonical projective descent to $\gal{F}$ of the $\gal{F}$-invariant, irreducible representation of $\gal{F''}$ on $H^1(B_{\overline{F}}, \Q(\alpha)_{\lambda})$.\footnote{Under the homomorphism $\mc{R} \to \End_{\gal{F''}}(\oplus {}^{sigma} H^1(B, \Ql)$, the $\lambda_{\sigma}$ permute the factors in the direct sum; the projection via $\mc{R} \xrightarrow{\alpha} \Q(\alpha)$ collapses all the factors to a single copy with scalars extended to $\Q(\alpha)$, i.e. to $H^1(B, \Ql) \otimes_{\Q} \Q(\alpha)$. This implies the claim about projective descents.} So, $\rho^M_{\lambda}|_{\gal{F''}}$ and $\rho^B_{\ell}|_{\gal{F''}} \cong r_{spin} \circ  \tilde{\rho}_{\ell}|_{\gal{F''}}$ are isomorphic up to twist, hence $r_{spin} \circ \tilde{\rho}_{\ell}$ and $\rho^M_{\lambda}$ are twist-equivalent as $\gal{F}$-representations. The representation $r_{spin} \colon \mr{GSpin}(V_{\Q}) \to \mr{GL}(W_{\Q})$ is the identity on the center, so after identifying $\rho^M_{\lambda}$ to a representation on $W_{\lambda} \otimes \Q(\alpha)_{\lambda}$, we see that it factors through $\mr{GSpin}(V_{\ell} \otimes \Q(\alpha)_{\lambda})$ as a lift of $\rho^V_{\ell} \otimes \Q(\alpha)_{\lambda}$. The required motivic lift is the representation of $\mc{G}_{F, \Q(\alpha)}$ corresponding to $H^1(M)$;\footnote{For more details on how to check this carefully, see Corollary \ref{transcor}.} and the various $\rho^M_{\lambda}$ form a compatible system because they are formed from Tate modules of abelian varieties.
\endproof

\subsection{Arithmetic descent: preliminary reduction}\label{prelimred}
Now we proceed to a more general argument, making first some preliminary reductions to the analogous lifting problem for the transcendental lattice. We must invoke Andr\'{e}'s work on the Tate 
%and Mumford-Tate [I don't understand the proof, and at least for the `generic' case I don't use this]
conjecture for $X$.\index{t}{Tate conjecture for hyperk\"{a}hler varieties}
\begin{thm}[see Theorem $1.6.1$ of \cite{andre:hyperkaehler}]
Let $(X, \eta)$ be a polarized variety over a number field\footnote{Faltings's theorem works for finitely-generated extensions of $\Q$, so this does as well.} $F$ satisfying $A_k$, $B_k$, and $B_k^+$. Then:
\begin{itemize}
\item $Prim^{2k}(X_{\overline{F}}, \Ql)(k)$ is a semi-simple $\gal{F}$-representation 
\item the Galois invariants $Prim^{2k}(X_{\overline{F}}, \Ql)(k)^{\gal{F}}$ are all $\Ql$-linear combinations of algebraic cycles;
%\item the inclusion $\overline{\rho_{\ell}(\gal{F})}^{Zar} \into \mc{G}^{V}_F \otimes_{\Q} \Ql$ of the algebraic monodromy group into the $\ell$-adic motivic Galois group of $V= Prim^{2k}(X)(k)$ is an isomorphism on connected components of the identity. 
\end{itemize}
\end{thm}
The Tate conjecture for $H^{2k}(X)$ then implies:
\begin{lemma}
There is an Artin motive $Alg$ over $F$ whose Betti realization is the subspace of $V_{\Q}$ spanned by algebraic cycle classes, and whose $\ell$-adic representation is the (Artin) $\gal{F}$-representation on $V_{\ell}^{\gal{F'}}$ for any $F'/F$ large enough (and Galois) that all of these cycle classes are defined over $F'$. The transcendental lattice $T$ likewise descends to an object of $\mc{M}_F$. In particular, there is an orthogonal decomposition
\[
V_{\ell}= Alg_{\ell} \oplus T_{\ell}
\]
of $\gal{F}$-representations (not merely $\gal{F'}$-representations).\index{s}{$Alg$: space of algebraic cycle classes as an object of $\mc{M}_F$}\index{s}{$T$: transcendental lattice as an object of $\mc{M}_F$}
\end{lemma}
\proof
Giving an Artin motive over $F$ is equivalent to giving a representation of $\gal{F}$ on a (finite-dimensional) $\Q$-vector space. In our case, the space ($\Q$-span) of cycles for homological equivalence 
\[
Z^k_{hom}(X_{F'}) \into H^{2k}(X_{\overline{F}}, \Ql)(k),
\]
or rather its intersection with $V_{\ell}$, does the trick. Since the object $Prim^{2k}(X)(k)$ of $\mc{M}_F$ is polarized, we can define $T$ in $\mc{M}_F$ as the orthogonal complement of $Alg$.
\endproof
We will enlarge $F'$ as in the Lemma, so that all algebraic classes in $Prim^{2k}(X_{\CC}, \Q)(k)$ are already defined over $F'$, and so that the motivic group $\mc{G}^T_{F'}$ is connected (see Lemma \ref{connected}). Note that $T_{\Q}$ is an orthogonal Hodge structure of type $(1, -1), (0,0), (-1, 1)$, with $h^{1, -1}=1$, so the Kuga-Satake construction applies to it as well (see Variant $4.1.5$ of \cite{andre:hyperkaehler}).  Since the Kuga-Satake variety associated to $V_{\Q}$ is simply an isogeny power of that associated to $T_{\Q}$, the latter, to be denoted $A(T)$, also descends to some finite extension $F'/F$.\index{s}{$A(T)$}

We introduce a little more notation. After extending scalars to a sufficiently large field $E$ (omitted from the notation unless we want to emphasize it) we denote by $(r_{spin, V}, W_V)$, $(r_{spin, Alg}, W_{Alg})$, and $(r_{spin, T}, W_T)$ the spin representations of these three spin groups; as before, by `the' spin representation in the $D_n$ case we will mean the direct sum of the two half-spin representations.\index{s}{$r_{spin, V}$, $r_{spin, Alg}$, $r_{spin, T}$: homomorphisms giving the spin representations associated to the various orthogonal spaces $V$, $Alg$, $T$}\index{s}{$W_V$, $W_{Alg}$, $W_T$: spin representations associated to the various orthogonal spaces $V$, $Alg$, $T$} 
\begin{lemma}\label{transclemma}
Suppose that $\rho^V_{\ell}$ factors through $\mr{SO}(Alg_{\ell}) \times \mr{SO}(T_{\ell}) \into \mr{SO}(V_{\ell})$.\footnote{At worst, ensuring this requires making a quadratic extension of $F$.} If we have found lifts $\tilde{\rho}^{Alg}_{\ell}$ and $\tilde{\rho}^T_{\ell}$ to $\mr{GSpin}(Alg_{\ell} \otimes \Qlb)$ and $\mr{GSpin}(T_{\ell} \otimes \Qlb)$ such that $r_{spin, Alg} \circ \tilde{\rho}^{Alg}_{\ell}$ and $r_{spin, T} \circ \tilde{\rho}^T_{\ell}$ are motivic, then $\rho^V_{\ell}$ has a lift $\tilde{\rho}_{\ell} \colon \gal{F} \to \mr{GSpin}(V_{\ell} \otimes \Qlb)$ such that $r_{spin, V} \circ \tilde{\rho}_{\ell}$ is also motivic. If the individual lifts $r_{spin, Alg} \circ \tilde{\rho}^{Alg}_{\ell}$ and $r_{spin, T} \circ \tilde{\rho}^T_{\ell}$ belong to compatible systems of $\ell$-adic representations, then the same holds for $r_{spin, V} \circ \tilde{\rho}_{\ell}$.\footnote{Recall that objects of $\mc{M}_F$ are not in general known to give rise to compatible systems.}
\end{lemma}
\proof
The isomorphism of \textit{graded} algebras (for the graded tensor product $\hat{\otimes}$) 
\[
C(Alg_{\Q}) \hat{\otimes}  C(T_{\Q}) \xrightarrow{\sim} C(V_{\Q})
\] 
induces an inclusion $C^+(Alg_{\Q}) \otimes C^+(T_{\Q}) \into C^+(V_{\Q})$, and then a map (not injective) 
\[
C^+(Alg_{\Q})^\times \times C^+(T_{\Q})^\times \to C^+(V_{\Q})^\times,
\]
which in turn induces a commutative diagram
\[
\xymatrix{
\mr{GSpin}(Alg_{\Q}) \times \mr{GSpin}(T_{\Q}) \ar[r] \ar[d] & \mr{GSpin}(V_{\Q}) \ar[d] \\
\mr{SO}(Alg_{\Q}) \times \mr{SO}(T_{\Q}) \ar[r] & \mr{SO}(V_{\Q}).
}
\]
As long as $\rho^V_{\ell} \colon \gal{F} \to \mr{SO}(V_{\ell})$ factors through $\mr{SO}(Alg_{\ell}) \times \mr{SO}(T_{\ell})$, this shows that we can lift $Alg_{\ell}$ and $T_{\ell}$ in order to lift $V_{\ell}$. We next want to understand the restriction to $\mr{GSpin}(Alg_{E}) \times \mr{GSpin}(T_{E})$ of the spin representation of $\mr{GSpin}(V_{E})$; Just for this argument, we will ignore the similitude factor, i.e. work with weights of $\mr{Spin}$ rather than $\mr{GSpin}$. We can write bases of the character lattices of $\mr{SO}(Alg_E)$, $\mr{SO}(T_E)$, and $\mr{SO}(V_E)$ as, respectively, $\chi_1, \ldots, \chi_a$, $\chi_{a+1}, \ldots, \chi_{a+t}$, and 
\begin{align*}
&\chi_1, \ldots, \chi_{a+t} \quad \text{if either $\dim (T)$ or $\dim (Alg)$ is even;}\\
&\chi_1, \ldots, \chi_{a+t+1} \quad \text{if both $\dim (T)$ and $\dim (Alg)$ are odd.}
\end{align*}
The set of weights of the spin representation of $\mf{so}(Alg_{E})$ (and similarly for the other cases) is then all $2^a$ weights of the form\footnote{The uniform description of the weights in the even and odd cases results from taking the sum of the two half-spin representations.}
\[
\sum_{i=1}^a \frac{\pm \chi_i}{2}.
\]
In the case where at least one of $\dim(T)$ and $\dim(Alg)$ is even, we see that the weights of $W_{Alg} \boxtimes W_{T}$ are precisely those of $W_{V}|_{\mr{GSpin}(Alg_E) \times \mr{GSpin}(T_E)}$, and therefore $W_V \cong W_{Alg} \boxtimes W_T$ as $\mr{GSpin}(Alg_E) \times \mr{GSpin}(T_E)$-representations. When both $\dim(T)$ and $\dim(Alg)$ are odd, weight-counting gives $W_V \cong (W_{Alg} \boxtimes W_T)^{\oplus 2}$.

Thus, if we have found lifts $\tilde{\rho}^{Alg}_{\ell}$ and $\tilde{\rho}^T_{\ell}$ (to $\mr{GSpin}(Alg_{\ell} \otimes \Qlb)$ and $\mr{GSpin}(T_{\ell} \otimes \Qlb)$) such that $r_{spin, Alg} \circ \tilde{\rho}^{Alg}_{\ell}$ and $r_{spin, T} \circ \tilde{\rho}^T_{\ell}$ are motivic (respectively, belong to compatible systems of $\ell$-adic representations), then the resulting lift to $\mr{GSpin}(V_{\ell} \otimes \Qlb)$, in its spin representation, is a direct sum of tensor products of motivic Galois representations (respectively, Galois representations belonging to compatible systems), hence is motivic. 
\endproof
\begin{cor}\label{transcendentalreduction}
If $\det T_{\ell}=1$ as $\gal{F}$-representation, and if we can find a lift $\tilde{\rho}^T_{\ell}$ of $\rho^T_{\ell}$ such that $r_{spin, T} \circ \tilde{\rho}^T_{\ell}$ is motivic (respectively, belongs to a compatible system), then we can do the same for $\rho^V_{\ell}$. 
\end{cor}
\proof
Since $Alg_{\ell}$ is an Artin representation, Tate's vanishing result allows us to lift to an Artin representation $\tilde{\rho}^{Alg}_{\ell} \colon \gal{F} \to \mr{GSpin}(Alg_{\ell} \otimes \Qlb)$. The corollary follows.
\endproof
\begin{cor}\label{transcor}
As in Corollary \ref{transcendentalreduction}, assume that $\det T_{\ell}=1$, and that for some number field $E$, and place $\lambda \vert \ell$ of $E$, we can find a lift $\tilde{\rho}^T_{\lambda} \colon \gal{F} \to \mr{GSpin}(T_{\ell} \otimes_{\Ql} E_{\lambda})$ such that $r_{spin, T} \circ \tilde{\rho}^T_{\lambda}$ is the $\lambda$-adic realization of an object $M$ of $\mc{M}_{F, E}$ whose base-change to some $F'/F$ is one of the spin direct factors of the Kuga-Satake motive (with scalars extended to $E$) associated to $T_{\Q}$ (see Lemma \ref{spinAV}). Then possibly enlarging $E$ to a finite extension $E'$, there is a lifting of representations of the motivic Galois group $\mc{G}_{F, E'}$:
\[
\xymatrix{
& \mr{GSpin}(V_{\Q} \otimes E') \ar[d] \\
\mc{G}_{F, E'} \ar@{-->}[ru] \ar[r]_-{\rho^V} & \mr{SO}(V_{\Q} \otimes E').
}
\]
\end{cor}
\proof
The Artin representation $\rho^{Alg}_{\ell}$ is definable over $\Q$, and lifts to $\mr{GSpin}(Alg_{E_1})$ after making some finite extension $E_1/\Q$. The conclusion of the corollary will hold with $E'$ equal to the composite $E_1 E$. We check it using the same principle as in the proof of Proposition \ref{taniyamalift}: that is, we check `geometrically' (for the restriction to $\mc{G}_{\overline{F}, E'}$) and for the $\lambda$-adic realization. By Corollary \ref{potentialmotiviclift}, the motivic representation $\rho^{A(T)}$ of $\mc{G}_{\overline{F}}$, and by extension of $\mc{G}_{\overline{F}, E'}$, factors through $\mr{GSpin}(T_{E'})$ and lifts $\rho^T \colon \mc{G}_{\overline{F}, E'} \to \mr{SO}(T_{E'})$.  Since $H^1(A(T)) \otimes_{\Q} E'$ is just some number of copies\footnote{$2^t$ for $\dim T$ odd; $2^{t-1}$ for $\dim T$ even.} of $M$, the same is true of $\rho^M$. By assumption, the $\lambda$-adic realization $\rho^M_{\lambda}$ is just $r_{spin} \circ \tilde{\rho}^T_{\lambda}$, so this also factors through $\mr{GSpin}(T_{\lambda})$, lifting $\rho^T_{\lambda}$. Using the section $s_{\lambda} \colon \gal{F} \to \mc{G}_{F, E'}(E'_{\lambda})$, so that $\mc{G}_{F, E'}(E'_{\lambda})= s_{\lambda}(\gal{F}) \cdot \mc{G}_{\overline{F}, E'}(E'_{\lambda})$, we conclude as in Proposition \ref{taniyamalift} that $\rho^M$ lifts $\rho^T \otimes E'$.

Combining $\rho^M$ with the lift of $\rho^{Alg}_{\ell}$, as in Lemma \ref{transclemma}, we similarly find a lift to $\mr{GSpin}(V_{\Q} \otimes_{\Q} E')$ of our given $\rho^V \colon \mc{G}_{F, E'} \to \mr{SO}(V_{\Q} \otimes_{\Q} E')$.
\endproof
\subsection{Arithmetic descent: the generic case}\label{genericcase}
To summarize, we have reduced the problem of finding motivic lifts of the motive (over $F$) $Prim^{2k}(X)(k)$ to the corresponding problem for the transcendental lattice $T$. In this section we treat the `generic' case in which the Hodge structure $T_{\Q}$ has trivial endomorphism algebra, using a variant of Ribet's method (which will return in \S \ref{nongenericodd}) We isolate this case both to demonstrate a slightly different argument, and because the non-generic cases will require even deeper input, Andr\'{e}'s proof of the Mumford-Tate conjecture in this context (see Theorem \ref{andreMT}). The starting-point of the analysis of the motive $T$ is Zarhin's calculation in \cite{zarhin:K3hodge} of the Mumford-Tate group:
\begin{prop}[Zarhin]\label{zarhin}
Let $T_{\Q}$ be a $\Q$-Hodge structure with orthogonal polarization and Hodge numbers $h^{1, -1}=1$, $h^{0. 0}>0$, and $h^{p, q}=0$ if $|p-q|>2$. Moreover assume that $T_{\Q}$ contains no copies of the trivial $\Q$-Hodge structure. Then $E_T= \End_{\Q-Hodge}(T_{\Q})$ is a totally real or CM field, and $T_{\Q}$ is a simple $MT(T_{\Q})$-module. There is a non-degenerate $E_T$-hermitian\footnote{In the case of totally real $E_T$, this means symmetric.} pairing 
\[
\langle \cdot , \cdot \rangle \colon T_{\Q} \times T_{\Q} \to E_T
\]
such that
\[
MT(T_{\Q})= \Aut(T_{\Q}, \langle \cdot , \cdot \rangle_E) \subset \mr{SO}(T_{\Q}).
\]
\end{prop}
For the rest of this section, we assume that $\End_{\Q-Hodge}(T_{\Q})= \Q$; the Mumford-Tate group $MT(T_{\Q})$ is then the full $\mr{SO}(T_{\Q})$.
\begin{lemma}
Assume $\End_{\Q-Hodge}(T_{\Q})= \Q$. Then the Mumford-Tate group, and therefore the motivic Galois group, of the Kuga-Satake variety $A(T)$ is equal to $\mr{GSpin}(T_{\Q})$. Consequently, $\End^0(A(T))= C^+_{\Q} \cong C^+(T_{\Q})$.
\end{lemma}
\proof
Easy.
\endproof
We now choose a number field $E/\Q$ splitting $C^+(T_{\Q})$, and consider the decomposition in $AV^0_{F', E}$
\begin{align*}
&A(T) \sim  B(T)^{2^t}  \quad \text{if $\dim(T_{\Q})= 2t+1$;}\\
&A(T) \sim B(T)^{2^{t-1}} \sim \left( B_{+}(T) \times B_{-}(T) \right)^{2^{t-1}}  \quad \text{if $\dim(T_{\Q})= 2t$,}
\end{align*}
as in Lemma \ref{spinAV}.\index{s}{$B(T)$} We saw in Corollary \ref{potentialmotiviclift} that $\rho^{A(T)}$ factors through $\mr{GSpin}(T_{\Q})$ and lifts $\rho^T$; viewing $\rho^{A(T)}$ in $\mr{GSpin}(T_{\Q})$, we then have the relation $r_{spin} \circ (\rho^{A(T)} \otimes E)= \rho^{B(T)}$, taking the Betti realization of $B(T)$ to be our model for the spin representation. We let $B_0$ equal $B(T)$ in the odd case and $B_{+}(T)$ in the even case.\index{s}{$B_0$} Similarly, we let $r_0$ denote $r_{spin}$ or one of the half-spin representations (which we may assume corresponds to $B_+(T)$).\index{s}{$r_0$} We also for convenience fix an embedding $E \into \Qlb$.
\begin{lemma}\label{spinlift}
(Without any assumption on $\End_{\Q-Hodge}(T_{\Q})$) There exists a lift $\tilde{\rho}^T_{\ell} \colon \gal{F} \to \mr{GSpin}(T_{\ell} \otimes \Qlb)$ of $\rho^T_{\ell}$ and a finite extension $F''/F'$ such that 
\[
r_{spin} \circ \tilde{\rho}^T_{\ell}|_{\gal{F''}} \cong H^1(B(T)_{\overline{F}}, \Ql) \otimes_{E} \Qlb |_{\gal{F''}}.
\]
\end{lemma}
\proof
As in the arguments of \S \ref{generalGalois}, we first choose a lift $\gal{F} \to \mr{GSpin}(T_{\ell} \otimes \Qlb)$ with finite-order Clifford norm. In the root datum notation of \S \ref{spineg}, the Hodge-Tate cocharacters $\mu_{\tau}$ of $\rho^T_{\ell}$, for all $\tau \colon F \into \Qlb$, are (conjugate to) $\lambda_1$, and the finite-order Clifford norm lifts have Sen operators (up to conjugacy) corresponding to the `co-characters'\footnote{That is, in the the co-character group tensored with $\Q$.} $\tilde{\mu}_{\tau}= \lambda_1$. We can modify this initial lift by twisting by $\lambda_0 \circ \omega'$, where $\omega' \colon \gal{F} \to \Qlb^\times$ is a Galois character whose square differs from the inverse $\omega^{-1}$ of the cyclotomic character by a finite-order twist (such $\omega'$ exists for any $F$). This gives a new lift $\tilde{\rho}^T_{\ell}$ whose composition with the Clifford norm (which, recall, is $2\chi_0$ in our notation) is $\omega^{-1}$, up to a finite-order twist, and whose $\tau$-labeled Hodge-Tate weights in the spin representation are $1$ and $0$ with equal multiplicity (compare the argument of \S \ref{liftinghodge}). Then $r_{spin} \circ \tilde{\rho}^T_{\ell}|_{\gal{F'}}$ differs from $H^1(B(T)_{\overline{F}}, \Ql) \otimes_E \Qlb$ by a finite-order twist, hence they are isomorphic after some additional finite base-change $F''/F'$.
\endproof
\begin{lemma}\label{softdescent}
There is a factor $M$ of $\Res_{F''/F}(B_0)$ (viewed as an object of $AV^0_{F, E}$), having endomorphisms by a finite extension $E'/E$, and an embedding $E' \into \Qlb$ extending our fixed $E \into \Qlb$ such that the associated $\ell$-adic realization $M_{\ell}$\footnote{Via $E' \into \Qlb$: that is, $M_{\ell}:= H^1(M_{\overline{F}}, \Ql) \otimes_{E'} \Qlb$.} is isomorphic as $\Qlb[\gal{F}]$-representation to $r_0 \circ \tilde{\rho}^T_{\ell}$. 
\end{lemma}
\proof
By the ($E$-linear) Tate conjecture,
\[
\End_{AV^0_{F, E}} \left( \Res_{F''/F}(B_0)\right) \otimes_E \Qlb \cong \End_{\Qlb[\gal{F}]}\left( \Ind_{F''}^F (H^1(B_{0, \overline{F}}, \Ql) \otimes_E \Qlb) \right).
\]
The Galois representation being induced is $r_0 \circ \tilde{\rho}^T_{\ell}$, and inside this endomorphism ring we can consider
\[
\Hom_{\Qlb[\gal{F}]} \left( r_0 \circ \tilde{\rho}^T_{\ell}, \Ind^F_{F''}(r_0 \circ \tilde{\rho}^T_{\ell})\right),
\] 
which by Frobenius reciprocity is just
\[
\End_{\Qlb[\gal{F''}]}\left( r_0 \circ \tilde{\rho}^T_{\ell} \right)= \Qlb,
\]
since $\End_{AV^0_{F'', E}}(B_0)= E$. In other words, there is a unique $\Qlb$-line in $\End_{AV^0_{F, E}} \left( \Res_{F''/F}(B_0)\right) \otimes_E \Qlb$ consisting of projectors onto the $r_0 \circ \tilde{\rho}^T_{\ell}$-isotypic piece of the $\ell$-adic representation. Decomposing the semi-simple $E$-algebra $\End_{AV^0_{F, E}} \left( \Res_{F''/F}(B_0)\right)$ into simple factors, we see that this line lives in a unique simple component (tensored with $\Qlb$), which itself must be just a finite field extension $E'$ of $E$ (else the $r_0 \circ \tilde{\rho}^T_{\ell}$-isotypic piece would have multiplicity greater than $1$); it then corresponds to exactly one of the simple factors of $E' \otimes_E \Qlb$, i.e. a particular embedding $E' \into \Qlb$. We can therefore take $M$ to be the abelian variety corresponding to this factor $E'$; $M$ has complex multiplication by $E'$, and via this specified embedding $E' \into \Qlb$, the $\ell$-adic realization $M_{\ell}$ is isomorphic to $r_0 \circ \tilde{\rho}^T_{\ell}$.
\endproof
\begin{rmk}
This is at its core the same proof as given in \S \ref{explicitvariant}; the latter proof is probably more transparent, but this one is somewhat `softer.' I don't think it translates as well to the more general context of \S \ref{nongenericodd}, however.
\end{rmk}
\begin{cor}
Theorem \ref{motiviclift} holds for the motive $Prim^{2k}(X)(k)$ over $F$.
\end{cor}
\proof
When $\dim(T_{\Q})$ is odd, we are done, by the previous lemma and Corollaries \ref{transcendentalreduction} and \ref{transcor}. When $\dim(T_{\Q})$ is even, we take the output $M$ in $AV^0_{F, E'}$ of the previous lemma, view it in $\mc{M}_{F, E'}$, and form the twisted dual $M^\vee(-1)$. Here there are two cases: if $\dim T_{\Q}$ is not divisible by four, this object corresponds to the composition of the other half-spin representation with $\tilde{\rho}^T_{\ell}$,\footnote{The two half-spin representations of $\mr{GSpin}_{2n}$ have highest weights $-\chi_0 + \frac{1}{2} (\sum_{i=1}^{n-1} \chi_i +\chi_n)$ (for $r_0$) and $-\chi_0 + \frac{1}{2} (\sum_{i=1}^{n-1} \chi_i -\chi_n)$ (for the other half-spin representation), so the lowest weight of $r_0^{\vee} \otimes (-2\chi_0)$ is $-\chi_0- \frac{1}{2}(\sum_{i=1}^n \chi_i)$, which, when $n$ is even, is visibly the lowest weight of the other half-spin representation.} and we can now apply the earlier corollaries to $M \oplus M^\vee(-1)$. If $\dim T_{\Q}$ is divisible by four (so the half-spin representations are self-dual), apply Lemma \ref{softdescent} to the composition of $\tilde{\rho}^T_{\ell}$ with the \textit{other} half-spin representation as well; so instead of a single motive $M$, we now have two motives $M_+$ and $M_-$\footnote{Enlarge the coefficient fields of $M_+$ and $M_-$, viewed as subfields of $\Qlb$ via the respective embeddings produced by Lemma \ref{softdescent}, to some common over-field.} such that the $\ell$-adic realization $(M_+ \oplus M_-)_{\ell}$ is isomorphic to $r_{spin} \circ \tilde{\rho}^T_{\ell}$. Then as before, we may apply Corollary \ref{transcor}.
\endproof
\subsection{Non-generic cases: $\dim(T_{\Q})$ odd}\label{nongenericodd}
To study the non-generic case $\End_{\Q-Hodge}(T_{\Q}) \neq \Q$, we have to understand the $\ell$-adic algebraic monodromy groups of the representations $\rho^T_{\ell}$, i.e. the $\ell$-adic analogue of Zarhin's result. Andr\'{e} has announced (\cite[Theorem 1.6.1]{andre:hyperkaehler}) a proof of the Mumford-Tate conjecture in this context\index{t}{Mumford-Tate conjecture for hyperk\"{a}hler varieties}; the arguments of \cite[\S 7.4]{andre:hyperkaehler} are not quite complete, but will be completed in a forthcoming paper (\cite{moonen:MTsurfaces}) of Ben Moonen, in the course of generalizing some of Andr\'{e}'s work. We therefore accept as proven the following: 
\begin{thm}[Andr\'{e}-Moonen]\label{andreMT}
Let $(X, \eta)$ be a polarized variety over a number field\footnote{Or, again, a finitely-generated extension of $\Q$.} $F$ satisfying $A_k$, $B_k$, and $B_k^+$. Then the inclusion $\overline{\rho_{\ell}(\gal{F})}^{Zar} \into \mc{G}^{V}_F \otimes_{\Q} \Ql$ of the algebraic monodromy group into the $\ell$-adic motivic Galois group of $V= Prim^{2k}(X)(k)$ is an isomorphism on connected components of the identity. 
\end{thm}
Recall that over the field $F'$, we may assume the groups $\overline{\rho_{\ell}(\gal{F'})}^{Zar}$ and $\mc{G}^{V}_{F'}$ are connected, and therefore isomorphic. Combining Theorem $6.5.1$ of \cite{andre:hyperkaehler} with Andr\'{e}'s result that Hodge cycles on abelian varieties are motivated, and with Zarhin's description (\cite{zarhin:K3hodge}) of the Mumford-Tate group of the transcendental lattice $T^{2k}(X_{\CC}, \Q)(k)$, we obtain (see Corollary $1.5.2$ of \cite{andre:hyperkaehler}):
\begin{cor}\label{explicitmonodromy}
The semisimple $\Q$-algebra $E_T:= \End_{\mc{G}^{T}_{F'}}(T)$ is a totally real or CM field, and there is a natural $E_T$-hermitian pairing $\langle \cdot , \cdot \rangle_{E_T} \colon T \times T \to E_T$. The motivic group (which equals the Mumford-Tate group, and equals, after $\otimes_{\Q} \Ql$, the $\ell$-adic algebraic monodromy group) $\mc{G}_{F'}^{T}$ is then isomorphic to the full orthogonal ($E_T$ totally real) or unitary ($E_T$ CM) group 
\[
\Aut(T, \langle \cdot , \cdot \rangle_{E_T}) \into \mr{SO}(T_{\Q}).
\]
\end{cor}
Before continuing, we formulate a variant of Ribet's method with coefficients:
\begin{prop}\label{generalribet}
Let $F'/F$ be an extension of number fields, and let $E/\Q$ be a finite extension. Suppose we are given an object $B$ of $AV^0_{F', E}$ such that for some embedding $E \into \Qlb$, the associated $\ell$-adic realization $B_{\ell}= H^1(B_{\overline{F}}, \Ql) \otimes_E \Qlb$ satisfies the invariance condition ${}^{\sigma} B_{\ell} \cong B_{\ell}$ for all $\sigma \in \gal{F}$. Further assume that $\End_{AV^0_{\overline{F}, E}}(B_{\overline{F}})= E$.\footnote{This necessarily holds for $B$ in $AV^0_{F', E}$ as well.} Then there exists a finite extension $E'/E$, an extension $E' \into \Qlb$ of the embedding $E \into \Qlb$, a number field $F''/F'$, and an object $M$ of $AV^0_{F, E'}$ such that 
\[
\left( H^1(M_{\overline{F}}, \Ql) \otimes_{E'} \Qlb \right) |_{\gal{F''}} \cong B_{\ell}|_{\gal{F''}}.
\]
That is, $B$, up to twist, has a motivic descent to $F$.
\end{prop}
\proof
This is proven as on page \pageref{ribet}, using the $E$-linear variant of Faltings' theorem: 
\[
\Hom_{AV^0_{F', E}}\left({}^{\sigma}B, B\right) \otimes_E \Qlb \xrightarrow{\sim} \Hom_{\Qlb[\gal{F'}]}\left(H^1({}^{\sigma}B_{\overline{F}}, \Ql) \otimes_E \Qlb, H^1(B_{\overline{F}}, \Ql) \otimes_E \Qlb \right).
\]
%and by the invariance assumption the right-hand side is non-trivial. Since $\End_{AV^0_{F', E}}(B)= E$, the left-hand side must then be a one-dimensional $\Qlb$-vector space, so we deduce the existence of an $E$-linear isogeny $\mu_{\sigma} \colon {}^{\sigma}B \to B$, arranging as before that $\mu_{\sigma}$ be trivial for $\sigma \in \gal{F'}$, and that it depend only on the $\gal{F'}$-coset of $\sigma$. This time the obstruction to descent $c_B$ lives in 
%\[
%H^2(\gal{F}, \End_{AV^0_{F', E}}(B)^\times)= H^2(\gal{F}, E^\times),
%\]
%and regarding $\Qb$ inside $\Qlb$ (say via $\iota_{\ell}$) we obtain an embedding $E \into \Qb$, and then we can trivialize the two-cocycle $c_B \in Z^2(\gal{F}, \Qb^\times)$ via a continuous $1$-co-chain $\alpha \colon \gal{F} \to \Q(\alpha)^{\times} \subset \Qb^\times$. Letting $F''$ denote an extension over which $\alpha$ is locally constant, we can then use $\alpha$ to define a homomorphism of $E$-algebras
%\[
%\mc{R}:= \End_{AV^0_{F, E}} \left( \Res_{F''/F}(B) \right) \xrightarrow{\alpha} \Q(\alpha).
%\]
In the notation of the earlier proof, $M= \Res_{F''/F}(B) \otimes_{\mc{R}, \alpha} \Q(\alpha)$ is the required motivic descent, to an isogeny abelian variety over $F$ with $\Q(\alpha)$-multiplication (so in the statement of the proposition, $E' \into \Qlb$ is $\Q(\alpha) \subset \Qb \into \Qlb$, extending the initial $E \into \Qlb$).
\endproof
We now assume that $\dim_{\Q} T_{\Q}$ is \textit{odd}, say of the form $cd$ where $d=2d_0+1= [E_T: \Q]$. In particular, $E_T$ is totally real. Let $\tilde{\rho}^T_{\ell}$ and $B(T)$ be as in Lemma \ref{spinlift}, and replace $F'$ by a large enough extension to satisfy the conclusion of that lemma, and such that $(\overline{\rho^T_{\ell}(\gal{F'})})^{Zar}$ is connected. Recall that $B(T)$ is an object of $AV^0_{F', E}$ for some finite extension $E/\Q$ large enough to split $C^+(T_{\Q})$. Corollary \ref{explicitmonodromy} implies that 
\[
\overline{ (\rho^T_{\ell} \otimes \Qlb)(\gal{F'}) }^{Zar} \cong \prod_{1}^d \mr{SO}_c(\Qlb) \subset \mr{SO}_{cd}(\Qlb).
\]
Restricting the spin representation $W_{cd}$ of $\mf{so}_{cd}(\Qlb)$ to $\prod_1^d \mf{so}_c(\Qlb)$, we obtain (via a weight calculation as in Lemma \ref{transclemma}, and writing $W_c$ for the spin representation of $\mf{so}_c$)
\[
W_{cd}|_{\prod \mf{so}_c} \cong (\boxtimes_1^d W_c )^{2^{d_0}}.
\]
The Lie algebra of the lift $\tilde{\rho}^T_{\ell}$ is one copy of the additive group $\mf{g}_a$ times the Lie algebra of $\rho^T_{\ell} \otimes \Qlb$, and this $\mf{g}_a$ acts by scalars in the spin representation, so $r_{spin} \circ \tilde{\rho}^T_{\ell}|_{\gal{F'}}$ has an analogous decomposition as $2^{d_0}$ copies of some (absolutely, Lie) irreducible representation $W'$ of $\gal{F'}$. For a suitable enlargement $E'$ of $E$ (and extension $E' \into \Qlb$), and $F''$ of $F'$, we can realize $W'$ as $H^1(B(T)', \Ql) \otimes_{E'} \Qlb$ for an object $B(T)'$ of $AV^0_{F'', E'}$.\footnote{More precisely, we need the decomposition of $W_{cd}|_{\prod \mr{so}_c}$ to be defined over $E'$; the first claim then follows from Faltings. The extension $F''/F'$ is needed to decompose $B(T)_{\overline{F}}\otimes_E E'$ over some finite extension.} Then $\End_{AV^0_{F'', E'}}(B(T)')= E'$, and the isomorphism $r_{spin} \circ \tilde{\rho}^T_{\ell} \cong (W')^{2^{d_0}}$ implies $\gal{F}$-conjugation invariance of $W'$. Thus we can apply Ribet's method to deduce:
\begin{lemma}
There exists a finite extension $E''/E'$, an embedding $E'' \into \Qlb$ extending the given $E' \into \Qlb$, and an object $M$ of $AV^0_{F, E''}$ such that 
\[
\left( H^1(M_{\overline{F}}, \Ql) \otimes_{E''} \Qlb |_{\gal{F'''}} \right) \cong W'|_{\gal{F'''}}
\]
for some still further finite extension $F'''/F''$.
\end{lemma} 
Let $M_{\ell}$ denote the associated $\ell$-adic realization (via $E'' \into \Qlb$). Since $r_{spin} \circ \tilde{\rho}^T_{\ell}$ is Lie-isotypic, and $M_{\ell}$ is a descent to $\gal{F}$ of its unique (after finite restriction) Lie-irreducible constituent, Corollary \ref{lieisotypic} shows that there is an Artin representation $\omega$ of $\gal{F}$ such that
\[
r_{spin} \circ \tilde{\rho}^T_{\ell} \cong M_{\ell} \otimes \omega.
\]
Possibly enlarging the field of coefficients yet again, we deduce:
\begin{thm}\label{nongenericoddmotiviclift}
Suppose $\dim T_{\Q}$ is odd, and that $\det T_{\ell}=1$ (as $\gal{F}$-representation). Then there exists 
\begin{itemize}
\item a number field $\tilde{E}$ and an embedding $\tilde{E} \into \Qlb$; 
\item an object $\tilde{M}$ of $\mc{M}_{F, \tilde{E}}$ that is a tensor product of an Artin motive and (the image of) an object of $AV^0_{F, \tilde{E}}$;
\item and a lift $\tilde{\rho}^T_{\ell} \colon \gal{F} \to \mr{GSpin}(T_{\ell} \otimes \Qlb)$ of $\rho^T_{\ell}$;
\end{itemize}
such that $r_{spin} \circ \tilde{\rho}^T_{\ell}$ is isomorphic to the $\ell$-adic realization (via $\tilde{E} \into \Qlb$) of $\tilde{M}$. Moreover:
\begin{itemize}
\item $\tilde{\rho}^T_{\ell}$ lives in a weakly-compatible system of lifts;
\item $\tilde{\rho}^T_{\ell}$ actually arises from a lifting of representations of the motivic Galois group $\mc{G}_{F, \tilde{E}}$, as in Theorem \ref{motiviclift}.
\item The same conclusions hold with the motive $V= Prim^{2k}(X)(k)$ in place of its transcendental lattice $T$ (and, again, a possible enlargement of $\tilde{E}$).
\end{itemize}
\end{thm}
\proof
The number field $\tilde{E}$ is the composite (inside the ambient $\Qlb$) of the field $E''$ and the field needed to define the Artin representation $\omega$. To conclude the proof of the theorem, we make three observations:
\begin{itemize}
\item $\mc{M}_{F, \tilde{E}}$ is Tannakian (note that we already know that $M$\footnote{Precisely, $M \otimes_{E''} \tilde{E}$.} and $\omega$ are motivic);
\item $M$ and $\omega$ both give rise to compatible systems of $\ell$-adic representations;
\item Corollary \ref{transcor} applies to lift the representations of motivic Galois groups.
\end{itemize}
\endproof
\subsection{Non-generic cases: $\dim T_{\Q}$ even}\label{Teven}
We do not fully treat the case of even-rank transcendental lattice, but here give a couple examples, describing the `shape' of the Galois representations in light of Proposition \ref{lietensorart}.

First, continue to assume $E_T$ is totally real. Let $\dim T_{\Q}= 2n= cd$, with $d= [E_T:\Q]$. For any $N$, denote by $W_{2N, \pm}$ (for each choice of $\pm$) the two half-spin representations of $\mf{so}_{2N}$, and continue to write $W_{2N+1}$ for the spin representation of $\mf{so}_{2N+1}$.\index{s}{$W_n$, for some odd integer $n$}\index{s}{$W_{n, +}$, $W_{n, -}$, for some even integer $n$}
\begin{lemma}\label{sorescalc}
When $c$ is even, the restriction $W_{cd, +}|_{\prod \mf{so}_c}$ is given by
\[
W_{cd, +}|_{\prod_1^d \mf{so}_c} \cong \bigoplus_{\substack{\epsilon= (\epsilon_i) \in \{\pm\}^d \\ \prod \epsilon_i=1}} \boxtimes_{i=1}^d W_{c, \epsilon_i},
\]
where the indexing set ranges over all choices of signs with $-$ occurring an even number of times. This is a direct sum of distinct Lie-irreducible representations.

When $c$ is odd, so $d= 2d_0$ is even, 
\[
W_{cd, +}|_{\prod_1^d \mf{so}_c} \cong 2^{d_0-1} \boxtimes_1^d W_c, 
\]
a single Lie-irreducible representation occurring with multiplicity $2^{d_0-1}$.
\end{lemma}
Now, recall (Lemma \ref{spinlift}) that after a sufficient base-change $F'/F$, we can find an abelian variety $B_+(T)$ (with coefficients in a number field $E$, embedded in $\Qlb$; we may by extending scalars assume $E$ is large enough that the above decomposition of spin representations is defined over $E$), and a lift $\tilde{\rho}^T_{\ell}$ such that $r_+ \circ \tilde{\rho}^T_{\ell}|_{\gal{F'}}$ is isomorphic to $H^1(B_+(T), \Ql) \otimes_E \Qlb$. We assume $F'$ sufficiently large that this Galois representation is a sum of Lie-irreducible representations.
\begin{prop}\label{gettingtired}
Suppose $c$ is even. Then motivic lifting holds for $\rho^T$.
\end{prop}
\proof
Since $c$ is even, the previous lemma shows that $r_+ \circ \tilde{\rho}^T_{\ell}$ is Lie-multiplicity-free, hence is a direct sum of inductions of non-conjugate, Lie-irreducible Galois representations. If $\pi_0(\rho^T_{\ell})$ is trivial, in which case $\tilde{\rho}^T_{\ell}$ also has connected monodromy group,\footnote{Since it contains the center of $\mr{GSpin}$; the sort of example this avoids is $\rho_{f, \ell} \otimes \omega^{\frac{1-k}{2}}$, where $f$ is a classical modular form of odd weight $k$.} then no inductions occur in this decomposition, so each factor $\boxtimes_1^d W_{c, \epsilon_i}$ in Lemma \ref{sorescalc} corresponds to a Lie-irreducible factor of $r_+ \circ \tilde{\rho}^T_{\ell}$, as $\gal{F}$-representation. Each of these factors is, over $F'$, of the form
\[
H^1(B_{\epsilon}, \Ql) \otimes_{E} \Qlb,
\]
where $B_{\epsilon}$ is an object of $AV^0_{F', E}$ with endomorphism algebra just $E$ itself (the usual application of Faltings' theorem, using the Lie-multiplicity-free property, and the fact that the spin representation decomposition holds over $E$). By $\gal{F}$-invariance of each factor, we can apply Proposition \ref{generalribet} to produce an object $M_{\epsilon}$ of $AV^0_{F, E_{\epsilon}}$, for some finite extension $E_{\epsilon}/E$ inside $\Qlb$, with $M_{\epsilon, \ell}|_{\gal{F'}}$ isomorphic to a finite-order twist of $H^1(B_{\epsilon}, \Ql) \otimes_{E_{\epsilon}} \Qlb$. By Lie-irreducibility, some finite-order twist (a character of $\Gal(F'/F)$, in fact, and, again, we may have to enlarge $E_{\epsilon}$) $M_{\epsilon}'$ of $M_{\epsilon}$ has $\ell$-adic realization isomorphic to the corresponding factor of $r_+ \circ \tilde{\rho}^T_{\ell}$.  Inside the ambient $\Qlb$, we take the compositum $E'$ of the various $E_{\epsilon}$, extending scalars on each $M_{\epsilon}'$. Then the object $\oplus_{\epsilon} M_{\epsilon}'$ of $\mc{M}_{F, E'}$ satisfies the hypotheses of Corollary \ref{transcor}, so we deduce the existence of the desired motivic lift.

We sketch the case of non-connected monodromy. Take the (motivic) Lie-irreducible factors of $r_+ \circ \tilde{\rho}^T_{\ell}|_{\gal{F'}}$, and partition them into $\gal{F}$-orbits. Fix a representative $W_i$ of each orbit, and consider the stabilizer $\gal{F_i}$ in $\gal{F}$ of $W_i$. Arguing as in Proposition \ref{gettingtired}, we can apply Ribet's method to descend each $W_i$ to an object $M_i$ of $\mc{M}_{F_i, E_i}$ for some extension $E_i$ of $E$ inside $\Qlb$. Moreover, twisting yields an $M_i'$ whose $\ell$-adic realization is isomorphic to a factor of $r_+ \circ \tilde{\rho}^T_{\ell}$.\footnote{The isomorphisms $\Hom_{\gal{F'}}(M_1, V) = \Hom_{\gal{F_i}}(\Ind_{F'}^{F_i}(M_i), V)= \Hom_{\gal{F_i}}(M_i \otimes \Qlb[\Gal(F'/F_i)], V)$, with $V= r_+ \circ \tilde{\rho}^T_{\ell}$, imply this claim, using the fact that $M_i$ is Lie-irreducible and after finite restriction occurs with multiplicity one in $r_+ \circ \tilde{\rho}^T_{\ell}$.} Then we can induce (the representation of motivic Galois groups) from $F_i$ to $F$ to obtain our motives over $F$. This completes the case of even $c$ (and $E_T$ totally real).
\endproof

Having demonstrated the available techniques in a couple of quite different situations (namely, where the lifts range between the extremes of being Lie-multiplicity-free and Lie-isotypic), we stop here, remarking only that when $E_T$ is CM, the analysis must begin not from Lemma \ref{sorescalc} but from the restriction $W_{cd, +}|_{\prod \mf{gl}_c}$, the product ranging over pairs of complex-conjugate embeddings $E_T \into \Qlb$. This restriction is given by:
\begin{lemma}
Suppose $d=2d_0$ is even, with notation otherwise as above. We denote irreducible representations of $\mf{gl}_c$ by $W(r)$, where $W$ is an irreducible representation of $\mf{sl}_c$, and $(r)$ indicates that the restriction to the center $\mf{g}_a \subset \mf{gl}_c$ is multiplication by $r$. Then, letting $V_i$ denote the standard representation of the $i^{th}$ copy of $\mf{sl}_c$ 
\[
W_{cd, +}|_{\prod_1^{d_0} \mf{gl}_c} \cong \bigoplus_{\substack{i_1, \ldots, i_{d_0}:\\
\sum i_j \in 2\Z}} \wedge^{i_1}V_1^*(\frac{c}{2}-i_1) \boxtimes \ldots \boxtimes \wedge^{i_{d_0}}V_{d_0}^*(\frac{c}{2}-i_{d_0}).
\]
\end{lemma}
\begin{rmk}
\begin{itemize}
\item Note that the representations occurring here do not necessarily extend to representations of $\mr{GL}_c$, since $c$ may be odd; they do extend on the (connected) double-cover of $\mr{GL}_c$.
\item In particular, we see that when $E_T$ is CM, $r_+ \circ \tilde{\rho}^T_{\ell}$ is Lie-multiplicity-free. This suggests proceeding as in Proposition \ref{gettingtired}, although we will stop here.
\end{itemize}
\end{rmk}
\section{Towards a generalized Kuga-Satake theory}\label{conclusion}
\subsection{A conjecture}
It is fair to assume that one could establish a motivic lifting result for the remaining hyperk\"{a}hler cases. More important, these lifting results clamor for generalization. Motivated by the Fontaine-Mazur conjecture, Theorem \ref{anymonodromylift}, Proposition \ref{taniyamalift}, and Theorem \ref{nongenericoddmotiviclift}, we are led to the following much more ambitious conjecture:
\begin{conj}\label{generalizedKS}
Let $F$ and $E$ be number fields, and let $H' \onto H$ be a surjection, with central torus kernel, of linear algebraic groups over $E$. Suppose we are given a homomorphism $\rho \colon \mc{G}_{F, E} \to H$. Then if $F$ is imaginary, there is a finite extension $E'/E$ and a homomorphism $\tilde{\rho} \colon \mc{G}_{F, E'} \to H'_{E'}$ lifting $\rho \otimes_E E'$. If $F$ is totally real, then such a lift exists if and only if the Hodge number parity obstruction of Corollary \ref{fullmonodromytotreal} vanishes.
\end{conj}
To indicate the scope of this conjecture, let $X/F$ be any smooth projective variety, and consider for any $k \leq \dim X$ the motive $H^{2k}(X)(k)$ (or $Prim^{2k}(X)(k)$, having chosen an ample line bundle). This gives rise to an orthogonal representation of $\mc{G}_F$, and the conjecture in this case (for $\mr{GSpin} \to \mr{SO}$, or the variant with the full orthogonal group in place of $\mr{SO}$) amounts to a generalization of the Kuga-Satake construction to \textit{arbitrary} orthogonally-polarized motivic (over $F$) Hodge structures. For other choices of $H'$ and $H$ (for instance, $\mr{GL}_n \to \mr{PGL}_n$, where necessarily we will have coefficient field larger than $\Q$), the conjectured generalization is even more mysterious.\footnote{Compare Lemma \ref{notCM}.}

Note the role of Lemma \ref{spreadout} in building our confidence in this conjecture. Let $H' \to H$ be a morphism of groups, say over $\Qb$, as in the conjecture. One way of formulating the Fontaine-Mazur conjecture is that $\mc{G}_F \otimes_{\Q} \Qlb$ should be isomorphic to the Tannakian group for the Tannakian category $\Rep^{g, ss}_{\Qlb}(\gal{F})$ of semi-simple geometric $\Qlb$-representations of $\gal{F}$.  In particular, a geometric lift $\tilde{\rho}_{\ell} \colon \gal{F} \to H'(\Qlb)$ of $\rho_{\ell} \colon \gal{F} \to H(\Qlb)$ should arise from an algebraic homomorphism of $\Qlb$-groups $\mc{G}_F \otimes \Qlb \to H'_{\Qlb}$. Lemma \ref{spreadout} then tells us that if $\rho$ arises from some $\mc{G}_F \otimes \Qb \to H$, then we can find a lift $\tilde{\rho} \colon \mc{G}_{F} \otimes \Qb \to H'$ of homomorphisms of $\Qb$-groups.

Now that we know to look for such a thing, we conclude by giving one more example in which it is easy to construct a `generalized Kuga-Satake motive.'
\subsection{Motivic lifting: abelian varieties}\label{AVlifting}
If $A/\CC$ is an abelian surface, then $H^2(A, \Q) \cong \wedge^2 H^1(A, \Q)$ has Hodge numbers $h^{2, 0}=1$ and $h^{1, -1}=4$. The classical Kuga-Satake construction then associates to $A$ another abelian variety, $KS(A)$, with rational cohomology $C^+(H^2(A, \Q))$ (or the analogue with $Prim^2(A, \Q)$, if we have fixed a polarization), and a theorem of Morrison (\cite{morrison:KSabsurface}) asserts that $KS(A)$ is isogenous to $A^8$ (or $A^4$ if we work with $Prim^2(A, \Q)$). We now show that this construction can be generalized to abelian varieties of any dimension.

Let $F$ be any subfield of $\CC$, and let $A/F$ be an abelian variety of dimension $g$ with a fixed polarization. Consider the algebraic representation 
\[
\rho \colon \mc{G}_F \to \mr{GSp}(H^1(A_{\CC}, \Q)) \cong \mr{GSp}_{2g},
\]
or, for the theory over $\CC$, its restriction to $\mc{G}_{\CC}$. We will compose $\rho$ with a homomorphism $\mr{GSp}_{2g} \to \mr{GSpin}_N$ for suitable $N$ to produce the generalized Kuga-Satake motive. Let $W= H^1(A_{\CC}, \Q)$, and let $r_{e_1+e_2} \colon \mr{Sp}(V) \to \mr{GL}(V_{e_1+e_2})$, in the weight notation of \S \ref{spineg}, denote the irreducible representation of $\mr{Sp}(V)$ obtained as the complement of the trivial representation in $\wedge^2(W)$. This absolutely irreducible representation, defined over $\Q$, has image in $\mr{SO}(V_{e_1+e_2})$, where the quadratic form is induced from the pairing canonically induced on $\wedge^2(W)$, and which coincides with the induced polarization on $Prim^2(A_{\CC}, \Q)$. Since $\mr{Sp}(W)$ is simply-connected, there is a lift to an algebraic homomorphism $\tilde{r} \colon \mr{Sp}(W) \to \mr{Spin}(V_{e_1+e_2})$.\footnote{Over $\CC$ this follows from standard Lie theory, and such a homomorphism descends to $\Qb$ (see Lemma \ref{spreadout}); since the map of Lie algebras $\mf{sp}(W) \to \mf{so}(V_{e_1+e_2})$ is defined over $\Q$, we see that the map over $\Qb$ is $\gal{\Q}$-invariant, hence descends to a morphism over $\Q$.} The dimension of $V_{e_1+e_2}$ is $\binom{2g}{2}-1$, which is odd ($=2n+1$) if $g$ is even, and even ($=2n$) if $g$ is odd. In either case, we have the representation $r_{\varpi_n}$ of $\mr{Spin}(V_{e_1+e_2})_{\Qb}$ having highest weight $\varpi_n= \frac{\sum_{i=1}^n \chi_i}{2}$.\footnote{We use the common fundamental weight notation here. This is the spin representation in the odd case; one of the half-spin representations in the even case.}
\begin{lemma}\label{where'sminus1}
Let $c$ denote the non-trivial (central) element of the kernel of $\mr{Spin}(V_{e_1+e_2}) \to \mr{SO}(V_{e_1+e_2})$. If $g \equiv 2, 3 \pmod 4$, then $\tilde{r}(-1)=c$, and if $g \equiv 0, 1 \pmod 4$, then $\tilde{r}(-1)=1$.
\end{lemma}
\proof
In all cases, $r_{e_1+e_2}(-1)=1$, so $\tilde{r}(-1)$ equals either $1$ or $c$. Since $c$ acts as $-1$ in any of the spin representations ($c$ \textit{is} the element $-1$ of the Clifford algebra), it suffices to compute $r_{\varpi_n} \circ \tilde{r}(-1)$. The weights of $V_{e_1+e_2}$ are
\[
\{\pm(e_i+e_j)\}_{1 \leq i < j \leq g} \cup \{e_i-e_j\}_{i \neq j},
\]
except with one copy of the weight zero deleted (so that zero has multiplicity $g-1$ rather than $g$). It follows that $r_{\varpi_n} \circ \tilde{r}$ has a weight equal to 
\[
\frac{1}{2} \left( \sum_{1 \leq i < j \leq g} (e_i+e_j) + \sum_{1 \leq i < j \leq g}(e_i- e_j) \right)= \sum_{i=1}^g (g-i) e_i.
\]
In particular, $r_{\varpi_n} \circ \tilde{r}(-1)$ is multiplication by $(-1)^{g(g-1)/2}$, and the lemma follows.
\endproof
\begin{cor}\label{AVKS}
$\tilde{r} \colon \mr{Sp}(W) \to \mr{Spin}(V_{e_1+e_2})$ extends to an algebraic homomorphism $\mr{GSp}(W) \to \mr{GSpin}(V_{e_1+e_2})$. If $g \equiv 2, 3 \pmod 4$, then this map can be chosen so the Clifford norm coincides with the symplectic multiplier; if $g \equiv 0, 1 \pmod 4$, then this map can be chosen to factor through $\mr{Spin}(V_{e_1+e_2})$. The composition
\[
\mc{G}_F \xrightarrow{\rho} \mr{GSp}(W) \xrightarrow{\tilde{r}} \mr{GSpin}(V_{e_1+e_2}) \into \mr{GL}(C^+(V_{e_1+e_2}))
\]
defines the generalized Kuga-Satake lift of $A$.
\end{cor}
\proof
This follows immediately from Lemma \ref{where'sminus1} and the identifications:
\begin{align*}
 \frac{\mathbb{G}_m \times \mr{Sp}(W)}{\langle (-1, -1) \rangle} & \xrightarrow{\sim} \mr{GSp}(W) \\
 \frac{\mathbb{G}_m \times \mr{Spin}(V_{e_1+e_2})}{\langle (-1, c) \rangle} & \xrightarrow{\sim} \mr{GSpin}(V_{e_1+e_2}).
\end{align*}
When $g \equiv 2, 3 \pmod 4$, we take the map $\mathbb{G}_m \to \mathbb{G}_m$ to be the identity, and when $g \equiv 0, 1 \pmod 4$, we take it to be trivial.
\endproof
\begin{rmk}
\begin{itemize}
\item Repeating the above arguments with $\wedge^2(W)$ in place of $V_{e_1+e_2}$, we can similarly construct lifts $\mc{G}_F \to \mr{GSpin}(\wedge^2(W))$.
\item When $g=2$, this recovers the classical construction. In that case, the composition $r_{\varpi_n} \circ \tilde{r}$ is the identity, and, decomposing $C^+(V_{e_1+e_2})$ as $4$ copies of $r_{\varpi_n}$ (as $\mr{GSpin}(V_{e_1+e_2})$ representation), the identification (up to isogeny) $KS(A) \sim A^4$ is nearly a tautology. 
\end{itemize}
\end{rmk}
The motivic formalism now tells us that $r_{\varpi_n} \circ \tilde{r} \circ \rho$ is (the Betti realization of) an object of $\mc{M}_F$. Since the $\ell$-adic realizations of $\rho$ form a weakly\footnote{In fact, strictly.} compatible system, the same is true for the $\ell$-adic realizations of this Kuga-Satake motive. In this case, however, we can say more, and will realize this explicitly (and unconditionally) as a Grothendieck motive. The first step is to compute the plethysm $r_{\varpi_n} \circ \tilde{r}$; we will do this, but first mention an equivalent, structurally appealing plethysm.  First we treat the $D_n$ case, that is when $d= 2g= \dim H^1(X)$ satisfies $2n= \binom{d}{2}$; this amounts to $\dim X$ being even. Let $V= H^2(X)= \wedge^2 W$. We use the following (common) notation for fundamental weights of $D_n$:
\begin{itemize}
\item $\varpi_i= \chi_1+ \ldots+ \chi_i$, and $r_{\varpi_i}= \wedge^i V$, for $i=1, \ldots, n-2$;
\item $\varpi_{n-1}= \frac{\chi_1+ \ldots +\chi_{n-1}-\chi_n}{2}$, $\varpi_n= \frac{\sum_1^n \chi_i}{2}$, and $r_{\varpi_{n-1}}$, $r_{\varpi_n}$ are the two half-spin representations.
\end{itemize}
As representations of $\mf{so}_{2n}$, we have the following identity:
\[
\bigoplus_{i=0}^{2n} \wedge^i(V)= \left(r_{\varpi_{n-1}} \oplus r_{\varpi_n} \right)^{\otimes 2}.
\]
The plethysm problem that we expect to solve, then, is to describe a representation $r_{KS(X)}$ of $\mr{Sp}(W)$ such that
\[
\bigoplus_{i=0}^{2n} \wedge^i(\wedge^2 W)= (r_{KS(X)})^{\otimes 2}.
\]
(More ambitiously, we could attempt this with $\wedge^{2k}(W)$ instead.) Similarly, in the $B_n$ case, we have fundamental weights $\varpi_i= \sum_1^i \chi_i$ for $i=1, \ldots, n-1$, and $\varpi_n= \frac{\sum \chi_i}{2}$, corresponding, respectively to the wedge powers $\wedge^i V$ of the standard representation ($r_{\varpi_1}$) and the spin representation. As representations of $\mf{so}_{2n+1}$, we have the identity
\[
\bigoplus_{i=0}^n \wedge^i(V)= r_{\varpi_n}^{\otimes 2},
\]
so in this case we want a representation $r_{KS}$ of $\mr{Sp}(W)$ satisfying
\[
\bigoplus_{i=0}^n \wedge^i(\wedge^2(W))= r_{KS}^{\otimes 2}.
\]
\begin{prop}
Writing $\omega_i$ ($i=1, \ldots, g$) for the usual fundamental weights\footnote{$\omega_i= e_1+ \ldots+ e_i$ in the standard coordinate system} of $C_g$, then in all cases ($g$ odd or even) we have
\[
\wedge^\bullet (\wedge^2(W)) \cong 2^g (V_{\omega_1+ \ldots + \omega_{g-1}})^{\otimes 2},
\]
as $\mr{Sp}(W)$-representations. Here $\wedge^\bullet$ denotes the full exterior algebra. This is deduced from the above discussion and the two calculations:
\begin{itemize}
\item (even) As $\mr{Sp}(W)$-representations,
\[
(r_{\varpi_{n-1}} \oplus r_{\varpi_n})\circ \tilde{r} \cong 2^{g/2} V_{\omega_1 + \ldots+ \omega_{g-1}}.
\]
\item (odd) As $\mr{Sp}(W)$-representations,
\[
r_{\varpi_n} \circ \tilde{r} \cong 2^{(g-1)/2} V_{\omega_1+ \ldots+ \omega_{g-1}}.
\]
\end{itemize}
\end{prop}
\proof
We treat the even case, the odd case being essentially identical. There are $g^2-g$ non-zero weights in $\wedge^2(W)$, and the weight zero occurs with multiplicity $g$. The weights of $(r_{\varpi_{n-1}} \oplus r_{\varpi_n})\circ \tilde{r}$ are sums of plus or minus any $n= g^2-\frac{g}{2}$ weights of $\wedge^2(W)$, then the total divided by two. It follows that the highest weight of $(r_{\varpi_{n-1}} \oplus r_{\varpi_n})\circ \tilde{r}$ is, as previously computed, $\sum_{i=1}^{g-1} \omega_i$, but moreover that it occurs with multiplicity $2^{g/2}$ (here $\frac{g}{2}= n-(g^2-g)$; we can choose the weight $+0$ or $-0$ this many times). It follows that $(r_{\varpi_{n-1}} \oplus r_{\varpi_n})\circ \tilde{r}$ contains $2^{g/2} V_{\omega_1+ \ldots+\omega_{g-1}}$; by dimension count (see the following Lemma \ref{dimcount}), they are isomorphic. 
\endproof
\begin{lemma}\label{dimcount}
With the above notation, the dimension of the irreducible $\mr{Sp}(W)$-representation $V_{\omega_1+ \ldots+ \omega_{g-1}}$ is $2^{g(g-1)}$.
\end{lemma}
\proof
Simplifying the Weyl dimension formula, we find
\[
\dim V_{\omega_1+ \ldots + \omega_{g-1}}= 2^{|\Phi^+|} \prod_{1 \leq i \leq j \leq g} \frac{2g+1-(i+j)}{2g+2-(i+j)},
\]
where the number $|\Phi^+|$ of positive roots is $2^{g^2}$. The product telescopes: fix an $i$, and the corresponding product over $j$ is equal to $\frac{1}{2}$. The lemma follows.
\endproof
\begin{cor}
Let $F$ be any subfield of $\CC$, and let $X/F$ be an abelian variety of any dimension $g$, giving rise to a representation 
\[
\rho_{H^1(X)} \colon \mc{G}_F \to \mr{GSp}(H^1(X, \Q)).
\]
Then the Kuga-Satake lift (Corollary \ref{AVKS} and remark following) of the representation 
\[
\mc{G}_F \to \mr{SO}(H^2(X, \Q)(1))
\]
can be explicitly realized in the spin (or sum of half-spin) representation as $2^{\lfloor \frac{g}{2} \rfloor}$ copies of the composition $r_{\omega_1+ \ldots+ \omega_{g-1}} \circ \rho_{H^1(X)}$ with the highest weight $\omega_1+ \ldots +\omega_{g-1}$ representation of $\mr{Sp}(H^1(X, \Q))$.\footnote{Extended to $\mr{GSp}$ on the center by the prescription of Corollary \ref{AVKS}. Throughout this argument we will ignore these extra scalars to simplify the notation.}

Let us call this object of $\mc{M}_F$ the Kuga-Satake motive $KS(X)$. Then $KS(X)$ is in fact a Grothendieck motive, for either numerical or homological equivalence.\index{s}{$KS(X)$, for an abelian variety $X$}\index{t}{Kuga-Satake motive associated to an abelian variety}
\end{cor}
\proof
It remains to check that $KS(X)$ can be cut out by algebraic cycles. We start with the explicit description (due to Weyl; see \S $17.3$ of \cite{fulton-harris:reptheory}) of the representation $V_{\omega_1+\dots+\omega_{g-1}}$. From now on, abbreviate $\lambda= \sum_1^{g-1} \omega_i$, and $r= \sum_1^{g-1} i$ (the length of $\lambda$); in fact, what follows applies to an arbitrary partition $\lambda$. Then
\[
V_{\lambda}= W^{\langle r \rangle} \cap \mbb{S}_{\lambda}(W)
\]
as $\mr{Sp}(W)$-representation. We have to explain this notation:\index{s}{$\mbb{S}_{\lambda}$}\index{s}{$c_{\lambda}$}\index{s}{$W^{\langle r \rangle}$} $\mbb{S}_{\lambda}(W)$ denotes the Schur functor associated to the partition $\lambda$, which explicitly is equal to the image of the Young symmetrizer $c_{\lambda}$ acting on $W^{\otimes r}$; and $W^{\langle r \rangle}$ is the subspace of $W^{\otimes r}$ given by intersecting the kernels of all the contractions ($1 \leq p < q \leq r$)
\begin{align*}
c_{p, q} \colon W^{\otimes r} &\xrightarrow{} W^{\otimes (r-2)} \\
v_1 \otimes \ldots \otimes v_{r} &\mapsto \langle v_p, v_q \rangle \cdot v_1 \otimes \ldots \otimes \hat{v}_p \otimes \ldots \otimes \hat{v}_q \otimes \ldots \otimes v_{r},
\end{align*}
where $\langle \cdot, \cdot \rangle$ represents the symplectic form on $W$. The result will now follow from some basic facts about algebraic cycles and, crucially, the Lefschetz and K\"{u}nneth\footnote{We will write $\pi^i_X$ for the algebraic cycle on $X \times X$ inducing $H(X) \onto H^i(X) \into H(X)$.} Standard Conjectures for abelian varieties (due to Lieberman; for a proof, see \cite{kleiman:algcycles}). Fix a polarization of $X$, giving rise to a Lefschetz operator $L_X$. Recall Jannsen's fundamental result (\cite{jannsen:numerical}) that numerical motives form a semi-simple abelian category. This and the K\"{u}nneth Standard Conjecture for abelian varieties imply that the category of numerical motives generated by abelian varieties over $F$ is a semi-simple Tannakian category. We deduce that the following are (numerical) sub-motives of $H(X^r)$:
\begin{itemize}
\item $M= (X, \pi^1_X, 0)$, i.e. the object corresponding to $H^1(X)$;
\item The kernel of each contraction $H^1(X)^{\otimes r} \to H^1(X)^{\otimes (r-2)}$. For notational simplicity, take $p=1, q=2$. If we were working in cohomology, we would compute the kernel of $c_{1, 2}$ as
\[
\ker \left( H^1(X) \otimes H^1(X) \xrightarrow{\langle \cdot, \cdot \rangle} \Q(-1) \right) \otimes H^1(X)^{\otimes r-2},
\]
where recall that the polarization $\langle \cdot, \cdot \rangle$ is defined by $\langle w_1, w_2 \rangle= L_X^{d-1} (w_1 \cup w_2)$. Now, $W^{\otimes 2}= \Sym^2(W) \oplus \wedge^2(W)$, and cup-product kills $\Sym^2(W)$, while mapping $\wedge^2(W)$ isomorphically to $H^2(X)$. The kernel of $c_{1,2}$ is therefore isomorphic to
\[
\left( \Sym^2(H^1(X)) \oplus Prim^2(X) \right) \otimes H^1(X)^{\otimes r-2}.
\]
To realize the analogous numerical motive, we can define $\Sym^2(M)$ as 
\[
(X \times X, \frac{1}{2}(1+ (1  2))\cdot (\pi_X^1 \times \pi_X^1), 0).
\]
For the projector, we have taken the product of two commuting idempotents, the first one being the Young symmetrizer associated to the Schur functor $\Sym^2$. There is also a projector $p^2$ onto primitive cohomology (see $1.4.4$, $2.3$, and $2A11$ of \cite{kleiman:algcycles}), so $(X, p^2, 0)$ realizes $Prim^2(X)$ as a numerical motive. Finally, since direct sums and tensor products exist in our category, we can therefore describe the kernel of any $c_{p, q}$ as a numerical motive $M_{p,q}$. 
\item In general, the Young symmetrizer $c_{\lambda}$ is not quite an idempotent. Write $c_{\lambda}^2= x^{-1} c_{\lambda}$,\footnote{$x$ depends on $\lambda$, but we have fixed a $\lambda$.} making $x c_{\lambda}$ an idempotent algebraic correspondence on $X^r$. $c_{\lambda}$ commutes with the $r$-fold K\"{u}nneth projector $(\pi^1_X)^r$ (this is easily checked at the level of cohomological correspondences, so \textit{a fortiori} holds for numerical equivalence), so 
\[
(X^r, x c_{\lambda} \cdot (\pi_X^1)^r, 0)
\]
defines a numerical motive, denoted $\mbb{S}_{\lambda}(H^1(X))$.
\item We would like to conclude the proof by intersecting the objects $M_{p, q}$ ($p<q$) and $\mbb{S}_{\lambda}(H^1(X))$. This is possible for numerical motives since the category is abelian. Finally, since all of the cycles considered in the proof are cycles on (disjoint unions of) abelian varieties in characteristic zero, where numerical and homological equivalence coincide,\footnote{In the presence of the Standard Conjecture of Hodge type, the Lefschetz conjecture implies `num $=$ hom': see Proposition $5.1$ of \cite{kleiman:standardconjectures}.} we also deduce the existence of a homological motive corresponding to $KS(X)$.
\end{itemize}
\endproof
\begin{rmk}
A weight calculation yields the Hodge numbers of $KS(X)$.
\end{rmk}
\subsection{Coda}
Identifying among all rational Hodge structures those that are motivic is one of the fundamental problems of complex algebraic geometry, and this generalized Kuga-Satake theory would systematically construct new \textit{motivic} Hodge structures from old, in a way not achievable by simply playing the Tannakian game. I hope that investigation of these phenomena will provide a stimulating testing-ground for thinking about Hodge theory in non-classical weights.

It is also tempting to ask what should be true if we replace $F$ by other fields, especially finitely-generated subfields of $\CC$. Our motivic descent in the hyperk\"{a}hler case works as written, except for the critical absence of Tate's basic vanishing result; in other contexts, it might be hoped that a similar descent works, conditional on the relevant cases of the Tate conjecture. I suspect only a qualitative \textit{potential} lifting result will hold in this generality-- that would suffice to imply the analogous lifting conjecture for $\CC$ itself. But another clearly important question to ask is: what, if any, is the analogue of Tate's vanishing theorem when the number field $F$ is replaced by any field finitely-generated over $\Q$?

\backmatter
\Printindex{s}{Index of symbols}
\Printindex{t}{Index of terms and concepts}

\bibliographystyle{amsalpha}
\bibliography{biblio.bib}

\end{document}